\documentclass[11pt]{article}
\usepackage{amsfonts,amsmath,amssymb,graphicx,color}
\usepackage{amsthm}
\usepackage[normalem]{ulem}
\usepackage{psfrag}
\usepackage{pdflscape}
\usepackage{multirow}
\usepackage{blindtext} 
\usepackage{enumitem}
\graphicspath{ {./figures/} }
\usepackage{hyperref}
\hypersetup{colorlinks=true}
\usepackage{algorithm}
\usepackage{algpseudocode}
\usepackage{epstopdf}
\usepackage{framed}
\usepackage{footnote}
\usepackage{geometry}
\usepackage{pifont}
\usepackage{color}
\usepackage{todonotes}
\usepackage{booktabs}
\usepackage{amssymb}% http://ctan.org/pkg/amssymb
\usepackage{pifont}% http://ctan.org/pkg/pifont
\usepackage[numbers,sort&compress]{natbib}
\usepackage{siunitx}
\usepackage{lscape}

\newtheorem{thm}{Theorem}[section]

\numberwithin{equation}{section}

\newcommand{\R}{\mathbb{R}}

\newcommand{\El}{\mbox{${\rm evals}_{\mbox{{\rm \tiny LBFGS}}}$}}
\newcommand{\En}{\mbox{${\rm evals}_{\mbox{{\rm \tiny NEWUOA}}}$}}

\begin{document}
% generates the title
\title{On the Numerical Performance of Derivative-Free Optimization Methods Based on Finite-Difference Approximations
%\footnote{Technical Report: NU-07242018}
}

\author{   
       Hao-Jun Michael Shi\thanks{Department of Industrial Engineering and Management Sciences, Northwestern University. Shi was supported by the Office of Naval Research grant N00014-14-1-0313 P00003. Xuan was supported by the National Science Foundation grant DMS-1620022. Nocedal was supported by the Office of Naval Research grant N00014-14-1-0313 P00003, and by National Science Foundation grant DMS-1620022. \newline\url{hjmshi@u.northwestern.edu}, \url{qxuan@u.northwestern.edu}, \url{j-nocedal@northwestern.edu}}
       	   \and
       Melody Qiming Xuan\footnotemark[1]
       \and
       Figen Oztoprak\thanks{Artelys Corporation. \url{figen.topkaya@artelys.com}}
       \and
       Jorge Nocedal\footnotemark[1]}
%****

\maketitle

%%%%%%%%%%%%%
%%% ABSTRACT %%% 
\begin{abstract}{
The goal of this paper is to investigate an approach for derivative-free optimization that has not received sufficient attention in the literature and is yet one of the simplest to implement and parallelize. It consists of computing gradients of a smoothed approximation of the objective function (and constraints), and employing them within established codes. These gradient approximations are calculated by finite differences, with a differencing interval determined by the noise level in the functions and a bound on the second or third derivatives. It is assumed that noise level is known or can be estimated by means of difference tables or sampling.  The use of finite differences has been largely dismissed in the derivative-free optimization literature as too expensive in terms of function evaluations and/or as impractical when the objective function contains noise. The test results presented in this paper suggest that such views should be re-examined and that the finite-difference approach has much to be recommended. The tests compared  {\sc newuoa, dfo-ls} and {\sc cobyla} against the finite-difference approach on three classes of problems: 
general unconstrained problems, nonlinear least squares, and general nonlinear programs with equality constraints. }

\bigskip\noindent
\textbf{Keywords:} derivative-free optimization, noisy optimization, zeroth-order optimization, nonlinear optimization
\end{abstract}

\bigskip
% \jn{Jorge: blue}, \fo{~~~Figen, magenta}, \ms{~~~Michael: red}, \mx{~~~Melody:green}
\newpage
\tableofcontents

\newpage
%%%%%%%%%%%%%%%%
%%% INTRODUCTION %%%%

\section{Introduction}
\label{sec:intro}
\setcounter{equation}{0}
                    
The problem of minimizing a nonlinear objective function  when gradient information is not available has received much attention in the last three decades; see  \cite{conn2009introduction,larson2019derivative} and the references therein.  A variety of methods  have been developed for unconstrained optimization, and some of these methods have been extended to deal with constraints. The important benchmarking paper by Mor\'e and Wild \cite{more2009benchmarking} showed that traditional methods, such as the Nelder-Mead simplex method \cite{NeldMead65} and a leading pattern-search method \cite{gray2006algorithm},
are not competitive with the model-based trust-region approach pioneered by Powell \cite{powell2006newuoa} and developed concurrently by several other authors \cite{conn2009introduction}. The advantages of Powell's approach reported in \cite{more2009benchmarking} were observed for both smooth and nonsmooth problems, as well as noisy and noiseless objective functions. Rios and Sahinidis \cite{rios2013derivative} confirmed the findings of Mor\'e and Wild concerning the inefficiency of traditional methods based on Nelder-Mead or direct searches; examples of direct-search methods are given in
\cite{audet2006mesh,abramson2011nomad,lewis2000direct,Wrig96}. 
Based on these studies, we regard the model-based approach of Powell as a leading method for derivative-free optimization (DFO).

There is, however, an alternative approach for DFO that is perhaps the simplest, but has been largely neglected in the nonlinear optimization literature. It consists of estimating derivatives using finite differences, and using them within standard nonlinear optimization algorithms. Many of the papers in the DFO literature dismiss this approach at the outset as being too expensive in terms of functions evaluations and/or as impractical in the presence of noise. As a result of this prevalent view, the vast majority of papers on DFO methods do not present comparisons with a finite-difference approach. We believe that if such comparisons had been made, particularly in the noiseless setting, research in the field would have followed a different path.  Ironically, as pointed out by Berahas et al. \cite{berahas2019derivative}, the finite-difference  DFO approach is widely used by practitioners for noiseless problems,  often unwittingly,  as many established  optimization codes invoke a finite-difference option when derivatives are not be provided.  A disconnect between practice and research therefore occurred in the last two decades, and no systematic effort was made to bridge this gap. This paper builds upon Berahas et al. \cite{berahas2019derivative} and Nesterov and Spokoiny \cite{nesterov2017random}, and  aims to bring the finite-difference DFO approach to the forefront of research by illustrating its performance on a variety of unconstrained and constrained problems, with and without noise. Our numerical experiments highlight its strengths and weaknesses, and identify some open research questions.

As is well known, the use of finite differences is delicate in the presence of noise in the function, and the straightforward application of rules designed to deal with roundoff errors, such as selecting the finite-difference interval on the order of the square-root of machine precision, leads to inefficiencies or outright failure. However, for various types of noise, the estimation of derivatives can be placed on a solid theoretical footing. One can view the task as the computation of derivatives of a smoothed function \cite{nesterov2017random}, or as the computation of an estimator that minimizes a mean squared error \cite{more2012estimating}.  The early work by Gill et al. \cite{GillMurrSaunWrig83} discusses an adaptive approach for computing the difference interval in the presence of errors in the objective function, but there appears to have been no follow-up on the application of these techniques for derivative-free optimization. The paper that  influenced our work the most is by  Mor\'e and Wild \cite{more2012estimating}, who describe how to choose the differencing interval as a function of the noise level and a bound on the second (or third) derivative, 
so as to obtain nearly optimal gradient estimates.  The noise level, defined as the standard deviation of the noise, can be estimated by sampling (in the case of stochastic noise) or using a table of differences \cite{more2011estimating} (in the case of computational noise).

To test whether the finite-difference approach to DFO is competitive with established methods, for noisy and noiseless functions, we consider three classes of optimization problems: general unconstrained problems, nonlinear least-squares problems, and inequality-constrained nonlinear problems.  For benchmarking, we selected the following well-established DFO methods:  {\sc newuoa} \cite{powell2006newuoa} for general unconstrained problems,  {\sc dfo-ls} \cite{cartis2019improving} for nonlinear least squares, and {\sc cobyla} \cite{powell1994direct} for inequality constrained optimization. The finite-difference approach can be implemented in many standard codes for smooth deterministic optimization. We choose  {\sc l-bfgs} \cite{mybook} for general unconstrained problems, {\sc lmder} \cite{MoreGarbHill80} for nonlinear least squares, and {\sc knitro} \cite{ByrdNoceWalt06} for inequality-constrained optimization problems.
We do not present comparisons for general problems involving both equality and inequality constraints because we were not able to find an established DFO code of such generality that was sufficiently robust in our experiments  ({\sc cobyla} accepts only inequality constraints.)

Finite-difference-based methods for DFO enjoy two appealing features that motivate further research. They can easily exploit parallelism in the evaluation of functions during finite differencing, which could be critical in certain applications. In addition, finite-difference approximations can be incorporated into existing software for constrained and unconstrained optimization, sometimes rather easily. This obviates the need to redesign existing unconstrained DFO to handle more general problems, an effort that has 
taken two decades for interpolation-based trust-region methods \cite{conn1997convergence,conn1998derivative,ConnScheiVic07g,ConSchVi06,conn2009introduction,powell2009bobyqa,lincoa,powell2015fast}, and is yet incomplete.  The strategy often suggested of simply applying an unconstrained DFO method to a penalty or augmented Lagrangian reformulation will not yield an effective general-purpose method, as is known from modern research in nonlinear programming.

\subsection{Literature Review}
The book by Conn, Scheinberg and Vicente \cite{conn2009introduction} gives a thorough treatment of the interpolation-based trust-region approach for DFO. It presents foundational theoretical results as well as  detailed algorithmic descriptions. 
Larson, Menickelly and Wild \cite{larson2019derivative} present a comprehensive review of DFO methods as of 2018. Their survey covers deterministic and randomized methods, noisy and noiseless objective functions, problem structures such as least squares and empirical risk minimization, and various types of constraints. But they devote only two short paragraphs to the potential of finite differencing as a practical DFO method, and focus instead on complexity results for these methods. 
Audet and Hare \cite{audet2017derivative}  review direct-search and pattern-search methods, and describe a variety of practical applications solved with {\sc mads}, and {\sc nomad}. Rios and Sahinidis \cite{rios2013derivative} report extensive numerical experiments comparing many DFO codes. Neumaier \cite{neumaier2004complete} reviews methods endowed with convergence guarantees, with emphasis on global optimization.
Nesterov and Spokoiny \cite{nesterov2017random} propose a Gaussian smoothing approach for noisy DFO that involves finite-difference approximations, and establish complexity bounds for smooth,  nonsmooth, and stochastic objectives.  Kelley et al. \cite{kelley2011implicit,choi2000superlinear} propose a finite-difference BFGS method, called implicit filtering, designed for the case when noise can be diminished at any iteration, as needed. In that approach, the finite-difference interval decreases monotonically.
Kimiaei \cite{kimiaei2020line} propose a randomized method, called VSBBON that implements a noisy line search and employs quadratic models in   subspaces determined adaptively.

Optimization papers that study the choice of the finite-difference interval $h$ in the presence of errors (or noise) in the function include Gill et al. \cite{GillMurrSaunWrig83}, which is inspired by earlier work by Lyness \cite{lyness1977has}. Gill et al. \cite{GillMurrSaunWrig83} gives examples where their approach fails, and pays careful attention to the estimation of bounds on the second (or third) derivatives, since these bounds are needed to obtain an accurate estimate of $h$. Mor\'e and Wild \cite{more2012estimating} derive formulae for the optimal choice of $h$, and propose {\tt ECnoise}, a practical procedure for estimating stochastic or deterministic noise \cite{more2011estimating}. They also give some  attention to the estimation of the second derivative.

\subsection{Contributions of this Paper}
To our knowledge, this is the first systematic investigation into the empirical performance of finite-difference-based DFO methods relative to established techniques, across a range of problems, with and without noise in the functions. 
The main contributions of this paper can be summarized as follows. We use the acronym ``FD-DFO method'' for a method that employs some form of finite differences to approximate the gradient of the objective function and (possibly) the constraints.
\begin{itemize}
\item For {\em noiseless functions}, we found that the FD-DFO methods are at least as efficient if not superior to established methods, across all three categories of problems, without the use of sophisticated procedures for determining the finite-difference interval. Surprisingly, the carefully crafted {\sc newuoa} code is not more efficient, as measured by the number of function evaluations, than a finite-difference {\sc l-bfgs} method, and has  higher linear algebra cost. 

\item For {\em noisy functions}, we observed that {\sc newuoa} is more efficient and accurate than the finite-difference {\sc l-bfgs} method for unconstrained optimization, but not by a wide margin.  For least-squares problems, the performance of the recently developed {\sc dfo-ls}  code is comparable to that of the finite-difference version of the classical {\sc lmder} code. For inequality-constrained problems, a simple finite-difference version of {\sc knitro} has comparable performance with {\sc cobyla}.
\item The differencing formulas suggested by Mor\'e and Wild perform well in our tests, unless the bound on the second (or third) derivative required in these formulas is poor. We found that existing procedures for estimating these bounds are not robust, and improvements on this seemingly small algorithmic detail may allow FD-DFO methods to close the performance gap observed in the unconstrained setting.
 \end{itemize}
One striking observation from our study is that interpolation-based trust-region methods are more robust in the presence of noise than we expected. Although they are not simple algorithms (e.g., they require a geometry phase), they have some  appealing features. For example,  {\sc newuoa} and  {\sc dfo-ls}  do not require knowledge of the noise level in the objective function or estimates of derivatives, and yet performed reliably for most levels of noise, suggesting that the internal logic of the algorithm normally reacts correctly to the noise inherent in the problem (although rare failures were observed). We are not aware of studies that comment on this robustness, except possibly for \cite{more2009benchmarking}.

\subsection{Limitations of this Work}

For each problem class, we employed  only one established DFO code to benchmark the efficiency of the corresponding FD-DFO method. We found it essential to work with a small number of codes  to allow us to understand them well enough and ensure fairness of the tests. As more research is devoted to the development of FD-DFO methods, comparisons with other codes will be essential. As in most benchmarking studies, the standard disclaimer is in order: codes were tested with default options and overall performance may vary with other settings. As our focus was on scalable local optimization methods, we did not consider global optimization methods, such as Bayesian optimization, surrogate optimization, or evolutionary methods \cite{eriksson2019pysot,frazier2018tutorial,hansen2016cma}. Most Bayesian and surrogate optimization methods are not scalable to problems above 20 variables unless modified as in TuRBO; see \cite{frazier2018tutorial,eriksson2019scalable}.  Other approaches to global optimization include restarts, grid searches or surrogate models, enhancements that were not considered here; see \cite{st:gamsbaron:7.2.5} and \cite{neumaier2004complete}.

Nonsmooth problems were not considered in this study.  Based on the results by Lewis and Overton \cite{lewis2013nonsmooth} and Curtis et al. \cite{curtis2013adaptive}, a finite-difference implementation of BFGS (not L-BFGS) could prove a strong competitor to current DFO methods. However,  how to perform finite differencing robustly in the nonsmooth setting is still an open research question. 

Perhaps the most important limitation of this study is that it considers only one model of noise: additive uniformly-distributed bounded noise. We do not know how the methods tested here behave for anisotropic noise, unbounded noise, or noise based on other distributions. We focused on just one model of noise because this already raised some important algorithmic questions that need to be resolved to improve the peformance of FD-DFO methods.

\subsection{Organization of the Paper}
A  large number of experiments were performed in this study. We provide some of these numerical results in the Appendices \ref{app:lipschitz}--\ref{sec:ext_num}. In the main body of the paper, we display graphs or tables that attempt to accurately summarize the conclusions of the experiments.
The paper is organized into five sections.  Unconstrained problems are studied in section~\ref{ch:unconstrained}; nonlinear least-square problems in section~\ref{ch:squares}, and nonlinear optimization problems with general inequality constraints in section~\ref{ch:constrained}. Concluding remarks and some open questions are described in section~\ref{ch:final}.

\section{Unconstrained Optimization}
\label{ch:unconstrained}

In this section, we consider the solution of unconstrained optimization problems of the form
\begin{equation}  \label{uncprob}
	\min_{x \in \R^n} \phi(x),
\end{equation}
given only noisy evaluations, 
\begin{equation}  \label{fdef}
f(x) = \phi(x) + \epsilon(x).
\end{equation}
Here, $\epsilon(x)$ is a scalar that models  deterministic or stochastic noise and $\phi$ is assumed to be a smooth function. We compare the performance of {\sc newuoa} and {\sc l-bfgs} with finite-difference gradients on 73 \texttt{CUTEst} problems \cite{gould2015cutest}, varying the dimension of each problem up to $n \leq 300$ whenever possible. All experiments were run in double precision.

We chose {\sc newuoa} because, as mentioned above, it is regarded as one of the leading codes for unconstrained derivative-free optimization. Another well-known  model-based trust-region method is {\sc dfotr} \cite{bandeira2012computation}, which in our experience is often competitive with {\sc newuoa} in terms of function evaluations, but has much higher per-iteration cost.

In summary, the methods tested in our experiments are as follows.
\begin{itemize}
\item \texttt{NEWUOA}: A model-based trust-region derivative-free algorithm that forms quadratic models using interpolation of function values, combined with a minimum Frobenius-norm update of the Hessian approximation; see Powell \cite{powell2006newuoa}. We called the code in Python 3.7 through {\sc pdfo} developed by Ragonneau and Zhang \cite{PDFO}. We used the default settings in {\sc newuoa}, which in particular, set the number of interpolation points to $2n+1$.
\item \texttt{L-BFGS}: A finite-difference implementation of the limited-memory BFGS algorithm \cite{mybook} with a bisection Armijo-Wolfe line search described below. We use a memory of size $t = 10$. Forward- and central-difference options are tested. At intermediate trial points generated during the line search, the directional derivative is computed via forward- or central-differences along the direction of interest.
% \item \texttt{L-BFGS-E}: A finite-difference noise-tolerant limited-memory BFGS algorithm with a modified Armijo-Wolfe line search with lengthening; see Shi, et al. \cite{???}. We use $t = 10$.
\end{itemize}

We first consider the case when noise is not present and later study the effect of noise on the performance of the algorithms.

%%%%%%%%%%%%%%%%%%%%%%%%%%%%
\subsection{Experiments on Noiseless Functions}
In the first set of experiments, we let $\epsilon(x) \equiv 0$, so that the only errors in the function evaluations are due to machine precision, which we denote by  $\epsilon_{M}$. The approximate gradient $g(x) \in \R^n$ of the objective function $\phi(x)$, computed by finite differencing, is given as:
\begin{align}
    [g(x)]_i & = \frac{f(x + h_i e_i) - f(x)}{h_i}, & & h_i = \max\{1, |[x]_i|\} \sqrt{\epsilon_{M}}; & & \mbox{(forward differencing)} \label{fdg}\\
    [g(x)]_i & = \frac{f(x + h_i e_i) - f(x - h_i e_i)}{2h_i}, & & h_i = \max\{1, |[x]_i|\} \sqrt[3]{\epsilon_M}. & & \mbox{(central differencing)} \label{cdg}
\end{align}  
Approximating the gradient through  \eqref{fdg}, \eqref{cdg} therefore requires $n + 1$ and $2n$ function evaluations, respectively. 
The $\max\{1, |[x]_i|\}$ term is incorporated to handle the rounding error in $[x]_i$.

The directional derivative $D_p \phi(x) = \nabla \phi(x)^T p$ of the objective function $\phi$ along a direction $p \in \R^n$ is required within the Armijo-Wolfe line search employed by {\sc l-bfgs}. It is approximated as: 
\begin{align}   
    g(x; p) & = \frac{f(x + h p_u) - f(x)}{h} \|p\|, & & h = \sqrt{\epsilon_M}; \label{dnal1} & & \mbox{(forward differencing)} \\
    g(x; p) & = \frac{f(x + h p_u) - f(x - h p_u)}{2h} \|p\|, & & h = \sqrt[3]{\epsilon_M}; & & \mbox{(central differencing)}\label{dnal2}
\end{align}
where $p_u = p / \|p\|$ is the normalized direction. 
% \jn{[Any reason why the scale of $x$ was not taken into account for the directional derivative approximation?]} \ms{[It was not clear to me how that should be done.]}

These choices of $h$ can be improved by including contributions of the second and third derivatives, respectively, as discussed in the next subsection. However, we find that in the noiseless setting, and for our test functions, such a refinement is not needed to make the finite-difference {\sc l-bfgs} approach competitive. 
%since the gradient approximation is accurate enough 

The algorithms are terminated when either
\begin{equation}\label{eq:term crit}
    \phi(x_k) - \phi^* \leq \tau \cdot \max\{1, |\phi^*|\},
\end{equation}
with $\tau= 10^{-6}$, or when the limit of $500 \times n$ function evaluations is reached. The optimal value $\phi^*$ is determined by running BFGS  with {\em exact} gradients (provided by {\tt CUTEst} \cite{gould2015cutest}) until no more progress can be made on the function.  Complete numerical results are given in Tables \ref{tab:no-noise1}-\ref{tab:no-noise5} in Appendix \ref{app:smooth, det}.

In Figure \ref{fig:feval no noise}, we summarize the results using {\em log-ratio profiles} proposed by Morales \cite{Mora02}, which in this case report the quantity
\begin{equation}   \label{jl}
   \log_2 \left(\frac{\El}{\En}\right) ,
\end{equation}
where $\El$ and $\En$ denote the total number of function evaluations for {\sc newuoa} and {\sc l-bfgs} to satisfy \eqref{eq:term crit} or reach the maximum number of function evaluations. In the figures, the ratios \eqref{jl} are plotted in increasing order. Thus, the area of the shaded region gives an idea of the general success of a method. 

For fair comparison, the runs with the following outcomes were not included in Figure~\ref{fig:feval no noise}.
There are six instances where {\sc newuoa}'s 
%significantly outperforms {\sc l-bfgs} stemming from {\sc newuoa}'s 
initial sampling of $2 n + 1$ points immediately discovers a point whose function value satisfies \eqref{eq:term crit}, thirteen instances where {\sc l-bfgs} and {\sc newuoa} converged to different minimizers, and two instances where either one of the methods failed. 

Overall, we observe in these tests that forward-difference {\sc l-bfgs} outperforms {\sc newuoa} on the majority of problems in terms of function evaluations. This is perhaps surprising because {\sc newuoa} was designed to be parsimonious in terms of function evaluations, whereas finite-difference {\sc l-bfgs} requires $n$ function evaluations per iteration. It is also notable that this is achieved without any additional information (for example, the Lipschitz constant of the gradient for forward differencing) to squeeze out the best possible accuracy of the gradient in the finite-difference {\sc l-bfgs} approach. As expected, central-difference {\sc l-bfgs} requires 
%approximately double the number of 
significantly more function evaluations than the forward-difference option, and does not provide significant benefit in terms of solution accuracy for the majority of the problems. In particular, when run in double precision, forward-difference {\sc l-bfgs} is able to converge to the same tolerance that one would expect with analytical gradients.

\begin{figure}\label{fig:feval no noise}
    \centering
    \includegraphics[width=0.32\linewidth]{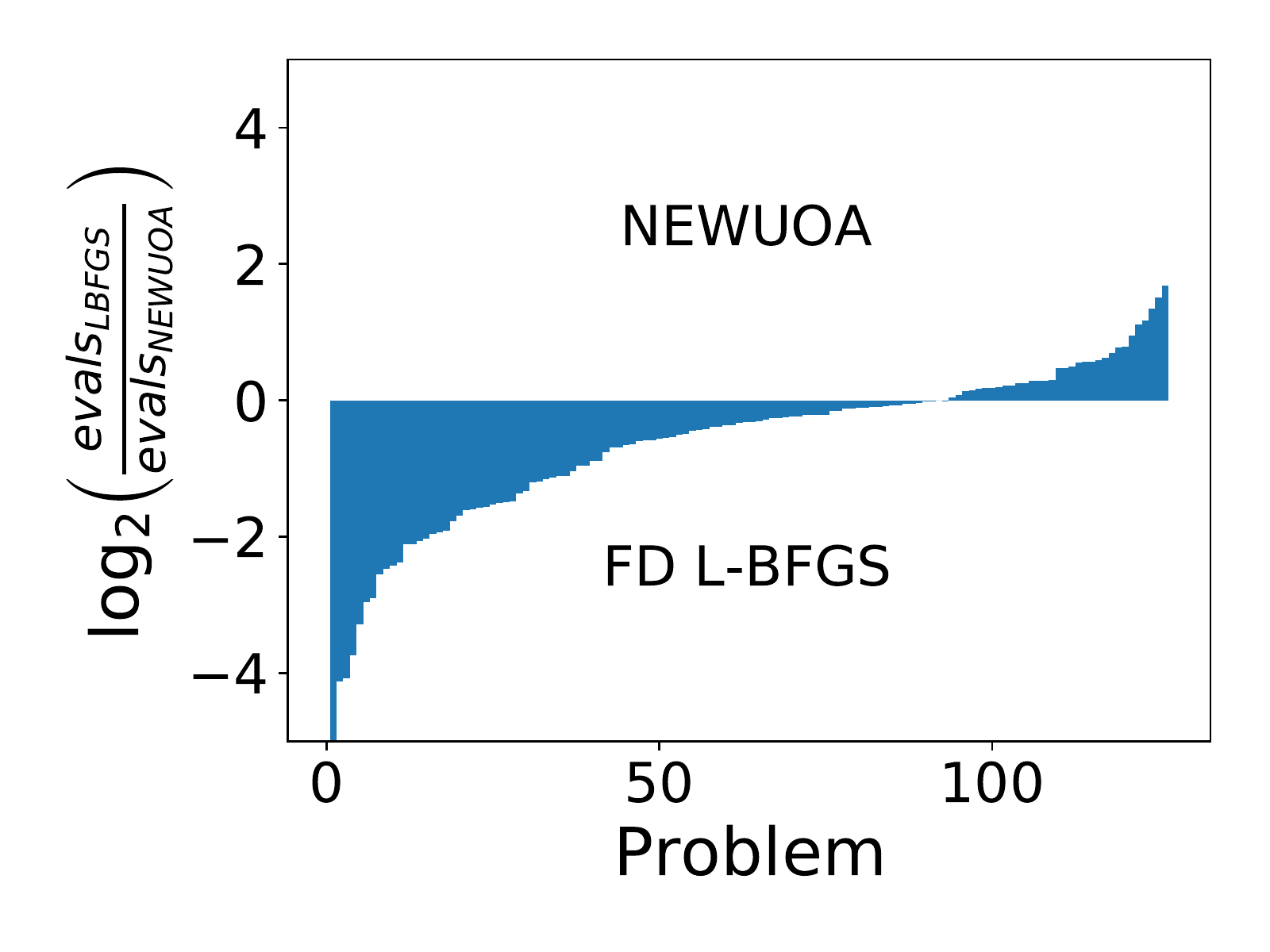}
    \includegraphics[width=0.32\linewidth]{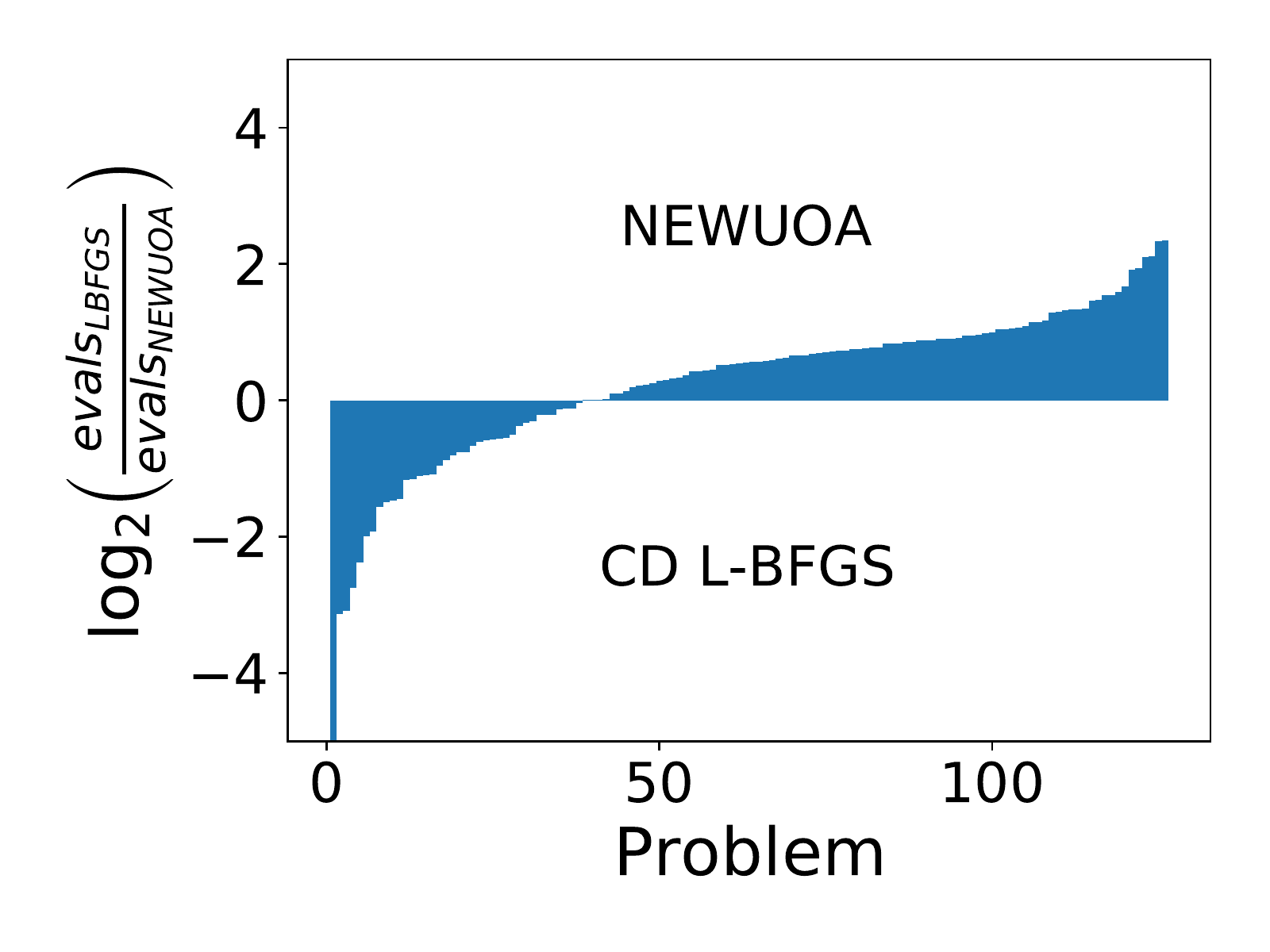}
    \caption{{\em Efficiency, Noiseless Case}. Log-ratio profiles for the total number of function evaluations to achieve \eqref{eq:term crit} with $\epsilon(x) = 0$. The left figure compares forward-difference {\sc l-bfgs} with {\sc newuoa}, and the right figure compares central-difference {\sc l-bfgs} with {\sc newuoa}.}
\end{figure}

In order to illustrate  iteration cost, we report in Table~\ref{tab:cpu time}  the CPU times for {\sc newuoa} and {\sc l-bfgs} for a representative problem as the number of variables $n$ increases. {\sc newuoa}'s execution time grows much faster than {\sc l-bfgs} because one iteration of {\sc newuoa} requires $O(n^2)$ flops, whereas the iteration cost of {\sc l-bfgs} is $4tn$ flops, where $t=10$ and all test functions are inexpensive to evaluate.
Across all the problems, we observe that when $n \approx 100$, {\sc newuoa} can  take at least $5-10$ times longer than {\sc l-bfgs} in terms of wall-clock time.

\begin{table}[htp]
    \centering
    \footnotesize
    \begin{tabular}{ c | c | c | c | c | c | c | c }
        \toprule
        $n$ & 10 & 20 & 30 & 50 & 90 & 100 & 500 \\
        \midrule
        NEWUOA & 9.61e-2 & 9.99e-2 & 1.21e-1 & 1.99e-1 & 5.83e-1 & 8.92e-1 & 2.42e2 \\
        FD L-BFGS & 1.13e-2 & 1.68e-2 & 2.26e-2 & 3.18e-2 & 5.50e-2 & 6.67e-2 & 4.67e-1 \\
        CD L-BFGS & 1.30e-2 & 2.09e-2 & 2.90e-2 & 4.28e-2 & 6.63e-2 & 8.56e-2 & 6.63e-1 \\
        \bottomrule
    \end{tabular}
    \caption{{\em Computing Time, Noiseless Case.} CPU time (in seconds) required for  {\sc newuoa} and {\sc l-bfgs} with forward (FD) and central (CD) differencing to satisfy \eqref{eq:term crit} on the \texttt{NONDIA} problem, as the dimension of the problem $n$ increases.}
    \label{tab:cpu time}
\end{table}

While {\sc newuoa} is an inherently sequential algorithm, finite-difference {\sc l-bfgs} offers ample opportunities for parallelism when computing the finite-difference approximation to the gradient, which can be distributed across multiple nodes, with the only sequential bottleneck arising in the function evaluations
within Armijo-Wolfe line search.
%evaluating the function within the Armijo-Wolfe line search. 
Given the ratio between the number of function evaluations needed for {\sc newuoa} versus finite-difference {\sc l-bfgs} reported in Tables~\ref{tab:no-noise1}-\ref{tab:no-noise5} in Appendix~\ref{app:smooth, det}, we hypothesize that even small amounts of parallelism across even a few processors may lead to significant speedup of {\sc l-bfgs}. This is an often overlooked benefit of finite-difference-based methods in the DFO literature, as we are not aware of implementations of model-based trust-region methods that benefit significantly from parallelism.

We conclude this subsection by noting that most DFO methods were developed and tested primarily in the noiseless setting \cite{conn2009introduction,larson2019derivative}.
 Our results suggest that most of these methods may not be competitive with the finite-difference {\sc l-bfgs} approach in the noiseless case. However, since no single set of experiments can conclusively establish such a claim, tests by other researchers should verify our observations. 
 %When noise is present, \jn{the relative performance of the two methods is different,} as we discuss next.
 
\subsection{Experiments on Noisy Functions}
\label{sub:noisy unc}

In this set of experiments, we synthetically inject uniform stochastic noise into the objective function. In particular, we  sample $\epsilon(x) \sim \sigma_f U(-\sqrt{3}, \sqrt{3})$ i.i.d. independent of $x$, where $\sigma_f \in \{10^{-1}, 10^{-3}, 10^{-5}, 10^{-7}\}$. (Thus, strictly speaking we should  write $\epsilon$, but we keep the notation $\epsilon(x)$ for future generality.)
By construction of $\epsilon(x)$, we have that $\sigma_f^2 = \mathbb{E}[\epsilon(x)^2]$. We refer to standard deviation of the noise $\sigma_f$ as the \textit{noise level}.  

We employ a more precise formula for the  finite-difference interval than in the noiseless setting --- one that depends both on the noise level $\sigma_f$ and on the curvature of the function. Specifically, we compute a different $h_i$ for each coordinate direction based on the following well-known result \cite{more2012estimating} that provides bounds on the mean-squared error of estimated derivatives.

\begin{thm}\label{thm:MSE}
Let $h_0 > 0$ and $x, p \in \R^n$, with $\|p\| = 1$, be given. Define the  interval $I_{\text{FD}} = \{x + tp : t \in [0, h_0]\}$. If $|D_p^2 \phi(y)| = |p^T \nabla^2 \phi(y) p| \leq L$ for all $y \in I_{\text{FD}}$, then for any $\tilde{h} \in (0, h_0]$, the following bound holds for forward differencing:
\begin{equation}  \label{mob1}
\sigma_{g, p}(\tilde{h})^2 \equiv \mathbb{E}\left[\left(g(x; p, \tilde{h}) - D_p \phi(x)\right)^2 \right] \leq \frac{L^2 \tilde{h}^2}{4} + \frac{2 \sigma_f^2}{\tilde{h}^2} ,
\end{equation}
where 
\begin{equation}
g(x; p, \tilde{h}) = \frac{f(x + \tilde{h} p) - f(x)}{\tilde{h}}.
\end{equation}
Similarly for central differencing, if the third derivative satisfies $|D_p^3 \phi(y)| \leq M$ for all $y \in I_{\text{CD}} = \{x \pm t p : t \in [0, h_0]\}$, then for any $\tilde h \in (0, h_0]$ we have 
\begin{equation}  \label{mob2}
\sigma_{g, p}(\tilde{h})^2 \equiv \mathbb{E}[(g(x; p, \tilde{h}) - D_p \phi(x))^2] \leq \frac{M^2 \tilde{h}^4}{36} + \frac{\sigma_f^2}{2 \tilde{h}^2},
\end{equation}
where $D_p^3 \phi(y)$ denotes the third directional derivative with respect to direction $p$ and
\begin{equation}
g(x; p, \tilde{h}) = \frac{f(x + \tilde{h} p) - f(x - \tilde{h} p)}{\tilde{h}}.
\end{equation}
\end{thm}
By minimizing the upper bounds in \eqref{mob1}, \eqref{mob2} with respect to $\tilde h$, one obtains the following expressions for each coordinate direction \cite{more2012estimating},
%Applying Theorem \ref{thm:MSE}, we can compute the approximate gradient $g(x) \in \R^n$ of the objective function $\phi(x)$ as suggested in Mor\'e and Wild \cite{more2012estimating}:
\begin{align}
    [g(x)]_i & = \frac{f(x + h_i e_i) - f(x)}{h_i}, & & h_i = \sqrt[4]{8} \sqrt{\frac{\sigma_f}{L_i}}; & & \mbox{(forward differencing)} \label{eq:fd formula}\\
    [g(x)]_i & = \frac{f(x + h_i e_i) - f(x - h_i e_i)}{2h_i}, & & h_i = \sqrt[3]{\frac{3 \sigma_f}{M_i}}, & & \mbox{(central differencing)} \label{eq:cd formula}
\end{align}
where $L_i$ and $M_i$ are bounds on the second and third derivative along the $i$-th coordinate direction $e_i$. The noise level of the derivative along an arbitrary direction $p \in \R^n$ can be shown to be
% \begin{align*}
%     \epsilon_g & = 2 \sqrt{L_2 \epsilon_f} \\
%     \epsilon_g & = \frac{3^{2/3}}{2} M^{1/3} \epsilon_f^{2/3}.
% \end{align*}
\begin{equation*}
    \sigma_{g, p} \leq \sqrt[4]{2} \sqrt{L \sigma_f} \qquad
    \sigma_{g, p} \leq \frac{\sqrt[6]{3}}{2} M^{1/3} \sigma_f^{2/3},
\end{equation*}
for forward and central differencing, respectively.
For the full gradient, the mean-squared error is given by
\begin{equation}
    \sigma_g^2 \equiv \mathbb{E}[\|g(x) - \nabla \phi(x)\|^2] = \sum_{i = 1}^n \sigma_{g, e_i}^2.
\end{equation}

With the formulas \eqref{eq:fd formula}, \eqref{eq:cd formula} in hand, we can now return to the practical implementation of the finite-difference {\sc l-bfgs} method. We estimate the constants $L_i$ in  \eqref{eq:fd formula} using a second-order difference. Given a direction $p \in \R^n$ where $\|p\| = 1$, we define
\[
   \Delta(t) = f(x + t p) - 2f(x) + f(x - t p),
\]
where $t$ is the second-order differencing interval. The Lipschitz constant $L$ along the direction $p$ can thus be approximated as $L \approx \Delta(t)/t^2$. However, as with the choice of $h$, we need to be careful in the selection of $t$, and for this purpose we employ an iterative technique proposed by Mor\'e and Wild \cite{more2012estimating}, whose goal is to find an interval $t$ that satisfies
\begin{align}
|\Delta(t)| & \geq \tau_1 \epsilon_f, & & \tau_1 \gg 1 \label{mw1} \\
|f(x \pm t p) - f(x)| & \leq \tau_2 \max\{|f(x)|, |f(x \pm t p)|\}, & & \tau_2 \in (0, 1) .\label{mw2}
\end{align} 
These  heuristic conditions aim to ensure that $t$ is neither too small nor too large; see \cite{more2012estimating} for a full description of the  Mor\'e-Wild (MW) technique. This technique cannot, however, be guaranteed to find a $t$ that satisfies \eqref{mw1}, \eqref{mw2}. To account for this, we developed the following procedure for estimating the constants $L_i$ for forward differencing.

\medskip
\noindent {\bf Procedure I.} {\em Adaptive Estimation of $L$}

1. At the first iteration of the {\sc l-bfgs} method, invoke the MW procedure to compute $t_i$, for $i = 1, ..., n$. If such a $t_i$ can be found to satisfy \eqref{mw1}, \eqref{mw2},  set $L_i = \max\{10^{-1}, |\Delta(t_i)| / t_i^2\}$ in \eqref{eq:fd formula}. Otherwise,  set it to $L_i = 10^{-1}$. Store the $L_i$ in a vector $\mathbf{L}$. 
%We use $L[i]$ whenever the $i$th component of the gradient is approximated and use $\|L\| / \sqrt{n}$ 
To calculate the directional derivative in the Armijo-Wolfe line search,  set the second-order differencing interval to $t=\|{\bf L} \| / \sqrt{n}$.

2.  If at any iteration of the {\sc l-bfgs} algorithm the line search returns a steplength $\alpha_k < 0.5$, then re-estimate the vector $\mathbf{L}$ by calling the MW procedure as in step 1 at the current iterate $x_k$.

\bigskip
When the gradient is sufficiently accurate, {\sc l-bfgs} generates well-scaled directions such that $\alpha_k=1$ is acceptable. Thus, the occurrence of $\alpha_k < 0.5$ is viewed as an indication that the curvature of the problem may have changed, and should be re-estimated.
The function evaluations performed in the estimation of $\mathbf{L}$ will be  accounted for in the numerical results presented below.

For central differencing, a practical procedure for estimating $M$ is more involved, and will be discussed in a forthcoming paper \cite{shi2021automatic}.  Here, we simply approximate $M$ using changes in the second derivative of the true objective:
\begin{equation} \label{gtn}
    M = \max\left\{10^{-1}, \frac{\left|p^T \left(\nabla^2 \phi\left(x + \hat{h} \frac{p}{\|p\|}\right) - \nabla^2 \phi(x) \right) p \right|}{\hat{h} \|p\|^2} \right\}
\end{equation}
where $\hat{h} = \sqrt{\epsilon_M}$. This is an idealized strategy, and its cost is not accounted for in the results below  since central differencing does not play a central role in this study.

%since there is no known method for estimating $L_3$ to our knowledge, we will use a theoretical estimate of $L_3$ without including the cost of evaluating $L_3$, as discussed further in the Appendix. For further discussion on the choice of $L_2$ and $L_3$ and known heuristics for estimating such quantities, please refer to Appendix \ref{app:lipschitz}.

% Three different variants of finite-difference L-BFGS are tested here \ms{(although one will be chosen as representative later)}:
% \begin{enumerate}
%     \item \texttt{L-BFGS}: The standard finite-difference limited-memory BFGS algorithm with a bisection Armijo-Wolfe line search, with memory $t = 10$. 
%     \item \texttt{L-BFGS-E}: A finite-difference noise-tolerant BFGS algorithm with a two-phase Armijo-Wolfe line search with lengthening, with memory $t = 10$.
%     \item \texttt{L-BFGS (Skips)}: A finite-difference BFGS algorithm with a two-phase Armijo-Wolfe line search with skipping, with memory $t = 10$. 
% \end{enumerate}
% Both L-BFGS-E and L-BFGS (Skips) rely on the noise control condition defined as
% \begin{equation}\label{eq:noise}
%     (g(x_k + \alpha_k p_k) - g(x_k))^T p_k \geq 2 (1 + c_3) \epsilon_g \|p_k\|.
% \end{equation}
% Here, we use $c_3 = 0.5$ for L-BFGS-E and $c_3 = 0$ for L-BFGS (Skips). L-BFGS-E will lengthen the gradient difference used in the L-BFGS update until \eqref{eq:noise}, while L-BFGS with skipping will skip the update if \eqref{eq:noise} does not hold.

In addition to the  choice of $h$ mentioned above, we have found it important to  modify the line search slightly to better handle noise. Following Shi et al. \cite{shi2020noise}, we relax the Armijo condition within the Armijo-Wolfe line search, as follows. Let the superscript $j$ denote the $j$th iteration of the line search performed from the iterate $x_k$. The modified Armijo condition is
\begin{equation}
    f(x_k + \alpha_k^j p_k) \leq
    \begin{cases}
    f(x_k) + c_1 \alpha_k^j g(x_k)^T p_k & \mbox{ if } j = 0, \, g(x_k)^T p_k < - \sigma_g \|p_k\| \\
    f(x_k) + c_1 \alpha_k^j g(x_k)^T p_k + 2 \sigma_f & \mbox{ if } j \geq 1,  \, g(x_k)^T p_k < - \sigma_g \|p_k\| \\
    f(x_k) & \mbox{ if } g(x_k)^T p_k \geq - \sigma_g \|p_k\|.
    \end{cases}
\end{equation}
Other than this modification to the line search, no other changes are made to the {\sc l-bfgs} method. We will refer to the resulting method simply as {\sc l-bfgs} to avoid introducing more acronyms. 

We perform tests to compare {\sc l-bfgs} with forward and central differencing
%, using this modified Armijo-Wolfe line search, 
against {\sc newuoa}, both in terms of the achievable solution accuracy and in terms of efficiency (as measured by function evaluations).

\subsubsection{Accuracy}

\begin{figure}[ht]
\centering
\includegraphics[width=0.32\textwidth]{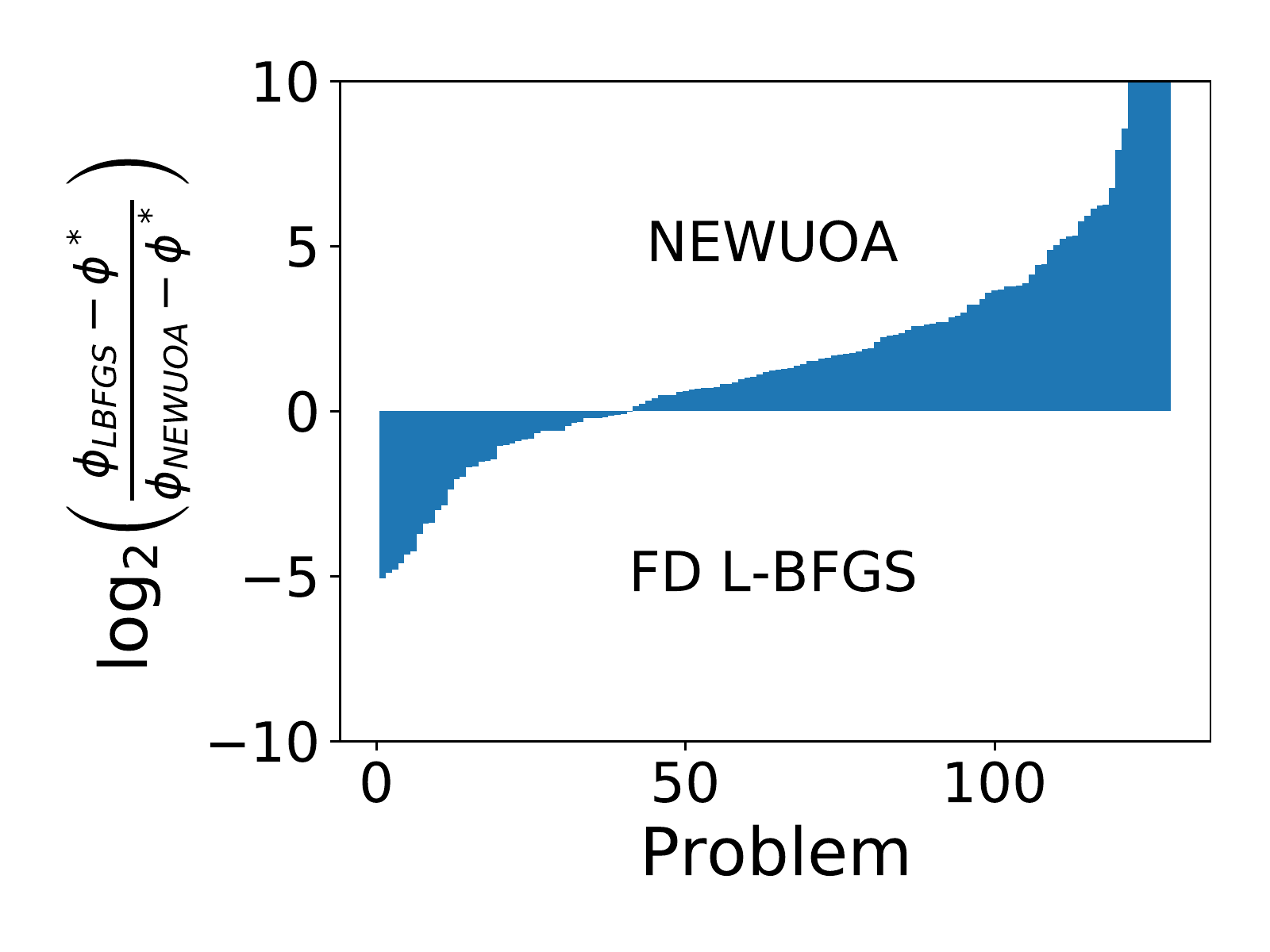}
\includegraphics[width=0.32\textwidth]{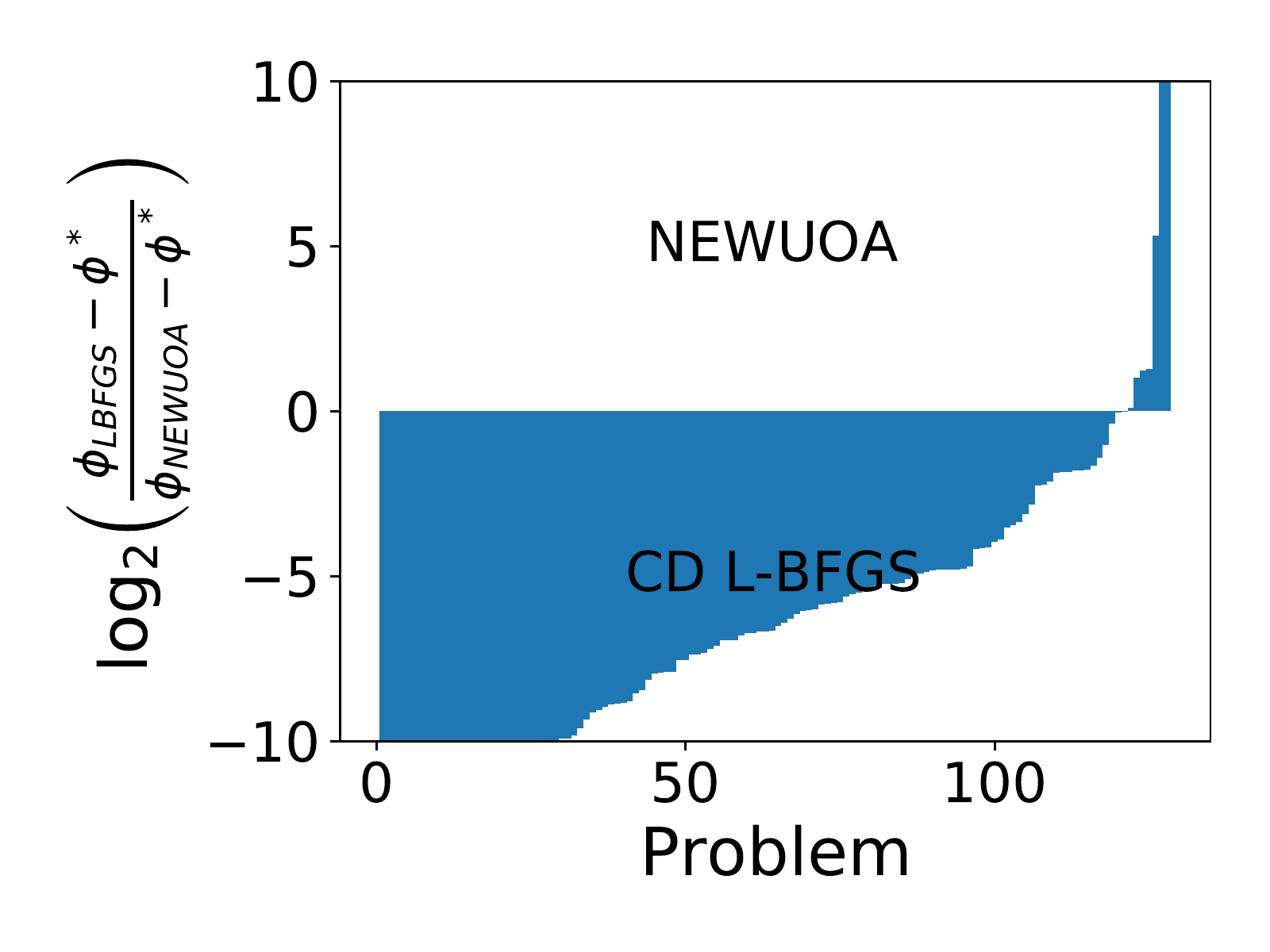}
\caption{{\em Accuracy, Noisy Case for $\sigma_f = 10^{-5}$}. Log-ratio optimality gap profiles comparing {\sc newuoa} against forward-difference {\sc l-bfgs} (left) and central-difference {\sc l-bfgs} (right). }
\label{fig:noisy opt}
\end{figure}

We  first compare the accuracy achieved by each algorithm, as measured by the  optimality gap $\phi(x_k) - \phi^*$. We do so by running {\sc newuoa} until $\rho = \rho_{\text{end}} = 10^{-6}$, and running the finite-difference {\sc l-bfgs} method until  the  objective function could not be improved over 5 consecutive iterations. In Figure \ref{fig:noisy opt}, we  report the log-ratio profile
\begin{equation}   \label{jl2}
   \log_2 \left(\frac{\phi_{\text{\sc lbfgs}} - \phi^*}{\phi_{\text{\sc newuoa}} - \phi^*} \right)
\end{equation}
for $\sigma_f=10^{-5}$, where $\phi_{\text{\sc lbfgs}}$, $ \phi_{\text{\sc newuoa}}$ denote the lowest objective achieved by each method.
(The results are representative of those obtained for $\sigma_f \in \{10^{-1}, 10^{-3}, 10^{-7}\}$.)
As in the noiseless case, we removed 4 problems where {\sc newuoa} terminates within the first $2n + 1$ function evaluations for all noise levels, as well as thirteen problems where the algorithms are known to converge to different minimizers without noise.

As seen in Figure \ref{fig:noisy opt}, {\sc newuoa} achieves higher accuracy in the solution than forward-difference {\sc l-bfgs}, while central-difference {\sc l-bfgs} yields  far better accuracy than both. It is not surprising that central differencing yields much higher accuracy than forward differencing since the noise level of their gradient approximations  is, respectively, $O(\sigma_f^{2/3})$ and $O(\sigma_f^{1/2})$.
%for forward differencing, the noise level of the gradient $\sigma_g$ is of order $\sigma_f^{1/2}$,while central differencing yields a noise level of the order $\sigma_f^{2/3}$. 
On the other hand, it is not straightforward to analyze the error contained in the gradient approximation constructed by {\sc newuoa}, and in turn its final accuracy. To try to shed some light into this question, we performed the following experiments. 

\begin{figure}[htp]
\centering
\includegraphics[width=0.32\textwidth]{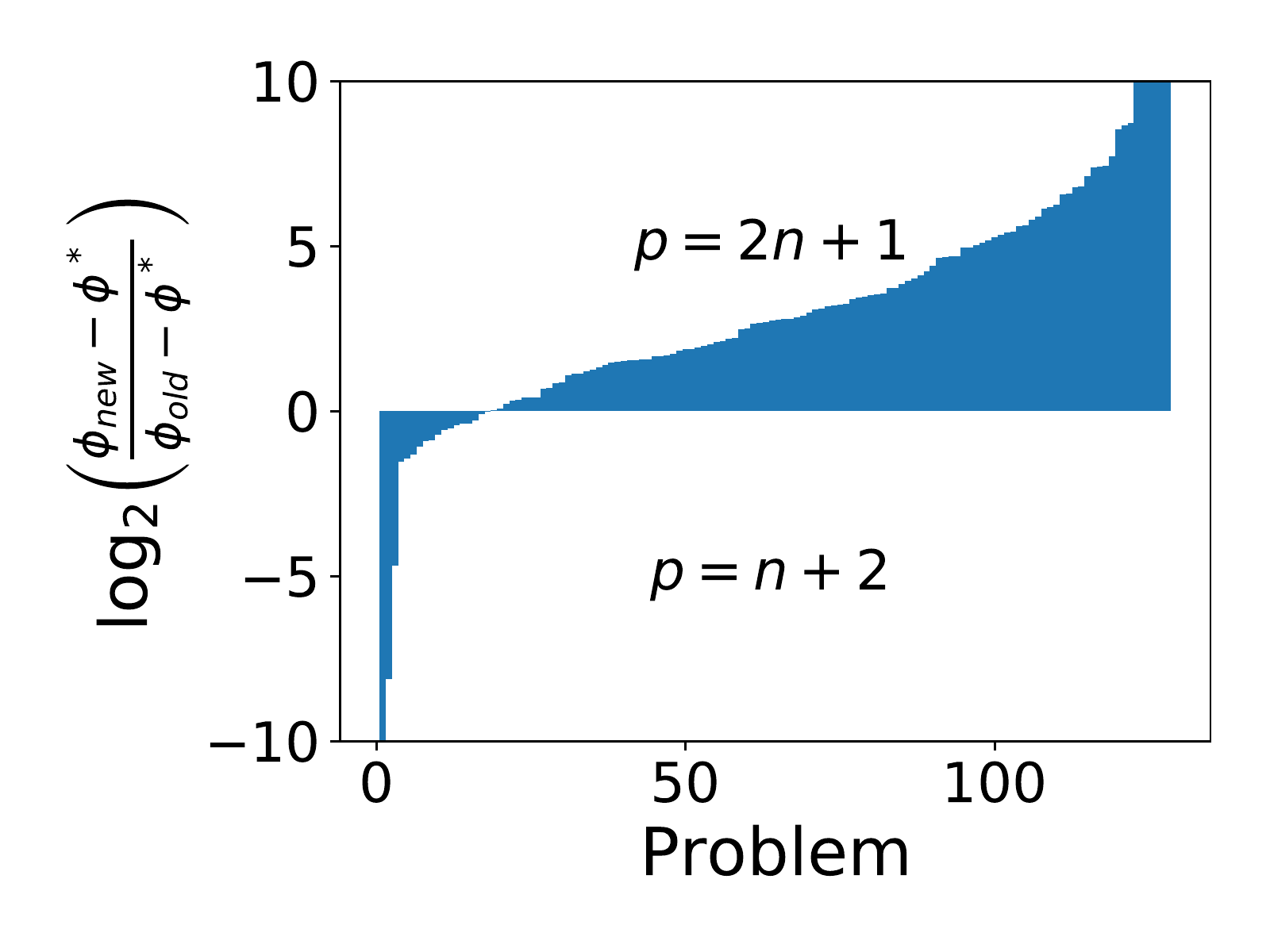}
\includegraphics[width=0.32\textwidth]{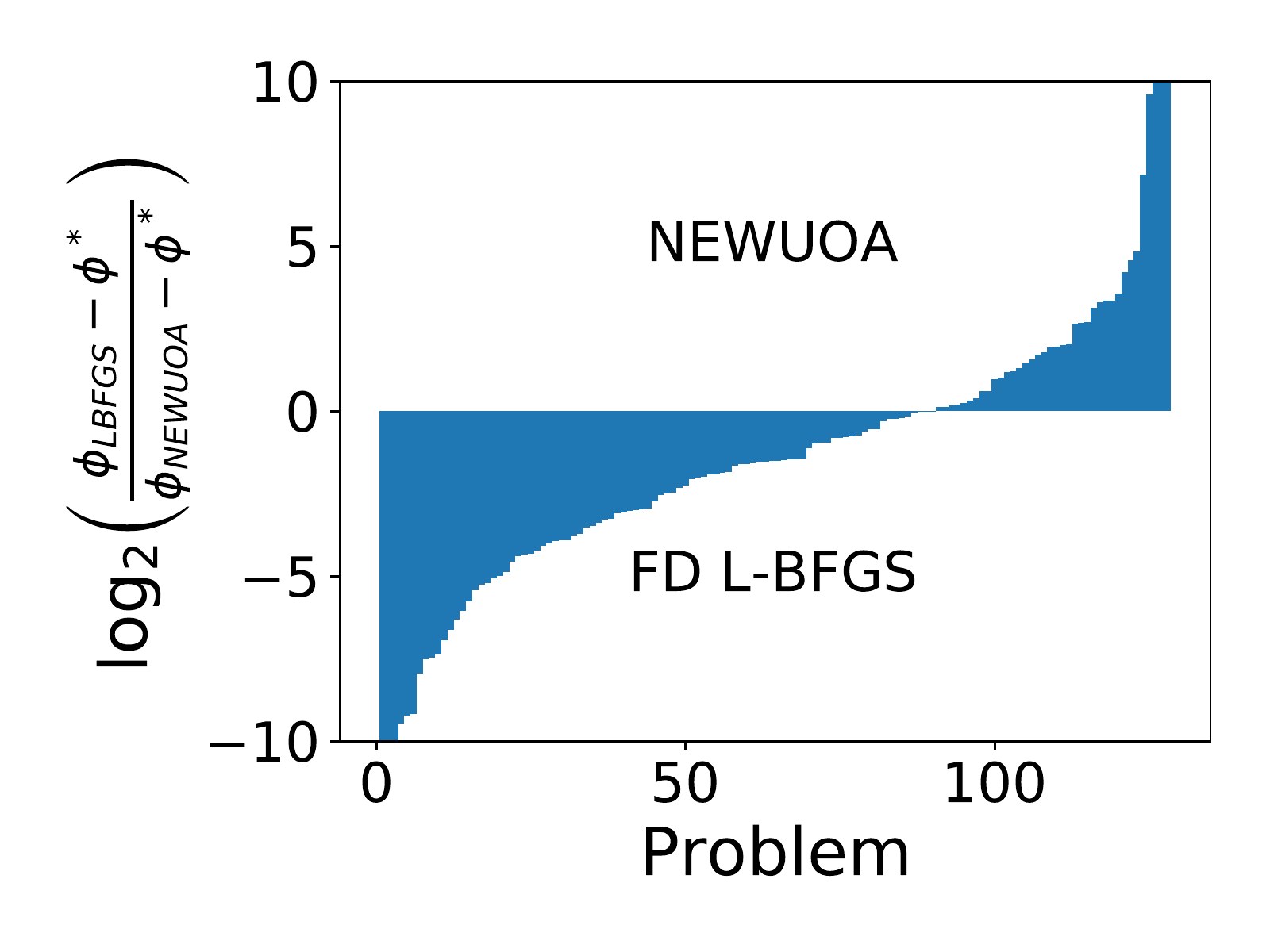}
\caption{{\em Accuracy for  Non-Default Version of {\sc newuoa}, Noisy Case for $\sigma_f = 10^{-5}$}. Log-ratio optimality gap profiles comparing {\sc newuoa} with $p = 2n + 1$ and $p = n + 2$ points (left), and {\sc newuoa} with $p = n + 2$ points against forward-difference {\sc l-bfgs} (right).}
\label{fig:n+2 obj}
\end{figure}

Recall that {\sc newuoa} by default uses $p = 2n + 1$ interpolation points when constructing the quadratic model of the objective, while forward differencing only uses $n + 1$ points. It is then natural to test the performance of {\sc newuoa} with only $p = n + 2$ interpolation points. The results in Figure \ref{fig:n+2 obj} show that {\sc newuoa} now lags behind finite-difference {\sc l-bfgs}, which may be due the fact that the quality of its gradient also depends on the poisedness of the interpolation set \cite{conn2009introduction,berahas2019linear,berahas2019theoretical}. Additional tests, given in Appendix~\ref{sec:newuoa}, suggest that the choice $p = 2n+1$ recommended by Powell strikes the right balance between accuracy in the solution and the speed of algorithm.

%When one reduces the number of interpolation points used in {\sc newuoa}, we observe a significant deterioration in its final accuracy. In particular, with a comparable number of interpolation points, forward difference {\sc l-bfgs} actually obtains solutions with comparable or higher accuracy than {\sc newuoa}. This is not surprising since the solution accuracy will depend on the quality of the gradient, which relies on the well-poisedness of the interpolation set; see \cite{berahas2019linear,berahas2019theoretical}. For more details, please refer to Appendix \ref{sec:newuoa}.

\begin{figure}[htp]
\centering
\includegraphics[width=0.32\textwidth]{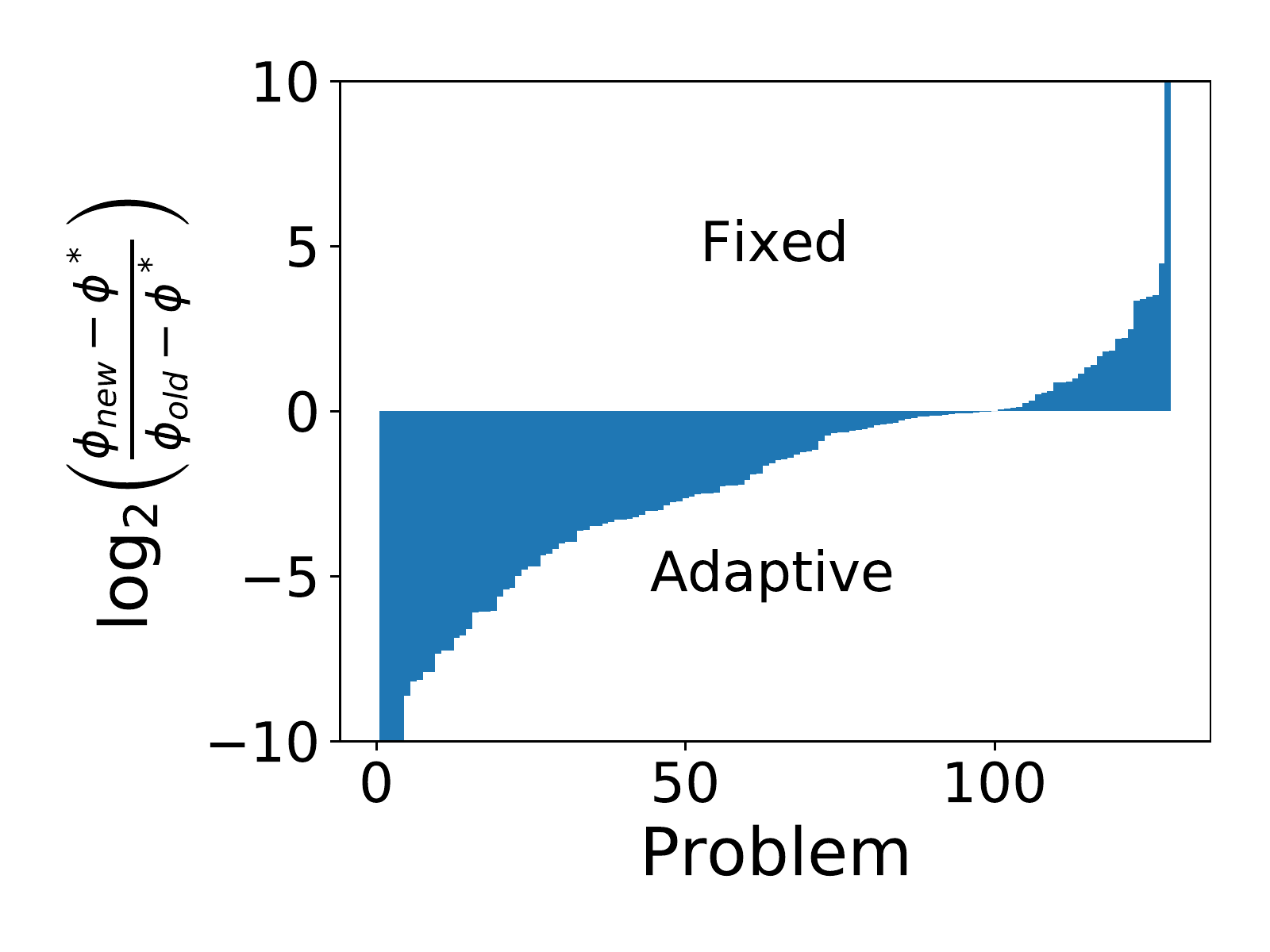}
\caption{{\em Accuracy of Fixed and Adaptive Lipschitz Estimation for {\sc l-bfgs}, Noisy Case for $\sigma_f = 10^{-5}$}. Log-ratio optimality gap profiles comparing fixed Lipschitz estimates against adaptive Lipschitz estimation for forward-difference {\sc l-bfgs}.}
\label{fig:adap vs fixed}
\end{figure}

It has been commonly suggested that estimating the Lipschitz constant $L_i$ along each coordinate at the initial point is sufficient for the entire run of a finite-difference method \cite{GillMurrSaunWrig83,more2012estimating}. This, in fact, is not the case in our experiments. We compare in Figure \ref{fig:adap vs fixed}  forward-difference {\sc l-bfgs} with the adaptive Lipschitz estimation given in Procedure~I against estimating $L$ only once at the beginning of the run. We see that, in terms of solution accuracy, fixing the Lipschitz constant $L$ at the start does not yield nearly as good accuracy as the adaptive procedure.
This is because that initial estimate is often not a proper bound on the second derivative in the later stages of the run, hence yielding a poor estimate of the finite-difference interval $h$, and consequently greater error in the gradient approximation. A variety of other strategies for estimating $L$ are discussed in  Appendix~\ref{app:lipschitz}.

\subsubsection{Efficiency}

\begin{figure}[ht]
\centering
\includegraphics[width=0.32\textwidth]{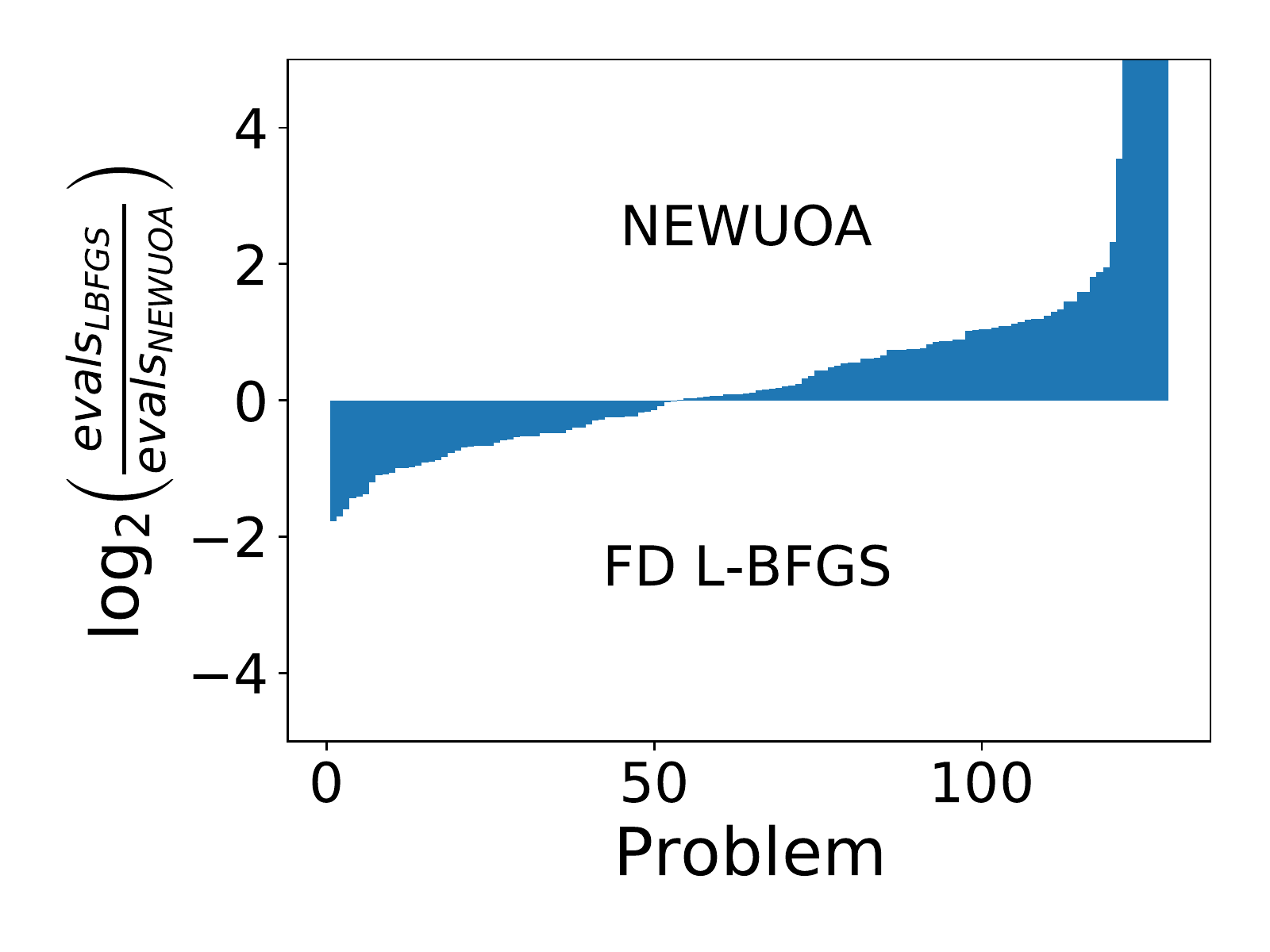}
\includegraphics[width=0.32\textwidth]{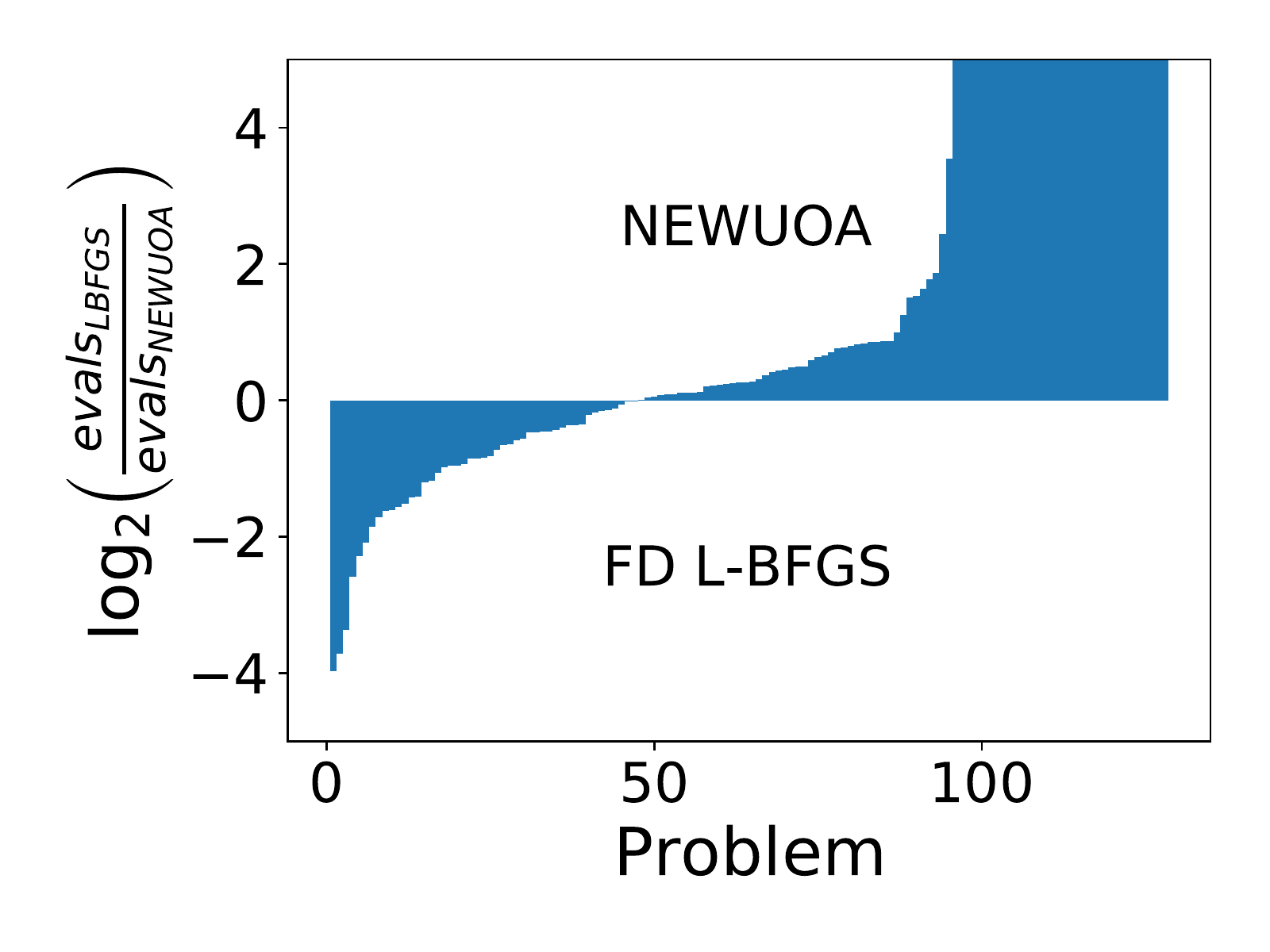}\\
\includegraphics[width=0.32\textwidth]{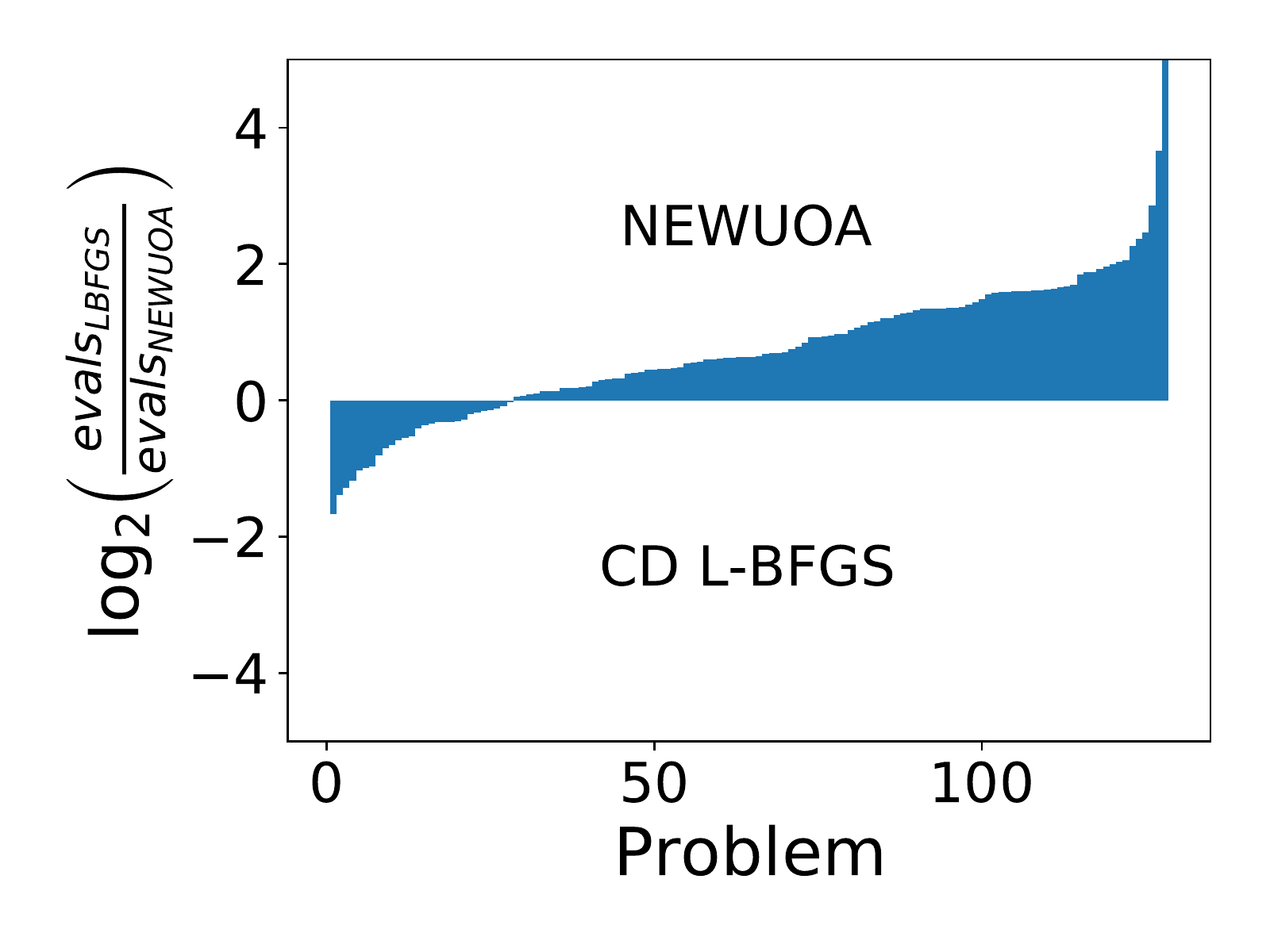}
\includegraphics[width=0.32\textwidth]{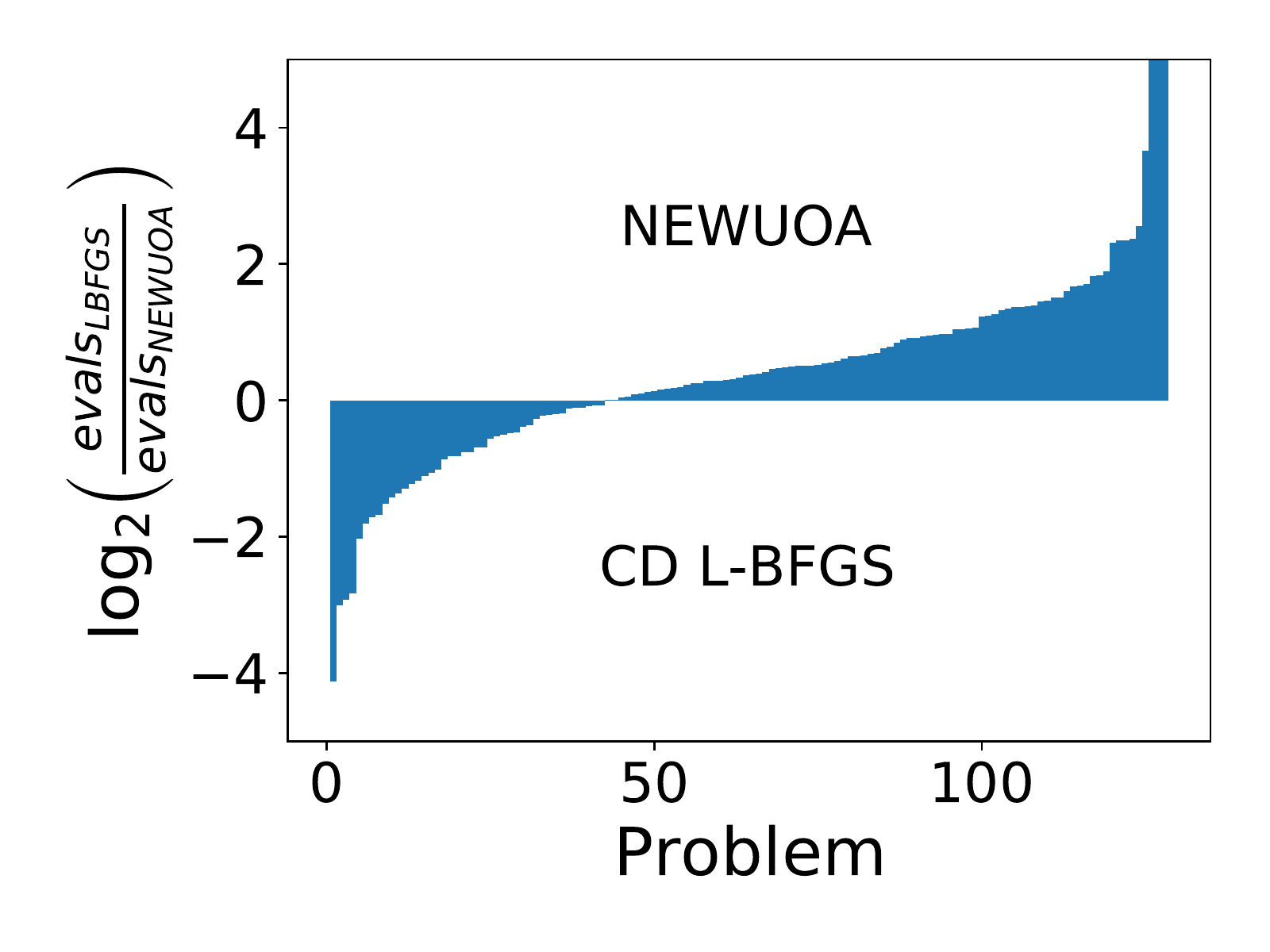}
\caption{{\em Efficiency, Noisy Case for $\sigma_f = 10^{-5}$}. Log-ratio profiles for the number of function evaluations to achieve \eqref{eq:term crit 2} for $\tau = 10^{-2}$ (left) and $10^{-6}$ (right),  comparing {\sc newuoa} against forward-difference {\sc l-bfgs} (top) and central-difference {\sc l-bfgs} (bottom). These plots are representative of the other noise levels $\sigma_f \in \{10^{-1}, 10^{-3}, 10^{-5}, 10^{-7}\}$.}
\label{fig:noisy feval fd l-bfgs}
\end{figure}

To compare the efficiency of different algorithms, we record the number of function evaluations required to reach a particular function value. We  use the best solution from {\sc newuoa} (denoted by $\phi_{new}$) as a baseline and report the number of function evaluations necessary to achieve
\begin{equation}\label{eq:term crit 2}
    \phi(x_k) - \phi_{new} \leq \tau \cdot (\phi(x_0) - \phi_{new}),
\end{equation}
for varied $\tau$. In Figure \ref{fig:noisy feval fd l-bfgs}, we report log-ratio profiles comparing the number of function evaluations necessary to satisfy \eqref{eq:term crit 2}, i.e.,
\begin{equation}   \label{jl2}
   \log_2 \left(\frac{\El}{\En}\right).
\end{equation}

\noindent
We observe from this figure that  {\sc newuoa} is more efficient than forward-difference {\sc l-bfgs}. The advantage is less significant for $\tau=10^{-2}$ but becomes  pronounced for $\tau=10^{-6}$, which is to be expected because, as mentioned above,  the forward-difference option struggles to achieve as high accuracy as {\sc newuoa}. Central-difference {\sc l-bfgs} is also less efficient overall than {\sc newuoa}, but becomes more competitive as $\tau$ is decreased, which is a product of its ability to converge to a higher accuracy than forward differencing. It is notable that {\sc newuoa} is able to deliver such strong performance without knowledge of the noise level of the function. The adjustment of the two trust-region radii in {\sc newuoa} seems to be quite effective: shrinking the radii fast enough to ensure steady progress, but not so fast as to create models dominated by noise. To our knowledge, there has been no in-depth study of the {\em practical} behavior of {\sc newuoa} (or similar codes) in the presence of noise;  we regard this as an interesting  research topic.

\begin{figure}[htp]
\centering
\includegraphics[width=0.32\textwidth]{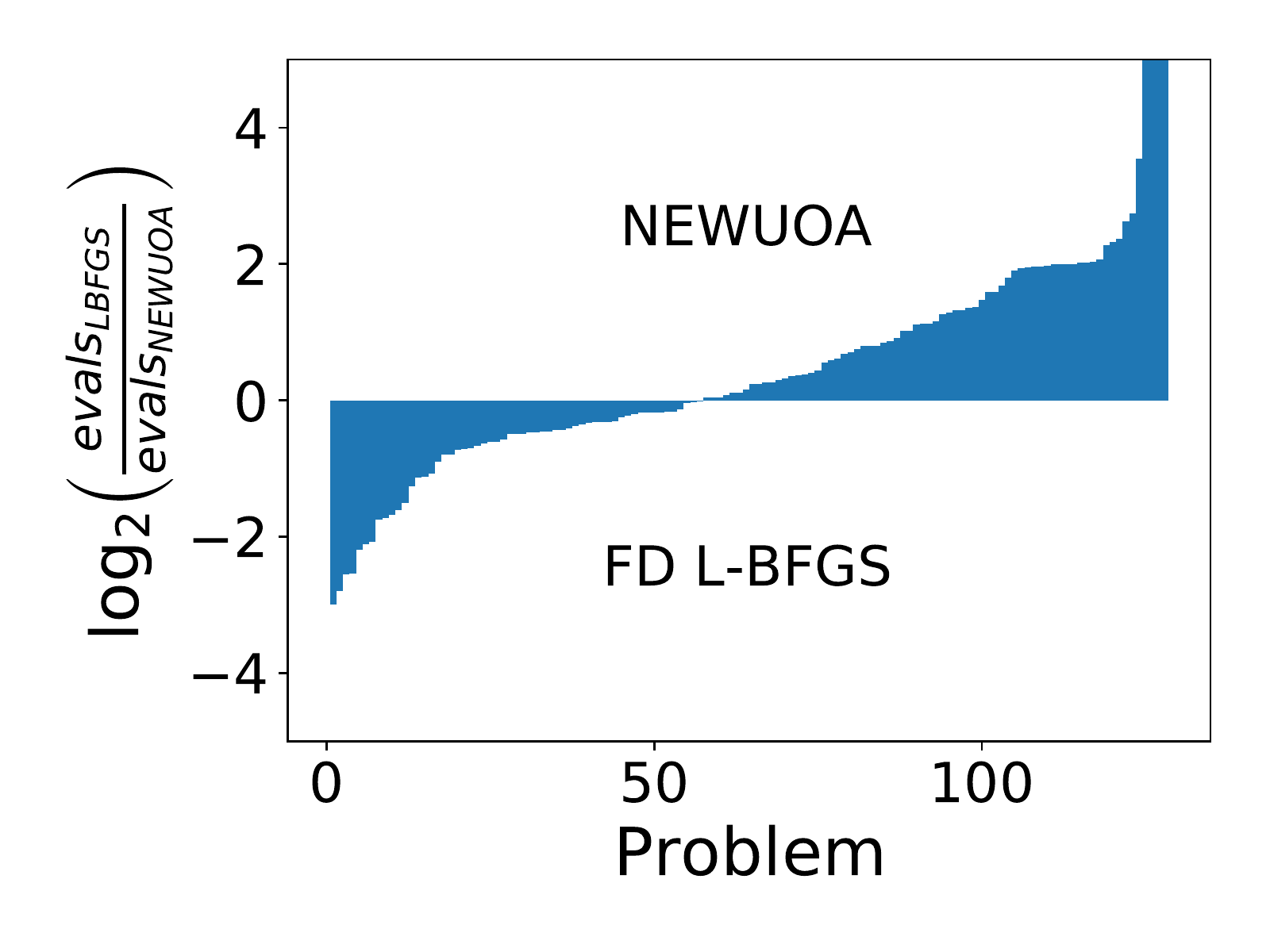}
\includegraphics[width=0.32\textwidth]{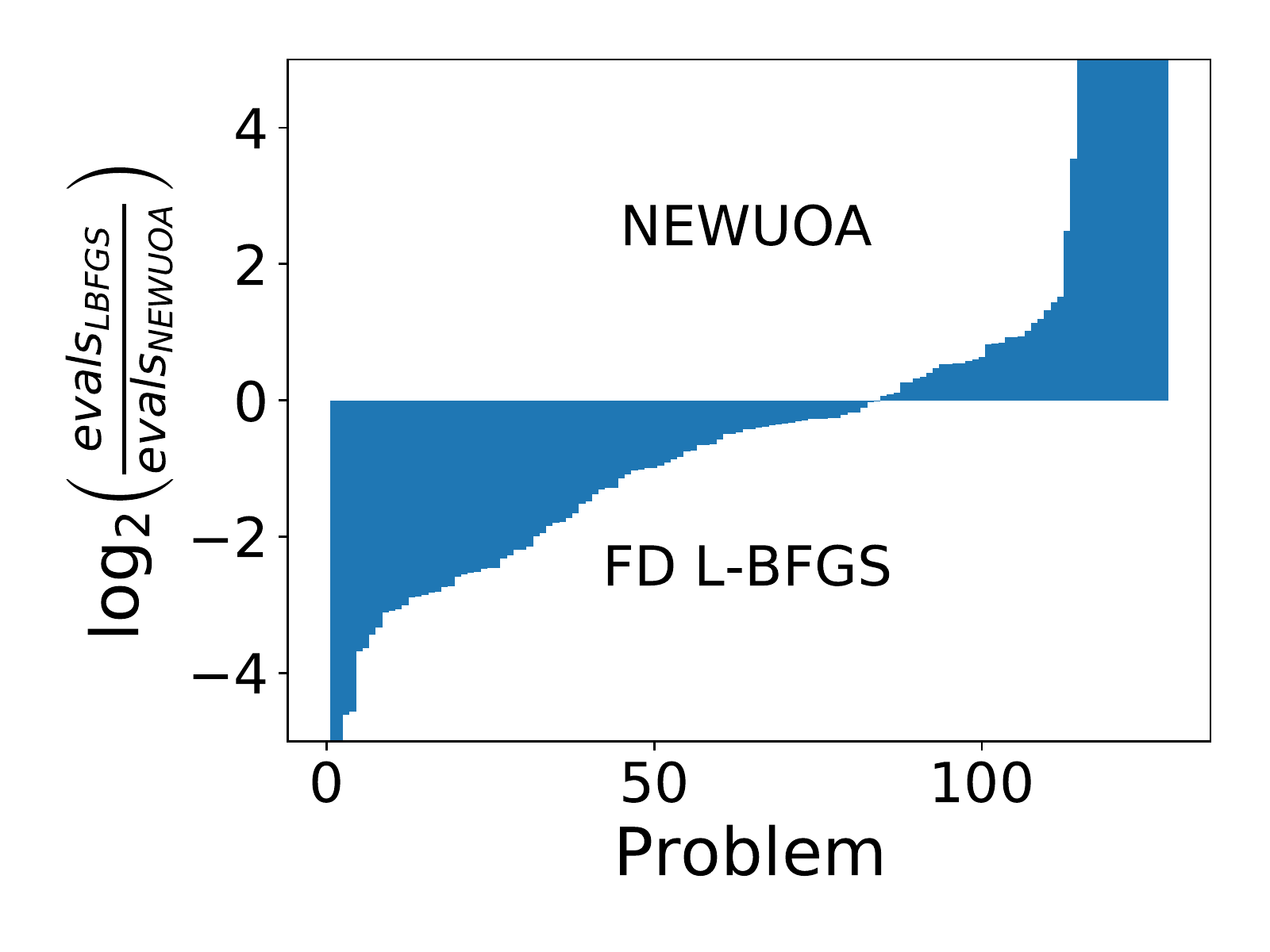}
\caption{{\em Efficiency for Non-Default Version of {\sc newuoa}, Noisy Case.} Log-ratio profiles for the number of function evaluations to achieve \eqref{eq:term crit 2} for $\tau = 10^{-2}$ (left) and $10^{-6}$ (right) for $\sigma_f = 10^{-5}$ comparing {\sc newuoa} with $p = n + 2$ interpolation points against forward-difference {\sc l-bfgs}.}
\label{fig:noisy feval n+2}
\end{figure}

We also compare the efficiency of {\sc newuoa} with only $p = n + 2$ interpolation points against forward-difference {\sc l-bfgs}; see Figure~\ref{fig:noisy feval n+2}. We observe that forward-difference {\sc l-bfgs} is  more competitive  in this case. More details are given in Appendix \ref{sec:newuoa}.

One conclusion from these results is that a more sophisticated implementation of finite differences is required to close the performance gap between {\sc l-bfgs} and  {\sc newuoa}. This requires a more accurate, yet affordable, mechanism for estimating the Lipschitz constants in the course of the run. The design of an adaptive procedure of this kind is the subject of a forthcoming paper \cite{shi2021automatic}.
It must be able to deal with some of the challenging situations described next.

\subsubsection{Commentary}  \label{potential}
It is well known in computational mathematics that finite-difference methods can be unreliable in the solution of differential equations when knowledge of higher-order derivative estimates is poor. Optimization, however, provides a more benign setting in that a poor choice of $h$ that leads to a bad step can, in many situations, be identified and corrective action can be taken. Gill et al. \cite{GillMurrSaunWrig83} propose an iterative procedure for estimating $h$, but we believe that more sophisticated strategies must be developed to improve the reliability of finite-difference DFO methods. The design of such techniques is outside the scope of this paper, but we now give one example illustrating the challenges that arise and some opportunities for addressing them.

Let us consider the use of forward differences for the solution of the \texttt{DENSCHNE} problem from the \texttt{CUTEst} collection, which is given by
\begin{equation*}
\min_{x \in \R^3} \phi(x) = x_1^2 + (x_2 + x_2^2)^2 + (-1 + e^{x_3})^2, ~~~ \text{with } x_0 = (2, 3, -8).
\end{equation*} 
Let us employ a different interval for each variable, and let us focus in particular on $x_3$. The choice of the differencing interval $h_3$ depends on an estimate $L$ of the second derivative. We have that
\begin{align*}
\frac{\partial \phi}{\partial x_3}(x) & = 2 e^{x_3} (e^{x_3} - 1), \qquad
\frac{\partial^2 \phi}{\partial x_3^2}(x)  = 2 e^{x_3} (2 e^{x_3} - 1).
\end{align*}
Let us define $L$ to be the value of the second derivative at $x_3=-8$.
This gives $L = \frac{2 e^8 - 4}{e^{16}} \approx 6.705 \times 10^{-5}$, which reflects the fact that $\phi$ is very flat at that point. Suppose that the noise level is $\sigma_f = 10^{-1}$, then from \eqref{eq:fd formula}, we have that $h \approx 20.539$, which is excessively long given that the function contains an exponential term. The forward-difference derivative approximation, $[g(x)]_3 \approx 3.791 \times 10^9$, is entirely wrong; the correct value is  $[\nabla \phi(x)]_3 = 2 e^{-8} (e^{-8} - 1) \approx -6.707 \times 10^{-4}$. Clearly, the problem is caused by the fact that the second derivative changes rapidly in the interval $[-8, -8 + 20.539]$ and we employed its lowest instead of its largest value. By greatly underestimating $L$ we computed an harmful differencing interval. (We should note in passing that {\sc newuoa} does not have any difficulties with the \texttt{DENSCHNE} problem.)

However, the difficulties arising in this problem can be prevented at two levels. First, it should be possible to design an automatic procedure similar to the one in Gill et al. \cite{GillMurrSaunWrig83} that would not allow for the creation of such a poor estimate of $h$. In particular, the procedure should estimate $L$ along an interval, not just in a pointwise fashion. Even if such a procedure is unable to take corrective action, monitoring the optimization step can identify difficulties. For the example given above, the fact that $g(x)$ is huge causes the algorithm to compute a very small steplength $\alpha_k$. This gives a warning signal that should cause an examination of the interval $h$, and would immediately reveal that the second derivative has varied rapidly in the differencing interval and that $L$ must be recomputed. 

It is easy to envision problems where the finite-difference interval $h$ is underestimated, in which case the effect of noise can be damaging. In addition, one must ensure that automatic procedures for estimating and correcting $L$ are not too costly. Taken, all together, these challenges motivate research into more reliable techniques for finite-difference estimation, {\em in the context of optimization.}

\section{Nonlinear Least Squares}
\label{ch:squares}
 In many unconstrained optimization problems,  the objective function has a nonlinear least-squares form. Therefore, it is important to pay particular attention to this problem structure in the derivative-free setting. We write the problem as
\begin{equation*}
	 \min_{x \in \mathbb{R}^n} \phi(x) =  \tfrac{1}{2}\Vert \gamma (x)\Vert^2 = \tfrac{1}{2}\sum_{i = 1}^m \gamma_i^2(x) ,
\end{equation*}
where  $\gamma : \mathbb{R}^n \rightarrow \mathbb{R}^m$ is a smooth function. We assume that the Jacobian matrix $[J(x)]_{ij} = \frac{\partial \gamma_i(x)}{\partial x_j}$ is not available, but the individual residual functions $\gamma_i(x)$ can be computed. More generally, the evaluation of the $\gamma_i$ may contain noise so that the observed residuals are given by
\[
        r_i(x) = \gamma_i(x) + \epsilon_i(x), \quad i=1, \ldots, m,
\] 
where $\epsilon_i(x)$ models noise as in \eqref{uncprob}. Thus, the minimization of the true objective function $\phi$ must be performed based on noisy observations $r_i(x)$ that define the observed objective function 
\begin{equation}
	 f(x)=  \tfrac{1}{2}\Vert r(x)\Vert^2 = \tfrac{1}{2}\sum_{i = 1}^m  r_i^2(x) .
\end{equation}
Since noise is incorporated into each residual function, the model of noise is different from additive noise model in the general unconstrained case discussed in the previous section. In particular, the function evaluation $f(x) = \tfrac{1}{2}\sum_{i = 1}^m \gamma_i^2(x) + 2 \epsilon_i(x) \gamma_i(x) + \epsilon_i^2(x)$ contains both multiplicative and additive components of noise.

Our goal is to study the viability of methods based on finite-difference approximations to the Jacobian. To this end, we employ a classical Levenberg-Marquardt trust-region method where the Jacobian is approximated by differencing, and perform tests comparing it against a state-of-the-art DFO code designed for nonlinear least-squares problems. The rationale behind the selection of codes used in our experiments is discussed next.

\medskip\noindent{\em Interpolation Based Trust Region Methods.}

Wild \cite{SWCHAP14} proposed an interpolation-based trust-region method that creates a quadratic model for each of the residual functions $\gamma_i$, and aggregates these models to obtain an overall model of the objective function. This method has been implemented in  {\sc pounders}
 \cite{SWCHAP14,munson2012tao}, which has proved to be significantly  faster than the standard approach that ignores the nonlinear least-squares structure and simply interpolates values of $f$ to construct a quadratic model. Zhang, Conn and Scheinberg \cite{zhang2010derivative} describe a similar method, implemented in   {\sc dfbols}, where each residual function is approximated by a quadratic model using $p\in [n+1, (n+1)(n+2)/2]$ interpolation points; the value $p=2n+1$ being recommended in practice. 
  
 More recently, Cartis and Roberts \cite{roberts2019derivative} proposed a Gauss-Newton type approach in which an approximation of the Jacobian $J(x_k)$ is computed by linear interpolation using $n+1$ function values at recently generated points $Y_k = \{x_k, y_1,...,y_n \}$. Here $x_k$ is the current iterate, which generally corresponds to the smallest objective value observed so far. The approximation $\hat{J}(x_k)$ of the true Jacobian $J(x_k)$ is thus obtained by solving the linear equations 
\begin{align}   \label{lins}
	r(x_k) + \hat{J}(x_k)(y_j - x_k) = r(y_j), \quad  j = 1,...,n.
\end{align}
The step of the algorithm $s_{k}$ is defined as an approximate solution of the trust region problem
\begin{align*}
	\min_{s \in \mathbb{R}^n} \ m_k(s) = f(x_k) + g_k^Ts + \tfrac{1}{2}s^TH_ks \qquad\mbox{s.t.}
	\ \ \|s\|_2 \leq \Delta_k,
\end{align*}
where $g_k = \hat{J}(x_k)^Tr(x_k)$ and $H_k = \hat{J}(x_k)^T \hat{J}(x_k)$. The new iterate is given by $x_{k+1}= x_k + s_k$ and the new set $Y_{k+1}$ is updated by removing a point in $Y_k$ and adding $x_{k+1}$.  As in all model-based interpolation methods, it is important to ensure that the interpolation points do not lie on a subspace of $\mathbb{R}^n$ (or close to a subspace). To this end, the algorithm contains a geometry improving technique that spaces out interpolation points, as needed. The implementation described in \cite{roberts2019derivative}, and referred to as {\sc dfo-gn}, is reported to be faster than  {\sc pounders}, and scales better with the dimension of the problem. An improved version of {\sc dfo-gn} is {\sc dfo-ls} \cite{cartis2019improving}, which provides a variety of options and heuristics to accelerate convergence and promote a more accurate solution. The numerical results reported by Roberts et al. \cite{cartis2019improving} indicate that {\sc dfo-ls} is a state-of-the-art code for DFO least squares, and therefore will be used in our benchmarking.

\medskip\noindent
\textit{Finite-Difference Gauss-Newton Method}.

One can employ finite differencing to estimate the Jacobian matrix $J(x)$ within any method for nonlinear least squares, and since this is a mature area, there are a number established solvers. We chose {\sc lmder} for our experiments, which is part of the {\tt MINPACK} package \cite{MoreGarbHill80} and is also available in 
the {\tt scipy} library. We did not employ {\sc lmdif}, the finite-difference version of {\sc lmder}, because it does not allow the use of different differencing intervals for each of the residual functions $r_i(x)$; we elaborate on this point below. Another code available in {\tt scipy} is  {\sc trf} \cite{branch1999subspace}, but our tests show that {\sc lmder} is slightly more efficient in terms of function evaluations, and tends to give higher accuracy in the solution. The code {\sc nls}, recently added to the {\tt Galahad} library \cite{galahad} would provide an interesting alternative. That method, however, includes a tensor to enhance the Gauss-Newton model, and since this may give it an advantage over {\sc dfo-ls}, we decided to employ the more traditional code {\sc lmder}.

 In summary, the solvers used in our tests are:
\begin{itemize}
	\item {\tt LMDER}: A derivative-based Levenberg-Marquardt trust-region algorithm from the {\tt MINPACK} software library \cite{MoreGarbHill80}, where the finite-difference module is user-supplied by us. We called the code in Python 3.7 through {\tt scipy} version 1.5.3, using the default parameter settings. 
	\item {\tt DFO-LS}: The most recent DFO software developed by Cartis et al. \cite{cartis2019improving} for  nonlinear least squares. This method uses linear interpolation to construct an approximation to the Jacobian matrix, which is then used in a Gauss-Newton-type method. We used version 1.0.2 in our experiments, with default settings except that the {\tt model.abs\_tol} parameter is set to 0 to avoid early termination %(unless specified otherwise). 
%\jn{[JN: Do we need the last parenthetical question?]}
\end{itemize}
 
The test problems in our experiments are those used by Mor\'e and Wild \cite{more2009benchmarking}, which have also been employed by Roberts and Cartis \cite{roberts2019derivative} and  Zhang et al. \cite{zhang2010derivative}. The 53 unconstrained problems in this test set include both zero and nonzero residual problems, with various starting points, and are all small dimensional, with $n \leq 12$. To measure efficiency, we regard $m$ evaluations of individual residual components, say $\{r_i(\cdot)\}_{i=1}^m$,  as one function evaluation. These $m$ evaluations may not necessarily be performed at the same point. We also assume that one can compute any component $r_i$ without a full evaluation of the vector $r(x)$. We terminate the algorithms when either: i) the maximum number of function evaluations ($500 \times n$) is reached, ii)  an optimality gap stopping condition, specified below, is triggered; or iii) the default termination criterion of the two codes is satisfied with tolerance of $10^{-8}$ (this controls the minimum allowed trust region radius).

%%%%KEEP THE FOLLOWING TEXT THAT HAS BEEN COMMENTED OUT
%\begin{verbatim}
%From Nick Gould. f Fortran is of any use to you, GALAHAD has a number
%of nonlinear least-squares packages. The latest is NLS
%that includes the tensor-Newton (higher order) method,
%while FILTRANE is an older edition, and of course
%ancient LANCELOT works pretty well on such problems.
%There is a Matlab interface (via the fortran) to NLS.
%And of course GALAHAD is fully hooked up to CUTEst,
%where the majority of the latest examples are from
%least-squares fitting (both academic and real-world
%from the RAL facilities).
%
%Full download and installation instructions are on
%
%https://urldefense.com/v3/__https://github.com/ralna/GALAHAD/wiki__;!!Dq0X2DkFhyF93HkjWTBQKhk!BZS0ZxVmy46wXpplRIzW3q_tDtV0pakJbqllBPH6nFUe5VxkJKgGeMka3JAiozlpKvP5Ag4$
%
%We also have a small-scale least-squares solver, RALFIT
%
%https://urldefense.com/v3/__https://github.com/ralna/RALFit__;!!Dq0X2DkFhyF93HkjWTBQKhk!BZS0ZxVmy46wXpplRIzW3q_tDtV0pakJbqllBPH6nFUe5VxkJKgGeMka3JAiozlp_LWzvg0$
%
%that is also in fortran, but with interfaces to C and Python
%\end{verbatim}

\subsection{Experiments on Noiseless Functions}
We first consider the deterministic case corresponding to $\epsilon_i(x) \equiv 0, \, \forall i$. For {\sc lmder}, we estimate the Jacobian $J(x)$ using forward differences. As before, let $\epsilon_M$ denote machine precision. We first evaluate $\{r(x + h_je_j) \}_{j=1}^n$, where 
\begin{equation*}
	h_j = \max\{1, |[x]_j|\}\sqrt{\epsilon_M},
\end{equation*}
and compute the Jacobian estimate as
\begin{align*}
	[\hat{J}(x)]_{ij} = \frac{r_i(x + h_je_j) - r_i(x)}{h_j}.
\end{align*}
This formula for $h_j$ does not include a term that approximates the size of the second derivative (c.f.\eqref{eq:fd formula}) because we observed that such a refinement is not crucial in our experiments with noiseless functions, and our goal is to identify the simplest formula for $h$ that makes a method competitive. 
%\jn{\sout{Therefore, the differencing interval is the same among all residuals and only depend on the coordinate directions. }}

As in the previous section, we consider log-ratio profiles to compare the efficiency and accuracy of {\sc lmder} and {\sc dfo-ls}. We record 
\begin{equation*}
	\log_2 \left(\frac{\text{evals}_\text{\sc dfols}}{\text{evals}_\text{\sc lmder}}\right) \quad \text{and} \quad \log_2 \left(\frac{\phi_{\text{\sc dfols}} - \phi^*}{\phi_{\text{\sc lmder}} - \phi^*}\right),
\end{equation*}
where  $\phi^*$ is obtained from \cite{roberts2019derivative}, $\phi_{\text{\sc dfols}},\,  \phi_{\text{\sc lmder}}$ denote the best function values achieved by the respective methods, and 
$\text{evals}_{\text{\sc dfols}}, \, \text{evals}_{\text{\sc lmder}}$ denote the number of function evaluations needed to satisfy the termination test
\begin{equation}\label{eq:term ls}
    \phi(x_k) - \phi^* \leq \tau \max\{1, |\phi^*|\},
\end{equation}
for various values of $\tau$. (We set $\text{evals}_{\text{\sc dfols}}$ and $\text{evals}_{\text{\sc lmder}}$ to be a very large number if the above condition is not achieved within the maximum number of function evaluations). 
The results comparing {\sc dfo-ls} and {\sc lmder} are summarized in Figures~\ref{fig:ls_det_accuracy} and \ref{fig:ls_det_efficiency}. The complete table of results is given in Appendix~\ref{app:ls det}.

 \begin{figure}[htp]  %Non-Noisy Case Morales Plot for solution quality
 \centering
	\includegraphics[width=0.32\textwidth]{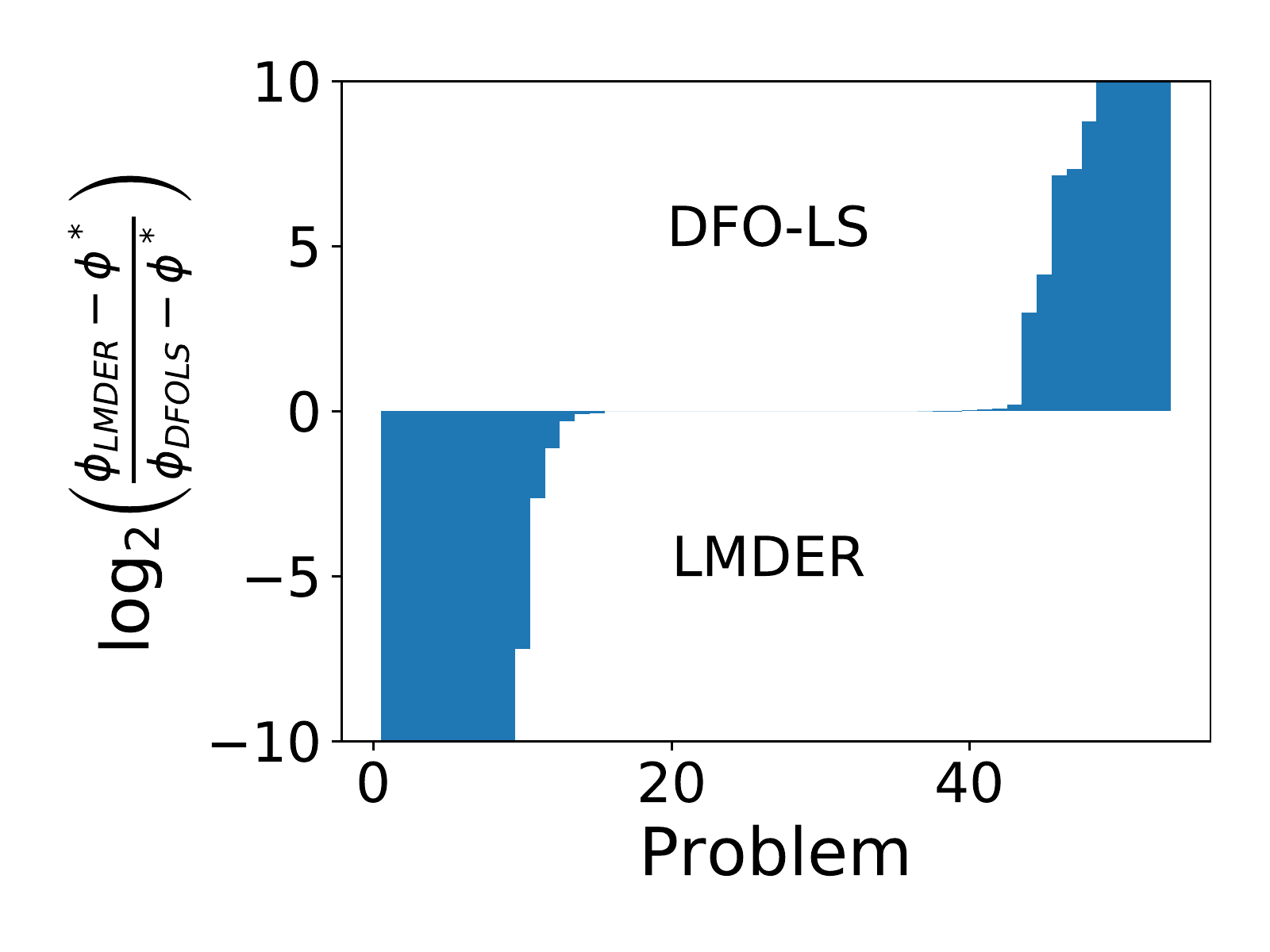}
    \caption{{\em Accuracy, Noiseless Case.} Log-ratio optimality gap profiles comparing {\sc dfo-ls} and {\sc lmder} for $\epsilon(x) = 0$. }\label{fig:ls_det_accuracy}
\end{figure}

\begin{figure}[ht] %Non-Noisy Case Morales Plot
    \centering
	\includegraphics[width=0.32\textwidth]{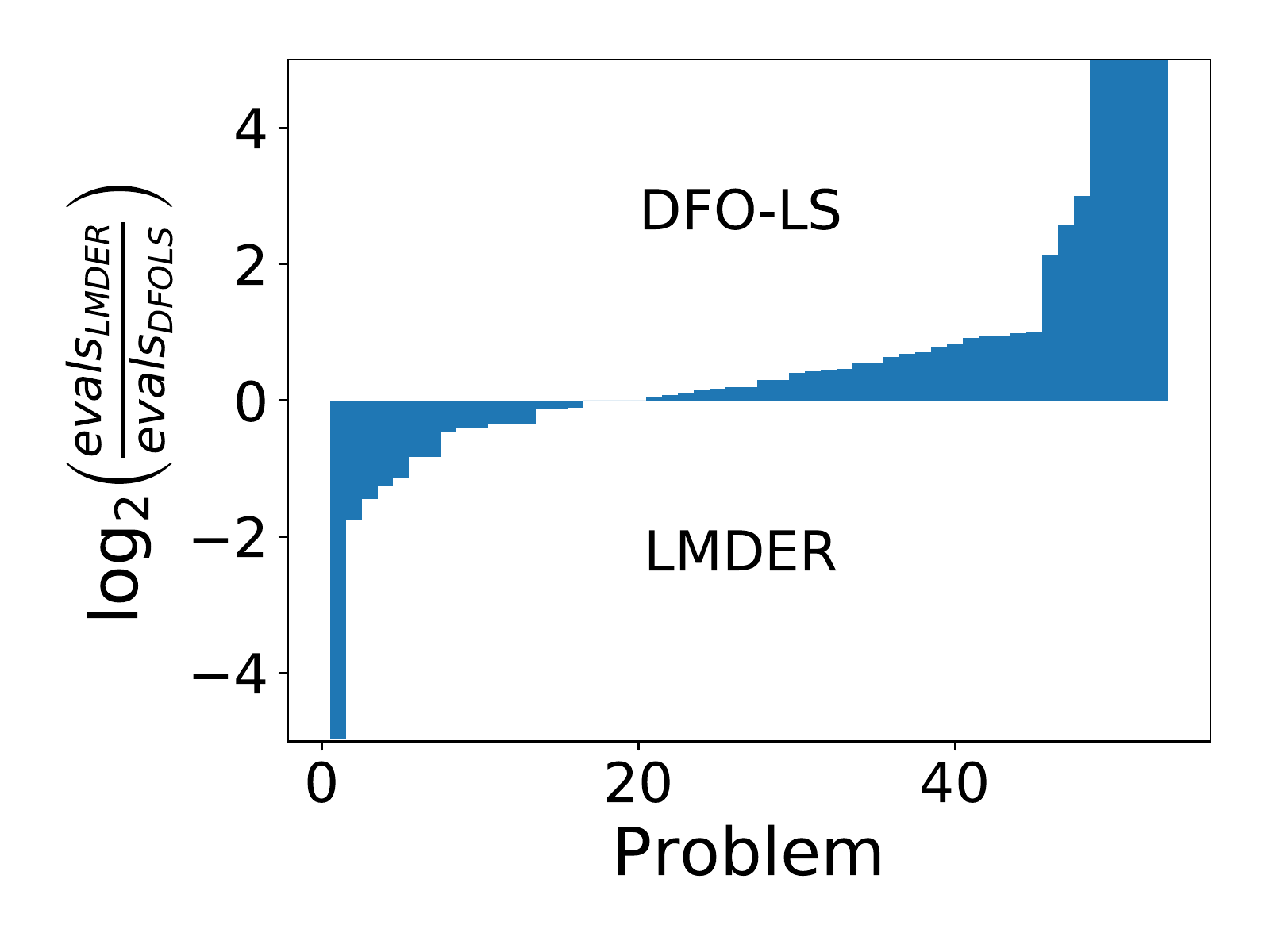}
	\includegraphics[width=0.32\textwidth]{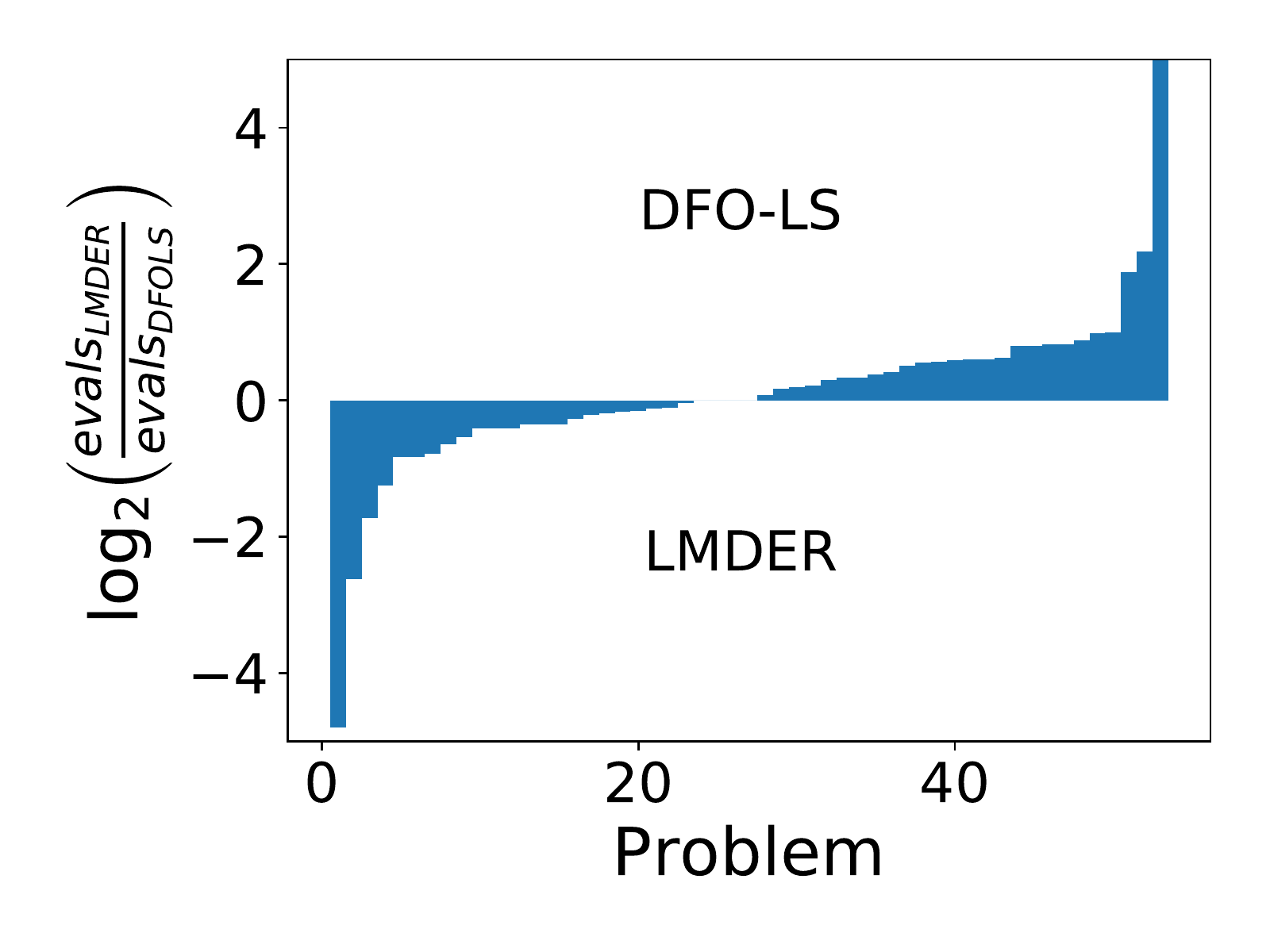}	
	\includegraphics[width=0.32\textwidth]{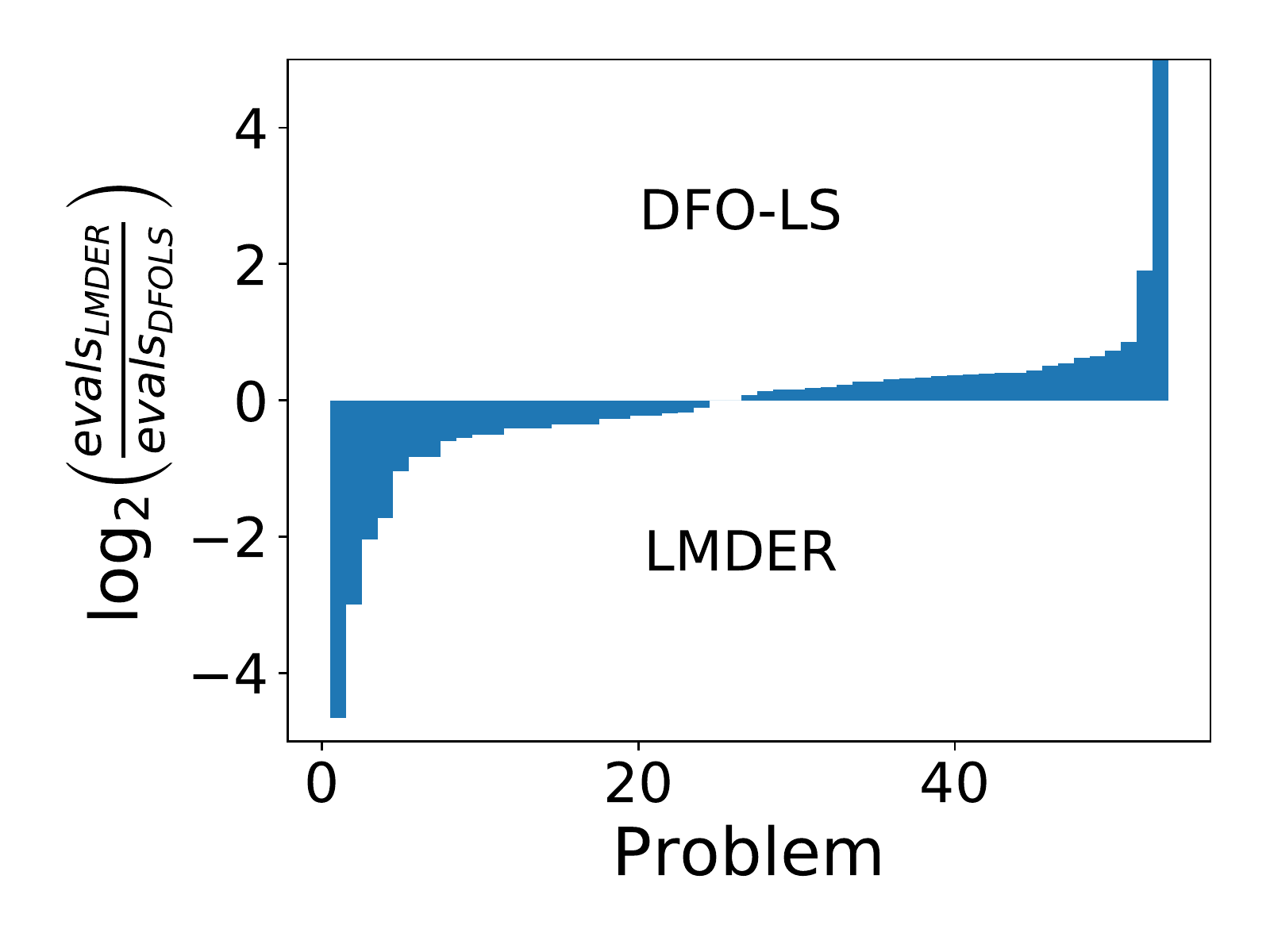}
    \caption{ {\em Efficiency, Noiseless Case}. Log-ratio profiles comparing {\sc dfo-ls} and {\sc lmder} for $\epsilon(x) = 0$. The figures measure number of function evaluations to satisfy \eqref{eq:term ls} for $\tau = 10^{-1}$ (left), $10^{-3}$ (middle), and $10^{-6}$ (right).  } 
    \label{fig:ls_det_efficiency}
\end{figure}

The performance of the two methods appears to be comparable since the area of the shaded regions is similar. This may be surprising since the Gauss-Newton approach seems to be a particularly effective way of designing an interpolation-based trust-region method, requiring only linear interpolation to yield useful second-order information.  Cartis and Roberts \cite{roberts2019derivative} dismiss finite-difference methods at the outset, and do not provide numerical comparisons with them. However, the tradeoffs of the two methods merit careful consideration.  {\sc dfo-ls} requires only 1 function evaluation per iteration, but its gradient approximation, $\hat{J}(x_k)^Tr(x_k)$, is  inaccurate until the iterates approach the solution and the trust region has shrunk. In contrast, the finite-difference Levenberg-Marquardt method in {\sc lmder} computes quite accurate gradients in this noiseless setting, requiring a much smaller number of iterations, but at a much higher cost per iteration in terms of function evaluations. The tradeoffs of the two methods appear to yield, in the end, similar performance, but we should note that  {\sc dfo-ls} is typically more efficient in the early stages of the optimization, as illustrated in Figure~\ref{fig:ls_det_efficiency} for the low tolerance level $\tau = 10^{-1}$. On the other hand, the finite-difference approach is more amenable to parallel execution, as mentioned in the previous section.

Let us now consider the linear algebra cost of the two methods. A typical iteration of {\sc dfo-ls} computes the LU factorization of the matrix $[y_1-x_k,...,y_n-x_k]$, and solves the interpolation system \eqref{lins} for $n$  different right hand sides, at a per-iteration cost of $O(mn^2)$ flops, assuming that $m > n$. ({\sc dfo-ls} offers a number of options, such as regression, which my involve a higher cost, but we did not invoke those options.) The linear algebra cost of the finite-difference version of {\sc lmder} described above is $O(mn)$ flops. Therefore, in terms of CPU time, {\sc lmder} is faster than {\sc dfo-ls} in our experiments involving inexpensive function evaluations. This is illustrated in Table~\ref{tab:ls cpu time} for a typical problem ({\tt BROWNAL}) with varying dimensions $n$ and $m$. 

\begin{table}[htp] 
    \centering
    \footnotesize
    \begin{tabular}{ c | c | c | c | c |  c | c | c}
        \toprule
        $(n,m)$ & (5,5) & (10,10) &  (30,30) & (50,50) & (100,100) & (300,300)& (500,500) \\
        \midrule
        {\sc lmder} & 0.002 & 0.003 &  0.01 & 0.017 &0.046& 0.302 &0.982 \\
        {\sc dfo-ls} & 0.127 & 0.059 &  0.762 & 0.104 & 0.34& 5.404 &24.085\\
        \bottomrule
    \end{tabular}
    \caption{{\em Computing Time, Noiseless Case}. CPU time (in seconds) required by {\sc dfo-ls} and {\sc lmder} with forward differencing to satisfy \eqref{eq:term ls} for $\tau = 10^{-6}$ on the \texttt{BROWNAL} function, as $n$ and $m$ increase. }
    \label{tab:ls cpu time}
\end{table}

\subsection{Experiments on Noisy Functions}
Let us assume that the noise model is the same across all residual functions, i.e., 
$\epsilon_i(x)$ are i.i.d. for all $i$. As in Section~\ref{sub:noisy unc}, we generate noise independently of $x$, following a uniform distribution, $\epsilon_i(x) \sim \sigma_f U(-\sqrt{3}, \sqrt{3})$, with noise levels $\sigma_f \in \{10^{-1}, 10^{-3}, 10^{-5}, 10^{-7} \}$. 

Differencing will be performed more precisely than in the noiseless case. Following \cite{more2011estimating}, the forward-difference approximation of the Jacobian is defined as
\begin{align}\label{ls fd}
	[\hat{J}(x_k)]_{ij} = \frac{r_i(x + h_{ij}  e_j) - r_i(x)}{h_{ij}}, \ \text{ where } h_{ij} = 8^{1/4}\left(\frac{\sigma_{f}}{L_{ij}} \right)^{1/2} ,
\end{align}
where $L_{ij}$ is a bound on $\vert e_j^T \nabla ^2 \gamma_i(x) e_j \vert$ within the interval $[x, x + h_{ij}e_j]$. To estimate $L_{ij}$ for every pair $(i,j)$, and at every iteration, would be impractical, and normally unnecessary. Several strategies can be designed to provide useful information at an acceptable cost.  For concreteness, we estimate $L_{ij}$ once at the beginning of the run of {\sc lmder}. 

More concretely, we compute a different $h_{ij}$ for every residual function $\gamma_i$ and each coordinate direction $e_j$ at the starting point, and keep $h_{ij}$ constant throughout the run of {\sc lmder}. To do so, we compute the Lipschitz estimates $\{L_{ij}\}$ for $i= 1,...,m, j = 1,...,n$  by applying the Mor\'e-Wild (MW) procedure \cite{more2012estimating} described in the previous section to estimate the Lipschitz constant for function $\gamma_i$ along coordinate directions $e_j$. If the MW procedure fails, we set $L_{ij}=1$.
The cost of computing the $L_{ij}$, in terms of function evaluations, is accounted for in the results reported below.

%In analogy with Section~\ref{sub:noisy unc}, we consider each residual function $\gamma_i: \mathbb{R}^n \rightarrow \mathbb{R}$ separately and estimate the Lipschitz constants individually, and use them to define differencing intervals to approximate $\nabla \gamma_i$.
%Specifically, we employ the Mor\'e-Wild procedure \cite{more2012estimating} described in the previous section to estimate the Lipschitz constant for function $\gamma_i$ along coordinate directions $e_j$, and denote it as $L_{ij}$.  (Thus, in contrast with the noiseless case where we computed a different $h$ for each coordinate direction, we now compute a different $h$ for each residual function $\gamma_i$ and each coordinate direction $e_j$.)  The Mor\'e-Wild procedure is invoked only at the beginning of a run of {\sc lmder}, and the Lipschitz estimates $\{L_{ij}\}$ for $i= 1,...,m, j = 1,...,n$ are kept fixed. (Whenever the procedure fails, we set $L_{ij}=1$). Other settings for second-derivative estimation are discussed in Section~\ref{further}. The cost, in terms of function evaluations, is accounted for in the results reported below.
 
In {\sc dfo-ls}, we did not employ restarts, and set the {\tt obj\_has\_noise} option to its default value {\tt False},  which also changes some trust-region-related parameters. We did so for two reasons. First, restarts introduce randomness, and as Cartis et al. \cite{cartis2019improving} observed, can lead the algorithm to a different minimizer, making comparisons difficult. In addition, restarts are designed to allow the algorithm to make further progress as it reaches the noise level of the function.
% \ms{\sout{and a similar benefit would be obtained by switching from forward to central differences in {\sc lmder}, a strategy that we explore in a separate paper  \cite{shi2021automatic}.} [MS: Are we sure that we want to promise this for the next paper? I'm not sure if we will investigate least squares then.]}

\begin{figure}[htp] \label{fig:ls_noise_accuracy} %Noisy Case Morales Plot - Accuracy
    \centering
	\includegraphics[width=0.32\linewidth]{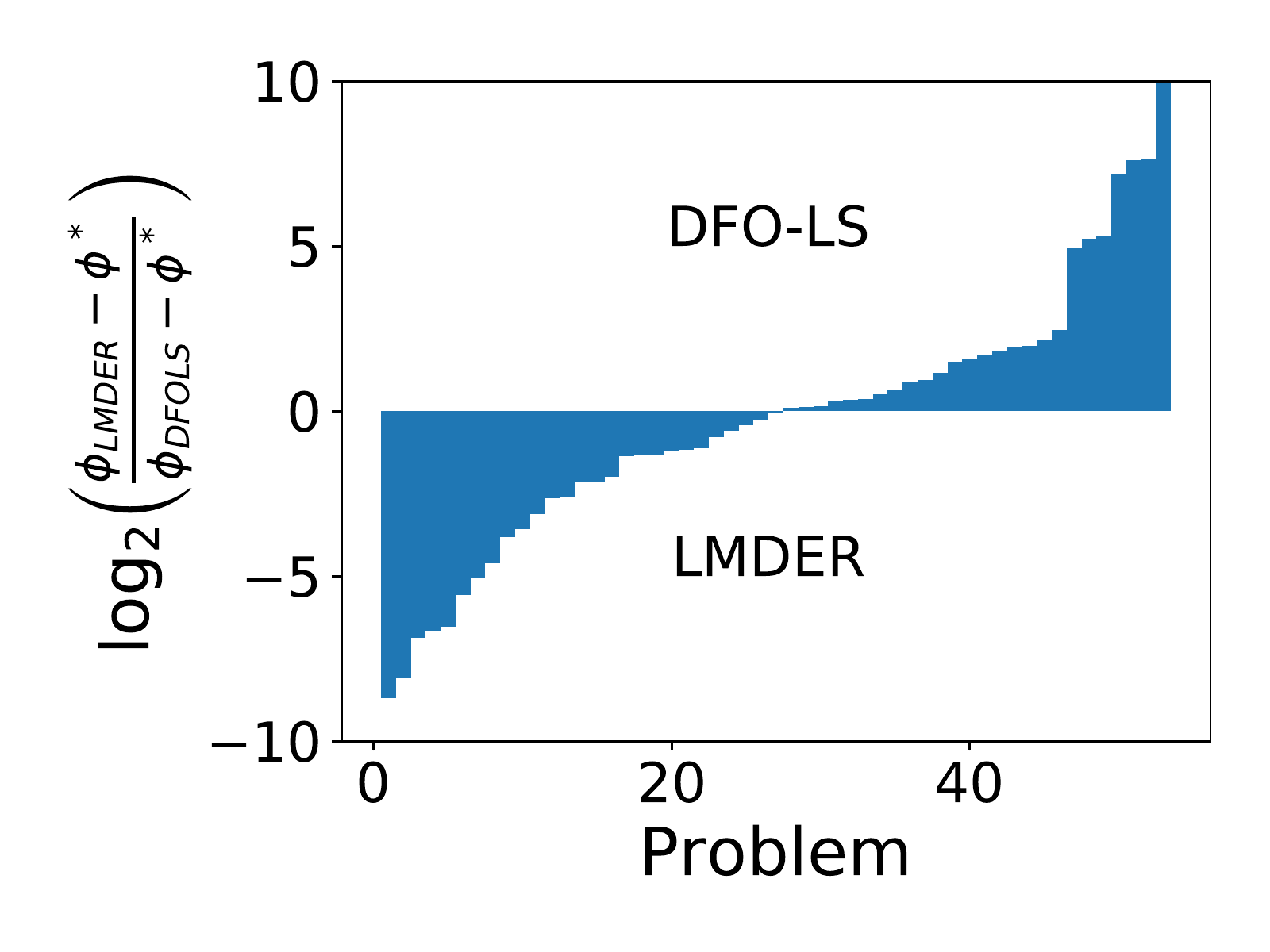}
	\includegraphics[width=0.32\linewidth]{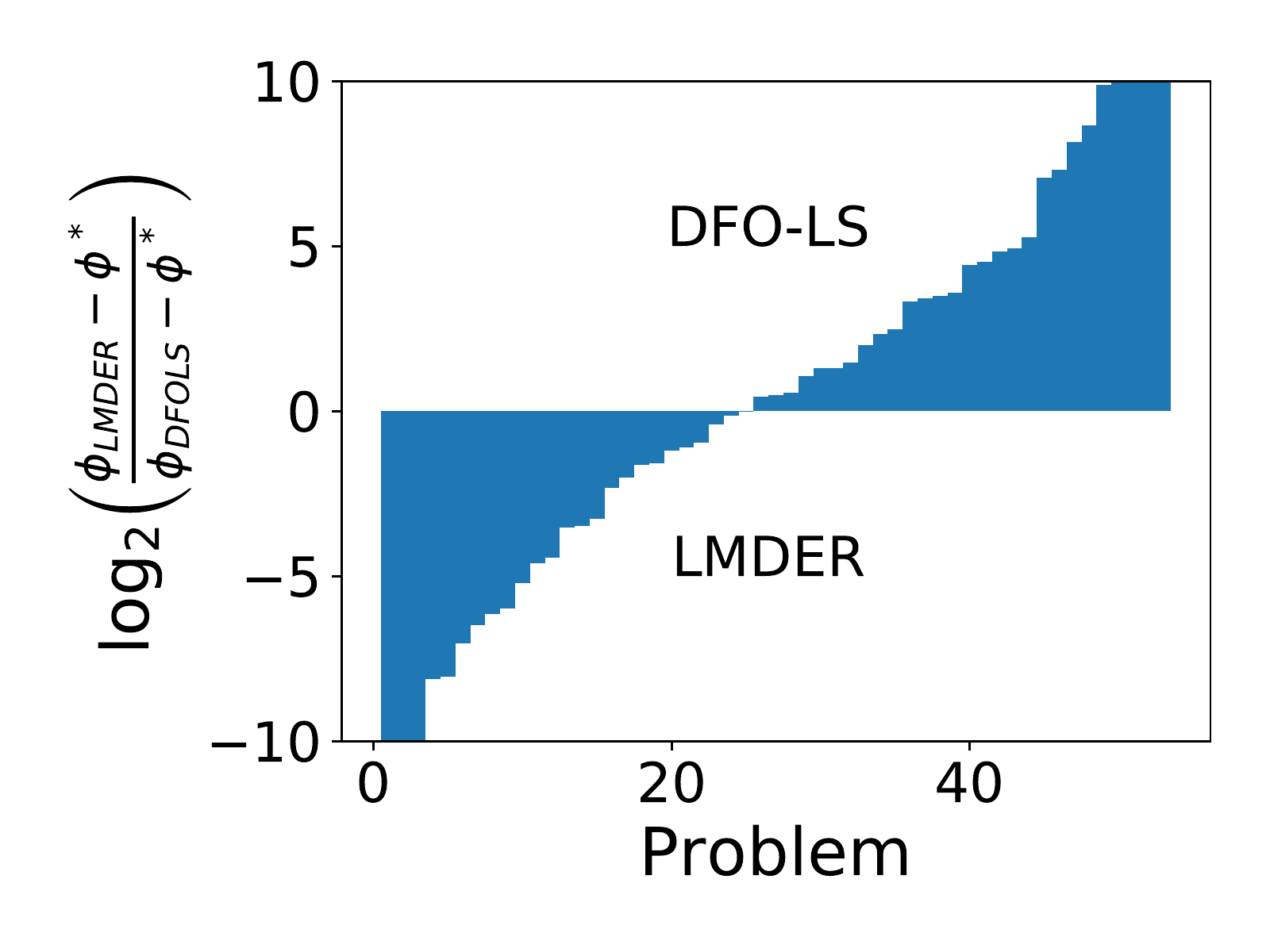}\\
	\includegraphics[width=0.32\linewidth]{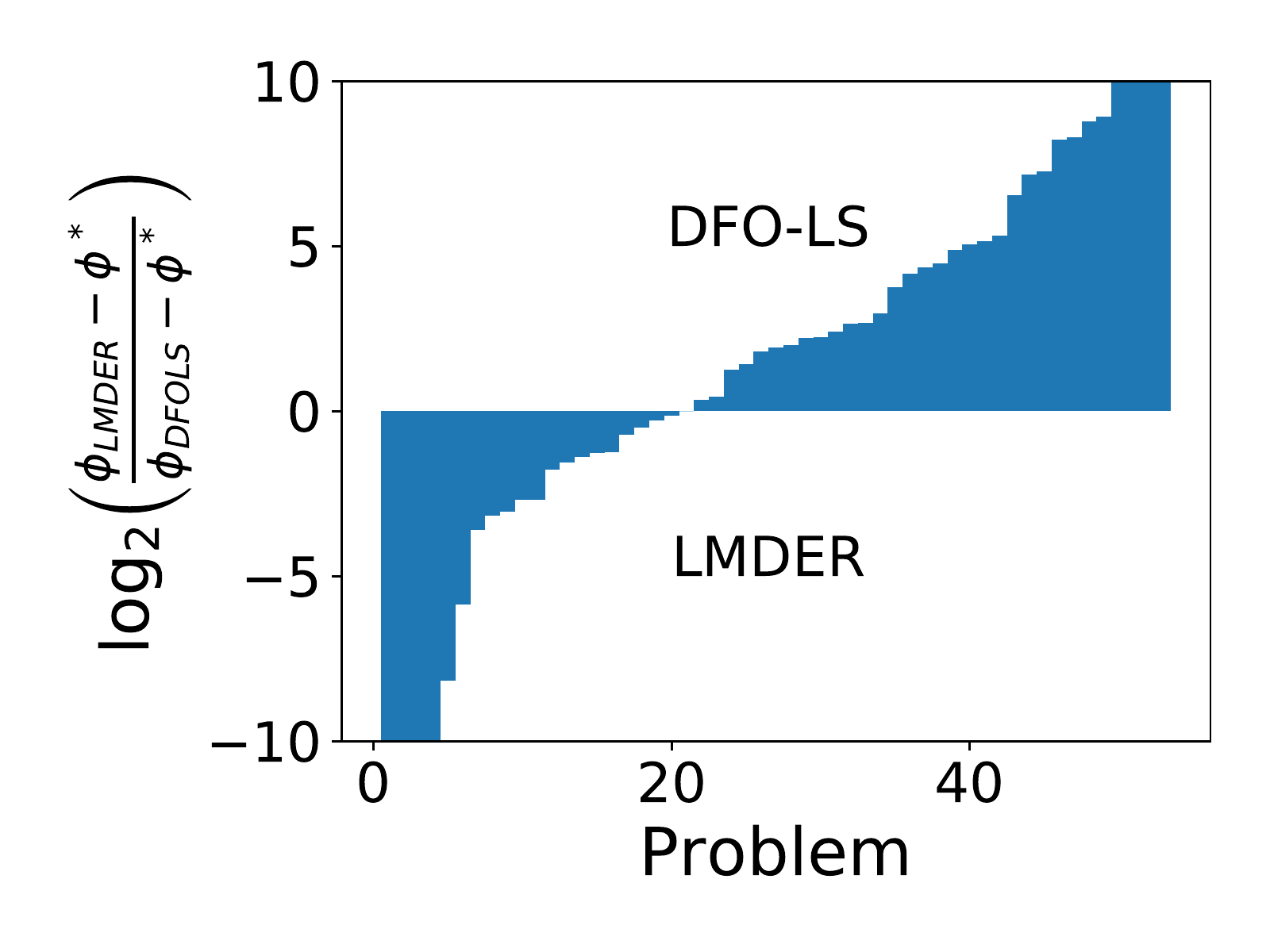}
	\includegraphics[width=0.32\linewidth]{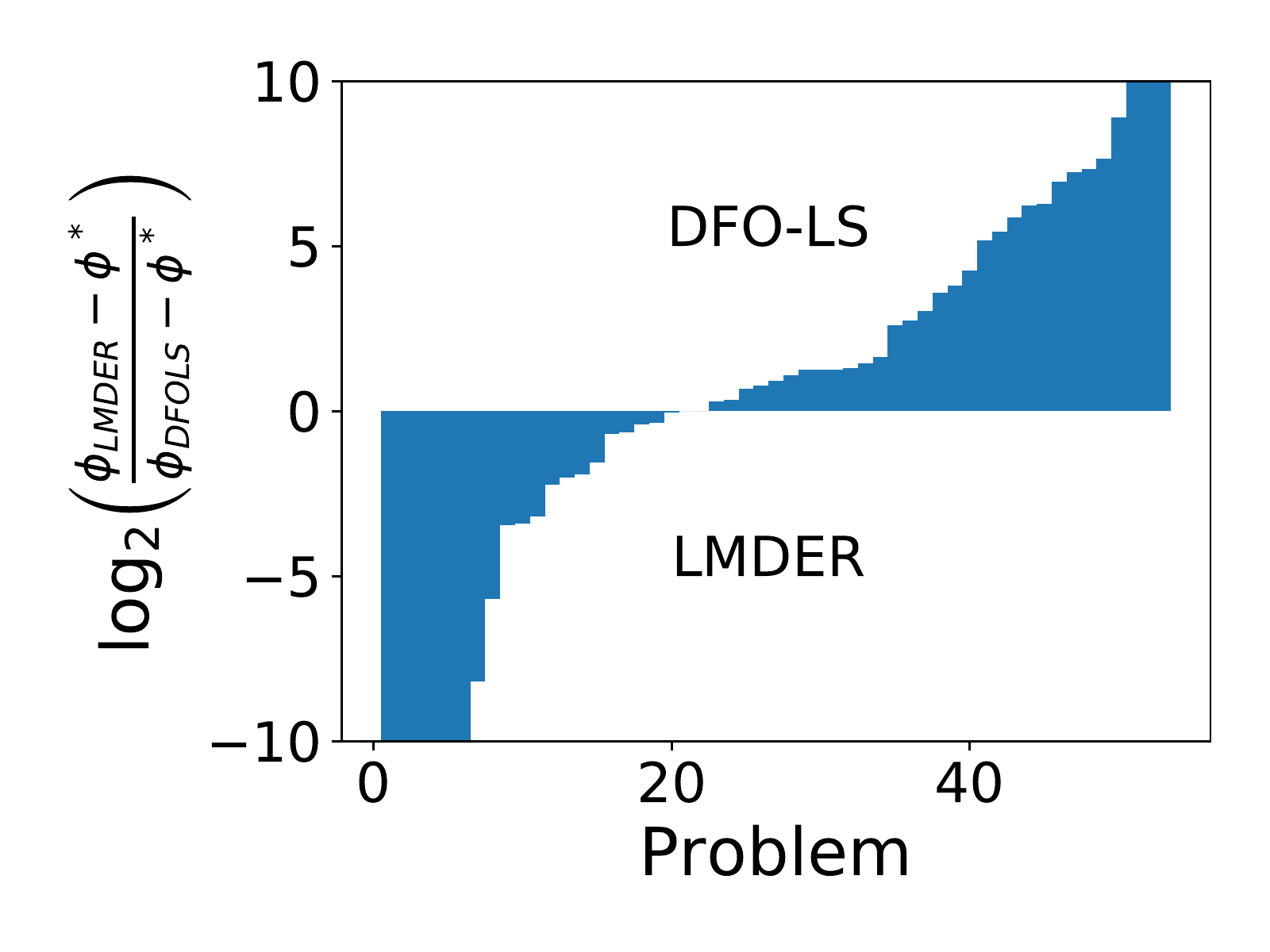}
    \caption{{\em Accuracy, Noisy Case.} Log-ratio optimality gap profiles  comparing {\sc dfo-ls} and {\sc lmder} for $\sigma_f = 10^{-1}$ (upper left), $10^{-3}$ (upper right), $10^{-5}$ (bottom left), $10^{-7}$ (bottom right). The bounds on the second derivatives are kept constant over the optimization process. }
%    \caption{{\em Accuracy Profiles; Nonlinear Least Squares with Noise.} Log-ratio profiles \eqref{eq: logratio accuracy noise} comparing {\sc dfo-ls} and {\sc lmder} for $\epsilon(x) = 10^{-1}(\text{upper left}),10^{-3}(\text{upper right}),10^{-5}(\text{bottom left}),10^{-7}(\text{bottom right})$. The bounds on the second derivatives are kept constant over the optimization process.}
    %
\end{figure}

\bigskip\noindent
{\em Accuracy.} 
We  compare the best optimality gap achieved by {\sc lmder} and {\sc dfo-ls}. We run the algorithms using their default parameters (except that we set {\tt model.abs\_tol = 0} for {\sc dfo-ls}) until no further progress could be made. In Figure~\ref{fig:ls_noise_accuracy}, we plot the log-ratio profiles
\begin{equation} \label{eq: logratio accuracy noise}
	\log_2 \left(\frac{\phi_{\text{\sc dfols}} - \phi^*}{\phi_{\text{\sc lmder}} - \phi^*}\right),
\end{equation}
for noise levels $\sigma_f \in \{10^{-1}, 10^{-3}, 10^{-5}, 10^{-7} \}$.

\begin{figure}[ht] %Noisy Case Morales Plot - Accuracy for theoretical Lipschitz estimation
    \centering
	\includegraphics[width=0.32\textwidth]{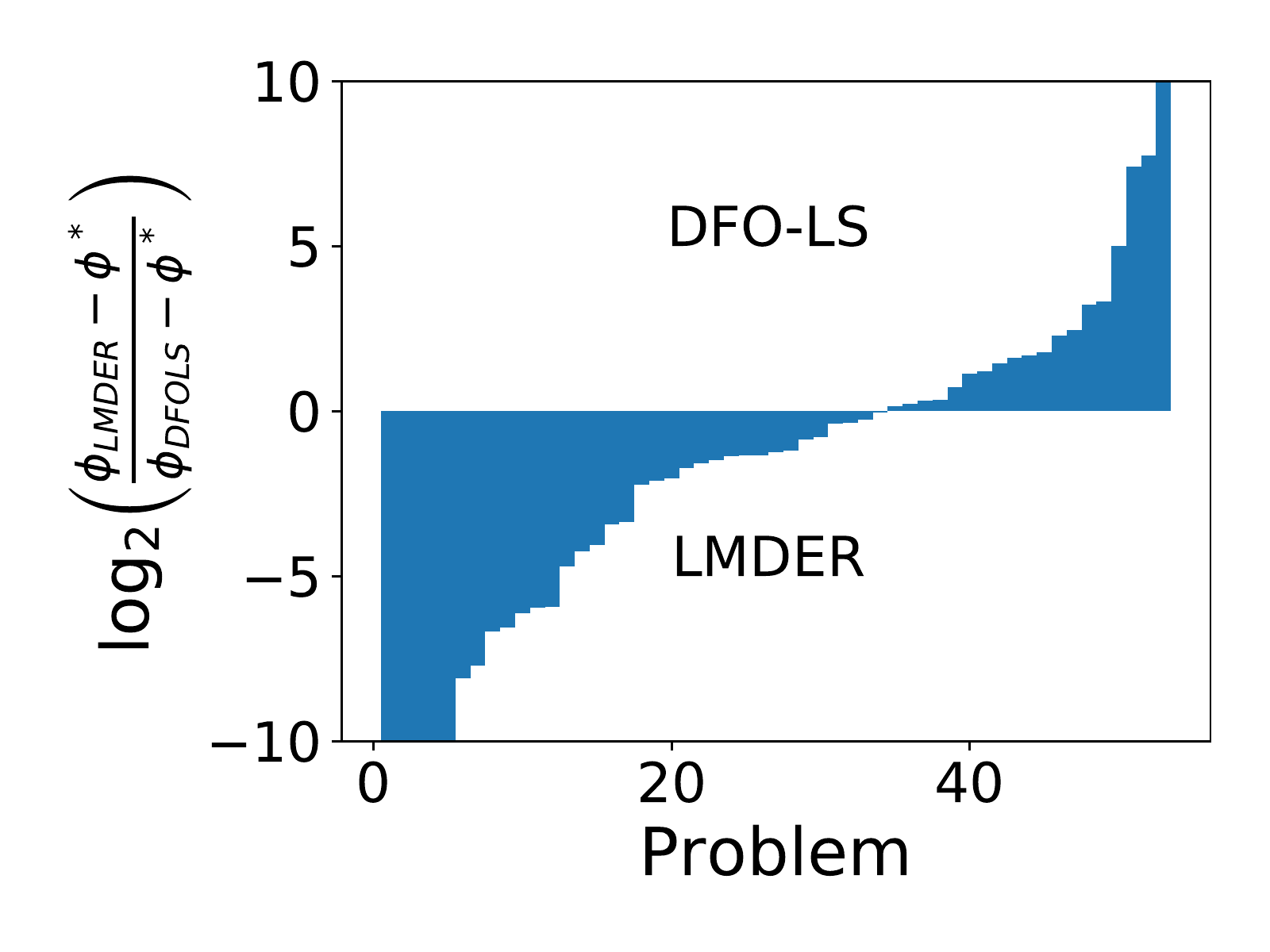}
	\includegraphics[width=0.32\textwidth]{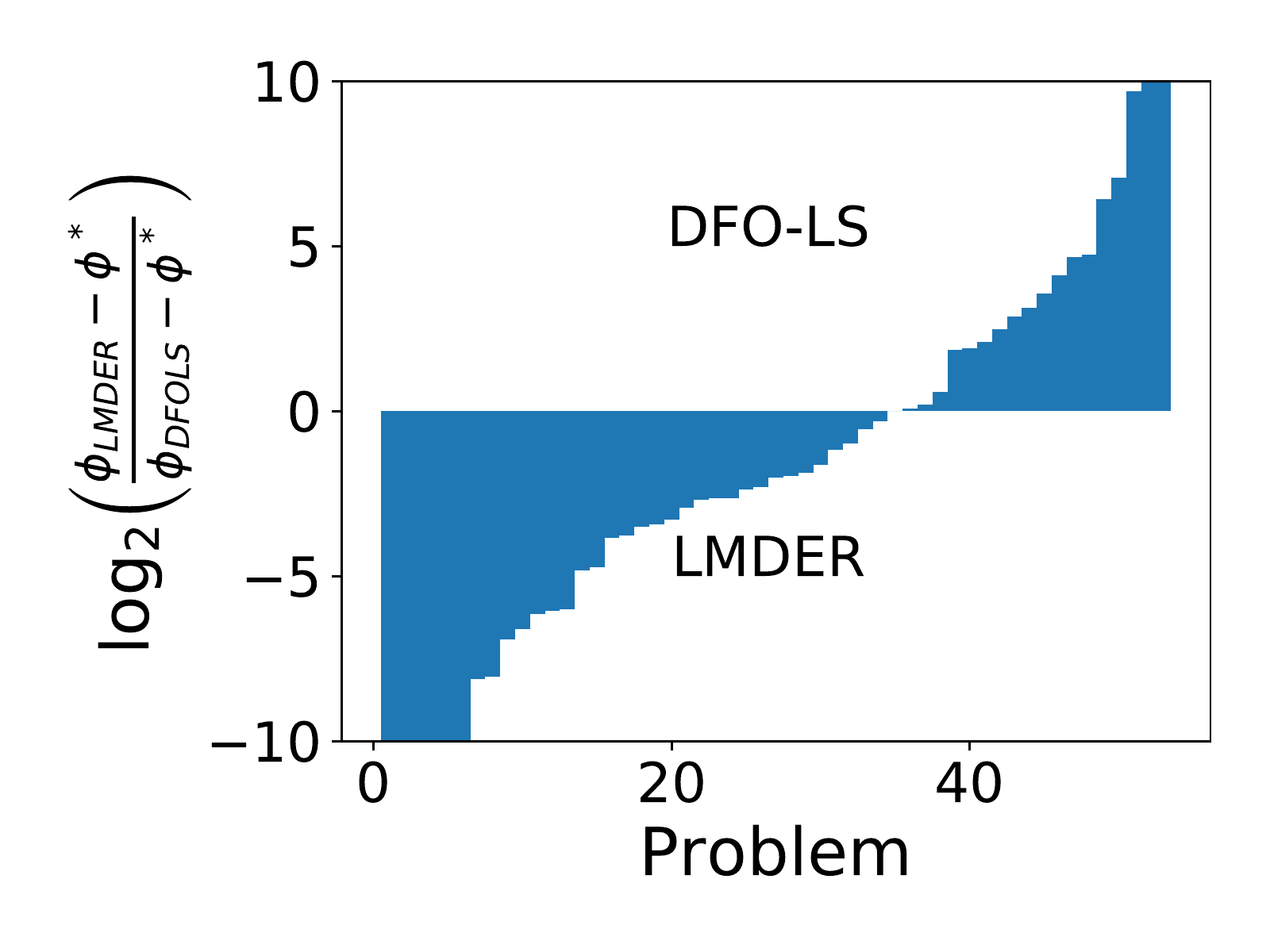}\\
	\includegraphics[width=0.32\textwidth]{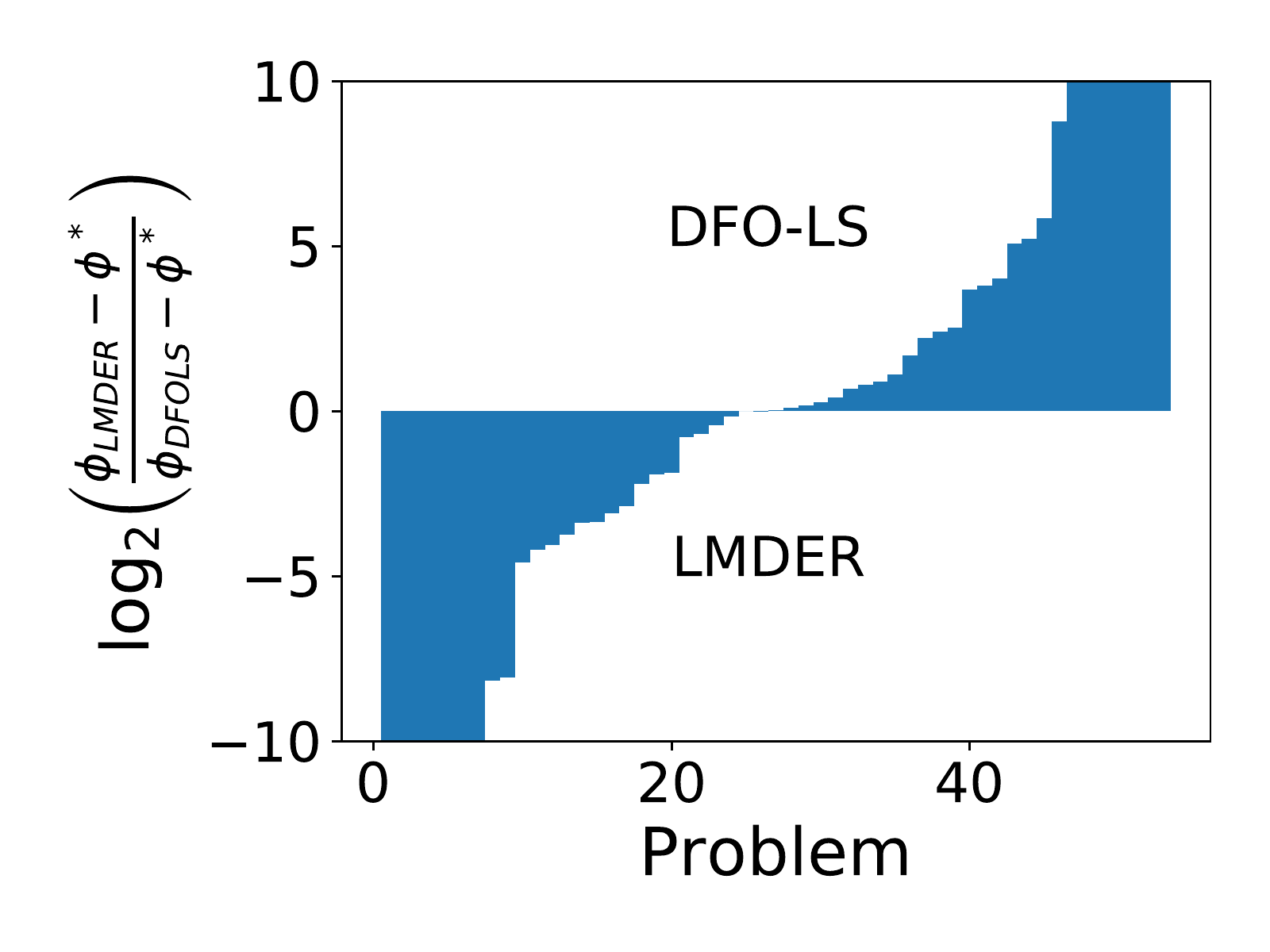}	
	\includegraphics[width=0.32\textwidth]{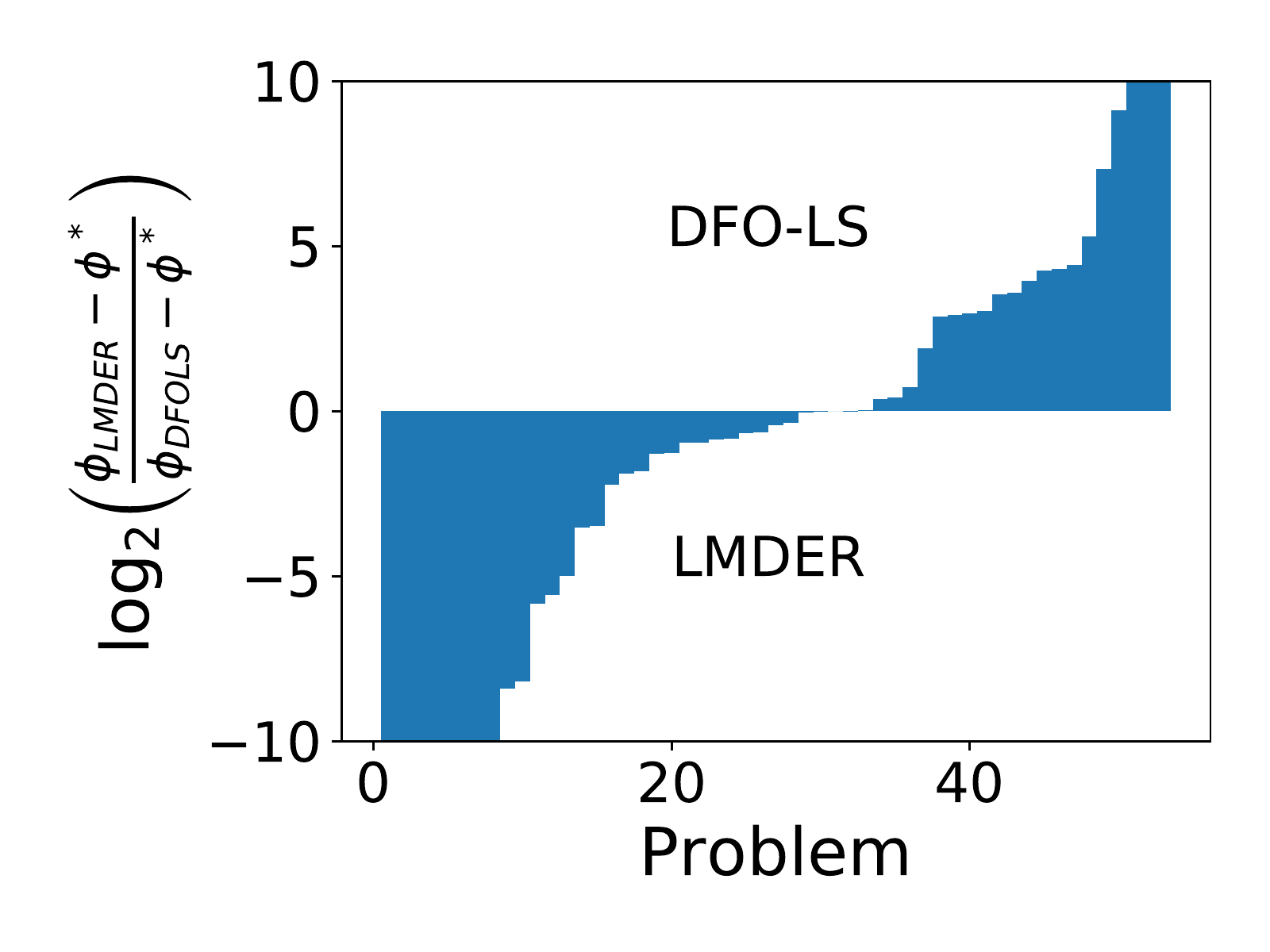}				
    \caption{ {\em Accuracy with Idealized Lipschitz Estimation, Noisy Case.} Log-ratio optimality gap profiles comparing {\sc dfo-ls} and {\sc lmder} for $\sigma_f = 10^{-1}$ (upper left), $10^{-3}$ (upper right), $10^{-5}$ (bottom left), $10^{-7}$ (bottom right). 
    }
    %The bound on the second derivative is reestimated at every iteration using$(\vert\gamma_i(x+he_j)+\gamma(x-he_j)-2\gamma(x)\vert)/h^2$ on the true function with differencing interval $h = \epsilon_M^{1/4}$} 
    \label{fig:ls_noise_accuracy_theory} 
\end{figure}

We observe that  {\sc dfo-ls} is  more accurate than {\sc lmder}, which points to the strengths of {\sc dfo-ls}, since it does not require knowledge of the noise level of the function in its internal logic.  
However, our implementation of {\sc lmder} is not sophisticated. As shown in the previous section, fixing the Lipschitz constant at the start of the finite-difference method is not always a good strategy, and one can ask whether a more sophisticated Lipschitz estimation strategy would  close the accuracy gap. To explore this question, we implemented an {\em idealized strategy} in which the Lipschitz estimation is performed as accurately as possible. For every $(i,j)$, $L_{ij}$ is computed at  every iteration, using  true function information, via the formula 
\begin{equation}
	L_{ij} = \frac{\vert \gamma_i(x+he_j)+\gamma_i(x-he_j)-2\gamma_i(x) \vert }{h^2}, \quad\mbox{with} \ h = (\epsilon_M)^{1/4}.
\end{equation}
We use this value in \eqref{ls fd}.  The results are given in Figure~\ref{fig:ls_noise_accuracy_theory}, and indicate that {\sc lmder} now achieves  higher accuracy than {\sc dfo-ls} on a majority of the problems across most noise levels. This highlights the importance of good Lipschitz estimates in the noisy settings, and suggests that the design of efficient strategies for doing so represent a promising research direction. 
%We expand on this topic in Section~\ref{further}.

\begin{figure}[htp]  %Noisy Case Morales Plot - Efficiency
    \centering
	\includegraphics[width=0.32\textwidth]{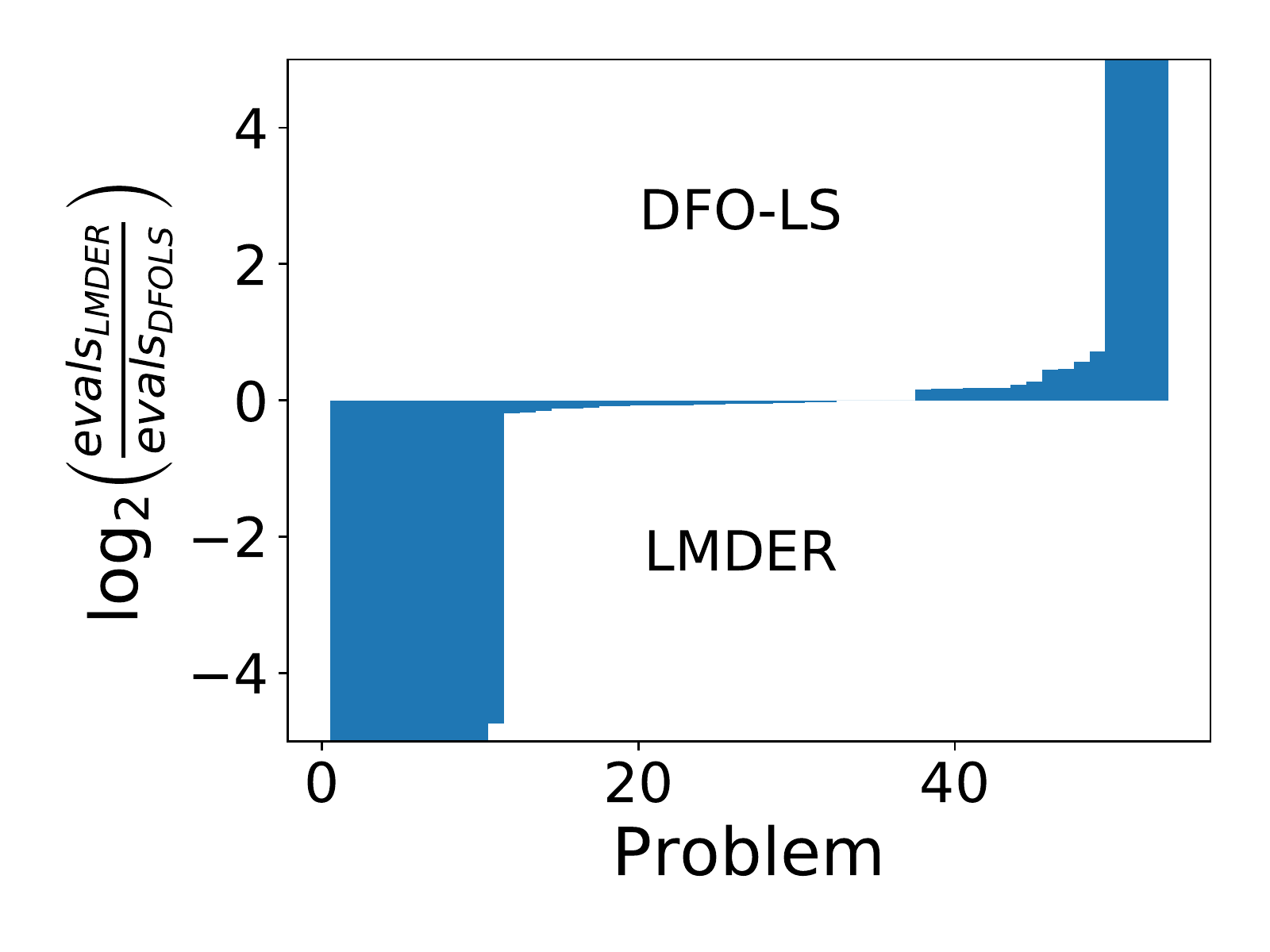}
	\includegraphics[width=0.32\textwidth]{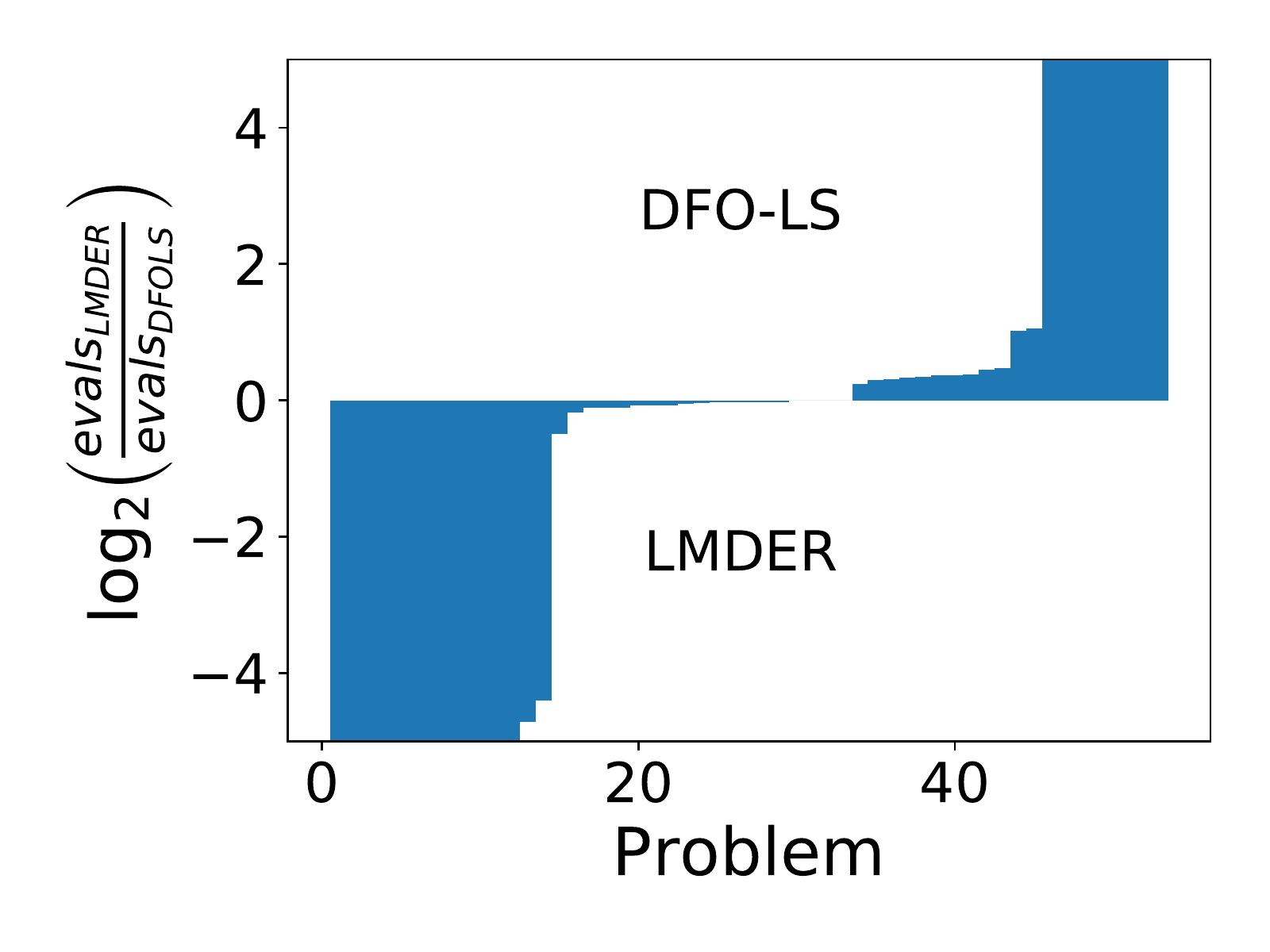} 	
	\includegraphics[width=0.32\textwidth]{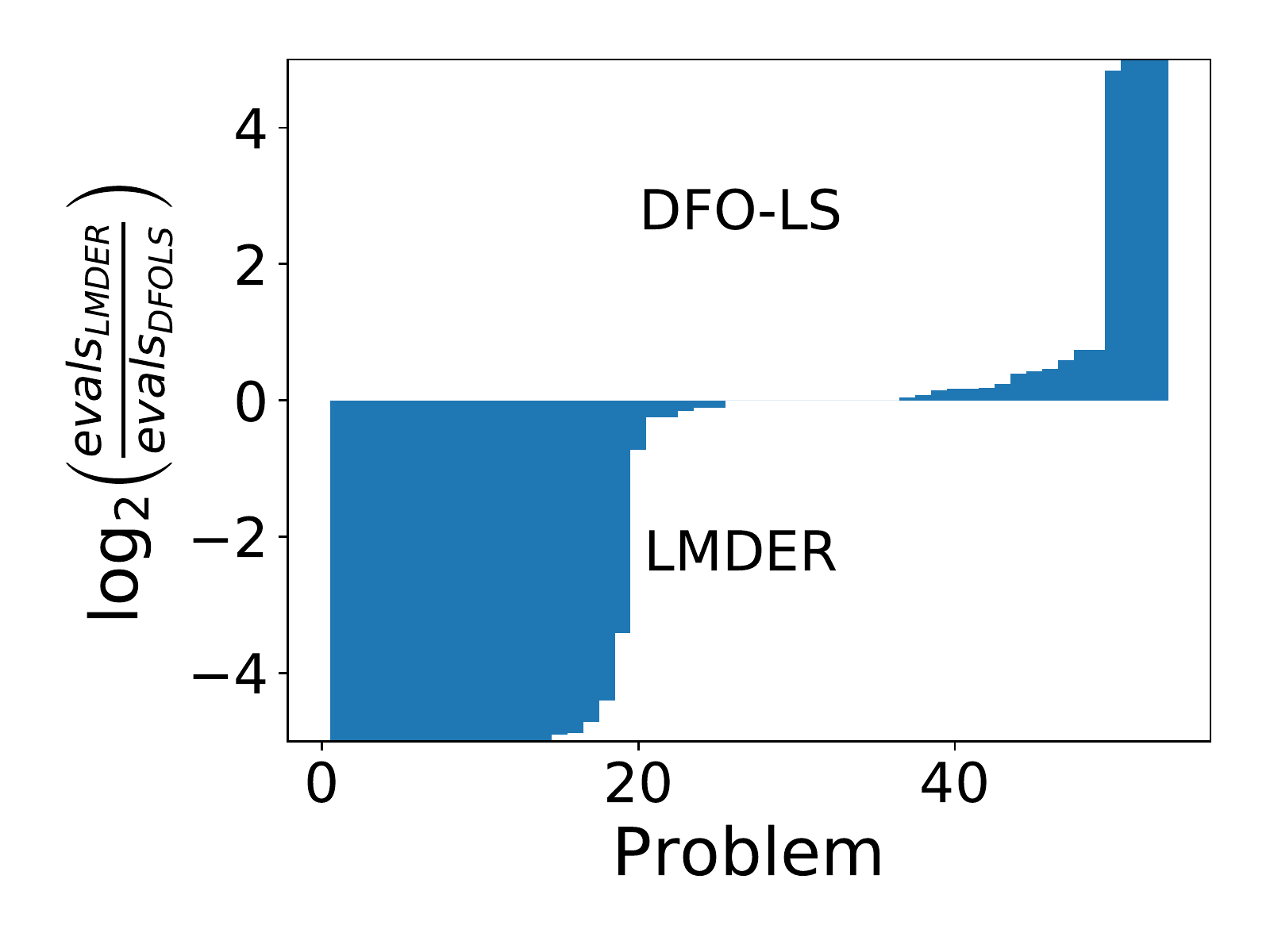} \\
	\includegraphics[width=0.32\textwidth]{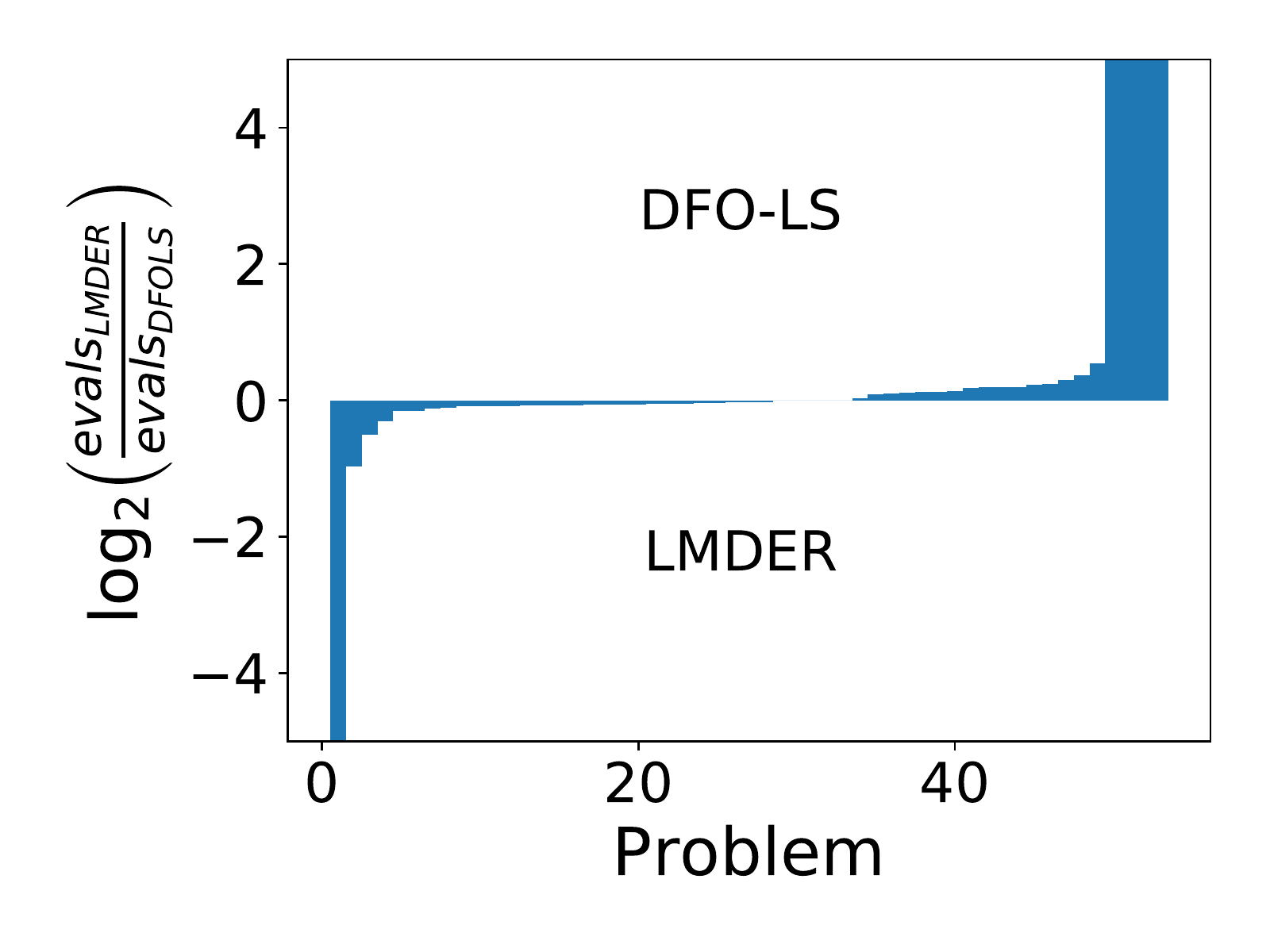}
	\includegraphics[width=0.32\textwidth]{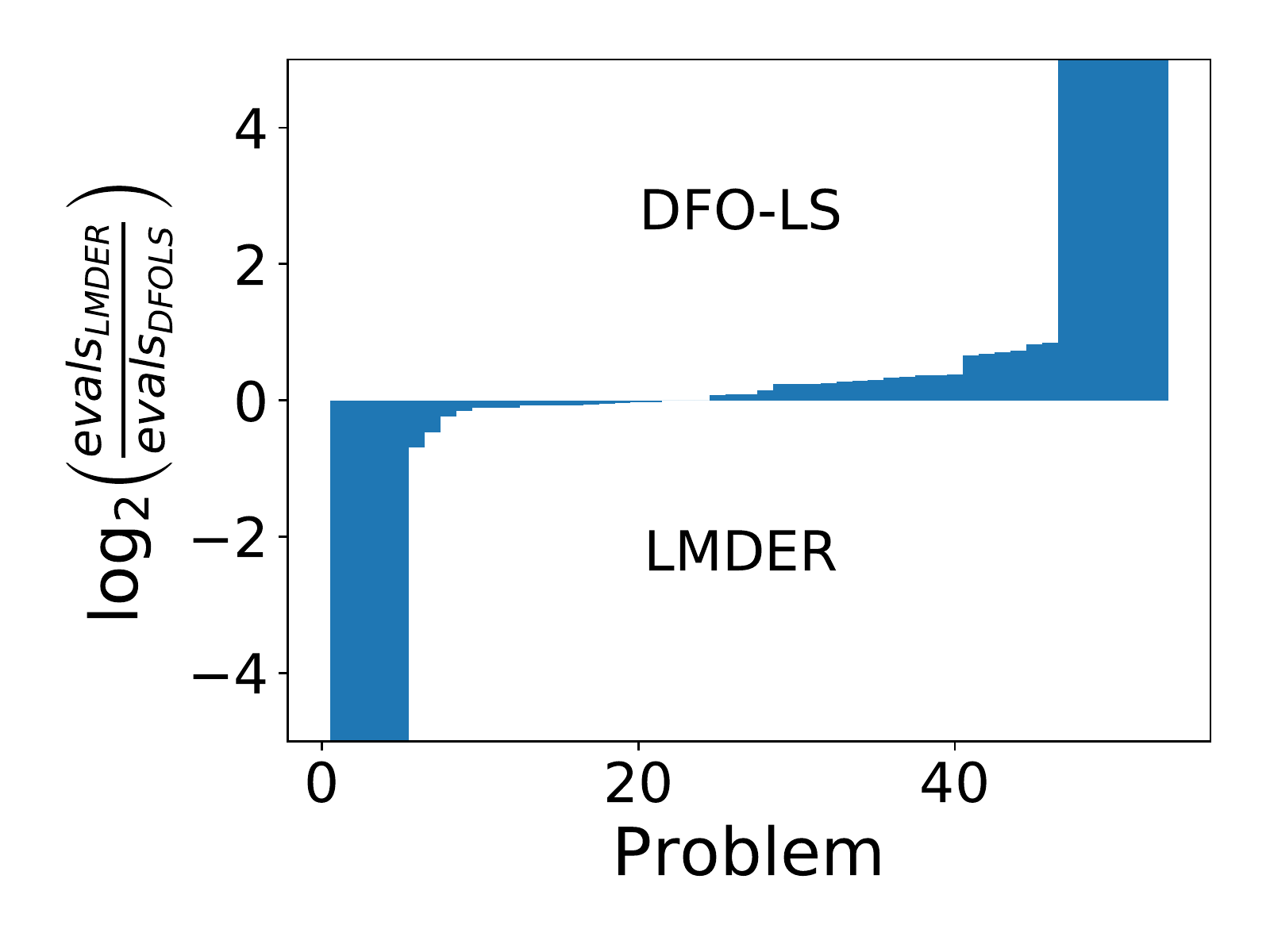}	
	\includegraphics[width=0.32\textwidth]{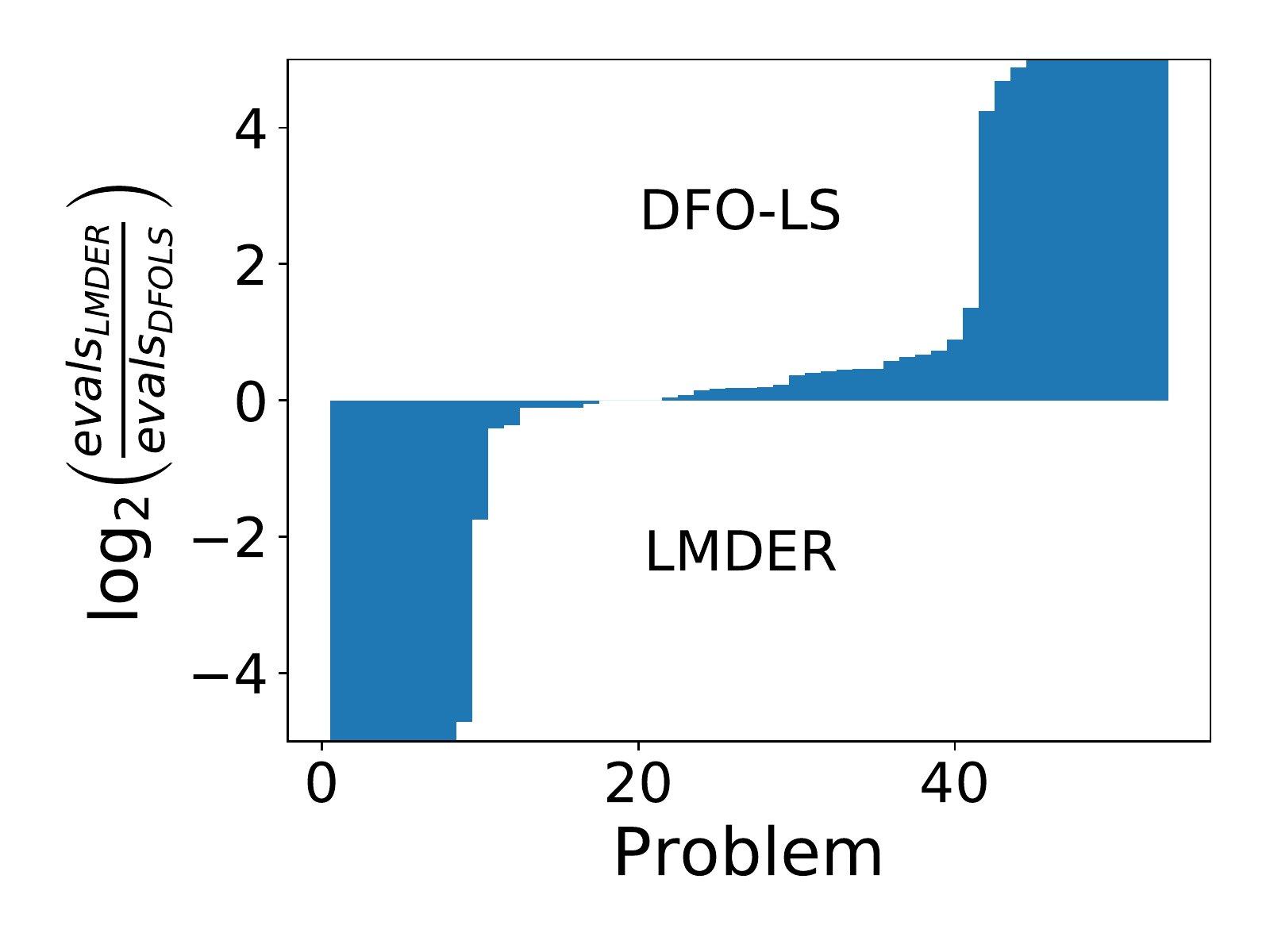}\\
	\includegraphics[width=0.32\textwidth]{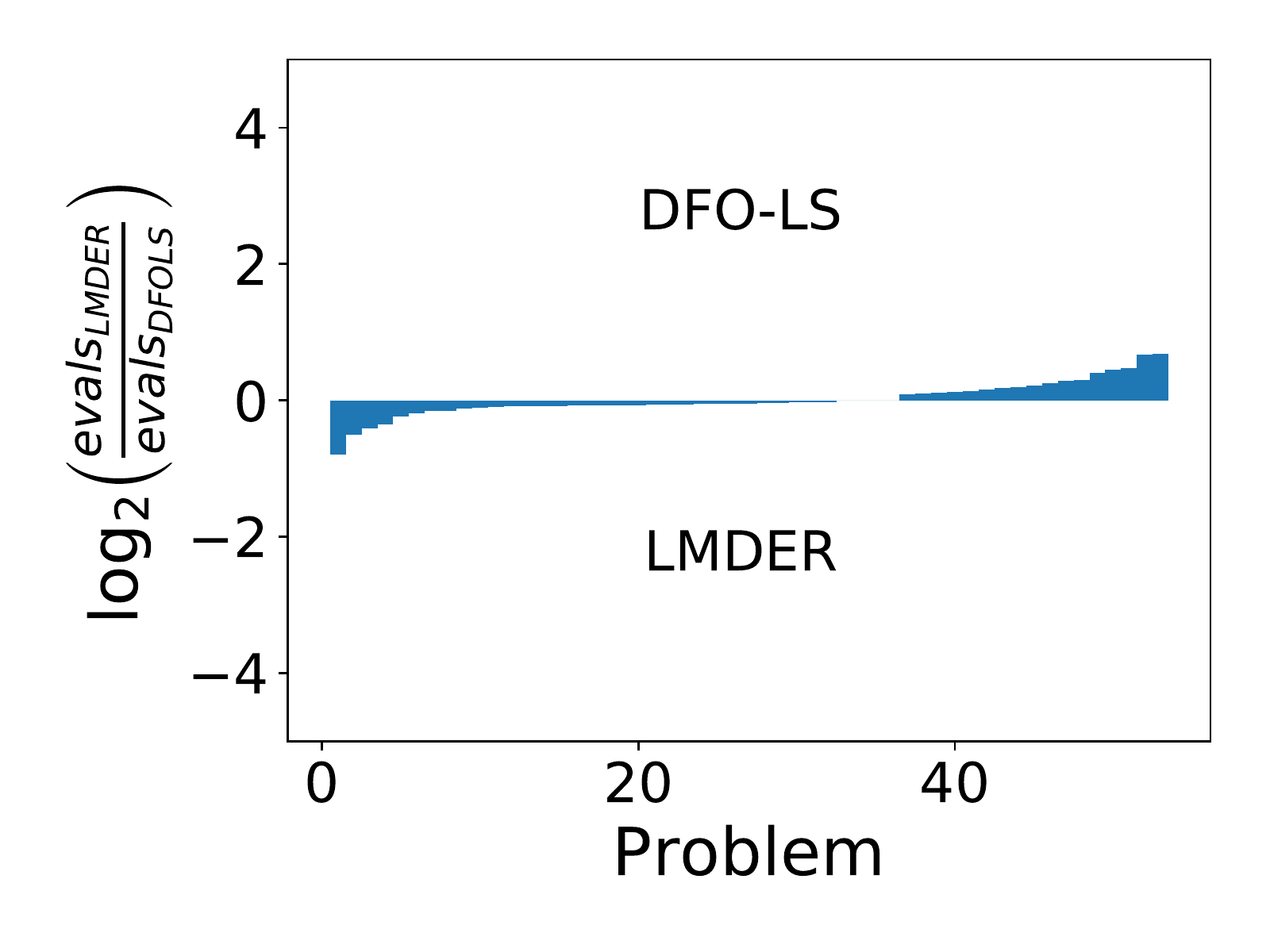}
	\includegraphics[width=0.32\textwidth]{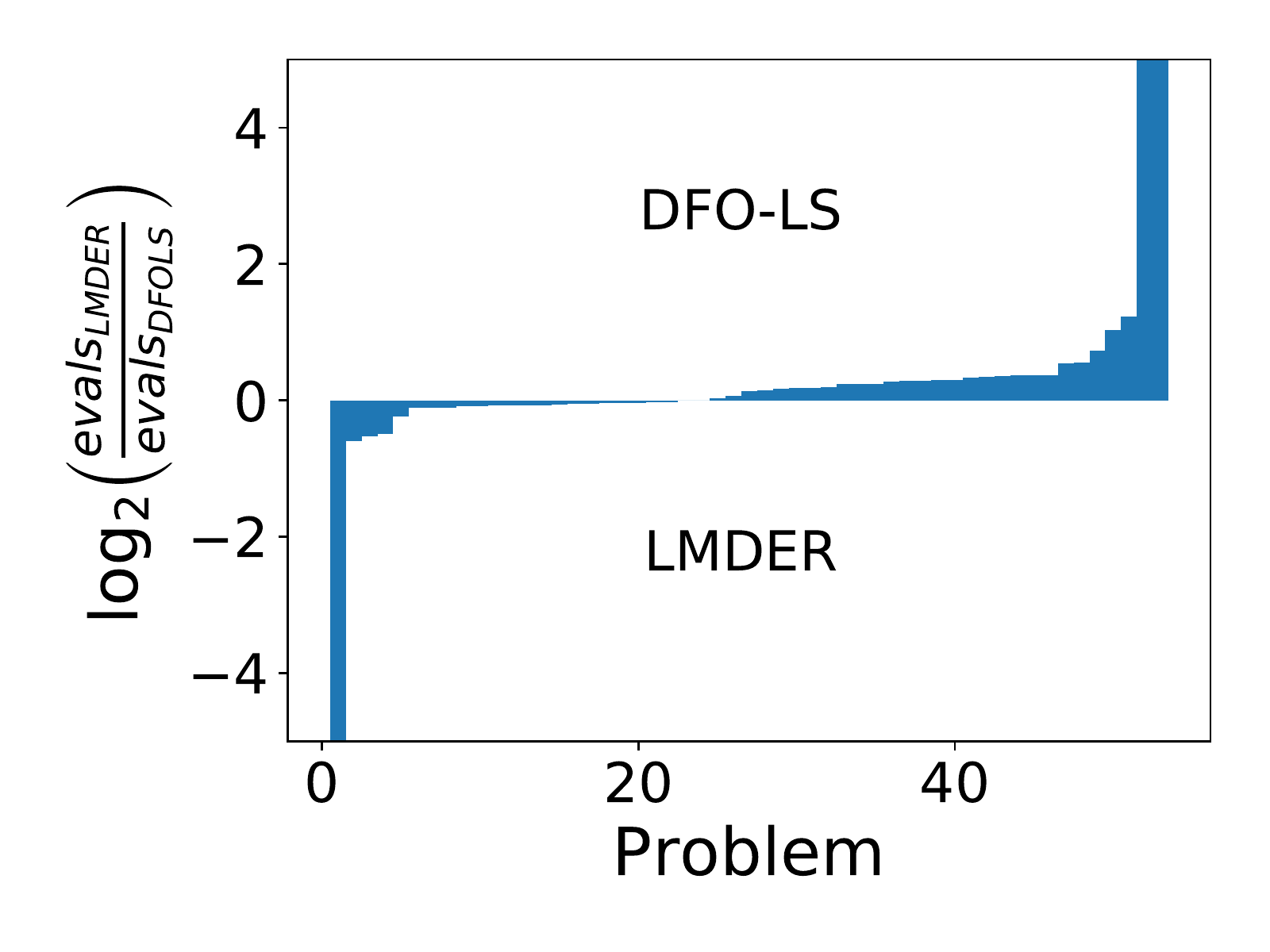}	
	\includegraphics[width=0.32\textwidth]{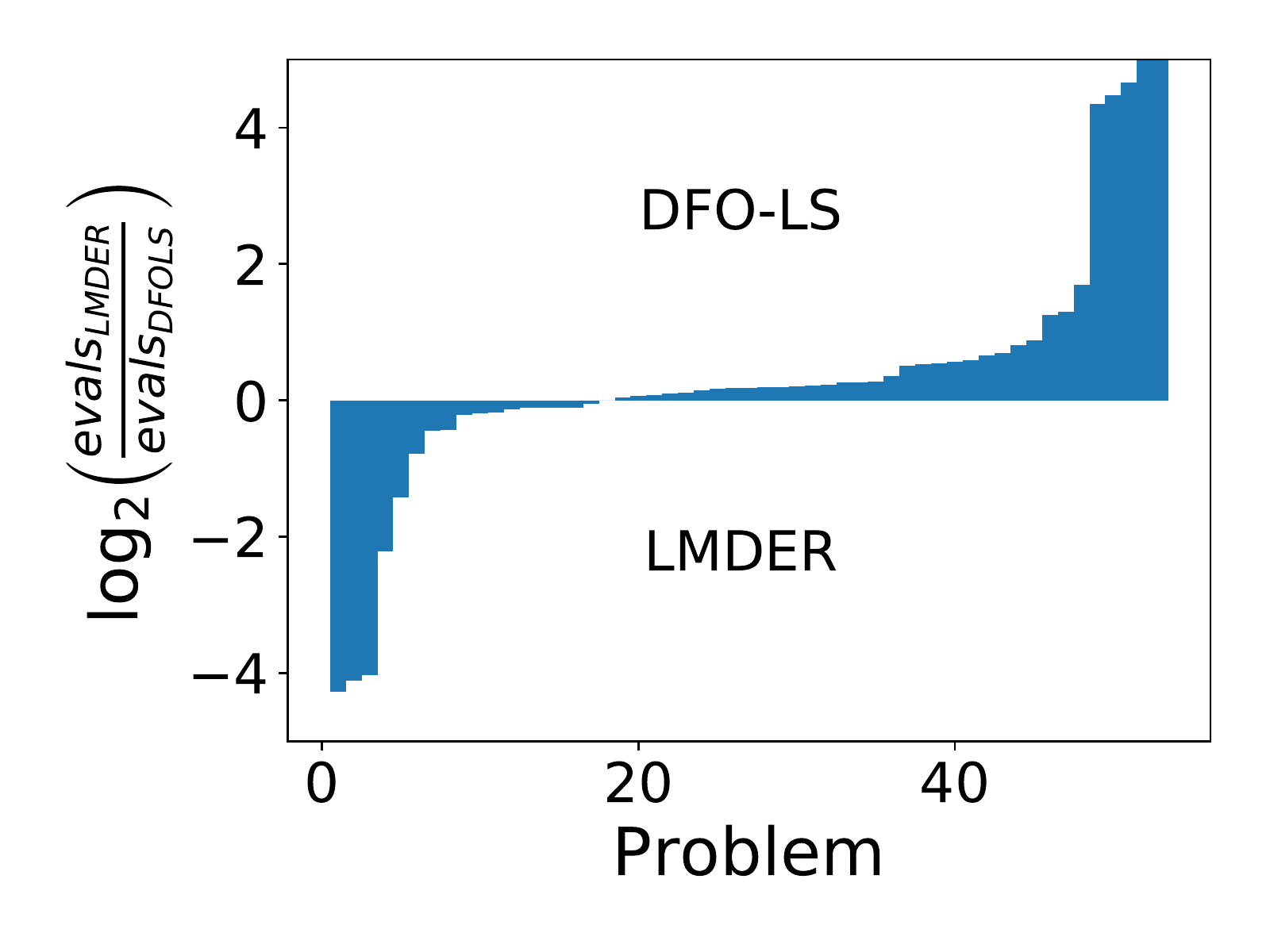}			
    \caption{{\em Efficiency, Noisy Case}. Log-ratio profiles comparing {\sc dfo-ls} and {\sc lmder} for $\sigma_f = 10^{-1}$ (top row), $10^{-3}$ (middle row), and $10^{-7}$ (bottom row). The figure measures the number of function evaluations to satisfy \eqref{eq:term ls} for $\tau = 10^{-1}$ (left column), $10^{-3}$ (middle column), and $10^{-6}$ (right column). Lipschitz constants were estimated only at the start of the {\sc lmder} run.}
	\label{fig:ls_noise_efficiency}
\end{figure}

\medskip\noindent{\em Efficiency.} 
To measure the efficiency of the algorithms in the noisy case, we record the number of function evaluations required to satisfy the termination condition
\begin{equation}\label{eq:term ls noise}
    \phi(x_k) - \tilde{\phi}^* \leq \tau (\phi(x_0) - \tilde{\phi}^*),
\end{equation}
where $\tilde{\phi}^*$, denotes the best objective value achieved by the two solvers, for a given noise level. To do so, both solvers were run until they could not make more progress.
The differencing interval in {\sc lmder} was computed as in the experiments measuring accuracy for the noisy setting, i.e., by employing only the Mor\'e-Wild procedure at the first iteration. 
In Figure~\ref{fig:ls_noise_efficiency}, we plot 
\begin{equation}  \label{skate}
	\log_2 \left(\frac{\text{evals}_\text{\sc dfols}}{\text{evals}_\text{\sc lmder}}\right)
\end{equation}
for three values of the tolerance parameter $\tau$ and for three levels of noise. We omit the plots for $\sigma_f = 10^{-5}$, which demonstrate similar performance. (When \eqref{eq:term ls noise} cannot be satisfied for a solver within the given budget of $500 \times n$ function evaluations, we set the corresponding value in \eqref{skate} to a very large number.) 

There is no clear winner among the two codes used in the experiments reported in Figure~\ref{fig:ls_noise_efficiency}. For low accuracy ($\tau=0.1$), {\sc lmder} appears to be more efficient, whereas the opposite is true for high accuracy ($\tau=10^{-6}$). We note again that {\sc dfo-ls} is able to handle different noise levels efficiently and reliably, without knowledge of the noise level or Lipschitz constants. On the other hand, the finite-difference approach is competitive even with a fairly coarse Lipschitz estimation procedure, and perhaps more importantly, it can be incorporated into existing codes (doing so in {\sc lmder} required little effort). 
In other words, in the finite-difference approach to derivative-free optimization, algorithms do not need to be constructed from scratch but can be built as adaptations of existing codes.

\subsubsection{Commentary}
The nonlinear least squares test problems employed in our experiments are generally not as difficult as some of the unconstrained optimization problems tested in the previous section. This should be taken into consideration when comparing the relative performance of methods in each setting.
 
  It is natural to question the necessity of estimating the Lipschitz constant for each component of each individual residual, as we did in our experiments, since this is affordable only if one can evaluate residual functions individually.
  %---and it requires the estimation of $m n$ derivatives. 
 One can envision problems for which a much simpler Lipschitz estimation suffices. For example, in data-fitting applications  all individual residual functions $\gamma_i$ may be similar in nature. In this setting, one could  use a single Lipschitz constant, say $L$, across all components and residual functions, especially when the variables are scaled prior to optimization. $L$ could be updated a few times in the course of the optimization process.
If the scale of the variables varies signficantly, one can compute Lipschitz constants $L_i$ for each component across all residual functions, requiring the estimation of $n$ Lipschitz constants.
%, and at the cost of $n + 1$ residual evaluations.

Another possibility is for the variables to be well scaled but the curvature of the residual functions $\gamma_i$ to vary signficantly.  In this case, one could estimate the root mean square of the absolute value of the second derivatives for each component; see Appendix \ref{app:lipschitz}. This approach notably only requires the estimation of a single constant for each residual function, yielding a total of $m$ Lipschitz constants. This permits the design of a more practical algorithm if the direct estimation of the root mean square is possible.

\section{Constrained Optimization}
\label{ch:constrained}

We now consider inequality-constrained nonlinear optimization problems of the form
\begin{equation} \label{consprob}
  \min_{x \in \mathbb{R}^n} \phi(x) \qquad \mbox{s.t.} \ \ \psi(x)\leq 0, \quad l \leq x \leq u,
\end{equation}
where $\psi: \mathbb{R}^n \rightarrow \mathbb{R}^m$ represents a set of $m$ linear or nonlinear constraints,  $l, u \in \mathbb{R}^n$, and $\phi$ and $\psi$ are twice continuously differentiable. We assume that the derivatives of $\phi$ and $\psi$ are not available, and more generally that we have access only to noisy function evaluations:
\begin{equation}  \label{connoise}
      f(x)= \phi(x) + \epsilon(x), \qquad c_i(x)= \psi_i(x) + \epsilon_i(x).
\end{equation} 
(We assume the same noise model for the objective $\epsilon$ and each of the constraint functions $\epsilon_i$, for simplicity.)
In this section, we compare the numerical performance of an established DFO method designed to solve problem \eqref{consprob} against a standard method for deterministic optimization that approximates the gradients of the objective function and constraints through finite differences.

\subsection{Choice of Solvers}
The solution of constrained optimization problems using derivative-free methods  has not been extensively studied in the literature; see the comprehensive review \cite{larson2019derivative}. There are few established codes for solving problem \eqref{consprob} and even fewer for problems that contain both equality and inequality constraints. 
To our knowledge, the best known software for solving problem \eqref{consprob} is {\sc cobyla}, developed by Powell \cite{powell1994direct}. The method implemented in that code constructs linear approximations to $\phi$ and $\psi$ at every iteration using function interpolation at points placed on a simplex in $\mathbb{R}^n$; it is designed to handle only inequality constraints. A more recent method by Powell, in the spirit of {\sc newuoa}, is {\sc lincoa} \cite{lincoa}. It implements an interpolation-based trust-region approach but it can handle only linear constraints. The {\sc nomad} package by Abramson et al. \cite{abramson2011nomad} is designed to solve general constrained optimization problems. It is a direct-search method that employs a progressive barrier approach to handle the constraints, and builds quadratic models (or other surrogate functions) as a guide.  
% In the author's  experience, PB is preferable to the other options on most, but not all, problems. One can try to solve constrained optimization problems by formulating a penalty and or augmented Lagrangian approach and employing unconstrained dfo methods, but these are unlikely to be competitive. 
In the numerical experiments reported in 
\cite{audet2018progressive},  {\sc cobyla} outperforms {\sc nomad}, and in \cite{augustin2014nowpac}, {\sc nomad} is not among the best performing methods. 
Based on these results, we chose {\sc cobyla} for our experiments.
%-PBTR presents experiments for (noiseless)comparison of PBTR with NOMAD and COBYLA: the results show that PBTR and COBYLA are comparable, and both outperform NOMAD as expected.
%-NOWPAC paper also includes (noiseless)comparisons between COBYLA and NOMAD: it seems that NOWPAC is more efficient than COBYLA in terms of function evaluations. In terms of accuracy, COBYLA and NOWPAC are competitive. The paper also includes experiments for noisy functions, but unfortunately no comparison is done between COBYLA and NOMAD. NOMAD is shown to be the least efficient.

There are many production-quality software packages for deterministic constrained optimization where we could implement the finite difference DFO approach. We chose {\sc knitro} because one of the algorithms it offers is a simple sequential quadratic programming (SQP) method that is close in spirit to {\sc cobyla}.  We did not to employ the interior-point methods offered by {\sc knitro}, which are known to be very powerful techniques for handling inequality constraints, because they may put {\sc cobyla} at an algorithmic disadvantage. In the same vein, we did not employ {\sc snopt} because it implements a sophisticated SQP method with many advanced features to improve efficiency and reliability. In short, we selected a simple nonlinear optimization method to more easily identify the strengths and weaknesses of the finite difference approach. 

 In summary, the codes tested are:

\begin{itemize}
\item {\tt COBYLA}. We ran the version of {\sc cobyla}  maintained in the  {\sc pdfo} package \cite{PDFO}.  We set the final trust region radius to $10^{-8}$ (\texttt{rhoend=1e-8}) to observe its asymptotic behavior, particularly in the noiseless case.  We ran {\sc pdfo} version 1.0, and called {\sc cobyla} via its Python interface ({Python 3.7.7}). 

\item {\tt  KNITRO}. We ran Artelys Knitro 12.2 with \texttt{alg=4} (an SQP algorithm), \texttt{gradopt=2} (forward differencing), and \texttt{hessopt=6} (L-BFGS). The choice of the finite difference interval $h$ is described below. In order to make the algorithm as close as possible to {\sc cobyla}, we set  the memory size of L-BFGS updating to its minimum value, $t=1$ (\texttt{lmsize=1}).  For consistency with {\sc cobyla},
we disabled the termination test based on the optimality error
by setting  \texttt{opttol=1e-16}, and \texttt{findiff\_terminate=0},  
and instead terminate when the computed step is less than $10^{-8}$ (by setting   \texttt{xtol=1e-8} and \texttt{xtol\_iters=1}). We called {\sc knitro} via its Python interface.
\end{itemize}

\subsection{Test Problems}

As in the previous sections, we used test problems from the \texttt{CUTEst} set \cite{gould2015cutest}, which were called through the Python interface, \texttt{PyCUTEst} version 1.0.  We recall from \eqref{consprob} that $n$ and $m$ refer to the number of variables and constraints, respectively (excluding bound constraints).  We first selected fixed-size problems that have at least one general nonlinear inequality and have no equality constraints, and for which $n\leq 100$ and $m\leq 100$.    
%We did not include problems with equality constraints to avoid any potential disadvantage for {\sc cobyla} because PDFO passes those to {\sc cobyla} as two inequalities. \fo{[I will verify that to make sure I recall correctly.]}  
The properties of the resulting 49 problems are listed in Table~\ref{tab:smallscaleprobs} in Appendix~\ref{app: cons property}.
We also selected 3 variable-size problems with at least one general nonlinear inequality and no equality constraints, and tested them with dimensions up to $n\leq 500$ and $m\leq500$. These problems are listed in Table~\ref{tab:varsizeprobs}.

%\item[C.] STYRENE, which simulates a styrene production process.  The objective is to maximize the net present value subject to several process and economic constraints. 

%\fo{Btw, I came across this very useful page by Sebastien Le Digabel (an earlier student of Prof.Audet)}
% https://www.gerad.ca/Sebastien.Le.Digabel/MTH8418/
%\fo{It includes links to DFO solvers and test problems; one of the problems is a constrained blackbox problem from production, which can be useful as an application-based testcase: https://github.com/bbopt/styrene}
%
%Problems with categorical variables? \jn{Try a few MINLP problems that Figen and Richard are familiar with and compare with NOMAD. [This could be a drag, so let's make sure that we are convinced that these tests will add to the paper]}

\bigskip

\subsection{Experiments on Noiseless Functions}   \label{cons-noiseless}

In the first set of experiments, we applied {\sc cobyla} and {\sc knitro} to solve the 49 small-scale, fixed-size CUTEst problems with exact function evaluations, i.e., with $\epsilon(x) = \epsilon_i(x) \equiv 0$ for all $i$. As in prior sections, the approximate gradient $g(x) \in \R^n$ of the objective function and Jacobian $\hat{J}(x) \in \R^{m \times n}$ of the constraints are evaluated by finite differencing:
\begin{align}
	[g(x)]_i & = \frac{f(x + h_i e_i) - f(x)}{h_i} \label{eq:approx g}\\
	[\hat{J}(x)]_{ij} & = \frac{c_i(x + h_j e_j) - c_i(x)}{h_j} \label{eq:approx J}
\end{align}
where
\begin{equation}
	h_i = \max\{1, |[x]_i|\} \sqrt{\epsilon_M}, \quad i=1,...,n.
\end{equation}
Both algorithms were stopped when the number of function evaluations exceed $500\max{(n,m)}$, or when the trust-region radius or steplength reaches its lower bound for {\sc cobyla} or {\sc knitro}, respectively.  Both solvers report  feasibility error as the max-norm of constraint violations.  

%The CPU time for KNITRO is reported by the solver; for COBYLA, we measure the CPU time using the \emph{time} package of Python.
%TODO: Measure the KNITRO CPU time with the same function.  FGN: I compared the two values on a few problems, they are identical.

Table~\ref{table:detA} reports the results of these tests.  A $^*$ indicates that the solver hits the limit of function evaluations.  We sorted the results in Table~\ref{table:detA} into four groups, separated by horizontal lines, according to the following outcomes:

(i) Both solvers converged to a feasible point with approximately the same objective function values.

(ii) The solvers converged to feasible points with different objective function values.

(iii) One of the solvers terminated at an infeasible point.

(iv) Both solvers terminated at infeasible points.

\medskip\noindent There were 32 problems associated with outcome (i). We can use them to safely compare the performance of the two solvers in terms of accuracy and efficiency.  (We comment on outcomes (ii)-(iv) in Appendix~\ref{app: cons det}.)

\begin{figure}[htp]
\centering
\includegraphics[width=0.32\textwidth]{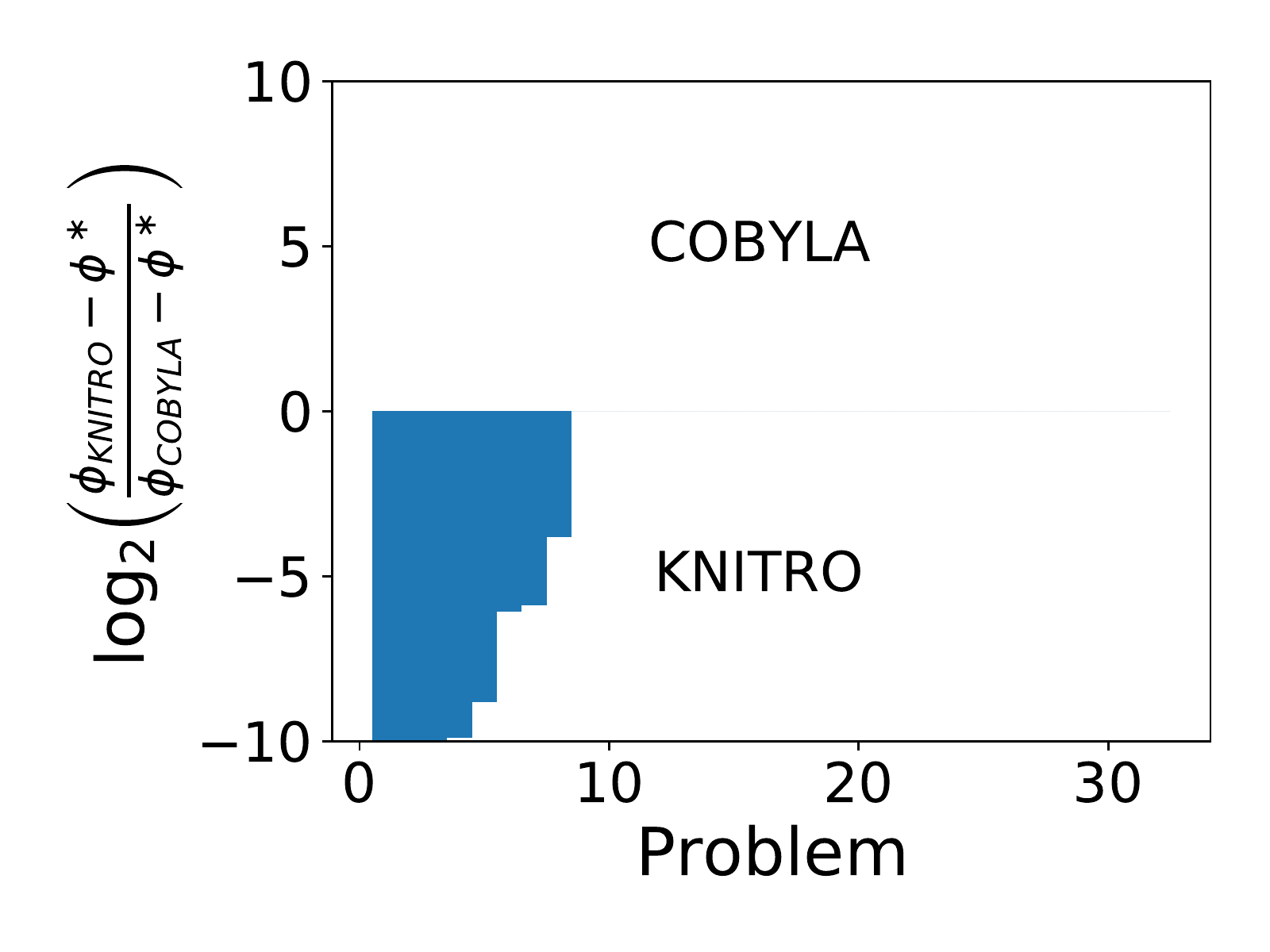}
\caption{{\em  Accuracy, Noiseless Case.} Log-ratio profiles comparing {\sc knitro} and {\sc cobyla} for $\epsilon(x) = \epsilon_i(x) = 0$. The figure plots the ratios \eqref{eight} for problems associated with outcome (i).}
\label{fig:accuracy_zero}
\end{figure}

Given that the two codes achieved feasibility in these 32 problems, we measure accuracy by comparing the best objective function value obtained by each solver with the value $\phi^\ast$ obtained by running {\sc knitro} with exact gradients until it could not make further progress.
%; we kept the remaining options the same so that it keeps running until the steplength becomes tiny.  
(In all the runs for determining $\phi^\ast$,  feasibility error was less than $10^{-10}$.) 
To measure accuracy, we report the ratios
\begin{equation} \label{eight}
\log_2\left(\frac{\max\{\phi_{\text{\tiny KNITRO}}-\phi^\ast,10^{-8}\}}{\max\{\phi_{\text{\tiny COBYLA}}-\phi^\ast,10^{-8}\}}\right).
\end{equation}
We compare accuracy up to eight digits because reporting, say, the ratio
$10^{-12}/10^{-15}$  would be misleading, given that the accuracy in the constraint violation could have the inverse ratio.
The results are presented in Figure~\ref{fig:accuracy_zero}, which shows that for most problems both solvers were able to achieve eight digits of accuracy; for the remaining problems, {\sc knitro} gave higher accuracy. 

%On the average, the number of function evaluations for COBYLA looks larger.  For a few problems with a large performance gap, we plot the progress of the two solvers as a function of the number of function evaluations to be able to make further comments.  A representative plot is provided in Figure~\ref{fig:progress}.

%\bigskip

%\footnotesize
%\begin{tabular}{r | l l l l | l l l l}
%\toprule
%\multicolumn{1}{c|}{ } &  \multicolumn{4}{c|}{KNITRO} & \multicolumn{4}{c}{COBYLA}\\
%\midrule
%problem	&	fval	&	evals	&	cputime	&	feaserr	&	fval	&	evals	&	cputime	&	feaserr	\\
%\midrule
%HS104	&	3.9512	&	718	&	0.282	&	0.00E+00	&	3.9512	&	3977	&	0.131	&	0.00E+00	\\
%HS92	&	1.3627	&	976	&	0.988	&	4.60E-17	&	1.3627	&	3000*	&	0.993	&	0.00E+00	\\
%HALDMADS	&	0.0001	&	113	&	0.073	&	4.44E-16	&	0.0001	&	669	&	0.341	&	0.00E+00	\\
%\bottomrule
%\end{tabular}
%\normalsize

%We observe that COBYLA achieves a close feasibility level and objective value to KNITRO within the same amount of function evalutions required by KNITRO up to its termination.  

\begin{figure}
	\centering
	\includegraphics[width=0.32\textwidth]{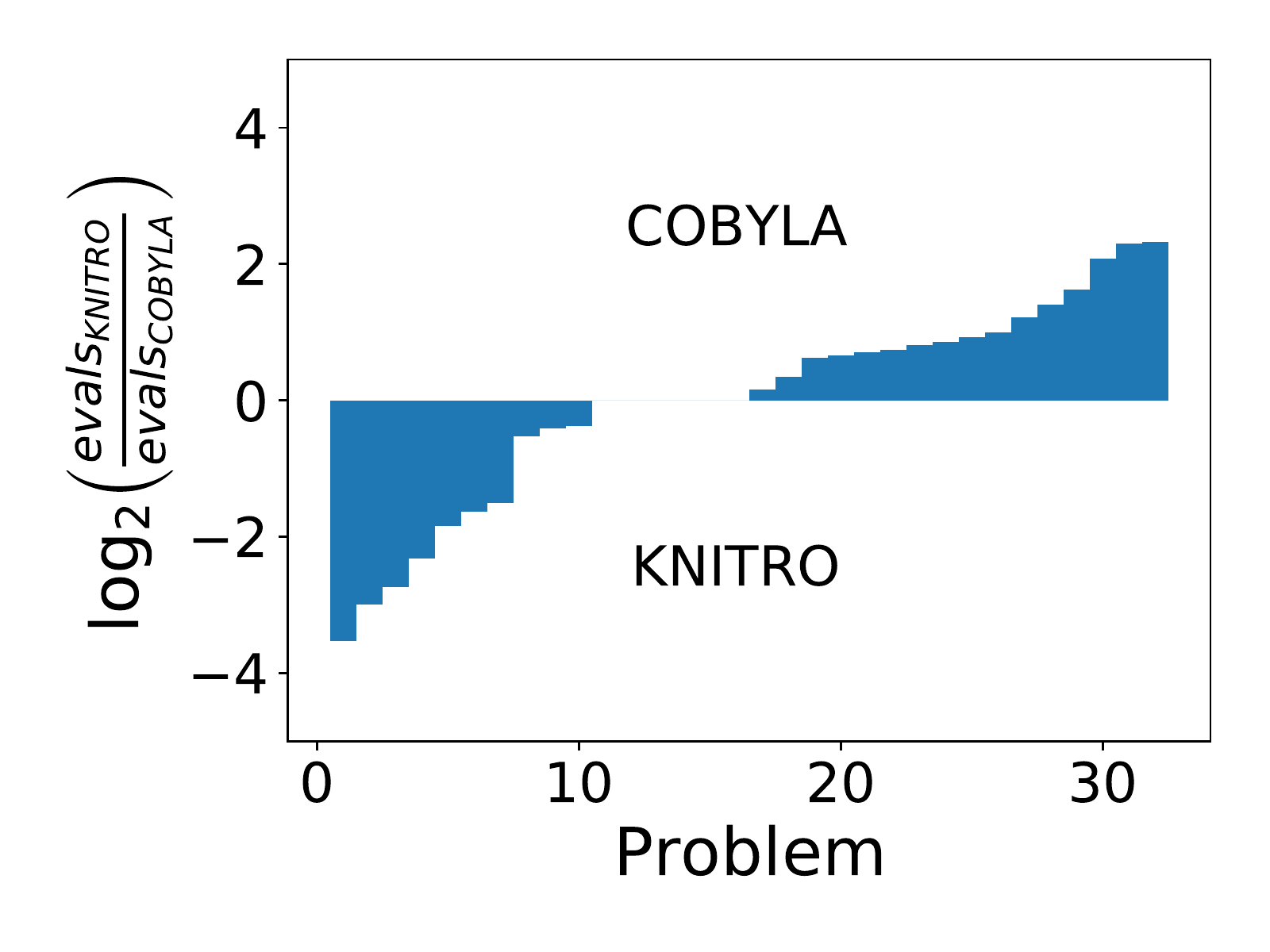}
	\includegraphics[width=0.32\textwidth]{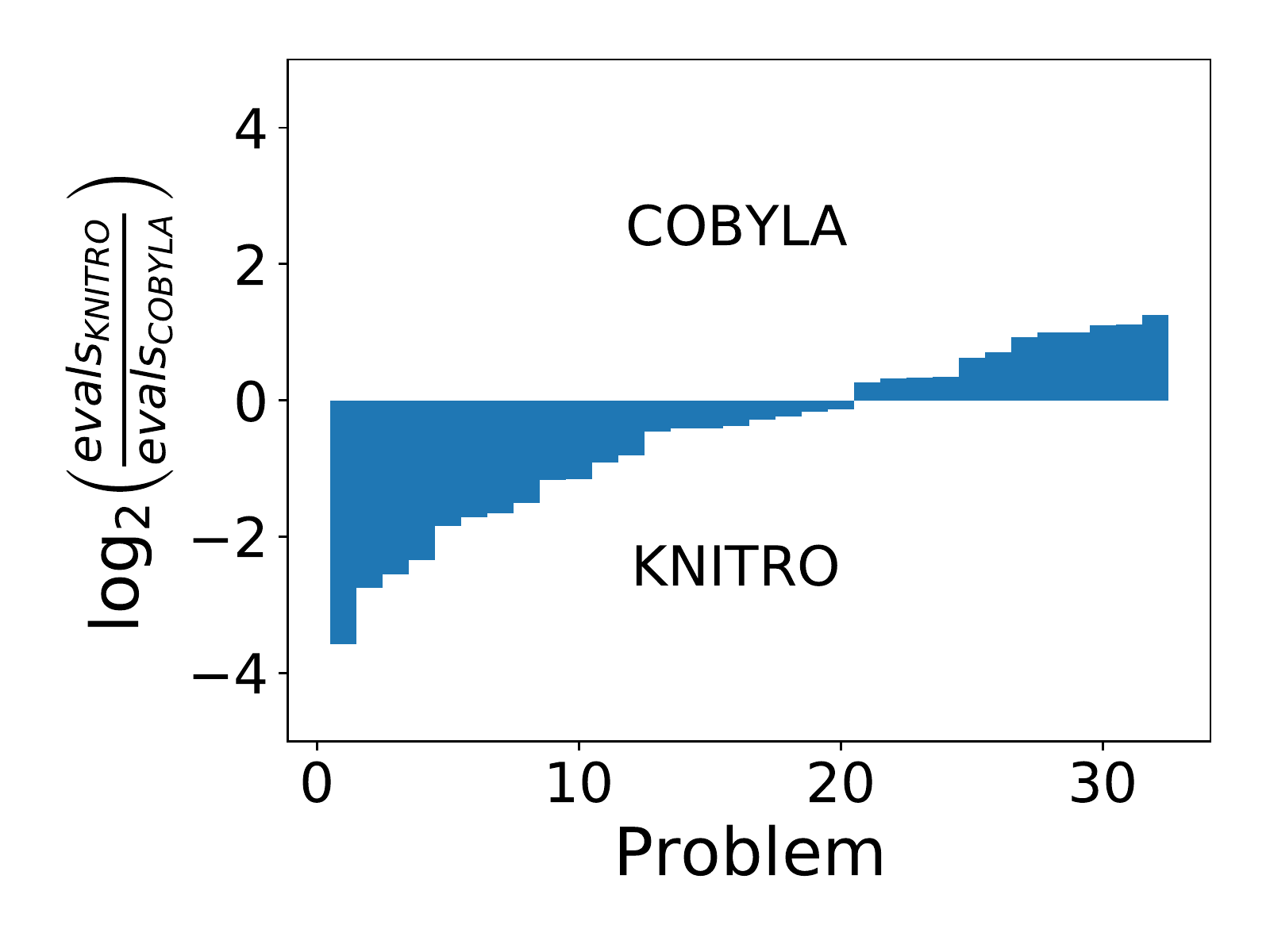}	
	\includegraphics[width=0.32\textwidth]{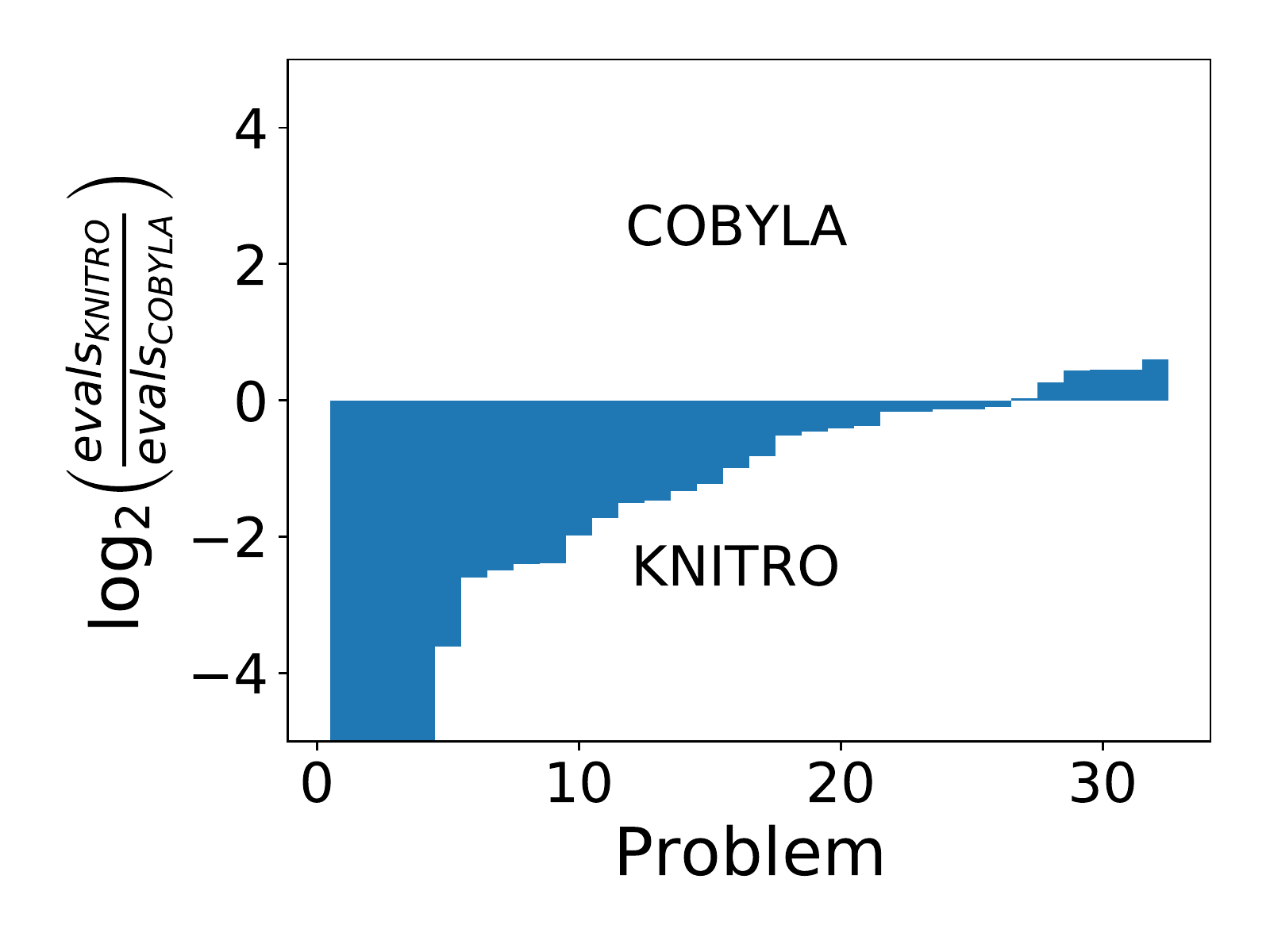}
    \caption{ {\em Efficiency, Noiseless Case.} Log-ratio profiles comparing {\sc knitro} and {\sc cobyla} for $\epsilon(x) = \epsilon_i(x) = 0$. The figures measure number of function evaluations to satisfy \eqref{eq:term} for $\tau = 10^{-1}$ (left), $10^{-3}$ (middle), and $10^{-6}$ (right) for outcome (i).} 
    \label{fig:efficiency_zero}
\end{figure}

To measure efficiency, we conducted a series of experiments using the same subset of 32 problems corresponding to outcome (i). We record the number of function evaluations,  $\text{evals}_{\text{\tiny KNITRO}}$ and $\text{evals}_{\text{\tiny COBYLA}}$, required by the two codes to satisfy the condition  
\begin{equation}
\label{eq:term}
\phi_k - \phi^\ast \leq \tau \cdot\max\{1, |\phi^\ast|\}
\end{equation}
for various values of $\tau$.  (If a solver fails to satisfy this test, we set the number of evaluations to a large value.)   Figure~\ref{fig:efficiency_zero} plots the ratios
\begin{equation} \label{plotthem}
\log_2\left(\frac{\text{evals}_{\text{\tiny KNITRO}}}{\text{evals}_{\text{\tiny COBYLA}}}\right)
\end{equation}
for $\tau \in \{10^{-1}, 10^{-3}, 10^{-6}\}$.  Table~\ref{table:detAtol} contains the complete set of results.  We observe that for low accuracy, {\sc cobyla} is slightly more efficient, while {\sc knitro} becomes significantly more efficient when high accuracy is required.

We should note that employing memory size $t= 1$ in L-BFGS updating yields a very weak quadratic model in the SQP method of {\sc knitro}. We experimented with a memory of size $t=10$ and observed that the performance of {\sc knitro} improved, particularly in the 
early iterations of the runs, but not dramatically;  see Figure~\ref{fig:efficiency_mem} in Appendix~\ref{app: cons det}.

%For HS67, it is likely that both solvers approach to the same solution but COBYLA runs out of its budget of evaluations before it can get sufficiently close.  For CRESC50 and CRESC4, the differences in the final function values are quite large; so, the solvers may not be converging to the same solution.  To check validity of all of these guesses, for each of the three problems we computed the optimality error (by using LS multipliers and exact gradeints) at the termination points of the two solvers.
%FGN: The least squares multipliers does not reflect the optimality error precisely, see the table below.
%\begin{center}
%\footnotesize
%\begin{tabular}{r | l l | l l }
%\toprule
%& \multicolumn{2}{|c|}{KNITRO} & \multicolumn{2}{c}{COBYLA}\\
%\midrule
%problem & optErrorLS & $\frac{\|x-x^\ast\|_\infty}{\max\{\|x^\ast\|_\infty, 1.0\}}$ & optErrorLS & $\frac{\|x-x^\ast\|_\infty}{\max\{\|x^\ast\|_\infty, 1.0\}}$ \\
%\midrule
%HS67 & 2.27E-3 & 1.31E-10 & 7.44E-2 & 1.89E-01\\
%CRESC4 & & & & \\
%CRESC50 & & & & \\ 
%\bottomrule
%\end{tabular}
%\normalsize
%\end{center}

\medskip\noindent
{\it Experiments on Variable-Dimension Problems without Noise.}
Next, we tested the two codes on the three problems of variable dimensions listed in Table~\ref{tab:varsizeprobs}, setting $\epsilon(x) = \epsilon_i(x) \equiv 0$. We impose a time limit of an hour; the runs exceeding the time limit are marked with a \emph{T}.  Table~\ref{cons cpu} displays the results. A clear picture emerges from these tests: {\sc knitro} is more reliable than {\sc cobyla} in terms of feasibility, tends to achieve a lower objective value, and can solve larger problems within the allotted time.  
%For problems  SVANBERG and READING with $N=500$, {\sc cobyla} is unable to trigger any of the stop tests within the one-hour time limit. 
We observed earlier that there is a perception in the DFO literature that finite differences require too many function evaluations, especially as the dimension of the problem increases. The results in Table \ref{cons cpu} do not support that concern.

  %Within the same time limit, KNITRO can solve ...-dimensional instances.  If we change to the interior point algorithm, it can even solve ... dimensional instances (TODO).

\bigskip
\begin{table}[htp]
\footnotesize
\begin{tabular}{r | l l l l | l l l l}
\toprule
\multicolumn{1}{c|}{ } &  \multicolumn{4}{c|}{KNITRO} & \multicolumn{4}{c}{COBYLA}\\
\midrule
problem	&	$\phi_k$	&	\#evals	&	CPU time	&	feaserr	&	$\phi_k$	&	\#evals	&	CPU time	&	feaserr	\\
\midrule
SVANBERGN10	&	15.7315	&	168	&	0.078	&	2.62E-14	&	26.0000	&	302	&	0.349	&	0.00E+00	\\
SVANBERGN50	&	82.5819	&	950	&	0.538	&	4.11E-15	&	136.0000	&	2140	&	0.941	&	0.00E+00	\\
SVANBERGN100	&	166.1972	&	1851	&	1.564	&	3.22E-15	&	273.5000	&	4715	&	9.754	&	0.00E+00	\\
SVANBERGN500	&	835.1869	&	10054	&	76.771	&	1.50E-14	&	-	&	-	&	T	&	-	\\
\midrule
READING4N2	&	-0.0723	&	20	&	0.014	&	0.00E+00	&	-0.0723	&	50	&	0.165	&	0.00E+00	\\
READING4N50	&	-0.2685	&	3556	&	1.767	&	8.22E-15	&	-0.0100	&	6099	&	0.123	&	3.61E+00	\\
READING4N100	&	-0.2799	&	9957	&	16.043	&	8.44E-15	&	0.0163	&	8118	&	22.609	&	5.40E-01	\\
READING4N500	&	-0.2893	&	88614	&	1198.368	&	6.00E-15	&	-	&	-	&	T	&	-	\\
\midrule
COSHFUNM3	&	-0.6614	&	379	&	0.152	&	3.77E-16	&	-0.6614	&	730	&	0.341	&	0.00E+00	\\
COSHFUNM8	&	-0.7708	&	4673	&	1.424	&	2.22E-16	&	-0.7708	&	$12500^*$	&	0.324	&	0.00E+00	\\
COSHFUNM14	&	-0.7732	&	$21500^*$	&	6.366	&	8.23E-13	&	-0.7731	&	$21500^*$	&	3.673	&	0.00E+00	\\
COSHFUNM20	&	-0.7733	&	$30500^*$	&	10.792	&	5.42E-09	&	-0.7731	&	$30500^*$	&	8.430	&	0.00E+00	\\
\bottomrule
\end{tabular}
\caption{{\em Variable-Size Problems, Noiseless Case.} Final objective value $\phi_k$, number of function evaluations (\#evals),  CPU time (in seconds), and final feasibility error (feaserr). A $^*$ indicates that the maximum number of function evaluations was reached, and \emph{T}  that the time limit (1 hour) was reached. }
\label{cons cpu}
\end{table}
\normalsize

\bigskip

\subsection{Experiments on Noisy Functions}

We now inject artificial noise in the evaluation of the objective and constraint functions.  We use the same noise model as in the unconstrained setting; i.e., uniform noise sampled i.i.d. given by 
\[
\epsilon(x), \epsilon_i(x) \sim \sigma_f U(-\sqrt{3}, \sqrt{3}),
\] 
where $\sigma_f \in \{10^{-1}, 10^{-3}, 10^{-5}, 10^{-7}\}$.
We use the same formulas for evaluating the finite-difference gradient \eqref{eq:approx g} and Jacobian \eqref{eq:approx J}, but with a different formula for the finite-difference interval: 
\[
h_i = \max\{1,|[x]_i|\}\sqrt{\sigma_f}, \quad i=1,...,n.
\]
We do not include a Lipschitz constant because this formula will suffice for our purposes.

\begin{figure}
\centering
	\includegraphics[width=0.32\linewidth]{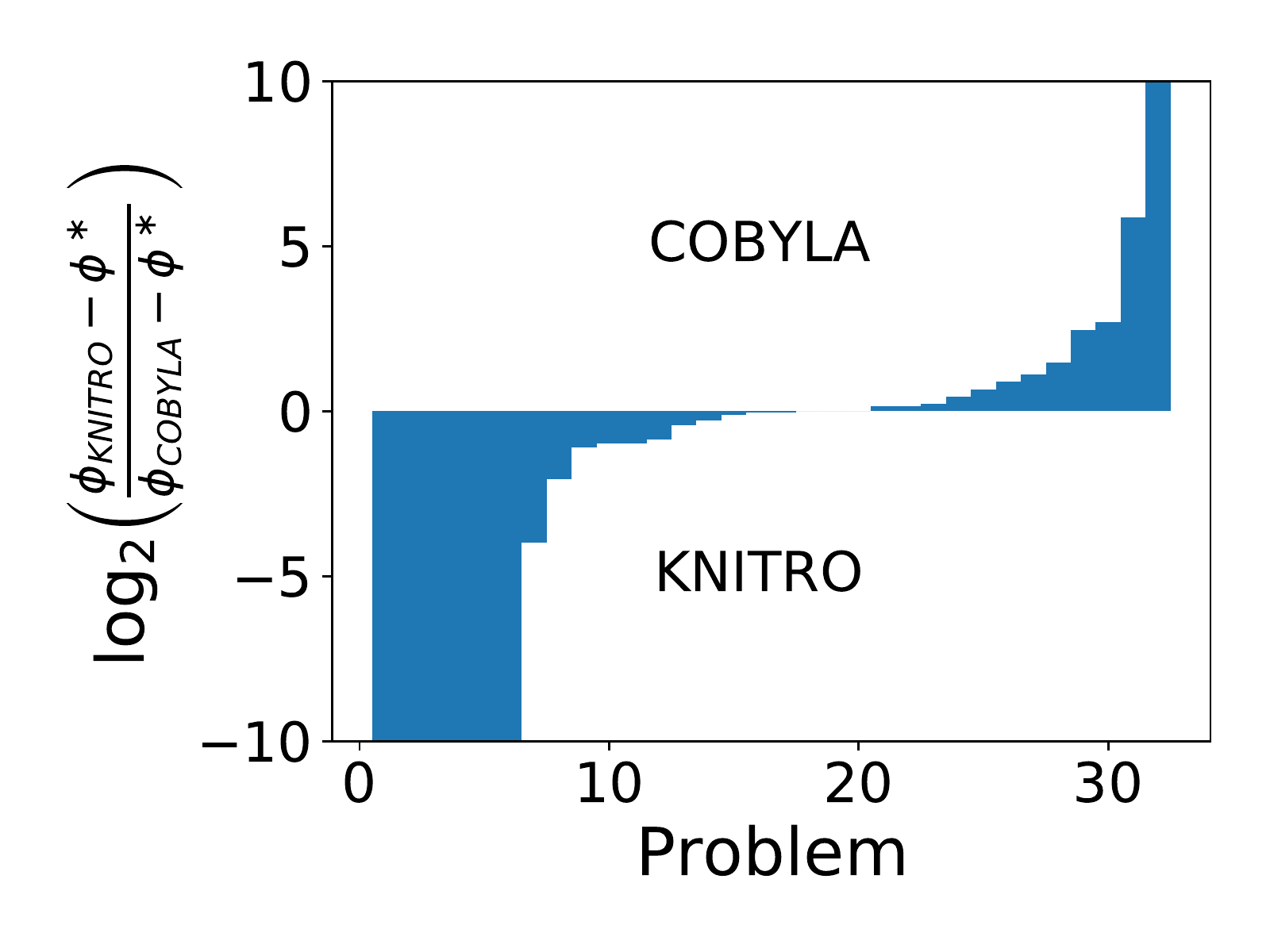}
	\includegraphics[width=0.32\linewidth]{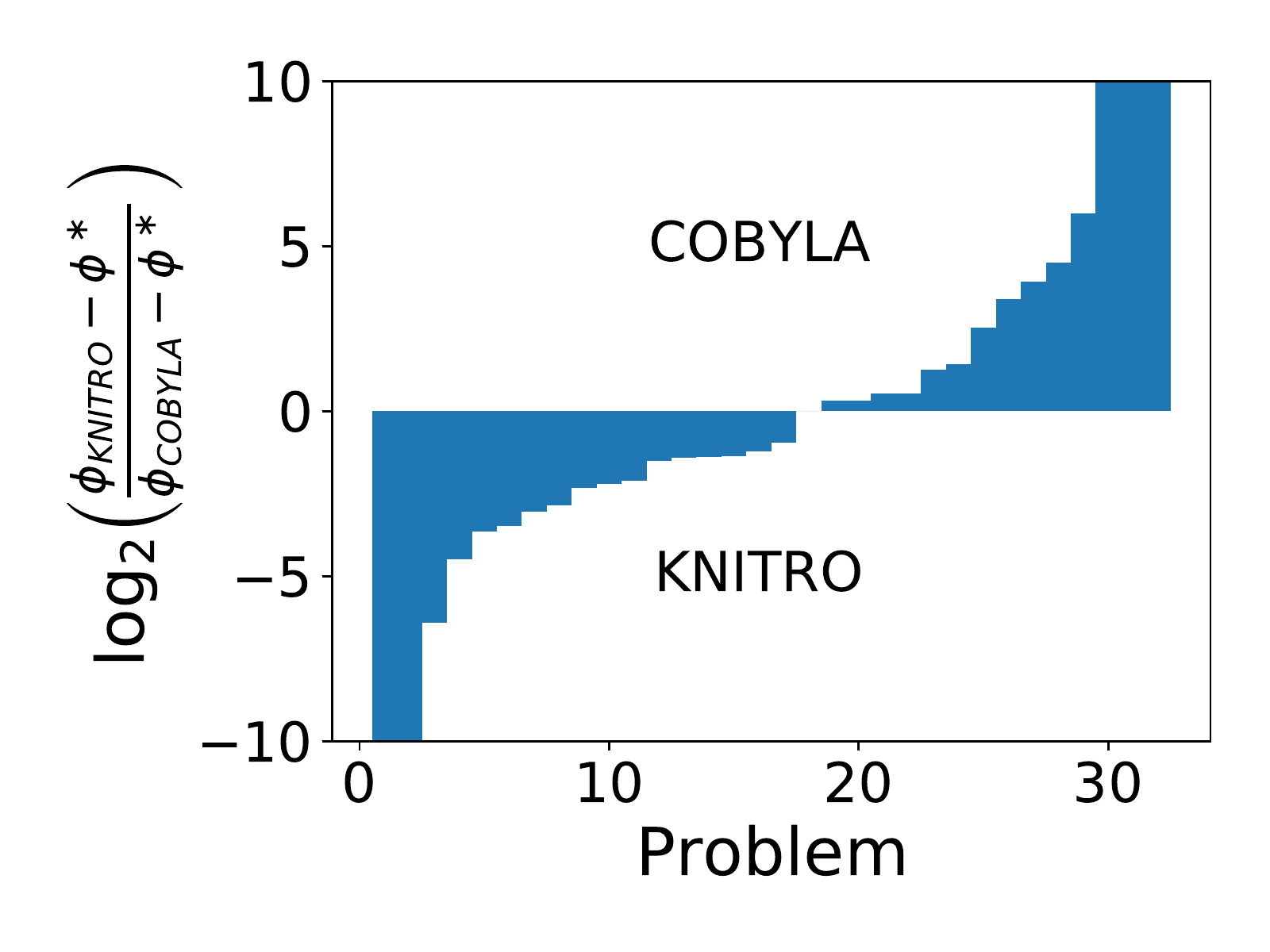}\\
	\includegraphics[width=0.32\linewidth]{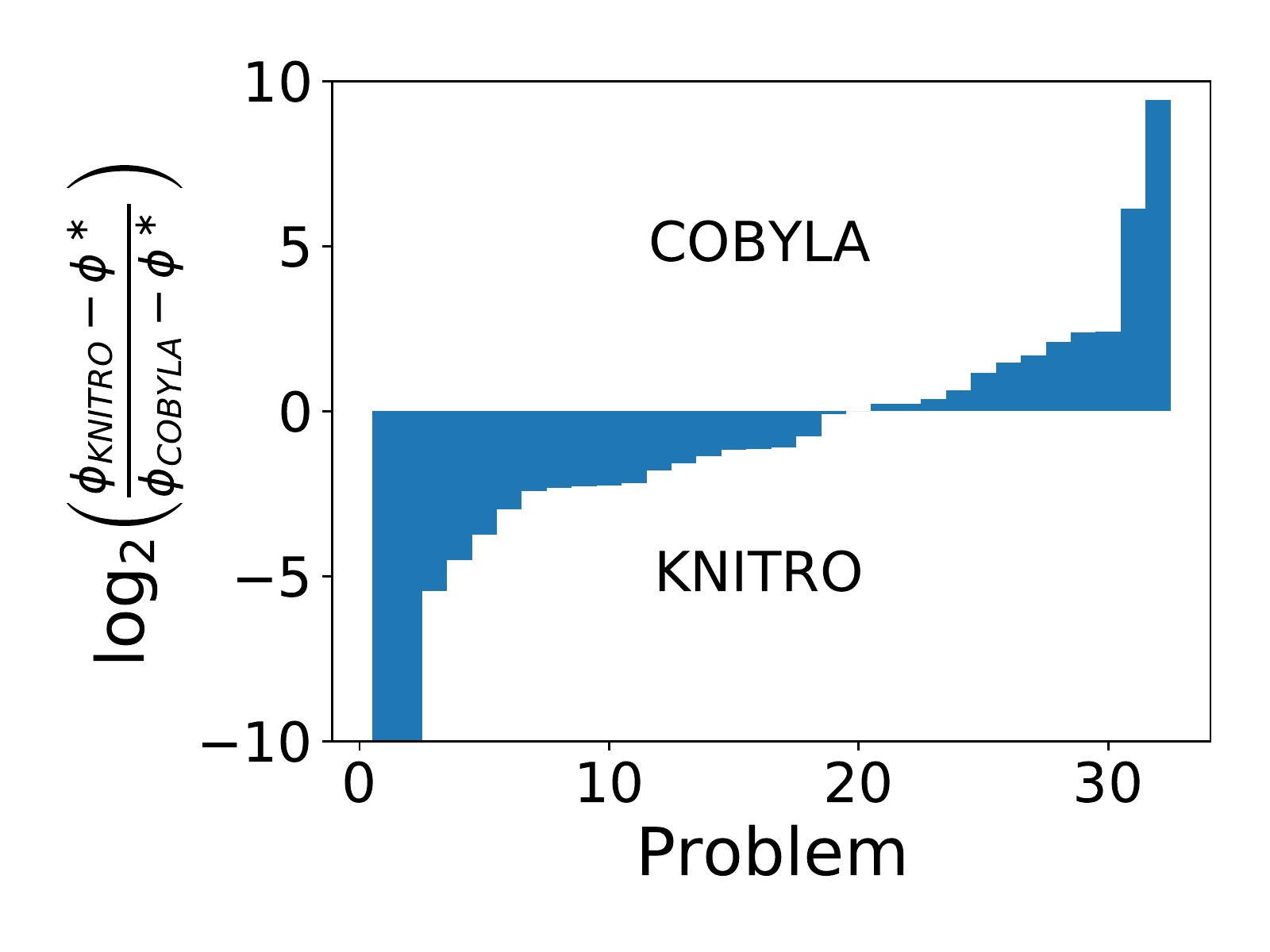}
	\includegraphics[width=0.32\linewidth]{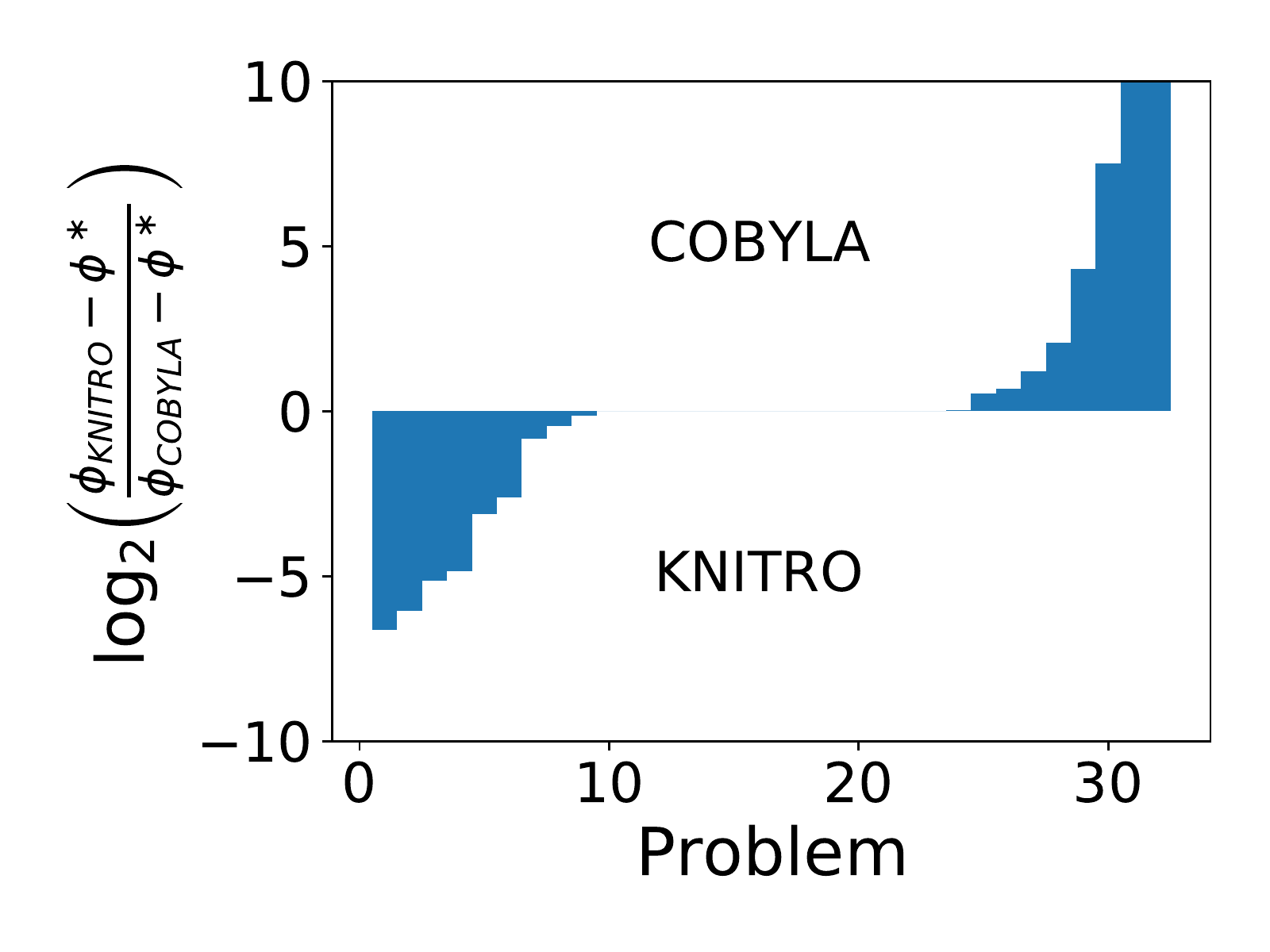}
    \caption{ {\em Accuracy, Noisy Case.} Log-ratio profiles comparing accuracy in the objective by {\sc knitro} and {\sc cobyla}, for $\sigma_f = 10^{-1}$ (upper left), $10^{-3}$ (upper right), $10^{-5}$ (bottom left), and $10^{-7}$ (bottom right).  
    %We compare objective values up to eight digits for problems with $\|\max\{\psi(x_k),0\}\|_\infty \leq \sqrt{3}\sqrt{\sigma_f}$ [\mx{$\sqrt{3}\sigma_f$?}].
    } 
    \label{fig:accuracy_eps}
\end{figure}

In the first experiment, we ran the two codes with their least stringent termination tests to observe the quality of the final solutions.  As  in the noiseless case, we set \texttt{rhoend=1e-8} for {\sc cobyla} and \texttt{xtol=1e-8} for {\sc knitro}, and impose a limit of $500\max\{n,m\}$ function evaluations. We tested the 32 CUTEst problems for which both solvers converged to the same feasible solution in the noiseless case (outcome (i) above).  In Figure~\ref{fig:accuracy_eps}, we compare the accuracy $\phi(x_k)-\phi_\ast$ in the true objective given by the code codes,  up to eight digits, as in \eqref{eight}.  If the feasibility violation is large compared to the noise level, i.e. $\|\max\{\psi(x_k),0\}\|_\infty \geq \sqrt{3}\sqrt{\sigma_f}$, we mark the corresponding run as a failure. 
We observe from Figure~\ref{fig:accuracy_eps} that the performance of the two solvers is comparable, with {\sc knitro} slightly more efficient for low accuracy. Interestingly, for the highest accuracy, $\sigma_f=10^{-7}$, the performance of the codes is remarkably close (note that the log ratio is nearly zero for about half of the problems). The complete set of results is given in Tables~\ref{table:noisyobj}-\ref{table:noisyfeas5} in Appendix \ref{sec:ext_num}.  
 
\begin{figure}
\centering
	\includegraphics[width=0.32\textwidth]{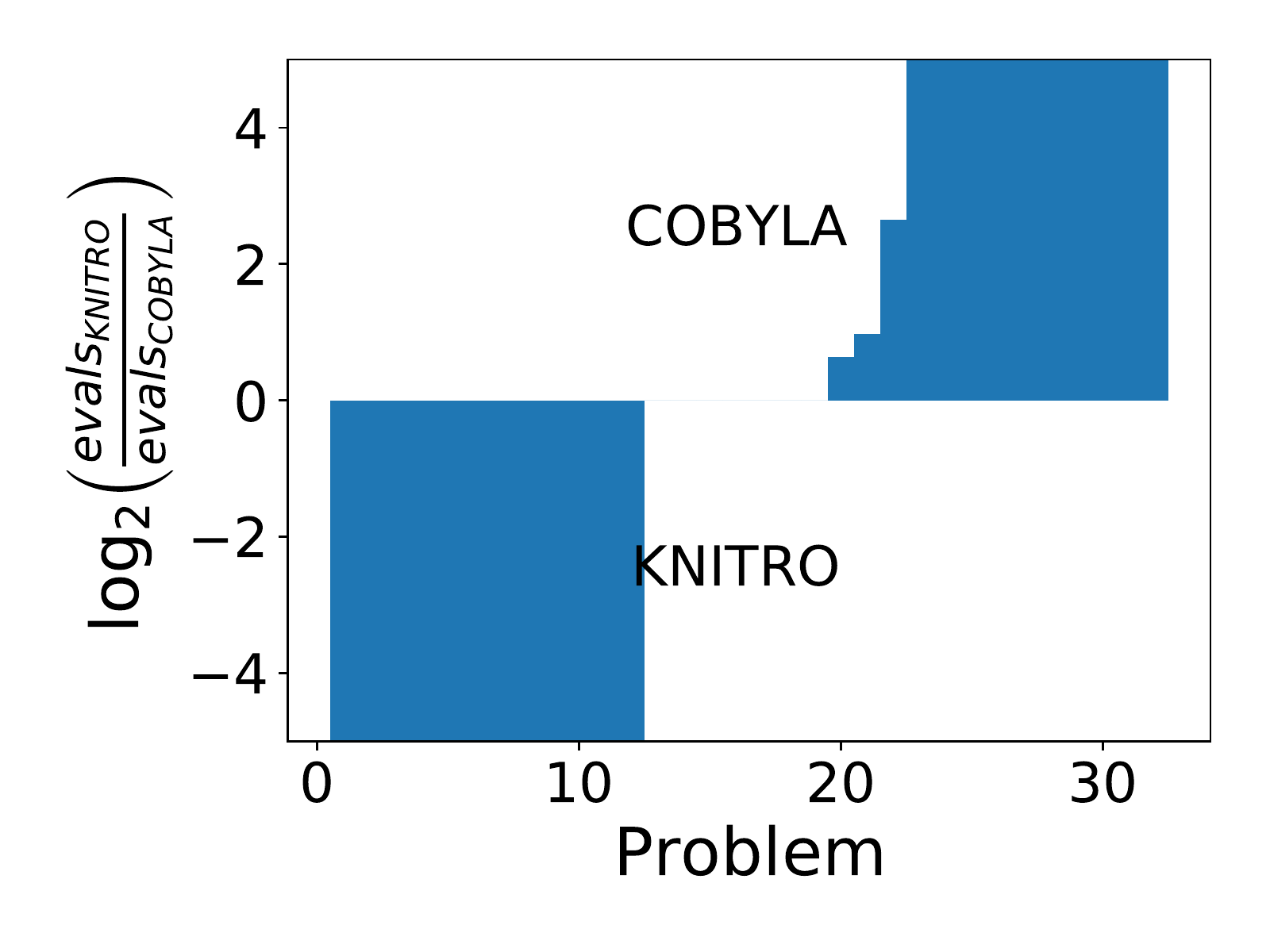}
	\includegraphics[width=0.32\textwidth]{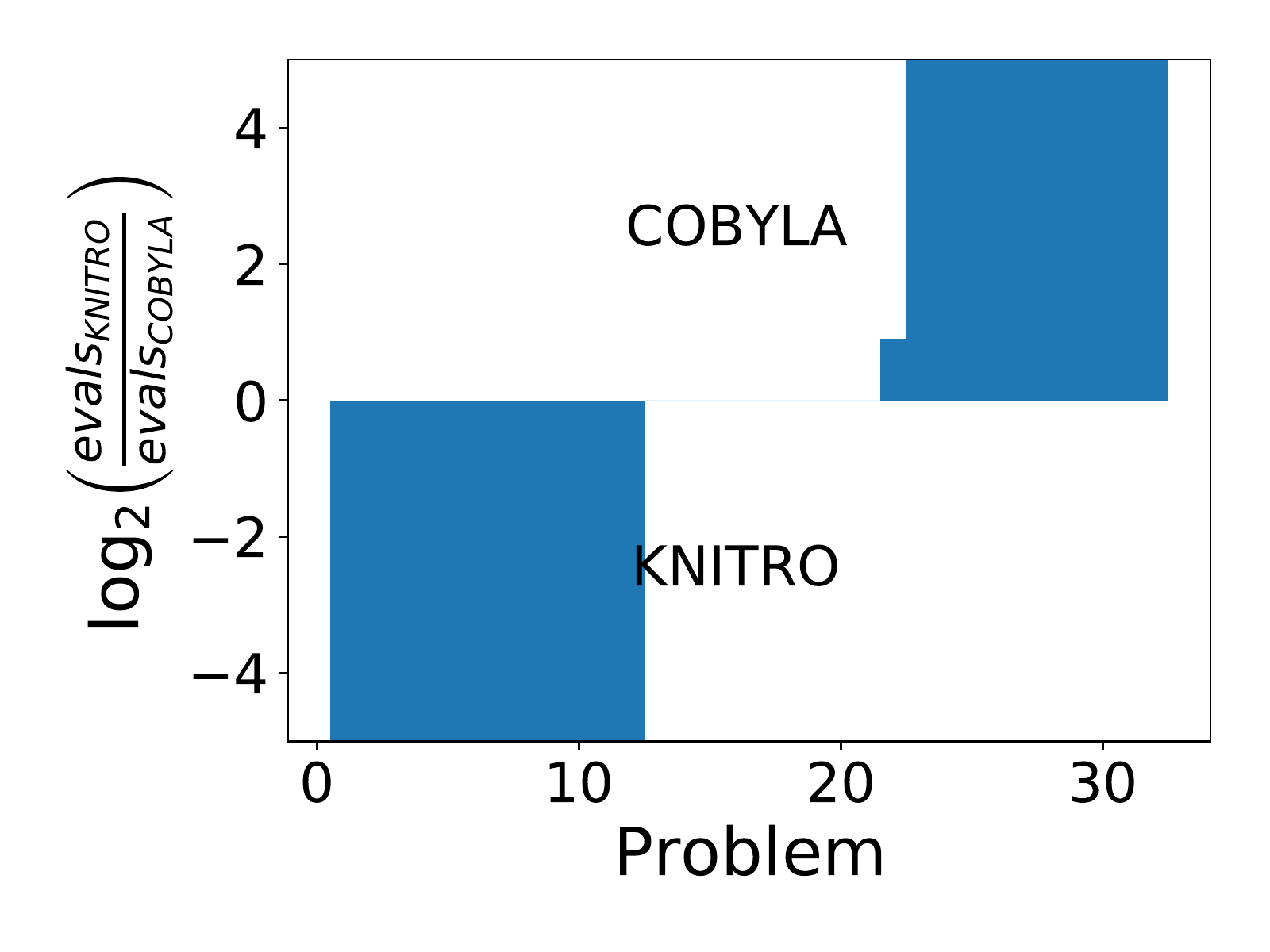}\\
	\includegraphics[width=0.32\textwidth]{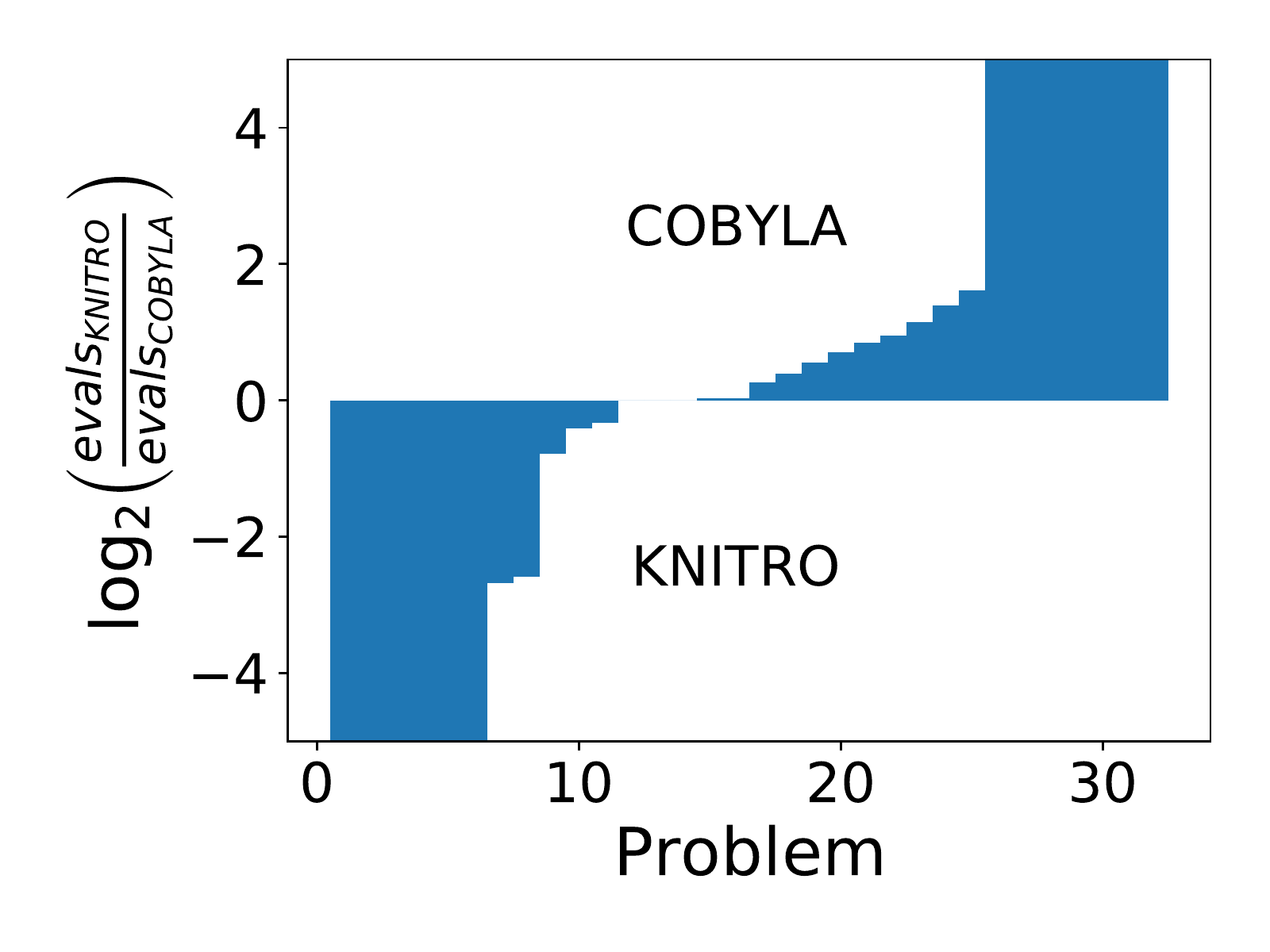}
	\includegraphics[width=0.32\textwidth]{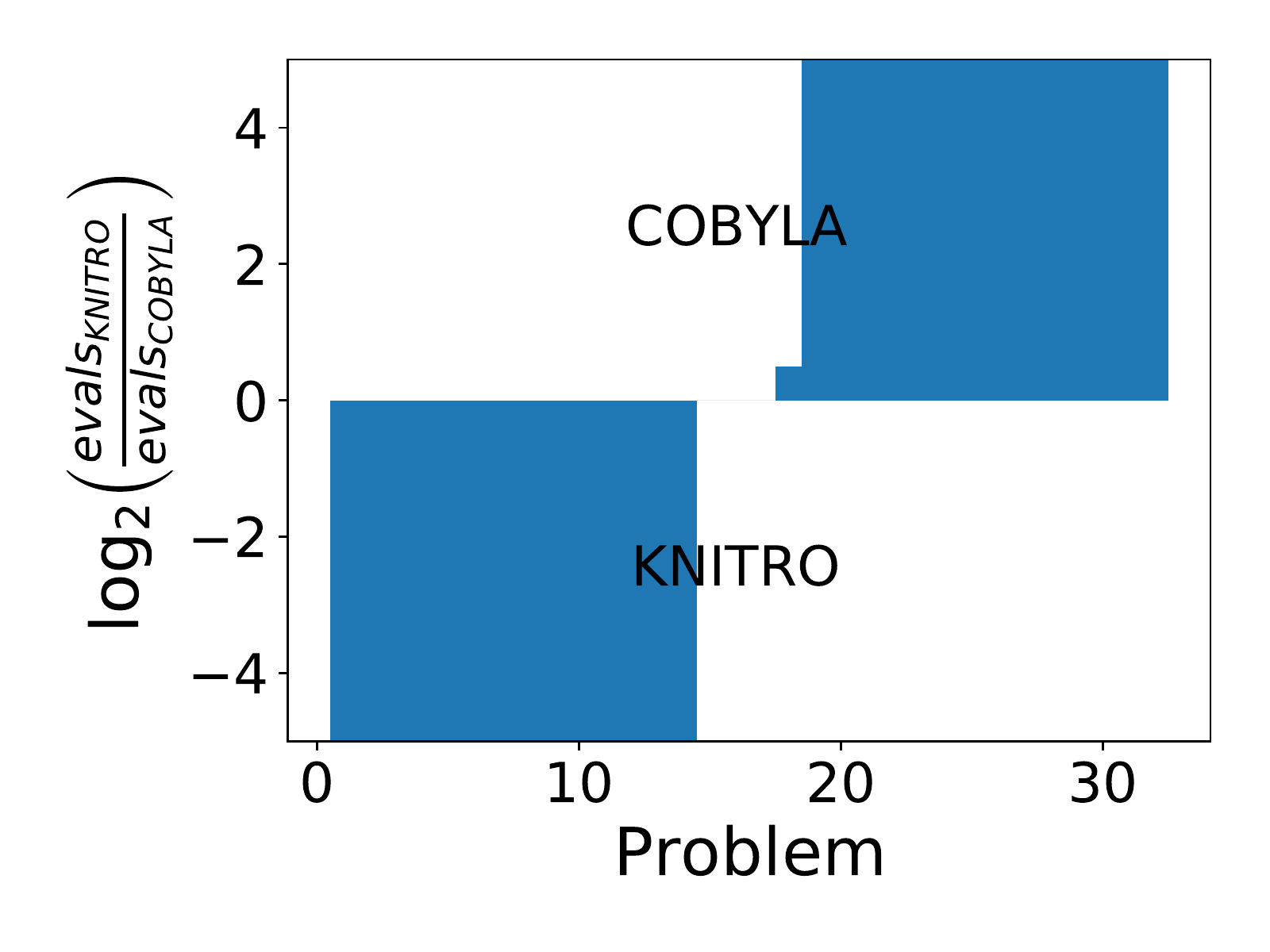}\\
	\includegraphics[width=0.32\textwidth]{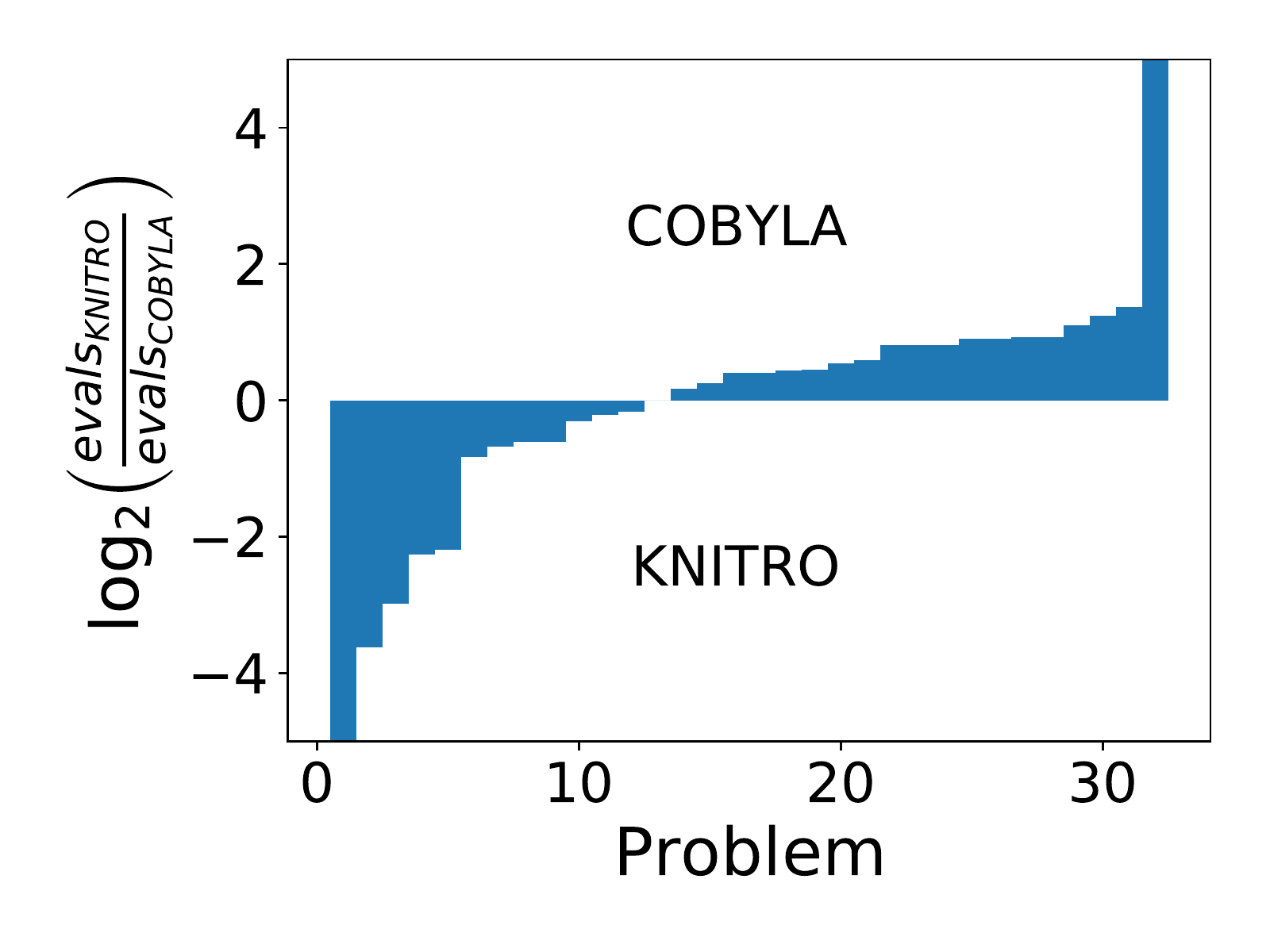}
	\includegraphics[width=0.32\textwidth]{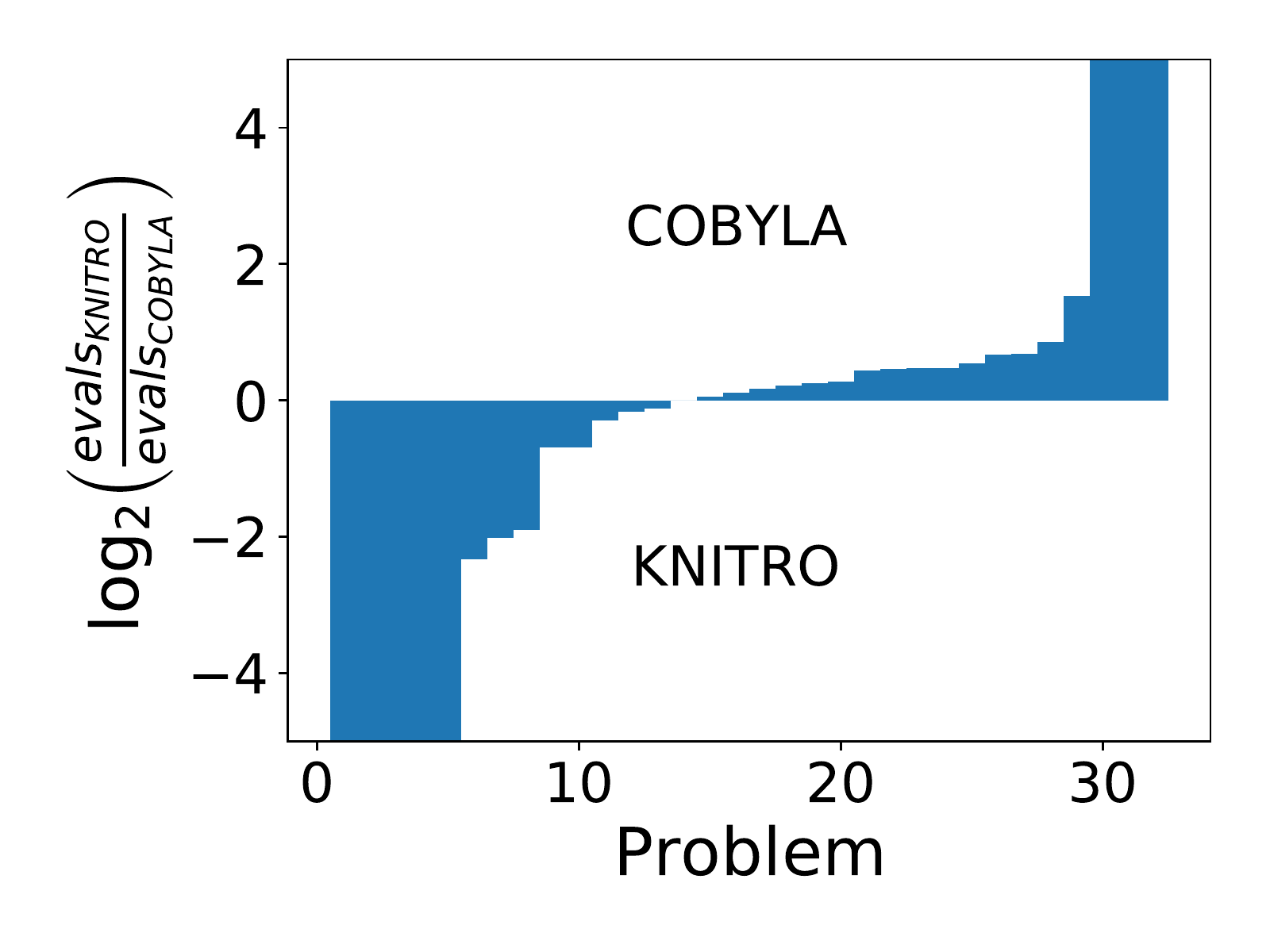}\\
    \caption{{\em Efficiency, Noisy Case}. Log-ratio profiles \eqref{plotthem} comparing  {\sc knitro} and {\sc cobyla} for $\sigma_f = 10^{-1}$ (top row), $10^{-3}$ (middle row),
    %10^{-5}(\text{row 3}),
    $10^{-7}$ (bottom row). The figure measures the number of function evaluations to satisfy \eqref{eq:term_eps} for $\tau = 10^{-2}$ (left) and $10^{-6}$ (right).}
	\label{fig:efficiency_eps}
\end{figure}

Next, we compare the efficiency of the two solvers for two levels of accuracy in the objective. Specifically, we record the number of function evaluations required to satisfy 
\begin{equation}
\label{eq:term_eps}
\phi(x_k)- \tilde \phi_\ast \leq \tau (\phi(x_0) - \tilde \phi_\ast) \quad\mbox{and}\quad \|\max\{\psi(x_k),0\}\|_\infty \leq \sqrt{3}\sqrt{\sigma_f},
\end{equation}      
for $\tau  \in \{10^{-2}, 10^{-6}\}$,
where $\tilde \phi_\ast$ is the minimum objective value obtained by the two codes.   Figure~\ref{fig:efficiency_eps} plots the ratios \eqref{plotthem}, which suggest that it is difficult to choose between the two codes. Therefore, in our tests on noisy problems, an unsophisticated code ({\sc knitro/sqp}) for deterministic nonlinear optimization, using a simple strategy for choosing the finite difference interval, is competitive with {\sc cobyla}, a method specifically designed for derivative-free optimization.

%\noindent
%\includegraphics[scale=0.5]{figures/constrained/feas_eps_1e-7}
%\includegraphics[scale=0.5]{figures/constrained/feas_eps_1e-5}\\
%\includegraphics[scale=0.5]{figures/constrained/feas_eps_1e-3}
%\includegraphics[scale=0.5]{figures/constrained/feas_eps_1e-1}

%\bigskip

%\noindent
%\includegraphics[scale=0.5]{figures/constrained/opt_eps_1e-7}
%\includegraphics[scale=0.5]{figures/constrained/opt_eps_1e-5}\\
%\includegraphics[scale=0.5]{figures/constrained/opt_eps_1e-3}
%\includegraphics[scale=0.5]{figures/constrained/opt_eps_1e-1}

\section{Final Remarks}
\label{ch:final}
Finite-difference approximations are widely employed within numerical analysis, particularly for solving ordinary and partial differential equations. The dangers and limitations of finite-difference methods are therefore also well-documented, such as the difficulty of approximating high-order derivatives accurately on unstructured grids or a failure to adapt to problems that require local mesh refinement.

Optimization, however, arguably provides a more benign setting than solving differential equations as errors do not accumulate; if a poor choice of the finite-difference interval yields a bad step, it can be detected and improved upon at a later point within the iteration. Two attractive features of finite-difference methods in optimization are the simplicity of incorporating them into existing nonlinear optimization solvers, and their ease of parallelization, in many situations.

Our empirical study has revealed that finite-difference methods can be made competitive, in many cases, against state-of-the-art derivative-free optimization methods. However, algorithmic development and analysis are still necessary to make them sufficiently robust for general-purpose optimization. In addition to more sophisticated procedures for computing the finite difference interval, research is needed in the design of globalization mechanisms such as line searches and trust regions in the noisy setting. We hope that our work draws attention to and provides an empirical basis for further investigation along this important line of inquiry. 

% The dangers and limitations of finite difference approximations have been well documented in the numerical analysis literature. Optimization problems, however, present a more benign setting in that errors do not always accumulate, and if they lead to poor steps, this can be detected and recovery mechanisms can be invoked; \cite{sbe,sbe}. We do not consider such techniques in this paper, but the results presented here provide clear guidelines to inform the developments of such methods. 

% In fact, several parts of the algorithms in {\sc knitro} such as line search are prone to failure with noisy evaluations.  Also, we are not aware of any analysis of {\sc cobyla} in the presence of noise.

\section*{Acknowledgements}

We are grateful to Richard Byrd, Oliver Zhuoran Liu, and Yuchen Xie for their initial feedback on this work. We also thank Philip Gill, Tammy Kolda, Arnold Neumaier, Michael Saunders, Katya Scheinberg, Luis Vicente and Stefan Wild for their correspondences that led to the design of the experiments in this work.

\bibliographystyle{plain}

\bibliography{references}

\appendix

\section{Numerical Investigation of Lipschitz Estimation}
\label{app:lipschitz}

\subsection{Investigation of Theoretical Lipschitz Estimates}
\label{app:theoretical lipschitz}

In Section \ref{sub:noisy unc}, we approximated the bound on the second and third derivative $L$ and $M$ (or the Lipschitz constant of the first and second derivative) along the interval $I = \{x \pm tp : t \in [0, h_0]\}$ for $h_0 > 0$ every time finite-differencing is performed. 

For forward-differencing, we employed the Mor\'e and Wild heuristic \cite{more2012estimating}. The heuristic estimates the bound on the second derivative of a univariate function with noise. Assume $\phi : \mathbb{R} \rightarrow \mathbb{R}$ is univariate. If we let $\Delta(t) = f(x + t) - 2f(x) + f(x - t)$, then $t > 0$ must satisfy two conditions:
\begin{align}
|\Delta(t)| & \geq \tau_1 \epsilon_f, & & \tau_1 \gg 1 \label{eq:more-wild error}\\
|f(x \pm t) - f(x)| & \leq \tau_2 \max\{|f(x)|, |f(x \pm t)|\}, & & \tau_2 \in (0, 1). \label{eq:more-wild equal}
\end{align} 
In practice, $\tau_1 = 100$ and $\tau_2 = 0.1$. The first condition ensures that $h$ is sufficiently large such that the second-order difference is not dominated by noise, while the second condition enforces that the difference is not dominated by a particular evaluation, so that there is ``equal'' contribution from each function evaluation in the finite-difference approximation. If conditions \eqref{eq:more-wild error} and \eqref{eq:more-wild equal} are satisfied, then we take $L = \max\{10^{-1}, |\Delta(t)| / t^2\}$.

For central differencing, we used a theoretical estimate based on knowledge of the true Hessian $\nabla^2 \phi(x)$:
\begin{equation} \label{gtn}
    M = \max\left\{10^{-1}, \frac{|p^T (\nabla^2 \phi(x + \tilde{h} \frac{p}{\|p\|}) - \nabla^2 \phi(x)) p|}{\tilde{h} \|p\|^2} \right\}
\end{equation}
where $\tilde{h} = \sqrt{\epsilon_M}$. Since this estimate is ideal, we do not include the cost of evaluating $M$ in the number of function evaluations. In both cases, the maximum is taken to ensure that the bound is strictly bounded away from $0$ so that \eqref{eq:fd formula} and \eqref{eq:cd formula} remain well-defined. 

One may ask whether or not a refined choice of $L$ or $M$ is necessary for finite-difference {\sc l-bfgs}. To show the impact of Lipschitz estimation on the performance of finite-difference methods, we focus on forward-difference {\sc l-bfgs} and compare nine different theoretical Lipschitz estimation schemes. The first five consider techniques where $L$ is fixed for the entire run based on information at the initial point. We test this because it is frequently claimed that using an initial estimate of $L$ is sufficient for the entire run; see \cite{GillMurrSaunWrig83,more2012estimating}. The first approach requires no additional information about the function, while the others incorporate information about the Hessian at the current point. The latter four techniques similarly incorporate information about the Hessian but re-estimate $L$ whenever the finite-difference gradient or directional derivative is evaluated. 

All of these techniques rely on the assumption that $|D_p^2 \phi(x)| \approx |D_p^2 \phi(\xi)|$ for some $\xi \in [x, x + h]$. As we will see, both of these approximations that are based on conventional wisdom are challenged in our experiments.

\begin{enumerate}
\item Fix $L = 1$. This requires no additional knowledge about the problem at no added cost.
\item Fix $L = \max\{10^{-1}, \|\nabla^2 \phi(x_0)\|_2\}$. Similar to fixing $L = 1$, but incorporates knowledge of the initial Hessian.
\item Fix $L = \max \left\{10^{-1}, \frac{1}{n} \sum_{i = 1}^n |[\nabla^2 \phi(x_0)]_{ii}| \right\}$. This can be obtained by choosing $h$ such that the bound on $\|g(x_0) - \nabla \phi(x_0)\|_1$ at the initial point is minimized. 
\item Fix $L = \max \left\{10^{-1}, \frac{1}{\sqrt{n}} \sqrt{\sum_{i = 1}^n [\nabla^2 \phi(x_0)]_{ii}^2} \right\}$. This can be obtained by choosing $h$ such that the bound on $\|g(x_0) - \nabla \phi(x_0)\|_2^2$ at the initial point is minimized.
\item Fix $L$ to a vector $(\max\{10^{-1}, |[\nabla^2 \phi(x_0)]_{ii}| \})_{i = 1}^n$ and use the $i$th component for estimating the $i$th component of $g(x)$. Uses $\|L\|_2 / \sqrt{n}$ when estimating $h$ for the directional derivative in the line search.
\item Evaluate $L = \max\{10^{-1}, \|\nabla^2 \phi(x)\|_2\}$ each time finite-differencing is performed. 
\item Evaluate $L = \max \left\{10^{-1}, \frac{1}{n} \sum_{i = 1}^n |[\nabla^2 \phi(x)]_{ii}| \right\}$ each time finite-differencing is performed.
\item Evaluate $L = \max \left\{10^{-1}, \frac{1}{\sqrt{n}} \sqrt{\sum_{i = 1}^n [\nabla^2 \phi(x)]_{ii}^2} \right\}$ each time finite-differencing is performed.
\item Evaluate $L = \max\left\{10^{-1}, |[\nabla^2 \phi(x)]_{ii}| \right\}$ when evaluating $[g(x)]_i$. When the directional derivative along $p \in \R^n$ is computed, evaluates $L = \max\left\{10^{-1}, \frac{|p^T \nabla^2 \phi(x) p|}{\|p\|_2^2} \right\}$ .
\end{enumerate}

Each method is terminated when it cannot make more progress over 5 consecutive iterations. The optimal value $\phi^*$ is obtained by running {\sc l-bfgs} on the non-noisy function until no more progress can be made.

To compare the solution quality between two algorithms, we will use log-ratio profiles as proposed in \cite{Mora02} over the optimality gaps for each algorithm. The log-ratio profiles report
\begin{equation}
   \log_2 \left(\frac{\phi_{\text{new}} - \phi^*}{\phi_{\text{old}} - \phi^*}\right)
\end{equation}
for each problem plotted in increasing order. The area of the shaded region is representative of the general success of the algorithm.

\begin{figure}[ht]
\centering
\includegraphics[width=0.32\textwidth]{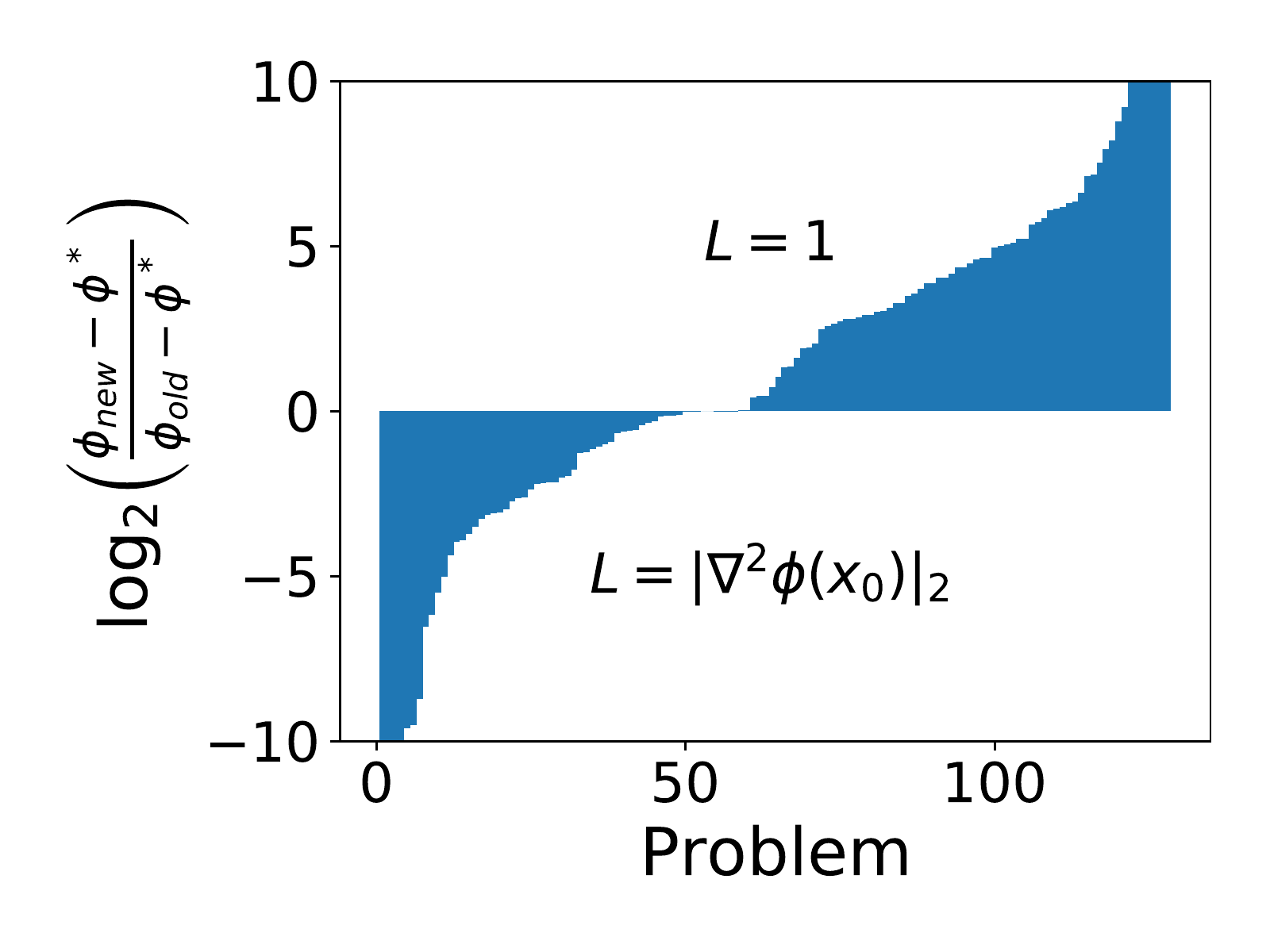}
\includegraphics[width=0.32\textwidth]{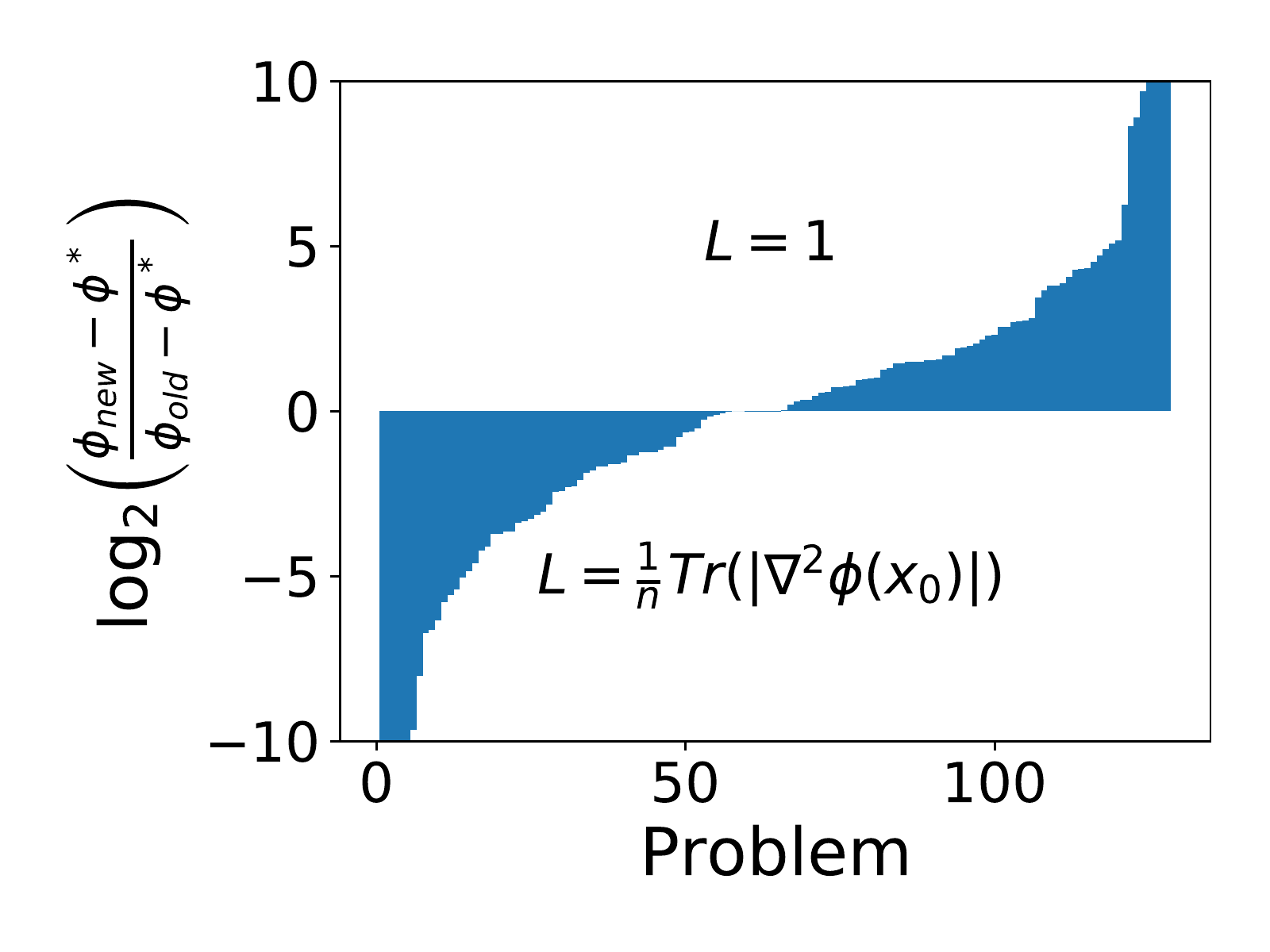}\\
\includegraphics[width=0.32\textwidth]{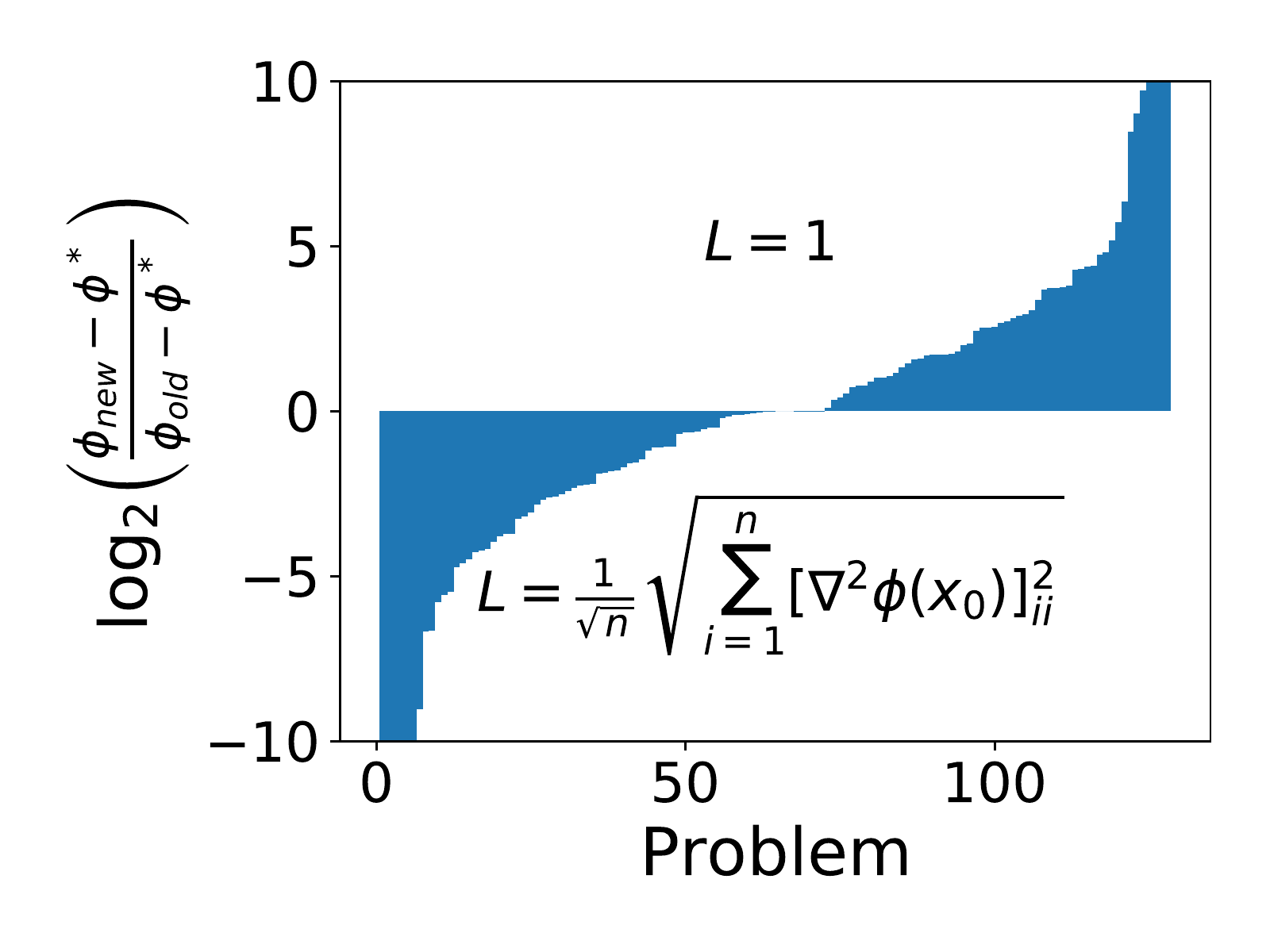}
\includegraphics[width=0.32\textwidth]{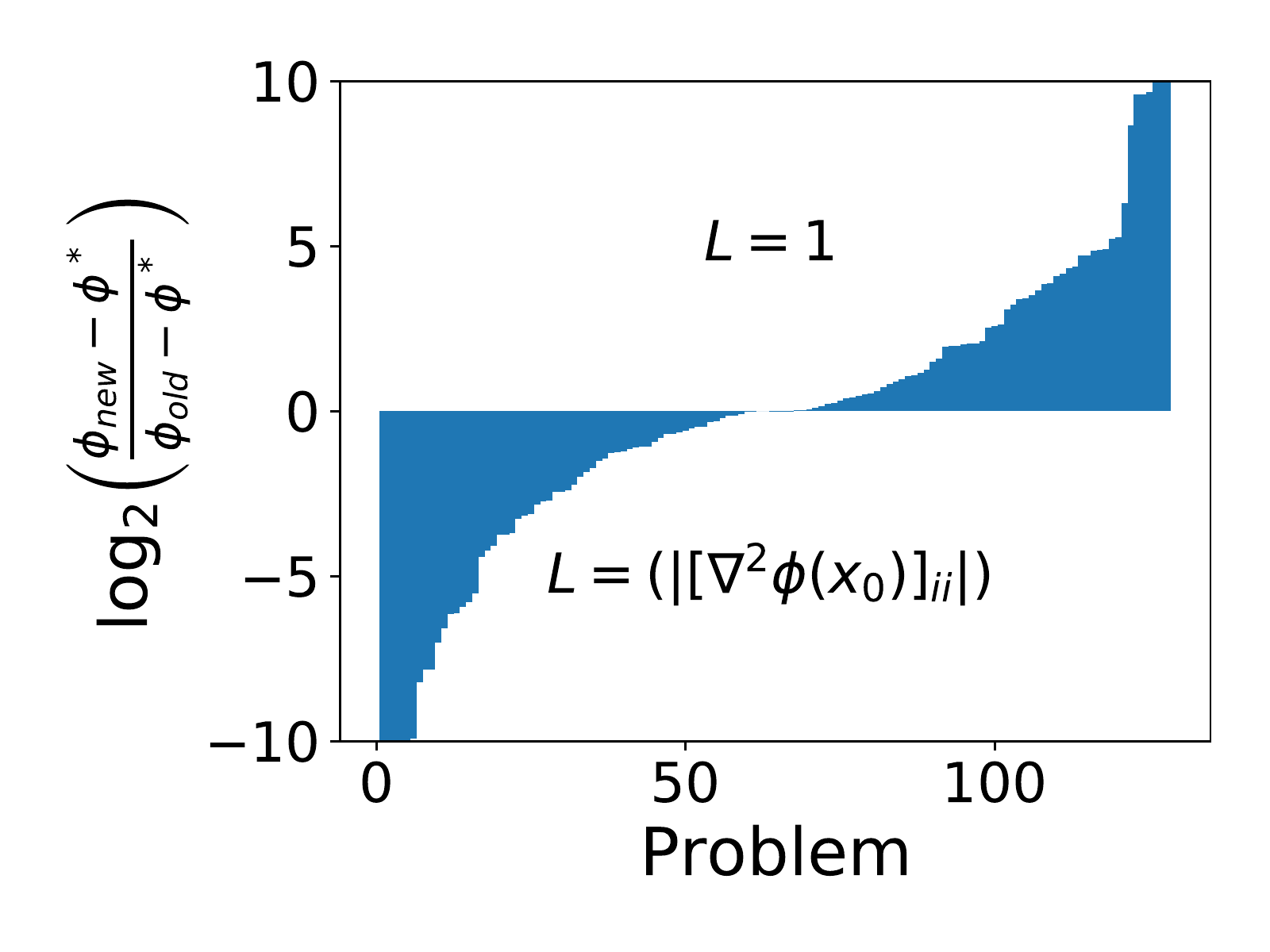}
\caption{{\em Accuracy, Noisy Case with $\sigma_f = 10^{-3}$}. Log-ratio optimality gap profiles comparing forward difference {\sc l-bfgs} with $L = 1$ and other fixed Lipschitz estimation schemes. The noise level is $\sigma_f = 10^{-3}$, but is representative for $\sigma_f \in \{10^{-1}, 10^{-3}, 10^{-5}, 10^{-7}\}$.}
\label{fig:L_1 vs fixed}
\end{figure}

One may ask whether or not using $L = 1$ is sufficient, without any additional knowledge of the second derivative. As seen in Figure \ref{fig:L_1 vs fixed}, when compared against the other fixed Lipschitz estimation schemes, the difference is not significant. In particular, the additional information from the Hessian does not yield significant benefits over setting $L = 1$ because the bound on the second derivative does not remain valid over the entire run of the algorithm. This may be surprising as it has been commonly suggested that it is sufficient to fix the Lipschitz estimate; see \cite{GillMurrSaunWrig83,more2012estimating}. In fact, we will see that it is more crucial for the algorithm to have a right estimate of $L$ during the later stages of the run than at the beginning of the run in order to achieve high accuracy.

\begin{figure}[ht]
\centering
\includegraphics[width=0.32\textwidth]{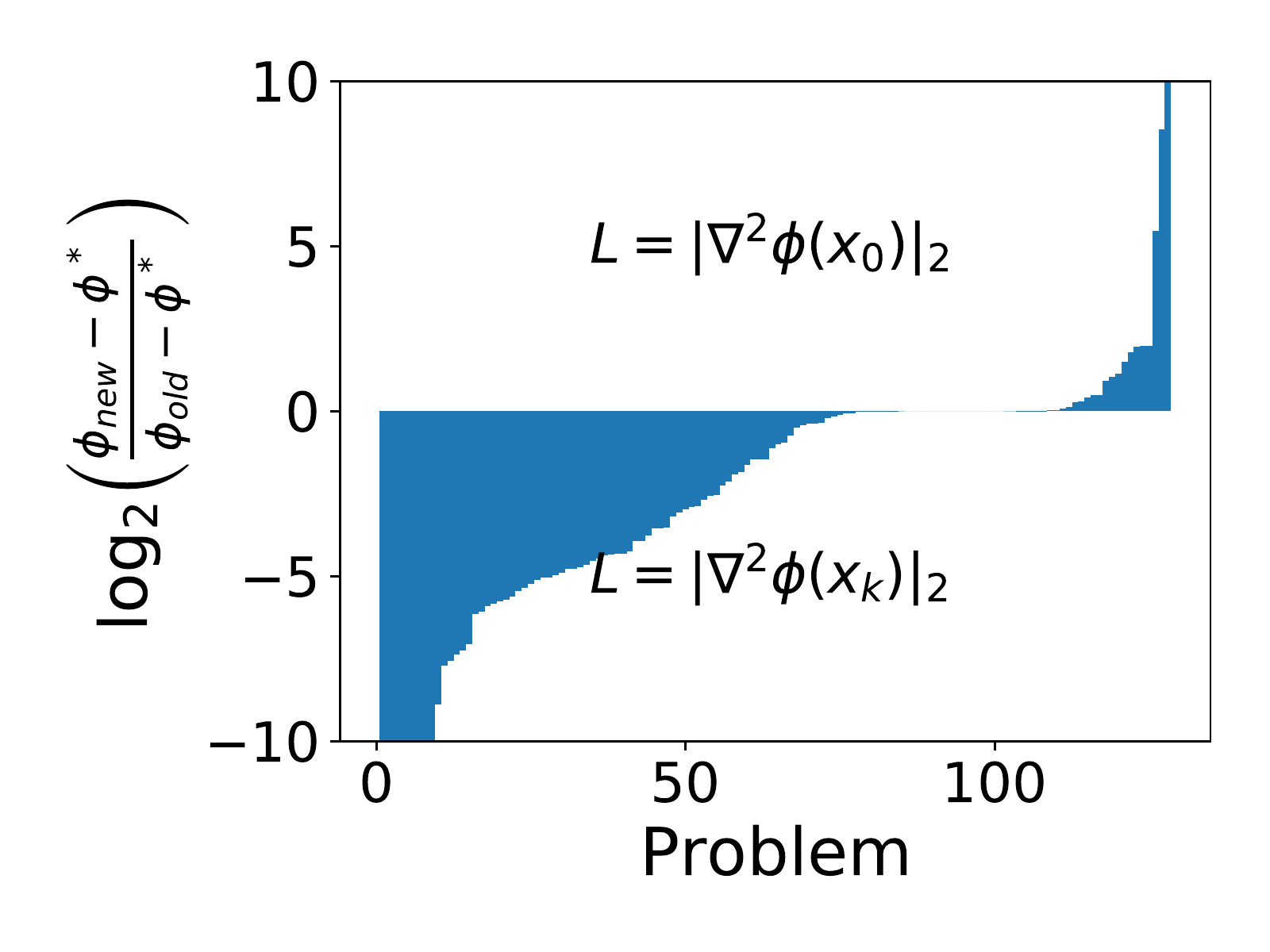}
\includegraphics[width=0.32\textwidth]{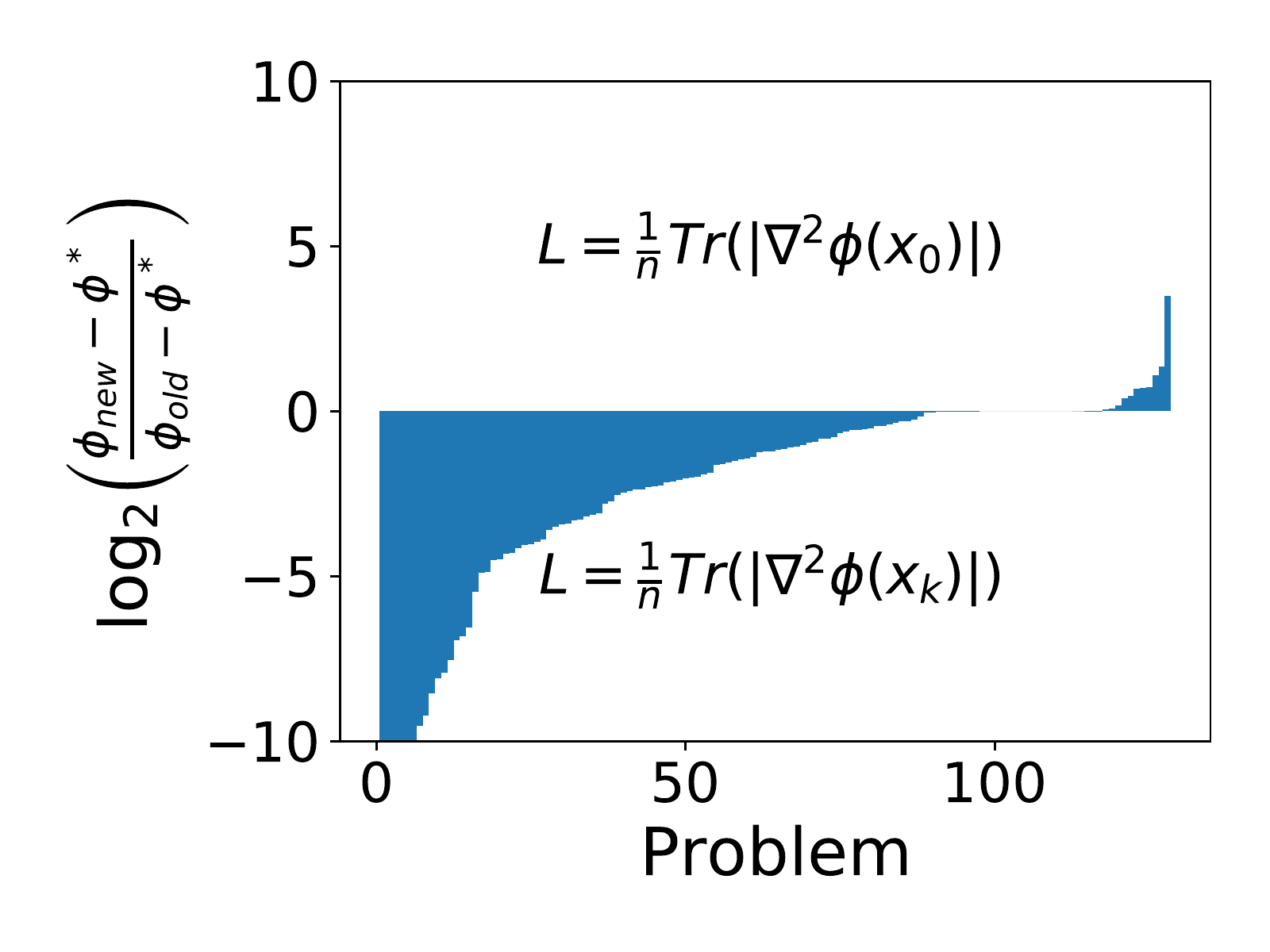}\\
\includegraphics[width=0.32\textwidth]{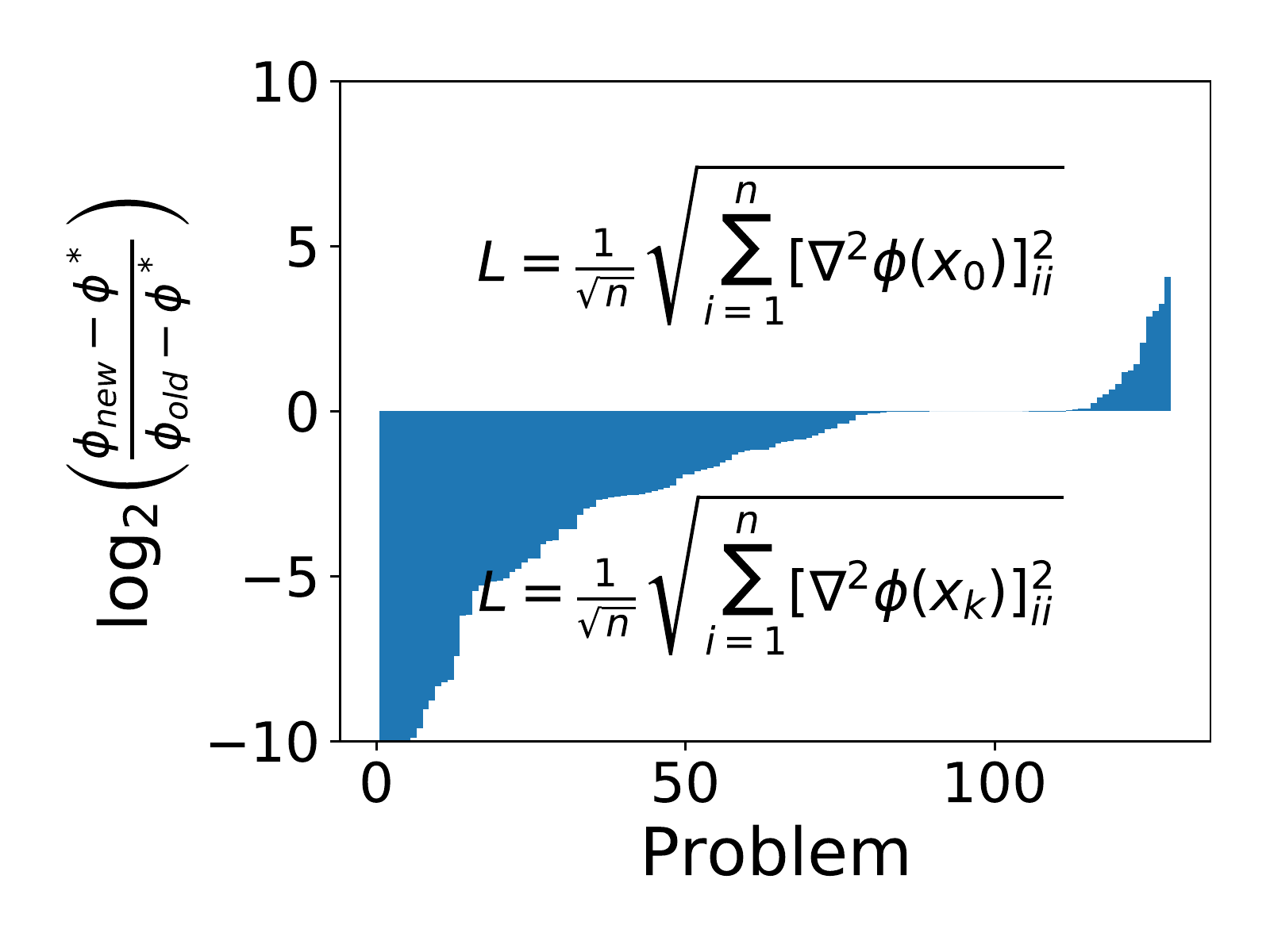}
\includegraphics[width=0.32\textwidth]{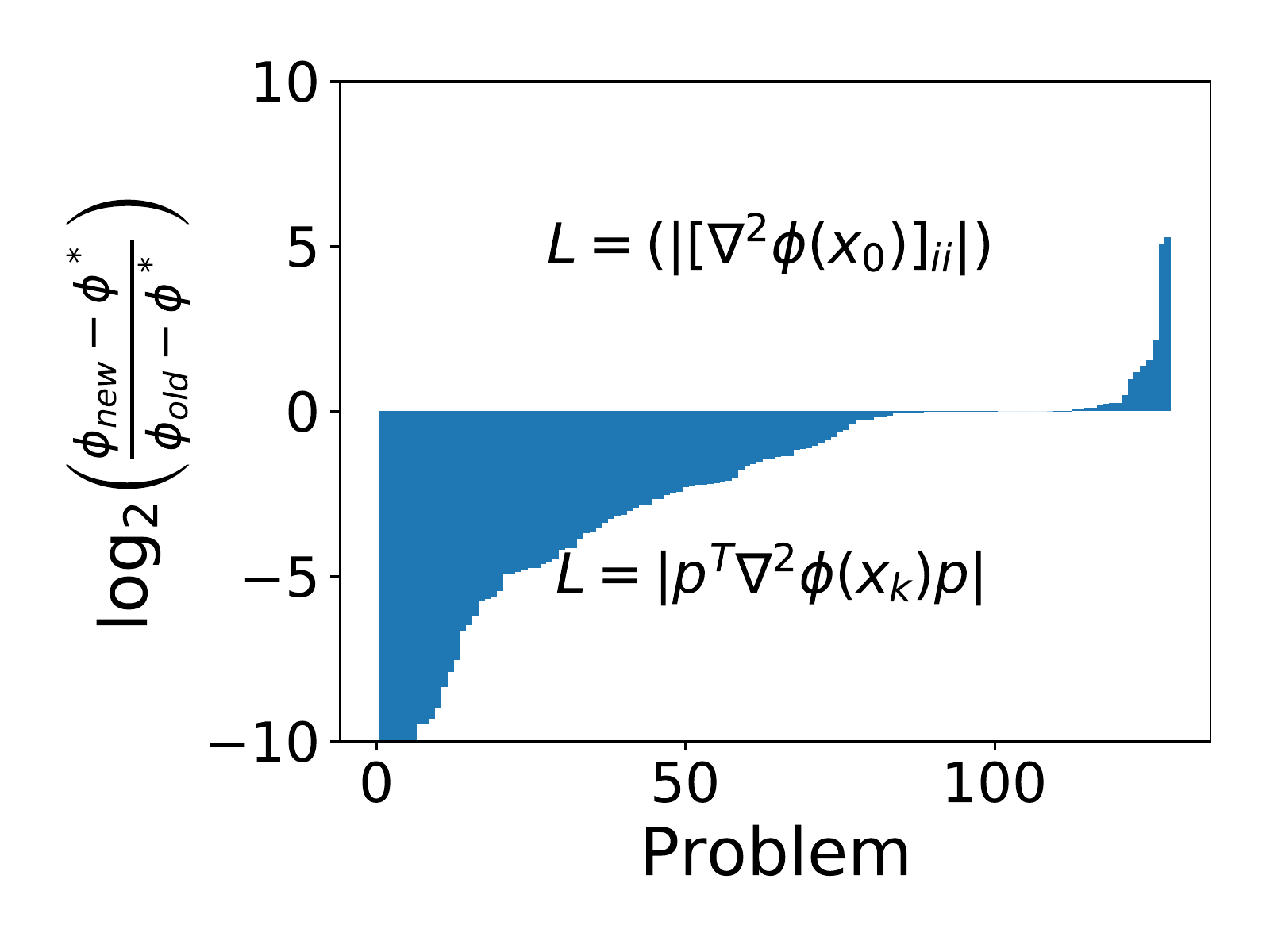}
\caption{{\em Accuracy, Noisy Case with $\sigma_f = 10^{-3}$}. Log-ratio optimality gap profiles comparing forward difference {\sc l-bfgs} with fixed and adaptive Lipschitz estimation schemes. The noise level is $\sigma_f = 10^{-3}$, but is representative for $\sigma_f \in \{10^{-1}, 10^{-3}, 10^{-5}, 10^{-7}\}$.}
\label{fig:fixed vs adaptive L}
\end{figure}

To further support why an adaptive $L$ is important to achieve high accuracy, we compare the fixed $L$ approaches against the adaptive $L$ in Figure \ref{fig:fixed vs adaptive L}.

As seen in Figure \ref{fig:fixed vs adaptive L}, we see that using Lipschitz estimation at every iteration yields much higher accuracy than the alternatives. Upon inspection of individual runs, one can observe many cases where fixing the Lipschitz constant is unstable and inadequate, potentially leading to poor gradient approximations that result in early stagnation of the algorithm. This suggests that it is imperative to adaptively re-estimate $L$ in order to reduce the noise in the gradient and squeeze the best possible accuracy out of forward difference {\sc l-bfgs}. 

\begin{figure}[ht]
\centering
\includegraphics[width=0.32\textwidth]{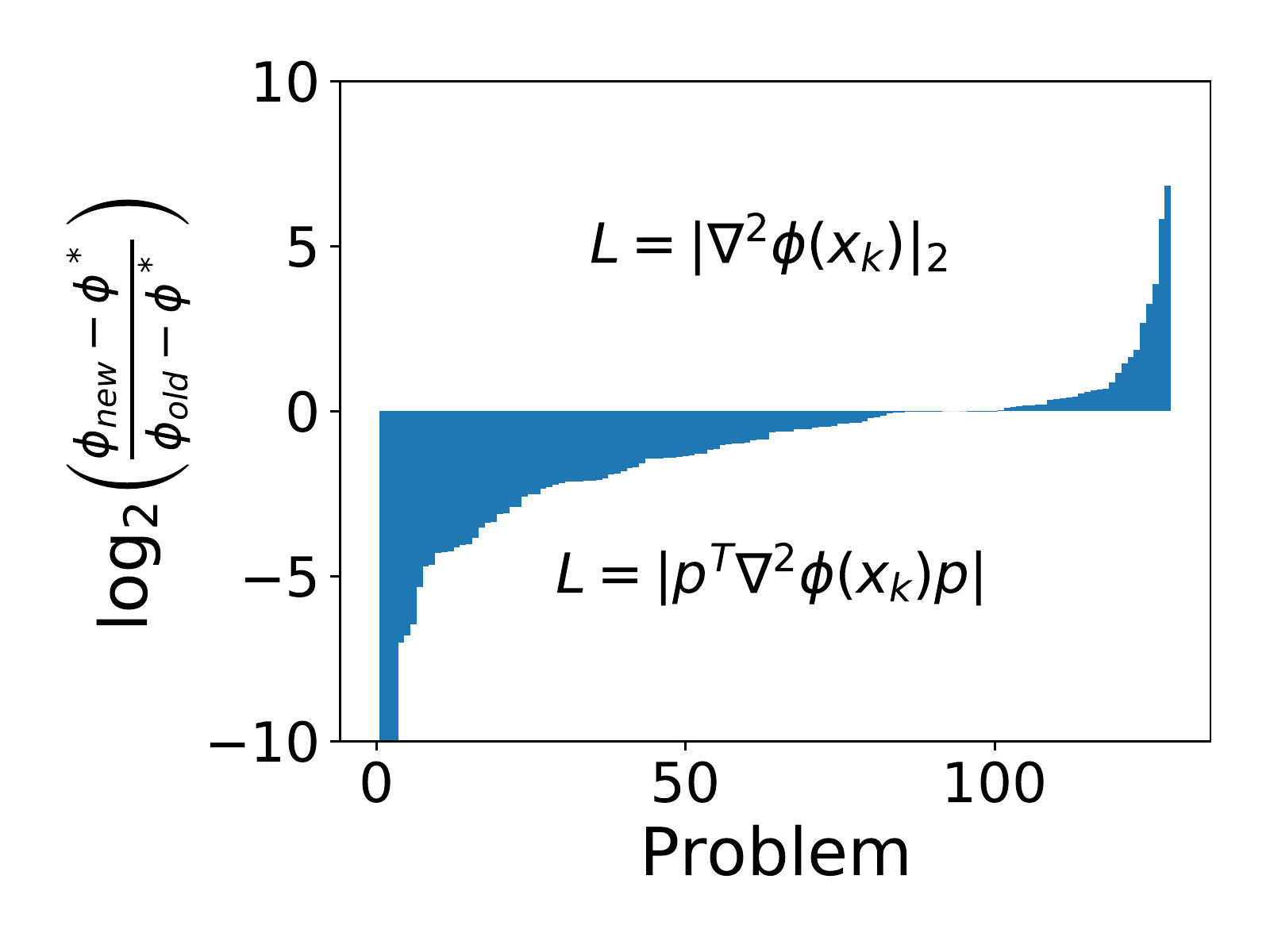}
\includegraphics[width=0.32\textwidth]{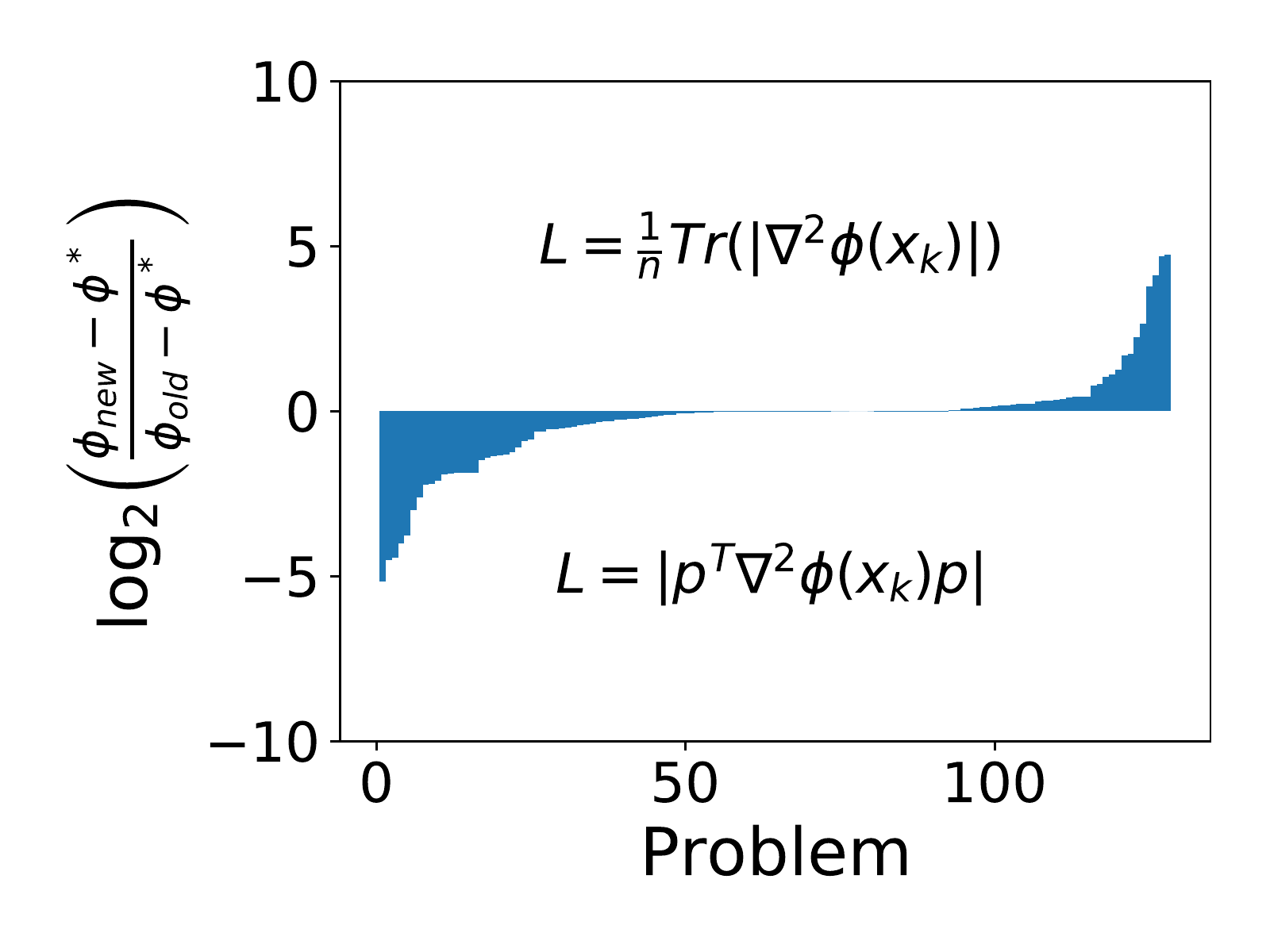}
\includegraphics[width=0.32\textwidth]{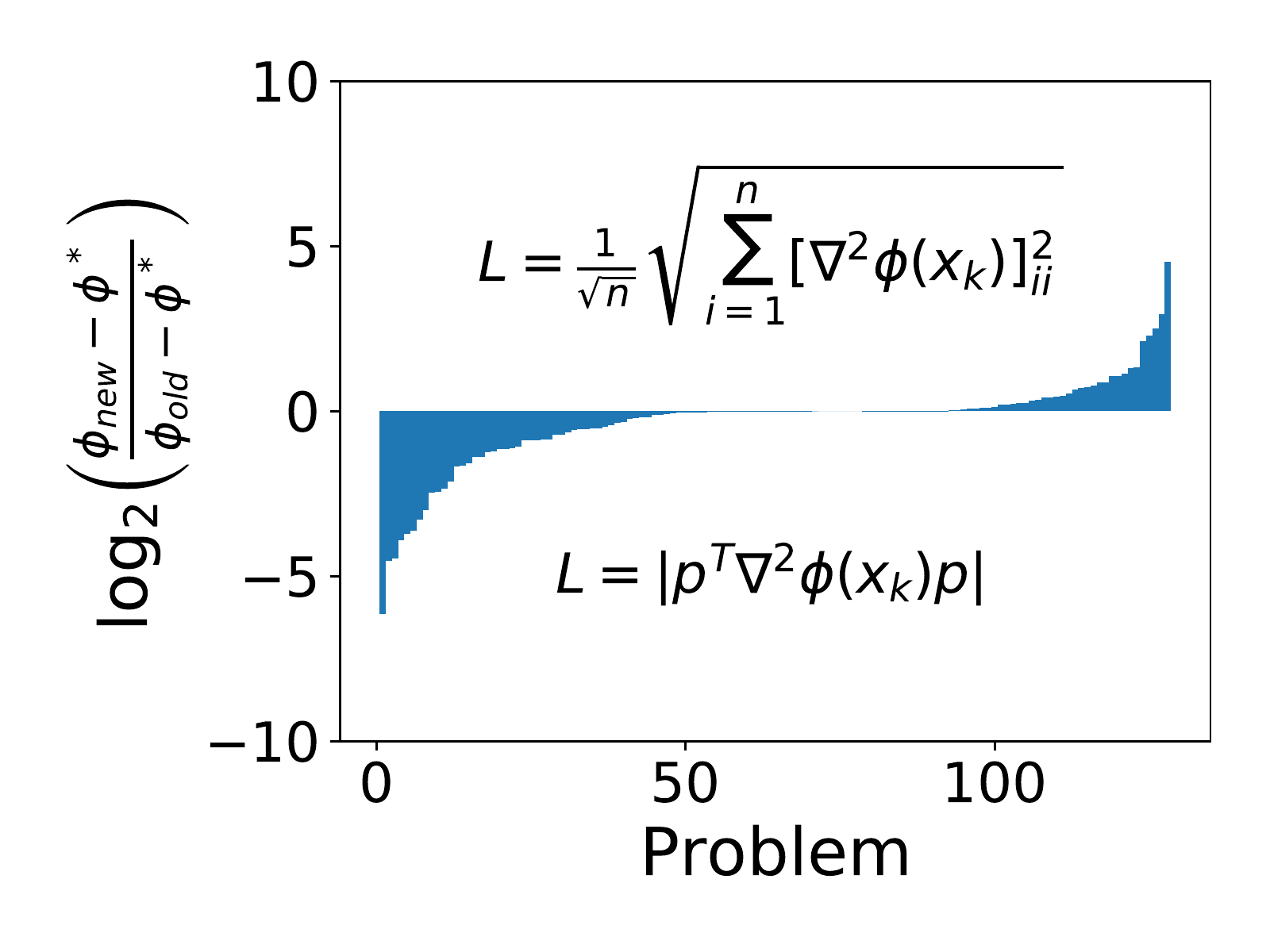}
\caption{{\em Accuracy, Noisy Case with $\sigma_f = 10^{-3}$}. Log-ratio optimality gap profiles comparing forward difference {\sc l-bfgs} with fixed and adaptive Lipschitz estimation schemes. The noise level is $\sigma_f = 10^{-3}$, but is representative for $\sigma_f \in \{10^{-1}, 10^{-3}, 10^{-5}, 10^{-7}\}$.}
\label{fig:adaptive vs component L}
\end{figure}

In addition, since estimating $L$ for each direction is most competitive out of the adaptive variants as seen in Figure \ref{fig:adaptive vs component L}, we present the results for a heuristic approach in our main work. However, our investigation reveals that the simple mean or root mean square can provide potentially cheaper alternatives if they can be estimated without knowledge of the componentwise Lipschitz constants.

To our knowledge, there are two weaknesses with the componentwise estimation approach. The first is that for a small subset of problems, it is prone to underestimate the Lipschitz constant, as discussed in Section \ref{ch:unconstrained}. The second weakness is that the approach, if performed at each iteration and approximated through finite-differencing, is far too expensive. We propose a few possible practical heuristics for handling this in the following section.

\subsection{Practical Heuristics for Lipschitz Estimation}

Two heuristics were proposed for estimating the Lipschitz constant. Mor\'e and Wild \cite{more2012estimating} proposed a simple heuristic for estimating the bound on the second derivative of a univariate function with noise, as described in Section \ref{ch:unconstrained}. Gill, et al. \cite{GillMurrSaunWrig83} proposed a similar procedure that instead enforces that the relative cancellation error lies within an interval
\begin{equation} \label{eq:GMSW heur}
\frac{4 \epsilon_f}{|\Delta(h)|} \in [0.001, 0.1].
\end{equation}
Note that the lower bound on the interval in \eqref{eq:GMSW heur} corresponds to \eqref{eq:more-wild error}. 

Both of these heuristics were proposed with particular initial guesses of $h$ and additional conditions to handle certain cases where the methods can fail. These heuristics were employed only at the beginning of the iteration for each variable component, then remains fixed for the entire run. 

The only work to our knowledge that employs re-estimation of the Lipschitz constant for finite-differencing derivative-free optimization is Berahas, et al. \cite{berahas2019derivative}. In his work, the Lipschitz constant is re-estimated whenever the line search fails and the recovery procedure is triggered. Similarly, we will use the Mor\'e and Wild heuristic to estimate the Lipschitz constant, but only re-estimate the Lipschitz constant when $\alpha < 0.5$ after the first iteration, as described in Procedure I in Section \ref{ch:unconstrained}. We compare two variants of the Lipschitz estimation procedure:

\begin{enumerate}
\item \texttt{Component MW}: We use the Mor\'e and Wild heuristic to estimate the Lipschitz constant with respect to each component. This is performed at the first iteration and subsequent iterations for line search methods when the line search from the prior iteration gives a steplength $\alpha_k < 0.5$. When estimating the directional derivative along the search direction $p_k$, we use the root mean square of the component-wise estimates. 
\item \texttt{Random MW}: We use the Mor\'e and Wild heuristic to estimate the Lipschitz constant along a random direction sampled uniformly from a sphere, i.e., $p \sim S(0, I)$. This is performed at the first iteration and subsequent iterations for line search methods when the line search from the prior iteration gives a steplength $\alpha_k < 0.5$.
\end{enumerate}

\begin{figure}[ht]
\centering
\includegraphics[width=0.32\textwidth]{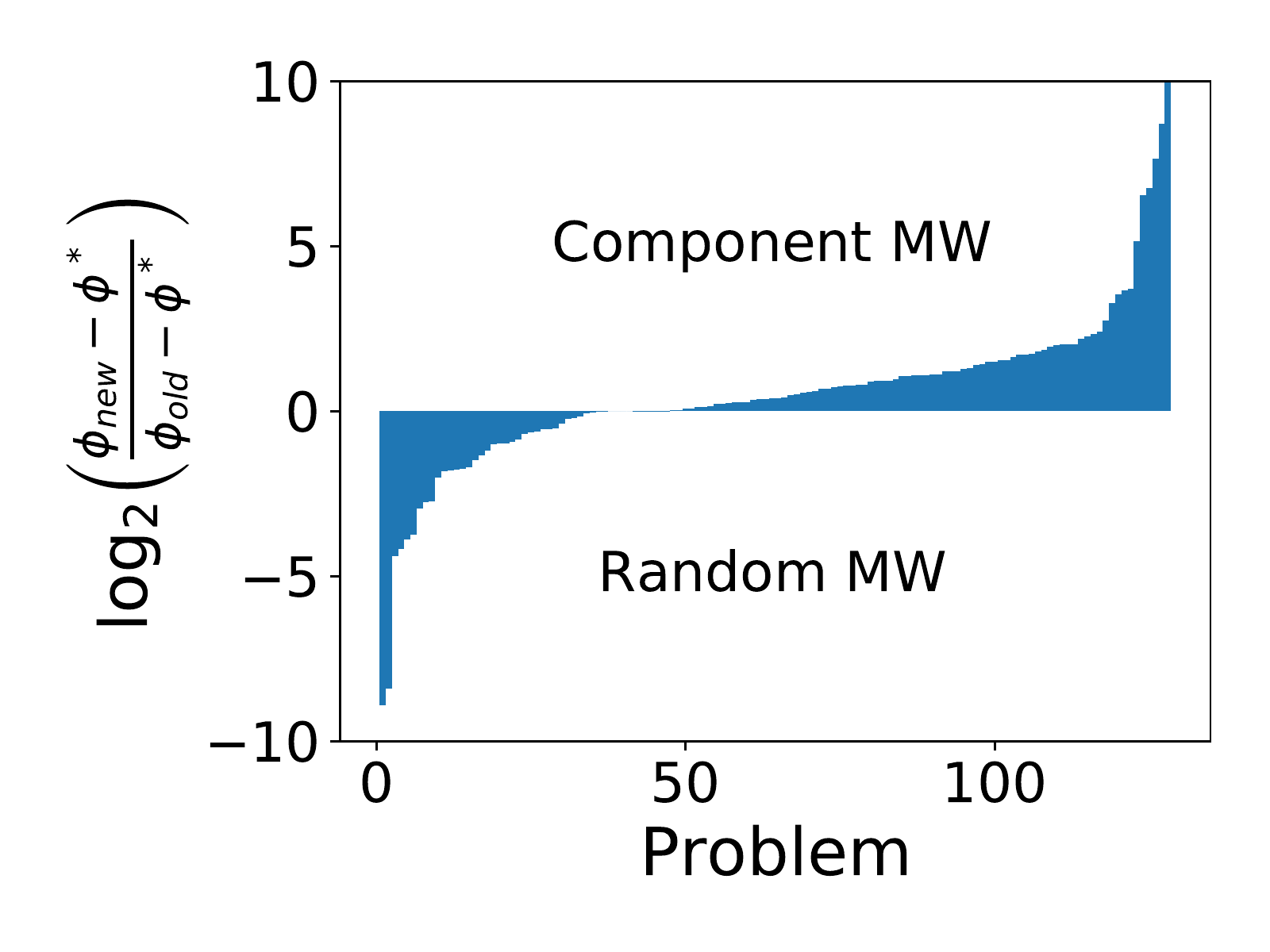}
\caption{{\em Accuracy, Noisy Case with $\sigma_f = 10^{-3}$.} Log-ratio optimality gap profiles comparing forward difference {\sc l-bfgs} with \texttt{component MW} and \texttt{random MW} Lipschitz estimation schemes.}
\label{fig:MW comparison}
\end{figure}

\begin{figure}[ht]
\centering
\includegraphics[width=0.32\textwidth]{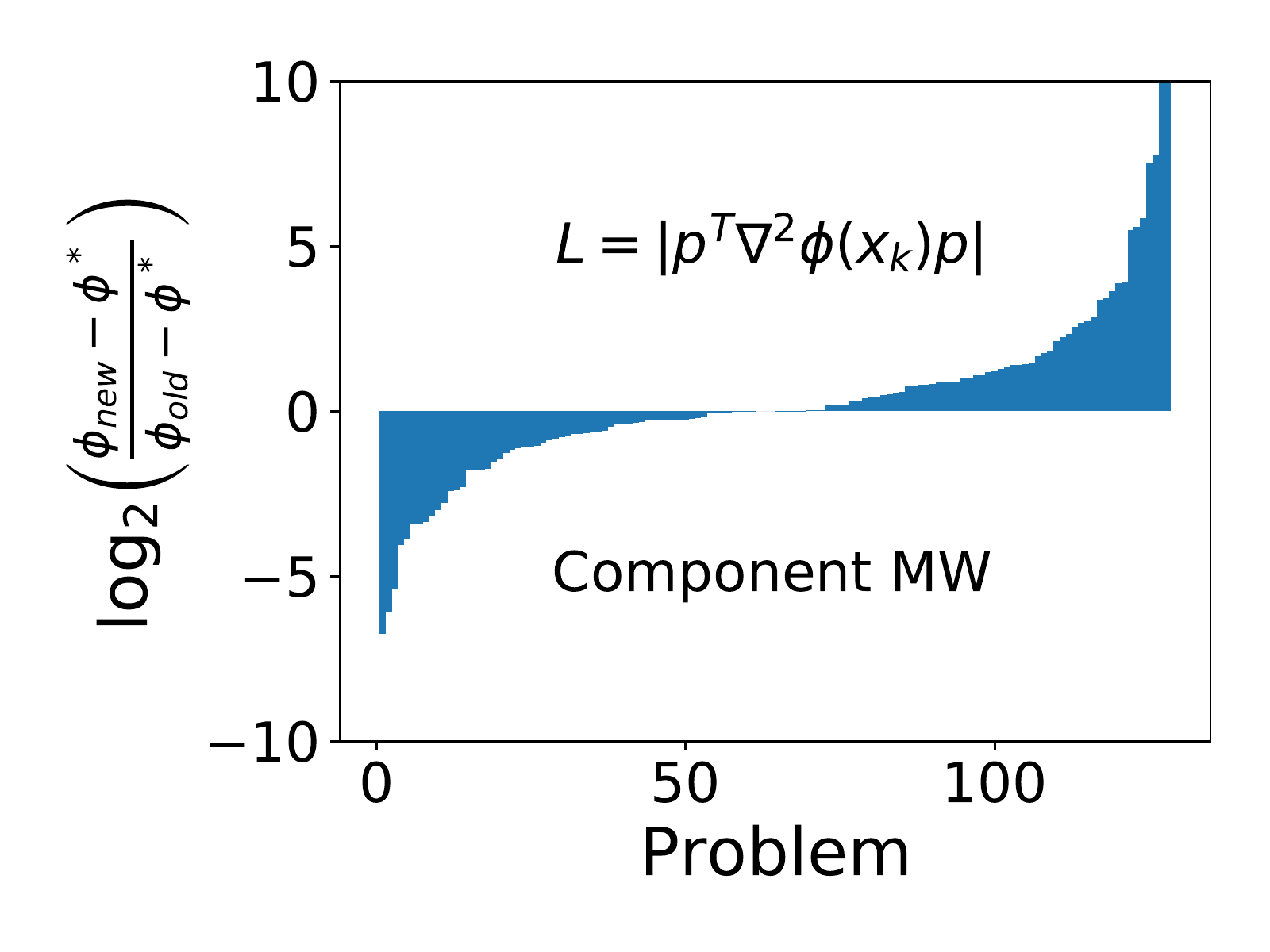}
\includegraphics[width=0.32\textwidth]{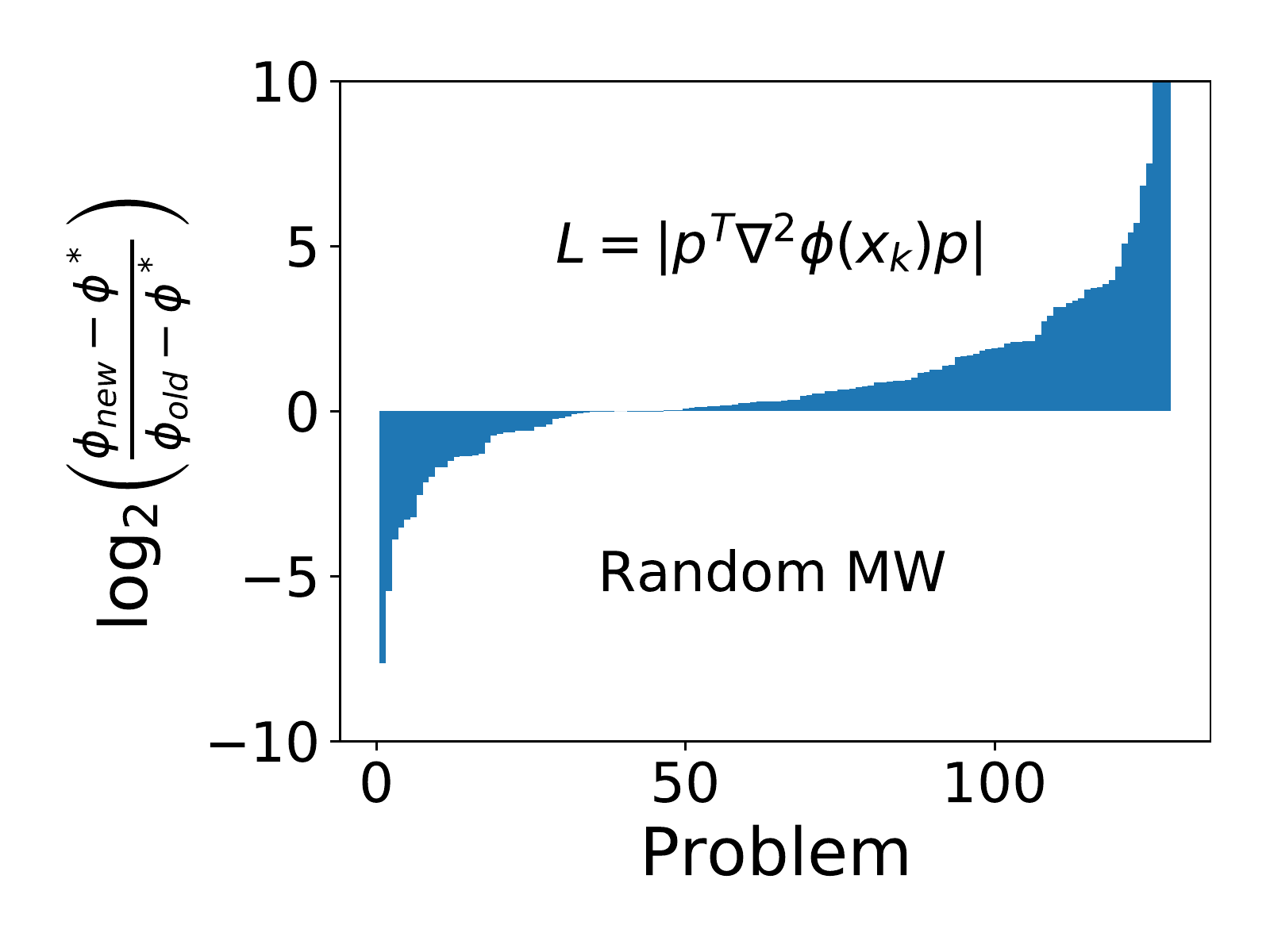}
\caption{{\em Accuracy, Noisy Case with $\sigma_f = 10^{-3}$}. Log-ratio optimality gap profiles comparing forward difference {\sc l-bfgs} with theoretical componentwise Lipschitz estimates and the Mor\'e and Wild Lipschitz estimation schemes.}
\label{fig:MW vs theory}
\end{figure}

\begin{figure}[ht]
\centering
\includegraphics[width=0.32\textwidth]{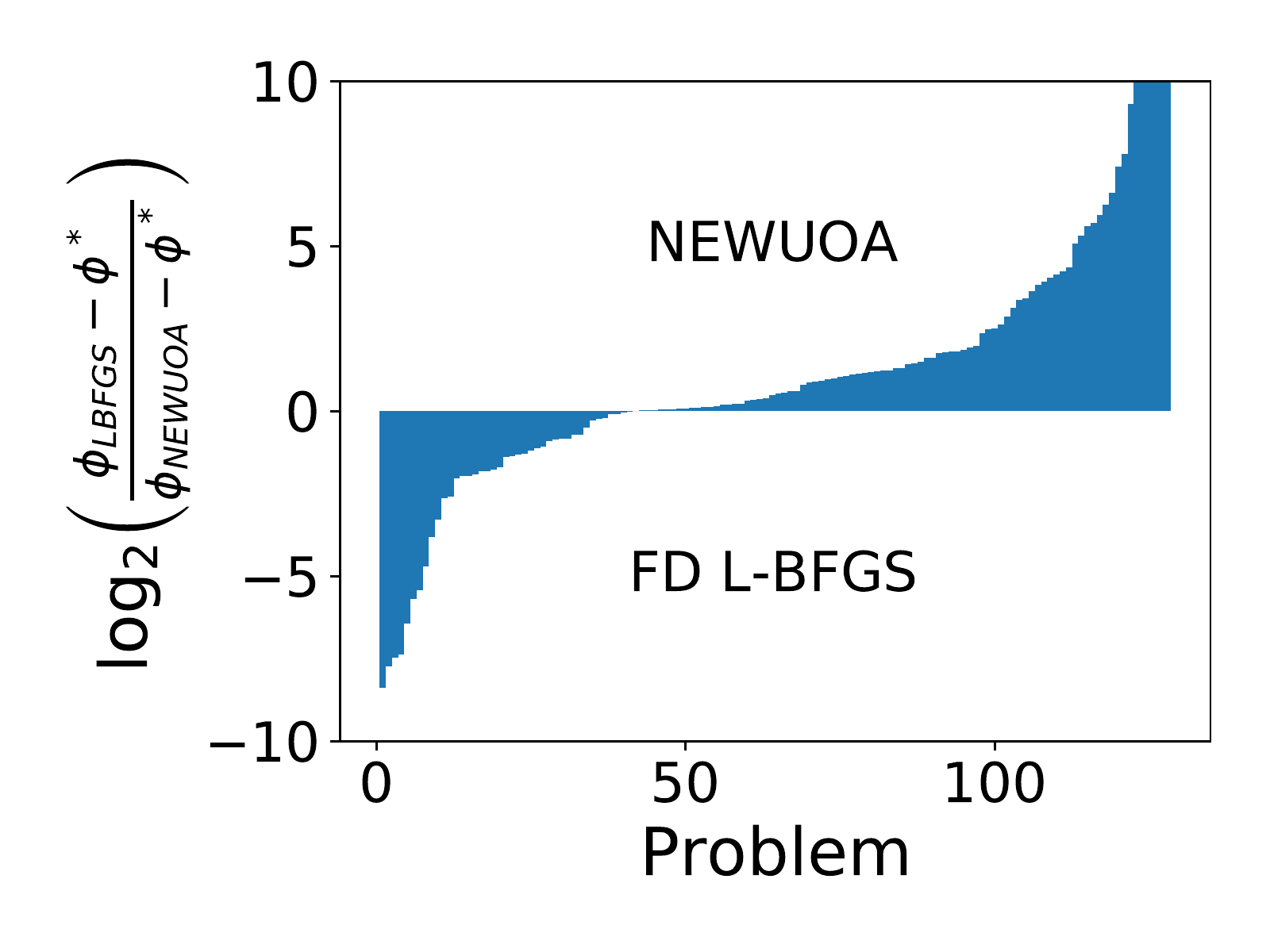}
\includegraphics[width=0.32\textwidth]{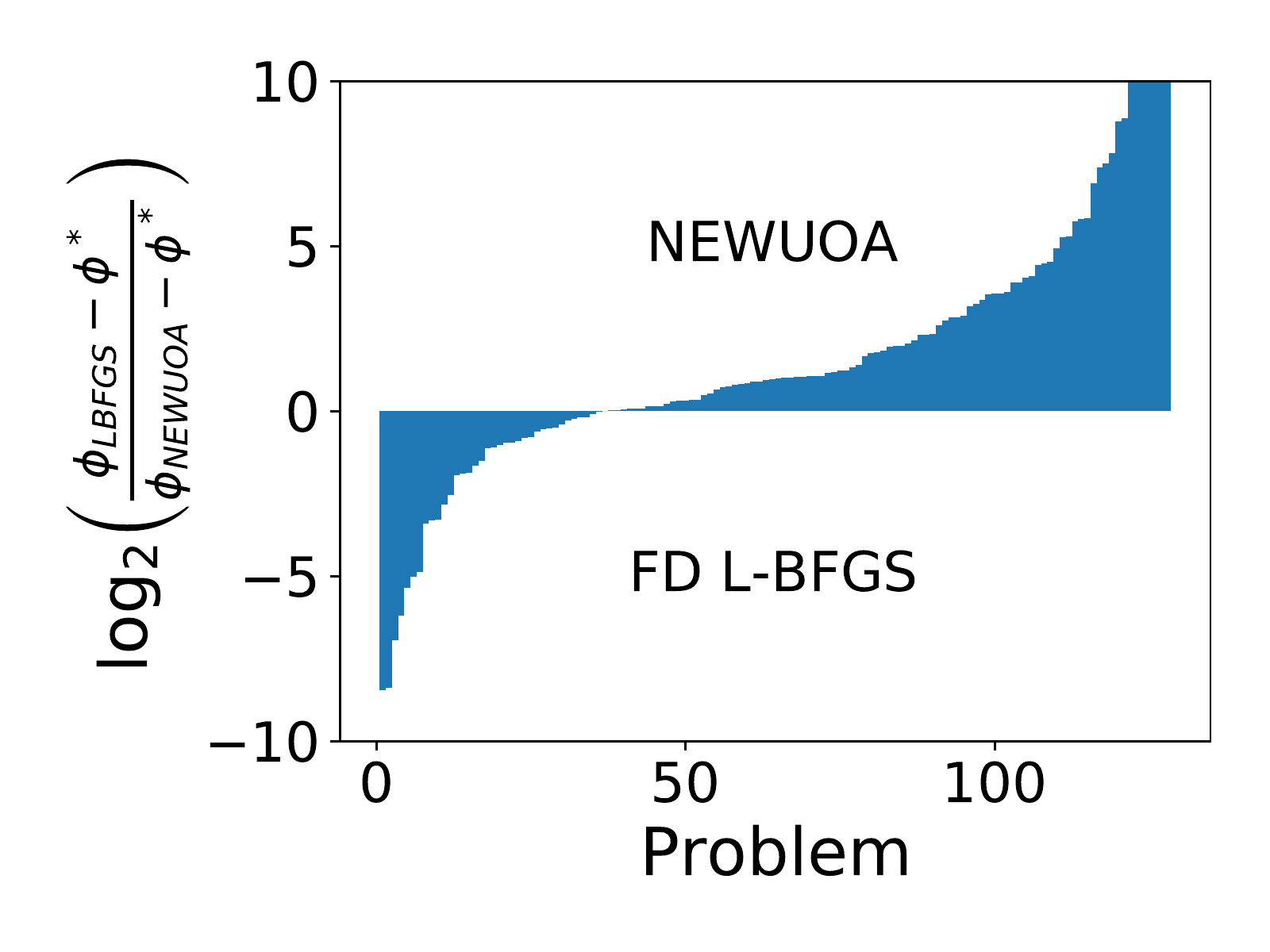}
\caption{{\em Accuracy, Noisy Case with $\sigma_f = 10^{-3}$}. Log-ratio optimality gap profiles comparing {\sc newuoa} against forward difference {\sc l-bfgs} with the Mor\'e and Wild Lipschitz estimation schemes. We compare {\sc l-bfgs} with \texttt{Component MW} (left) and \texttt{Random MW} (right).}
\label{fig:new vs MW}
\end{figure}

We compare both Mor\'e and Wild heuristics against {\sc newuoa} and theoretical componentwise Lipschitz estimation in Figures \ref{fig:MW comparison}, \ref{fig:MW vs theory}, and \ref{fig:new vs MW}. In general, \texttt{Component MW} gives a better solution than \texttt{Random MW}. When we compare these methods against {\sc newuoa} in Figure \ref{fig:new vs MW}, the performance of the methods are relatively the same as the theoretical Lipschitz estimates. 

\section{Investigation of Parameters for NEWUOA}
\label{sec:newuoa}

In this section, we empirically investigate the influence of the parameters for {\sc newuoa} for the noisy setting. Although these parameters have been optimized for the noiseless setting, it is not generally known how these parameters may impact its performance on noisy unconstrained problems. By default, {\sc newuoa} employs $p = 2n + 1$ interpolation points when constructing the quadratic model at each iteration with an initial trust region radius of $\rho_{\text{beg}} = 1$ and a final trust region radius of $\rho_{\text{end}} = 10^{-6}$.

\subsection{Number of Interpolation Points}

In Section \ref{ch:unconstrained}, we observed that {\sc newuoa} is able to converge to a better quality neighborhood than {\sc l-bfgs} with forward differencing, but inferior to central differencing. Since {\sc newuoa} by default employs $p = 2n + 1$ points, it is natural to both ask if: (1) decreasing the number of interpolation points would yield a less accurate quadratic model, with a potentially less accurate gradient; and (2) increasing the number of points used in the interpolation may improve the accuracy of the solution to be competitive with central differencing. 

\begin{figure}[ht]
\centering
\includegraphics[width=0.32\textwidth]{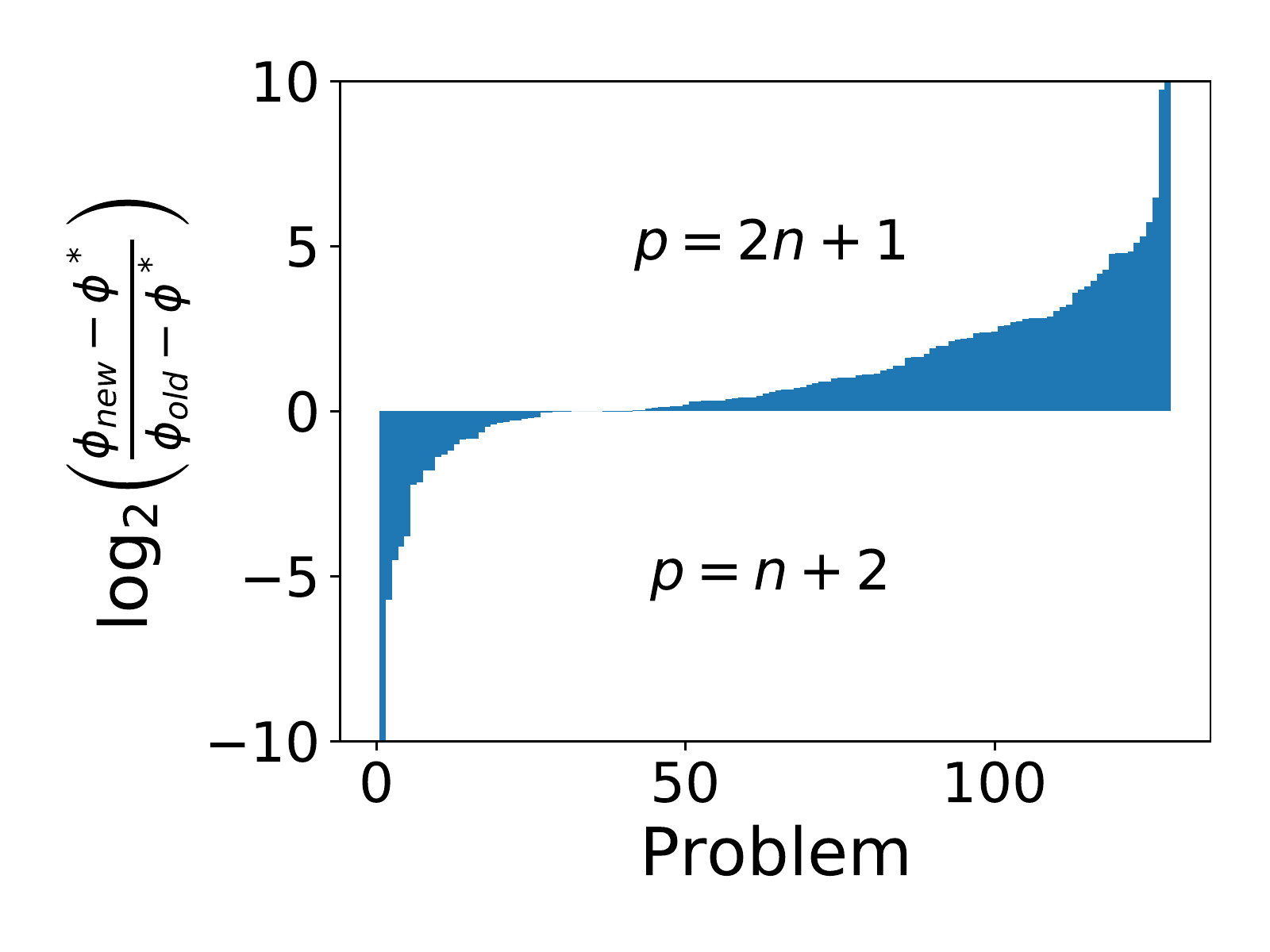}
\includegraphics[width=0.32\textwidth]{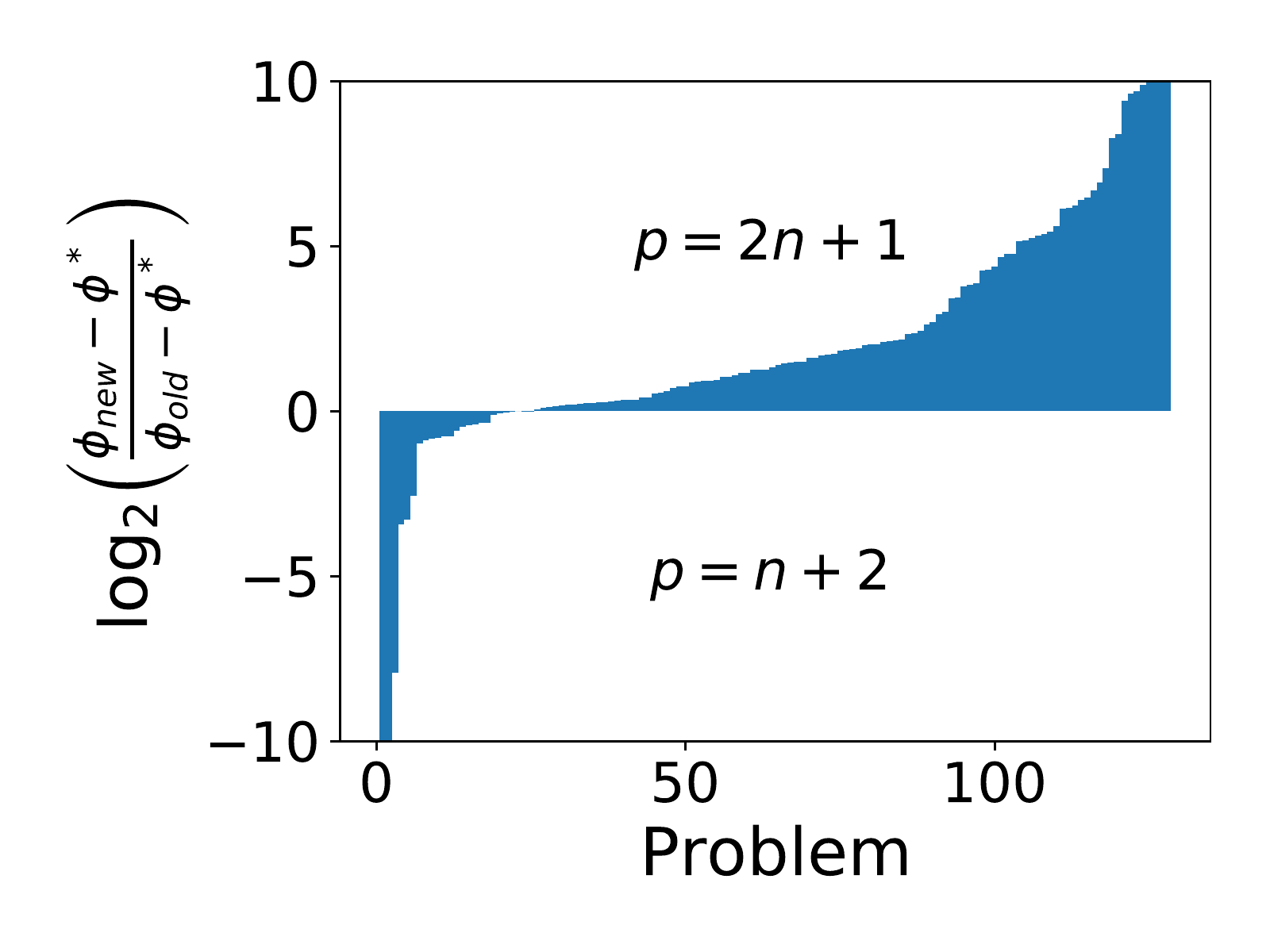}\\
\includegraphics[width=0.32\textwidth]{{figures/unconstrained/newuoa_tests/morales_obj_new_new_n+2_1e-05}.pdf}
\includegraphics[width=0.32\textwidth]{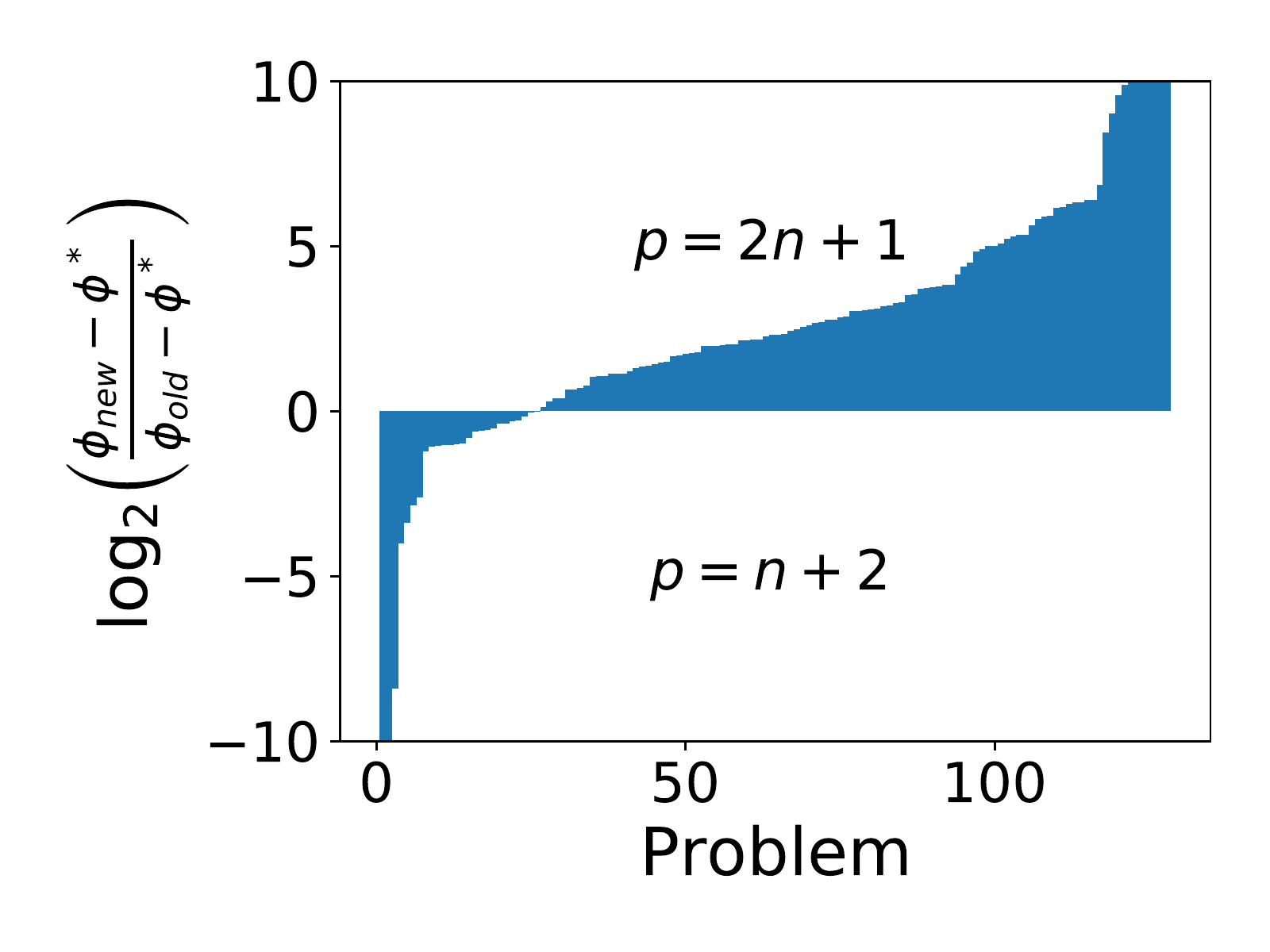}
\caption{Noisy Case. Log-ratio optimality gap profiles comparing {\sc newuoa} with $p = 2n + 1$ and $p = n + 2$ points. The noise levels are $\sigma_f = 10^{-1}$ (top left), $\sigma_f = 10^{-3}$ (top right), $10^{-5}$ (bottom left), and $10^{-7}$ (bottom right).}
\label{fig:2n+1 vs n+2 obj}
\end{figure}

\begin{figure}[ht]
\centering
\includegraphics[width=0.32\textwidth]{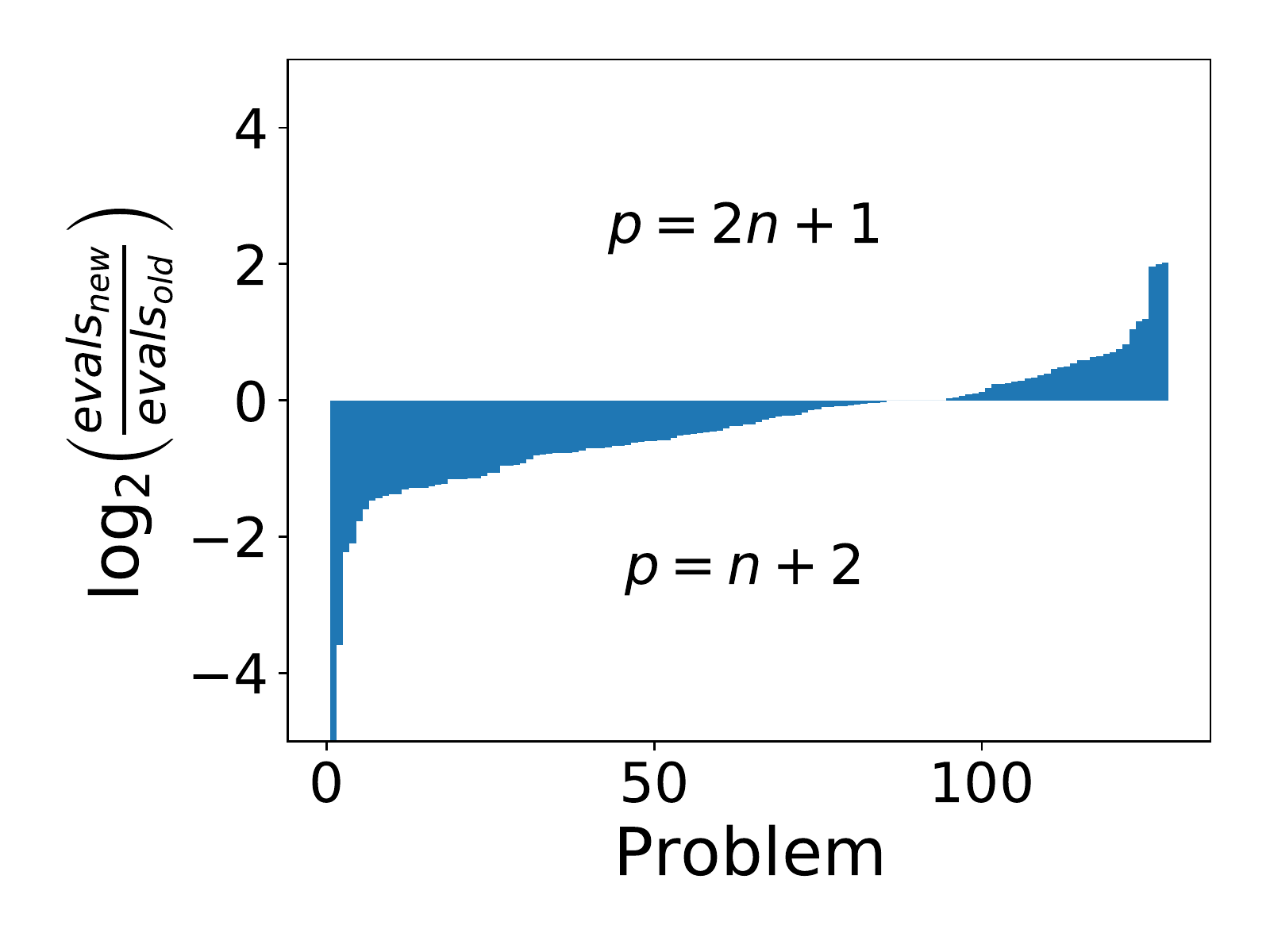}
\includegraphics[width=0.32\textwidth]{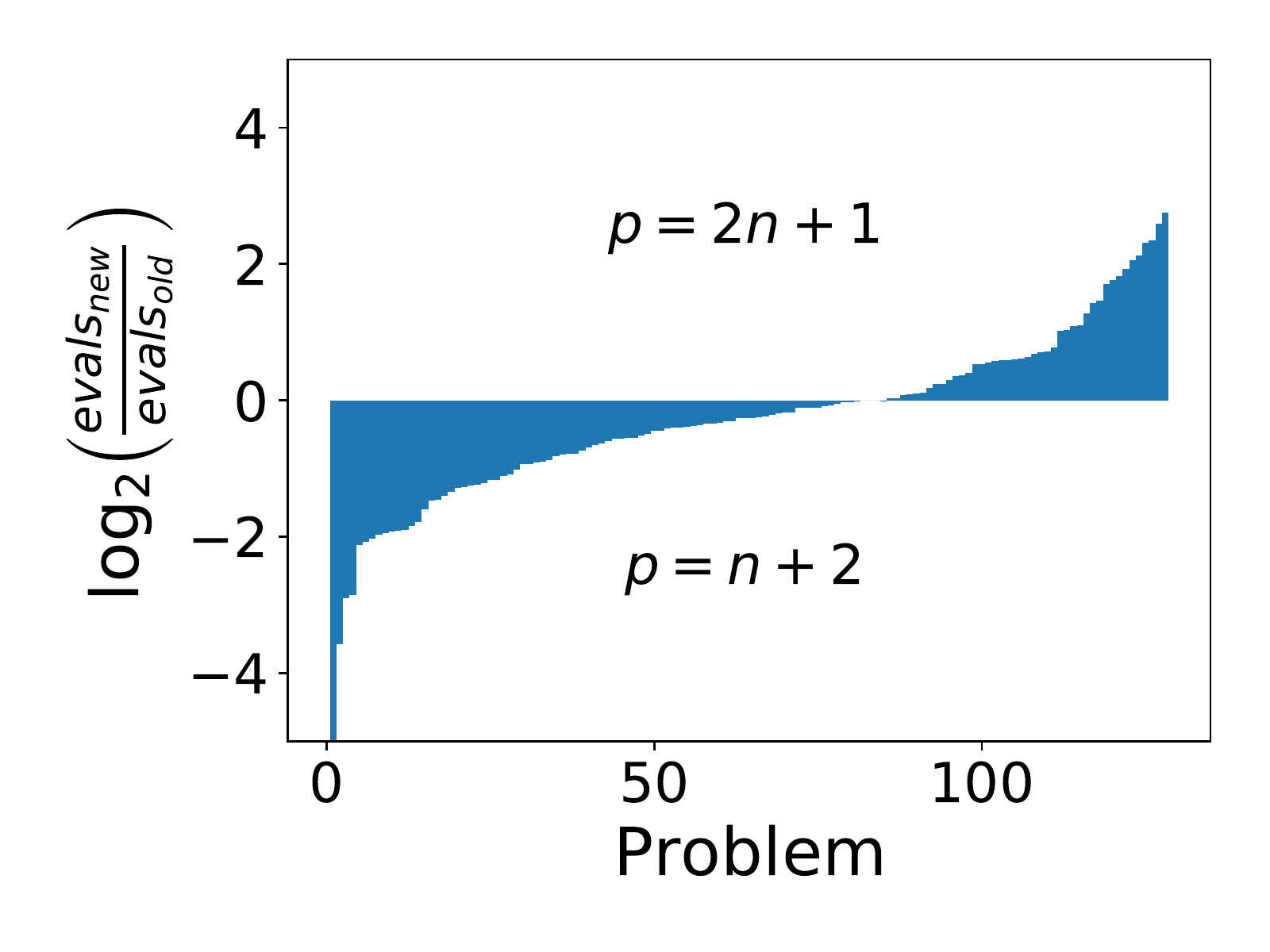}\\
\includegraphics[width=0.32\textwidth]{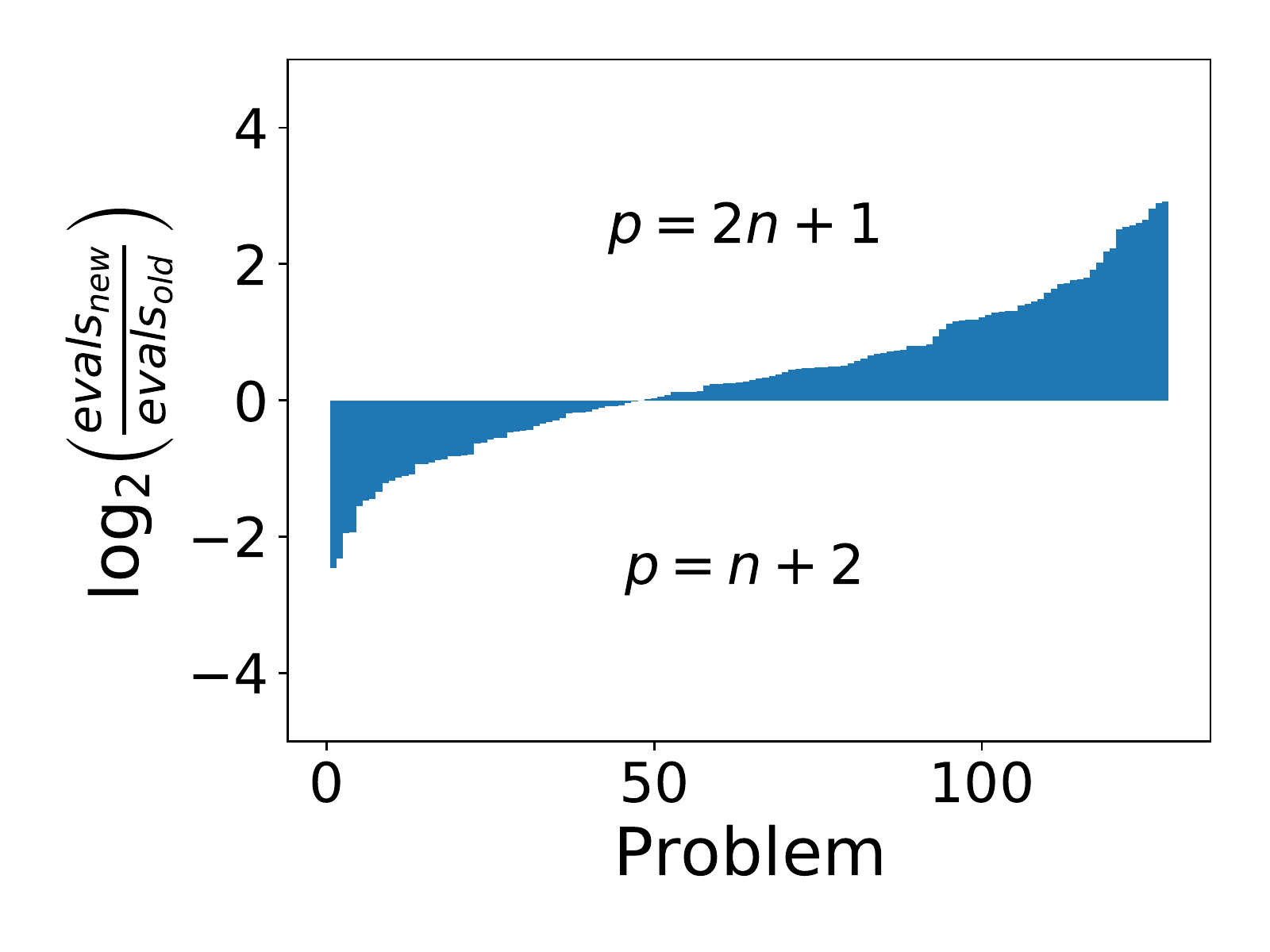}
\includegraphics[width=0.32\textwidth]{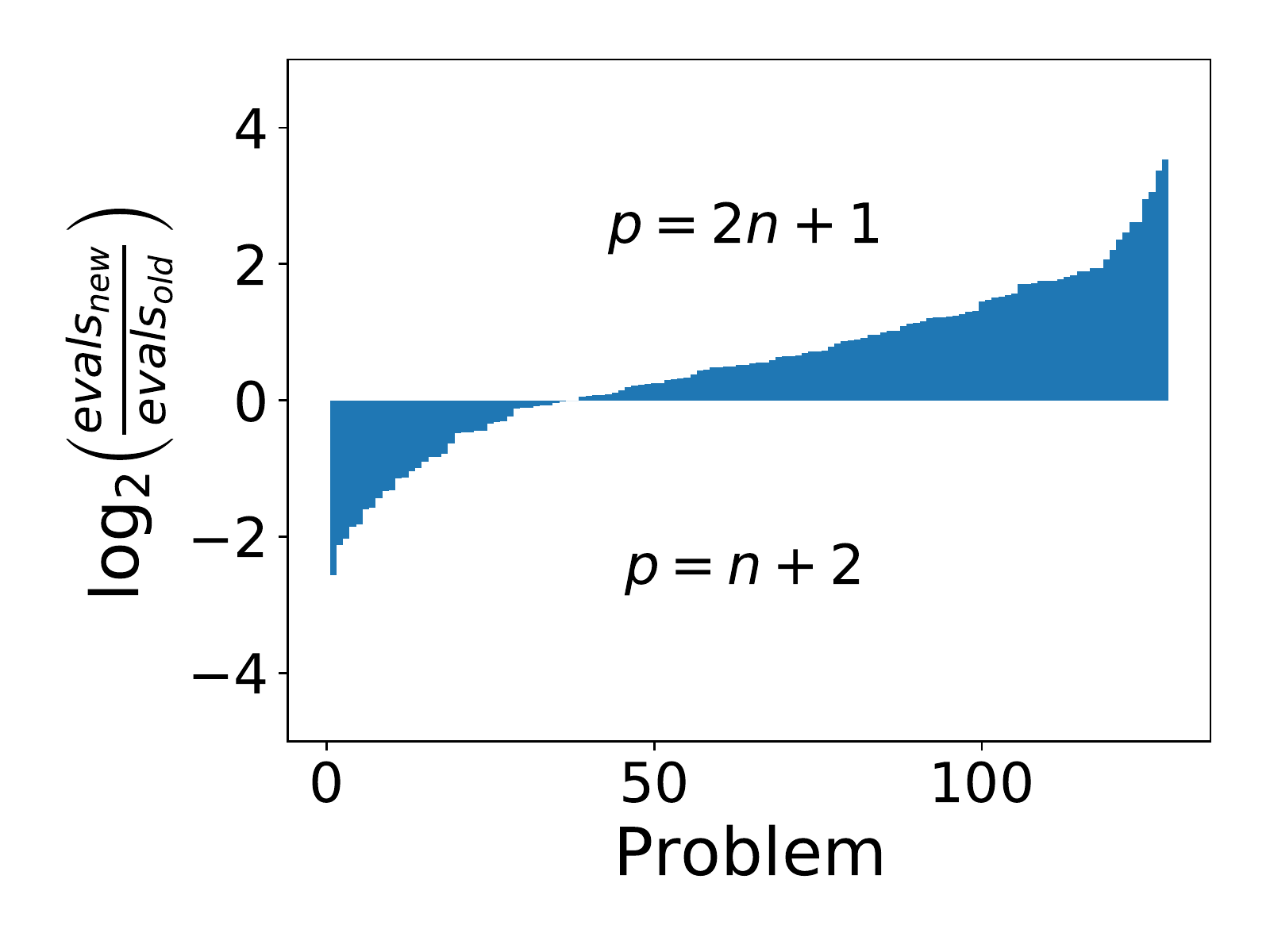}
\caption{Noisy Case. Log-ratio function evaluation profiles comparing {\sc newuoa} with $p = 2n + 1$ and $p = n + 2$ points. The noise levels are $\sigma_f = 10^{-1}$ (top left), $\sigma_f = 10^{-3}$ (top right), $10^{-5}$ (bottom left), and $10^{-7}$ (bottom right).}
\label{fig:2n+1 vs n+2 fevals}
\end{figure}

\begin{figure}[ht]
\centering
\includegraphics[width=0.32\textwidth]{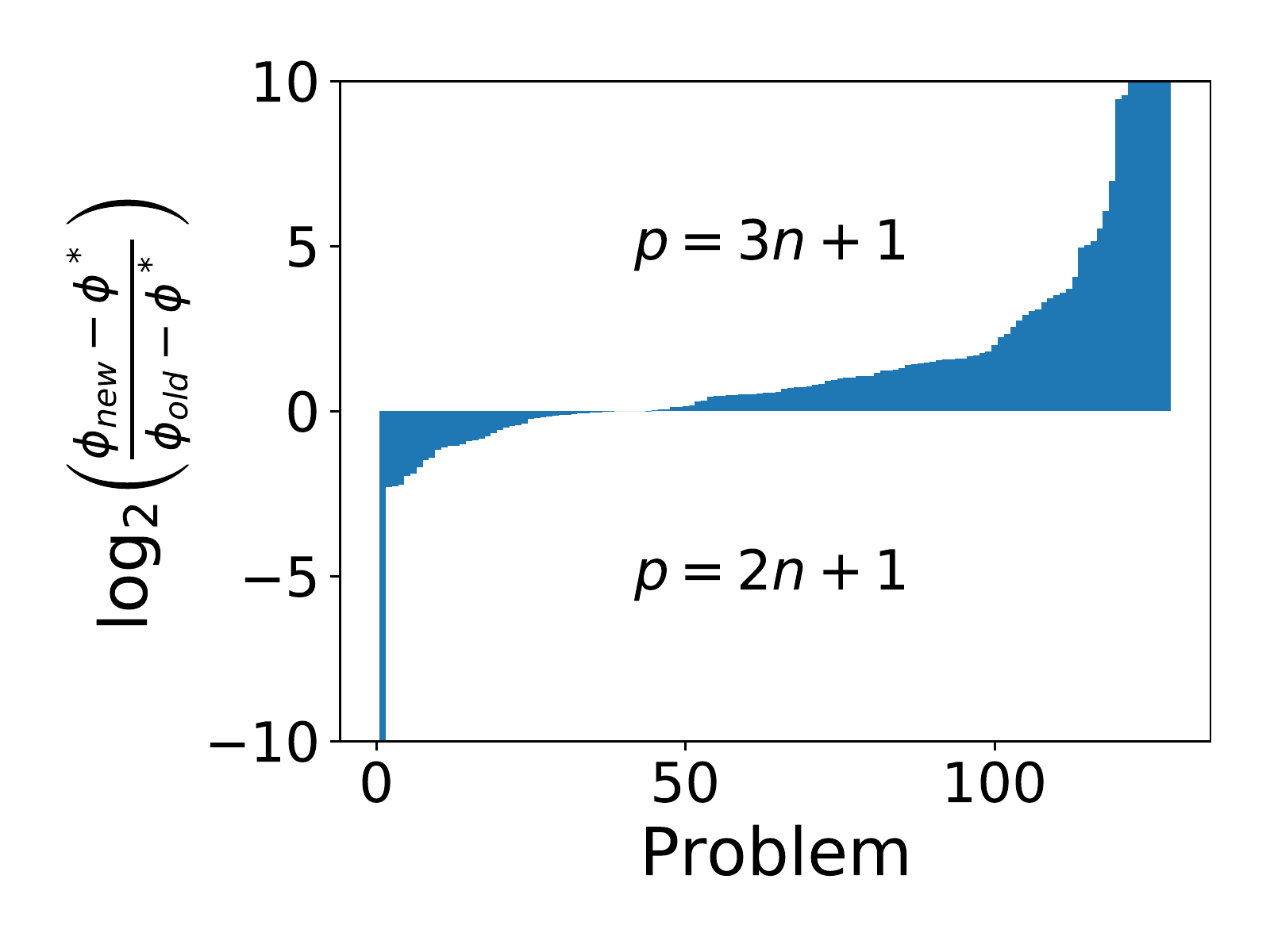}
\includegraphics[width=0.32\textwidth]{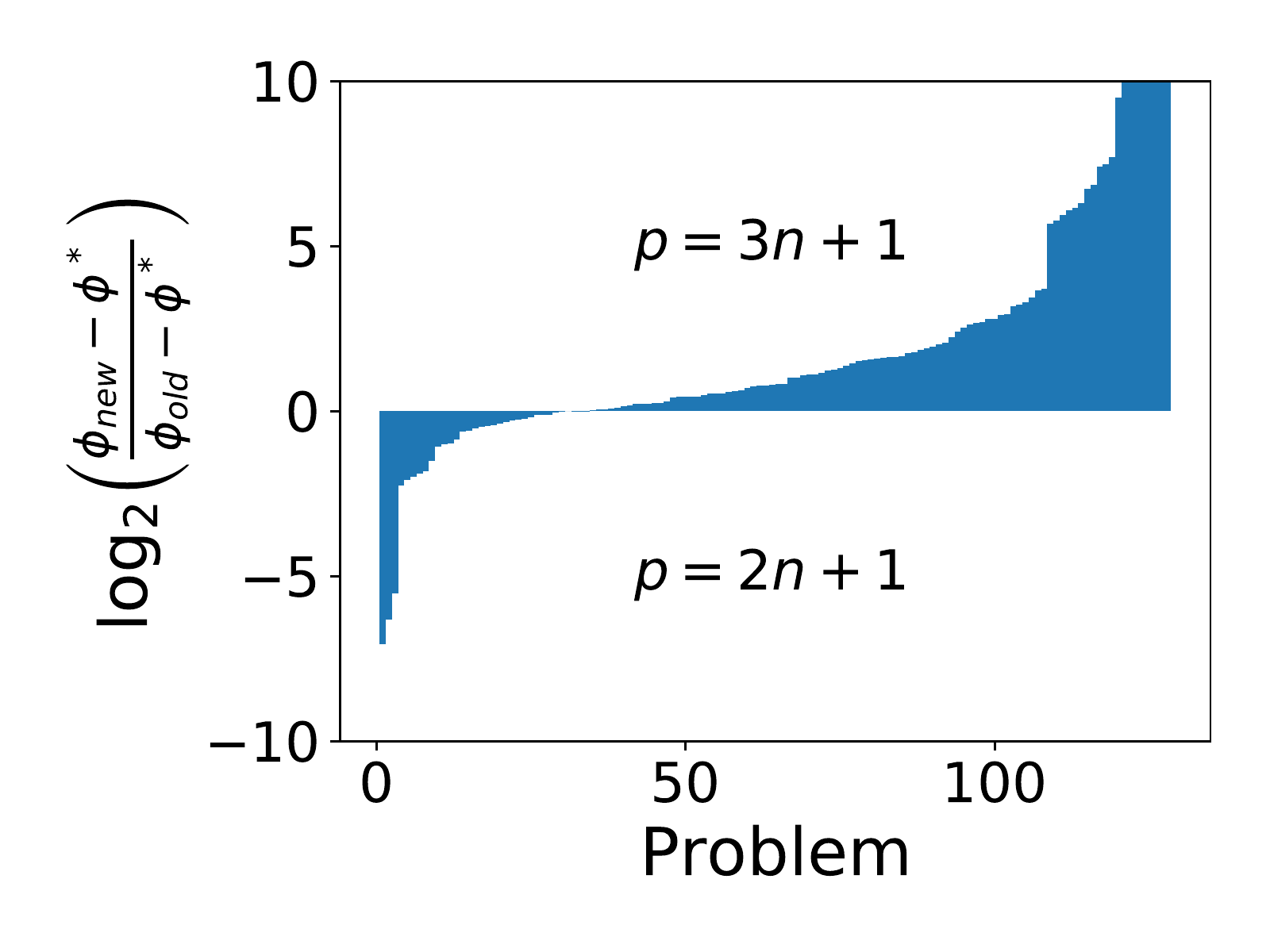}\\
\includegraphics[width=0.32\textwidth]{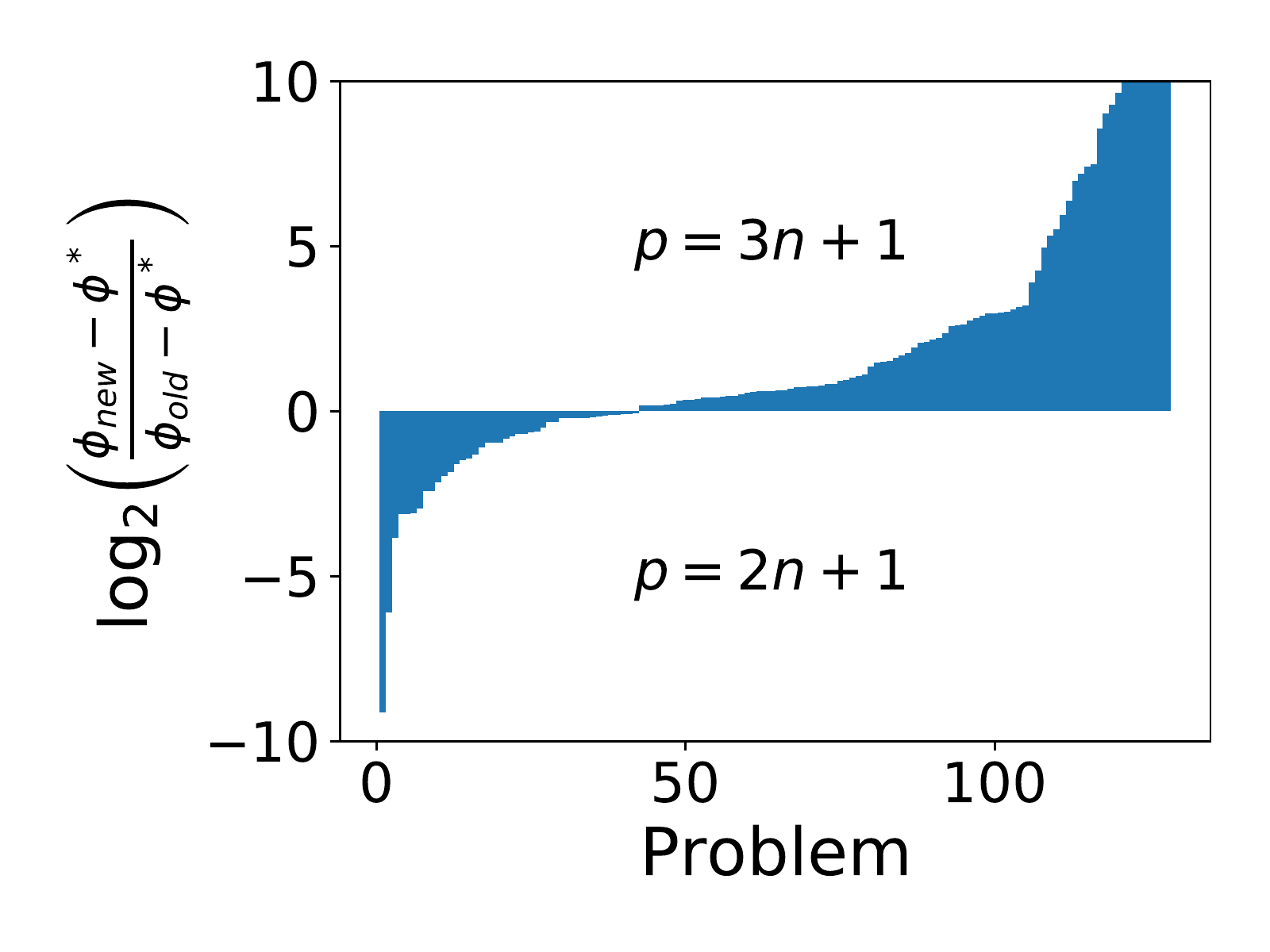}
\includegraphics[width=0.32\textwidth]{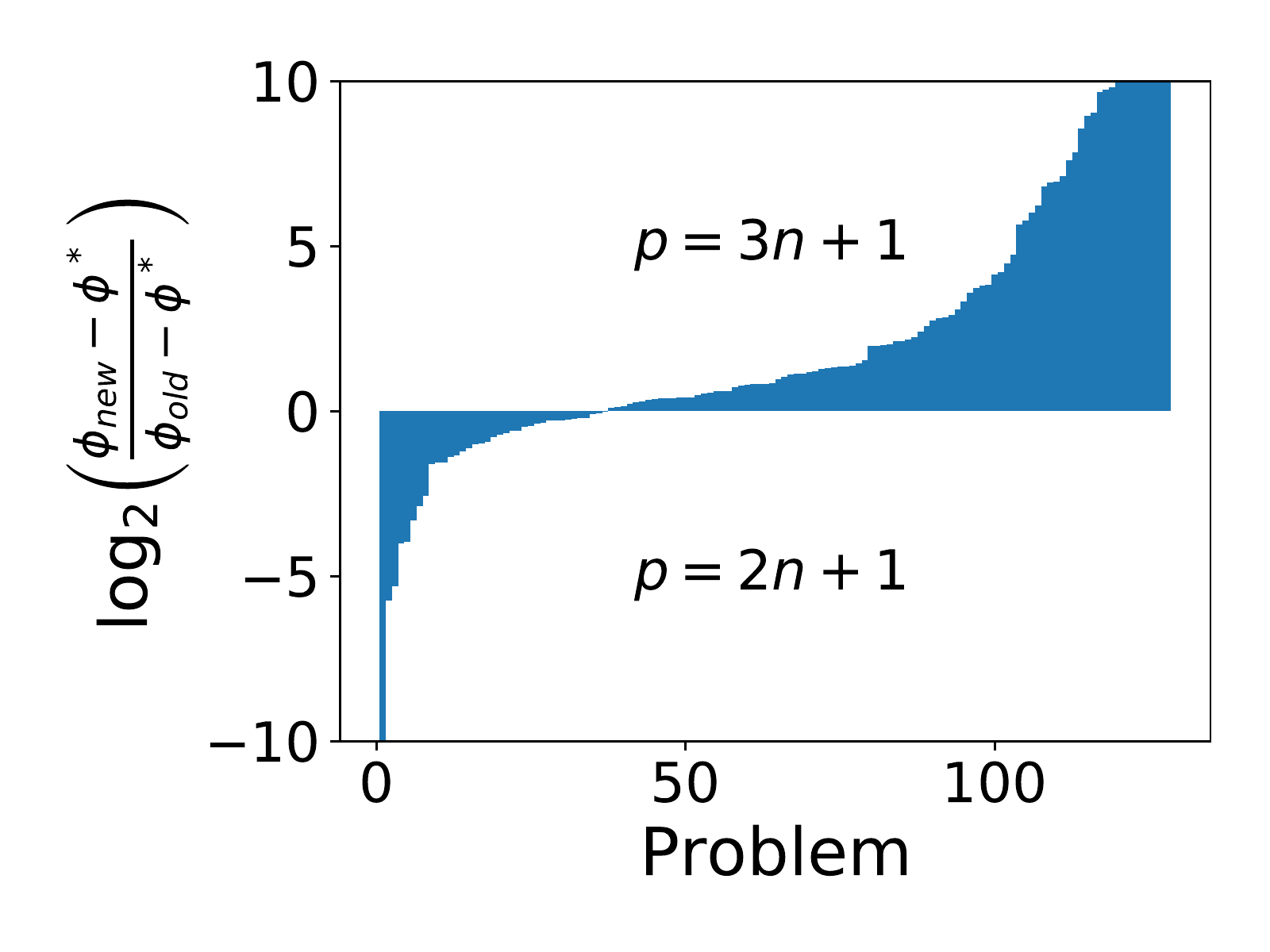}
\caption{Noisy Case. Log-ratio optimality gap profiles comparing {\sc newuoa} with $p = 3n + 1$ and $p = 2n + 1$ points. The noise levels are $\sigma_f = 10^{-1}$ (top left), $\sigma_f = 10^{-3}$ (top right), $10^{-5}$ (bottom left), and $10^{-7}$ (bottom right).}
\label{fig:3n+1 vs 2n+1 obj}
\end{figure}

\begin{figure}[ht]
\centering
\includegraphics[width=0.32\textwidth]{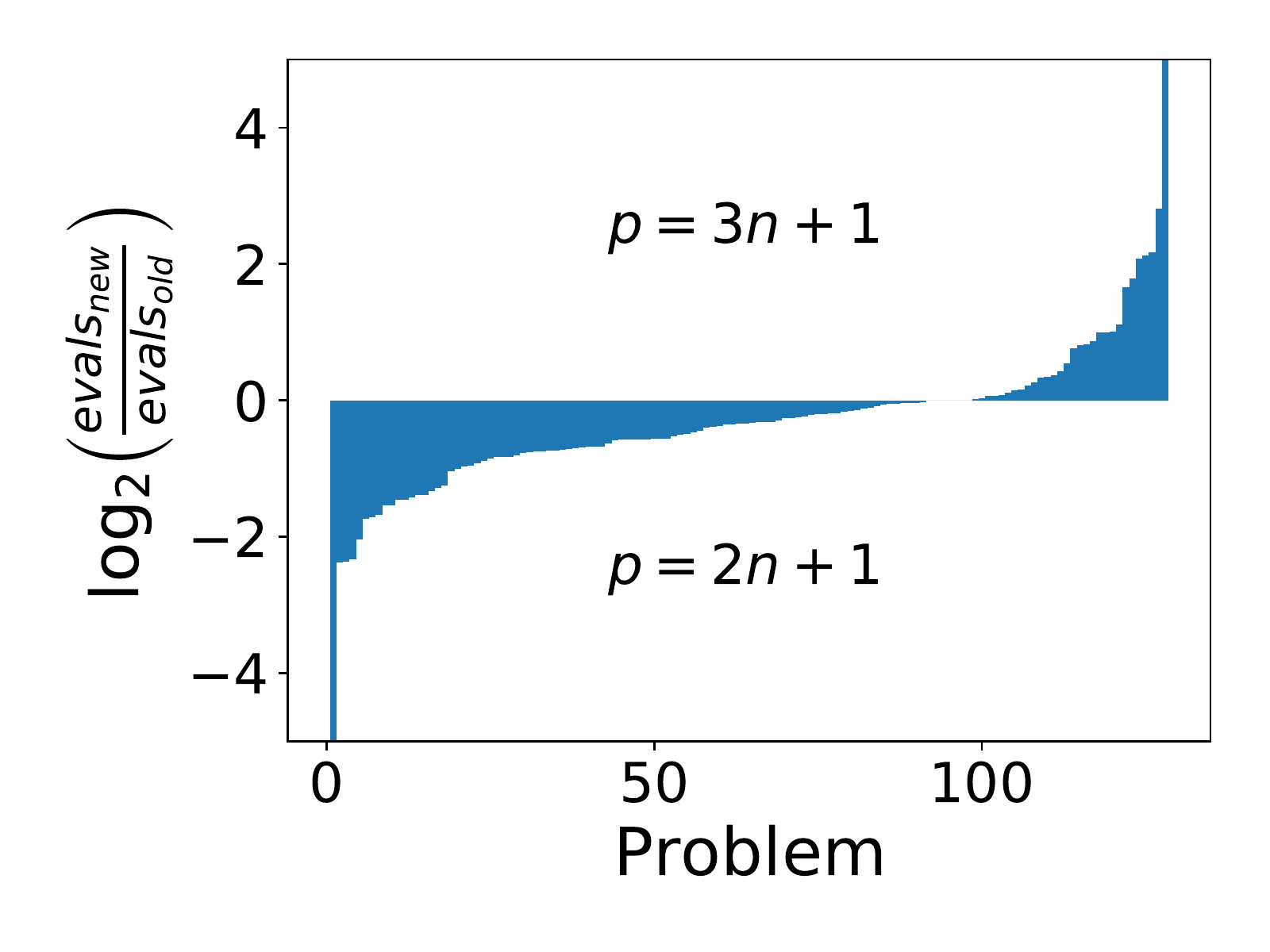}
\includegraphics[width=0.32\textwidth]{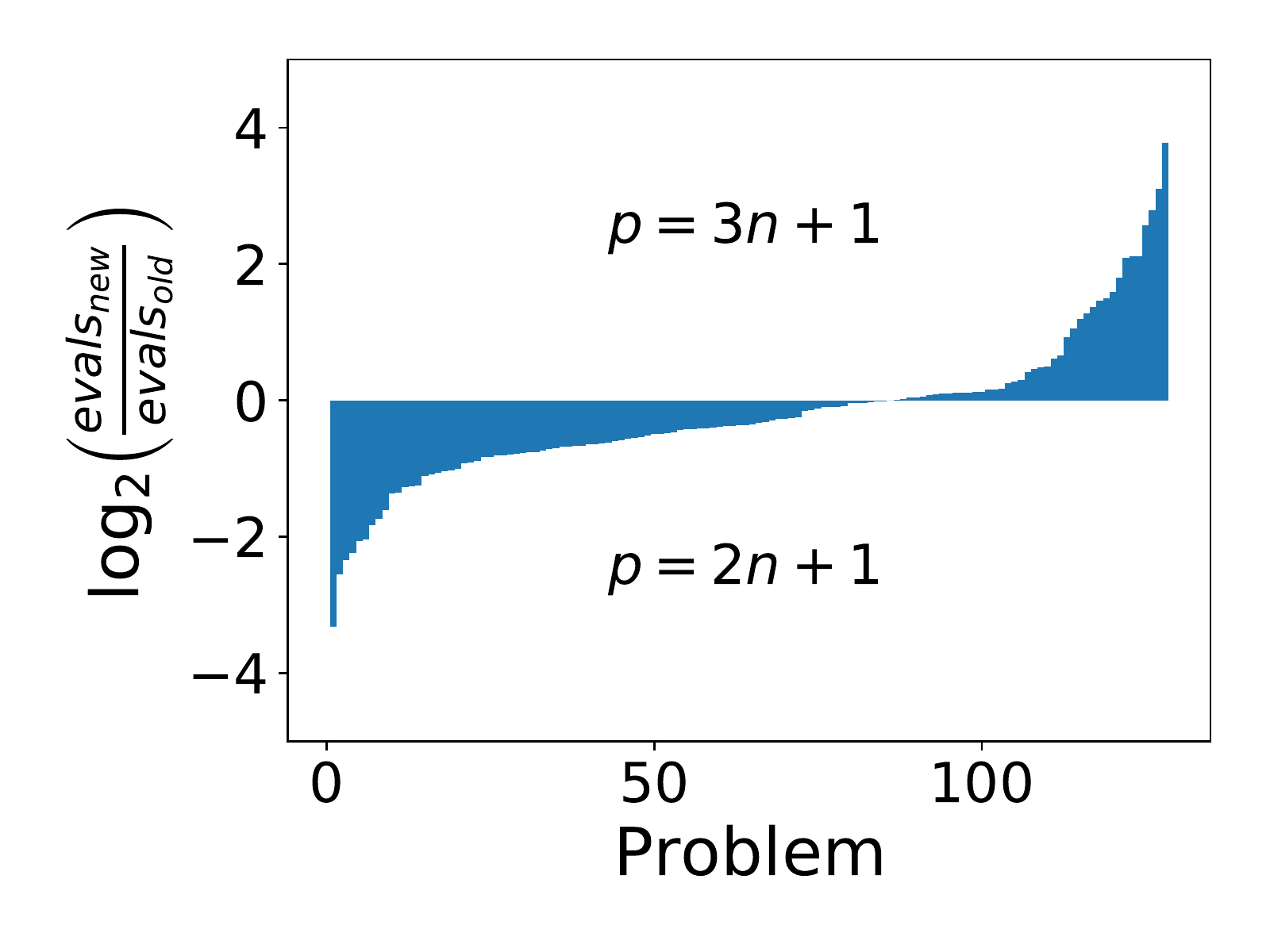} \\
\includegraphics[width=0.32\textwidth]{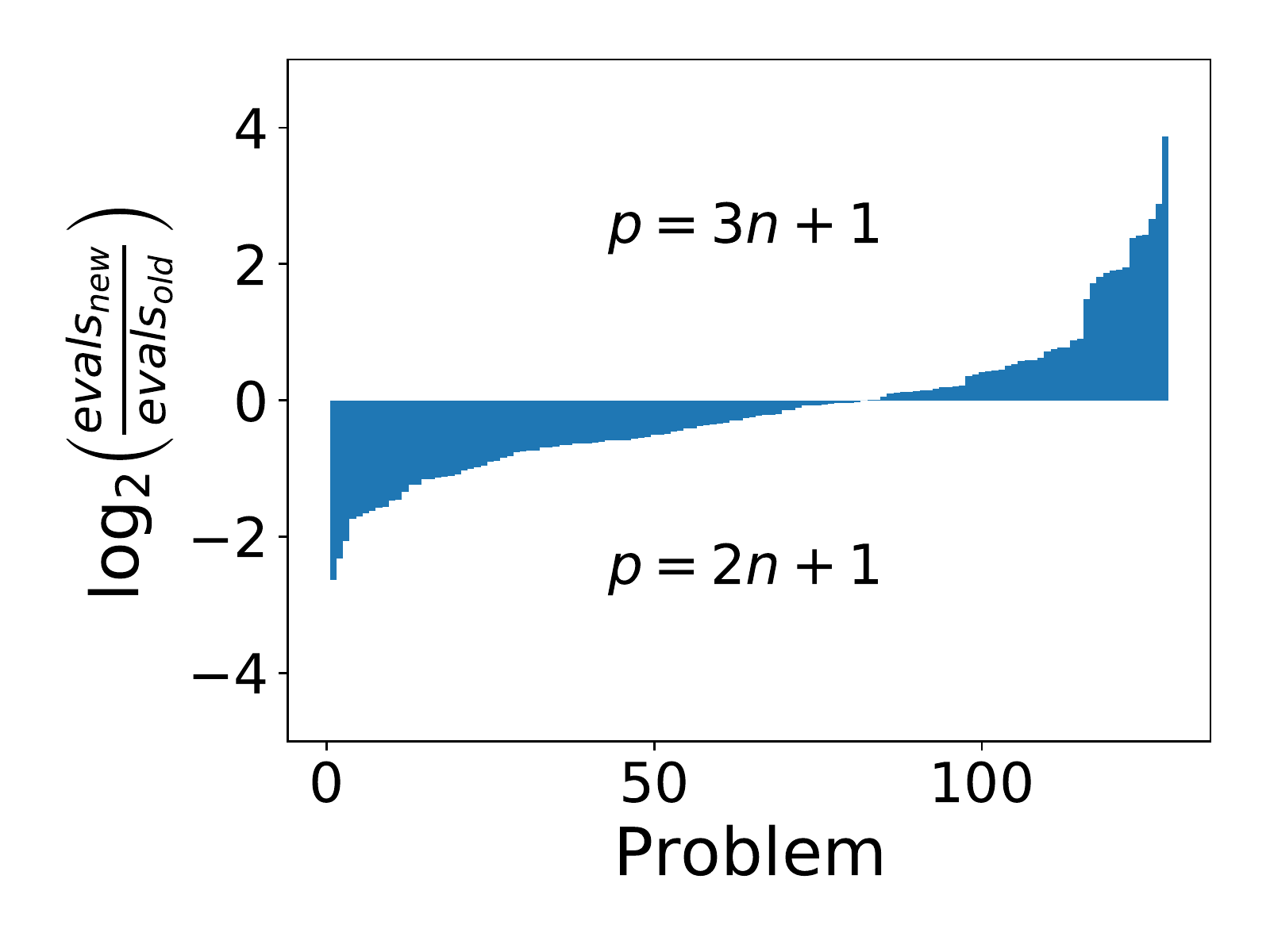}
\includegraphics[width=0.32\textwidth]{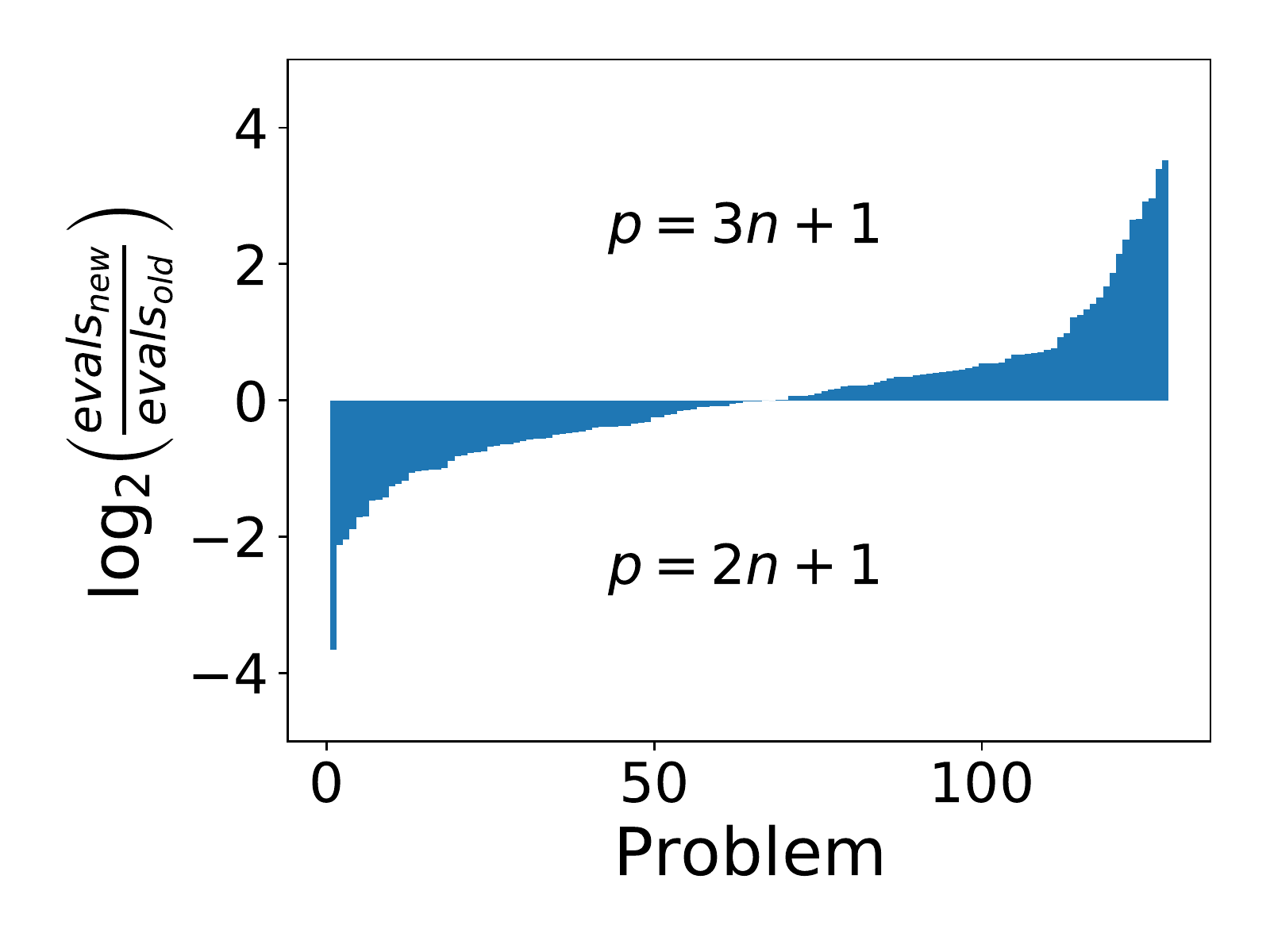}
\caption{Noisy Case. Log-ratio function evaluation profiles comparing {\sc newuoa} with $p = 3n + 1$ and $p = 2n + 1$ points. The noise levels are $\sigma_f = 10^{-1}$ (top left), $\sigma_f = 10^{-3}$ (top right), $10^{-5}$ (bottom left), and $10^{-7}$ (bottom right).}
\label{fig:3n+1 vs 2n+1 fevals}
\end{figure}

To do this, we run {\sc newuoa} with $p = n + 2$ and $p = \min\left\{3n + 1, \frac{(n + 1)(n + 2)}{2}\right\}$ interpolation points and compare their optimality gaps and number of function evaluations. In Figures \ref{fig:2n+1 vs n+2 obj} and \ref{fig:2n+1 vs n+2 fevals}, we report the log-ratio profiles for the objective function and function evaluations when comparing {\sc newuoa} with $p = n + 2$ and $p = 2n + 1$ points. We see that with higher noise levels, {\sc newuoa} with $p = n + 2$ tends to terminate earlier, while for lower noise levels, it is less efficient than {\sc newuoa} with $p = 2n + 1$ points. In terms of solution quality, {\sc newuoa} with $p = 2n + 1$ is able to converge to a higher accuracy in general than {\sc newuoa} with $p = n + 2$ points. This is similarly seen when we compare $p = 3n + 1$ and $p = 2n + 1$ points in Figures \ref{fig:3n+1 vs 2n+1 obj} and \ref{fig:3n+1 vs 2n+1 fevals}.

\begin{figure}[ht]
\centering
\includegraphics[width=0.32\textwidth]{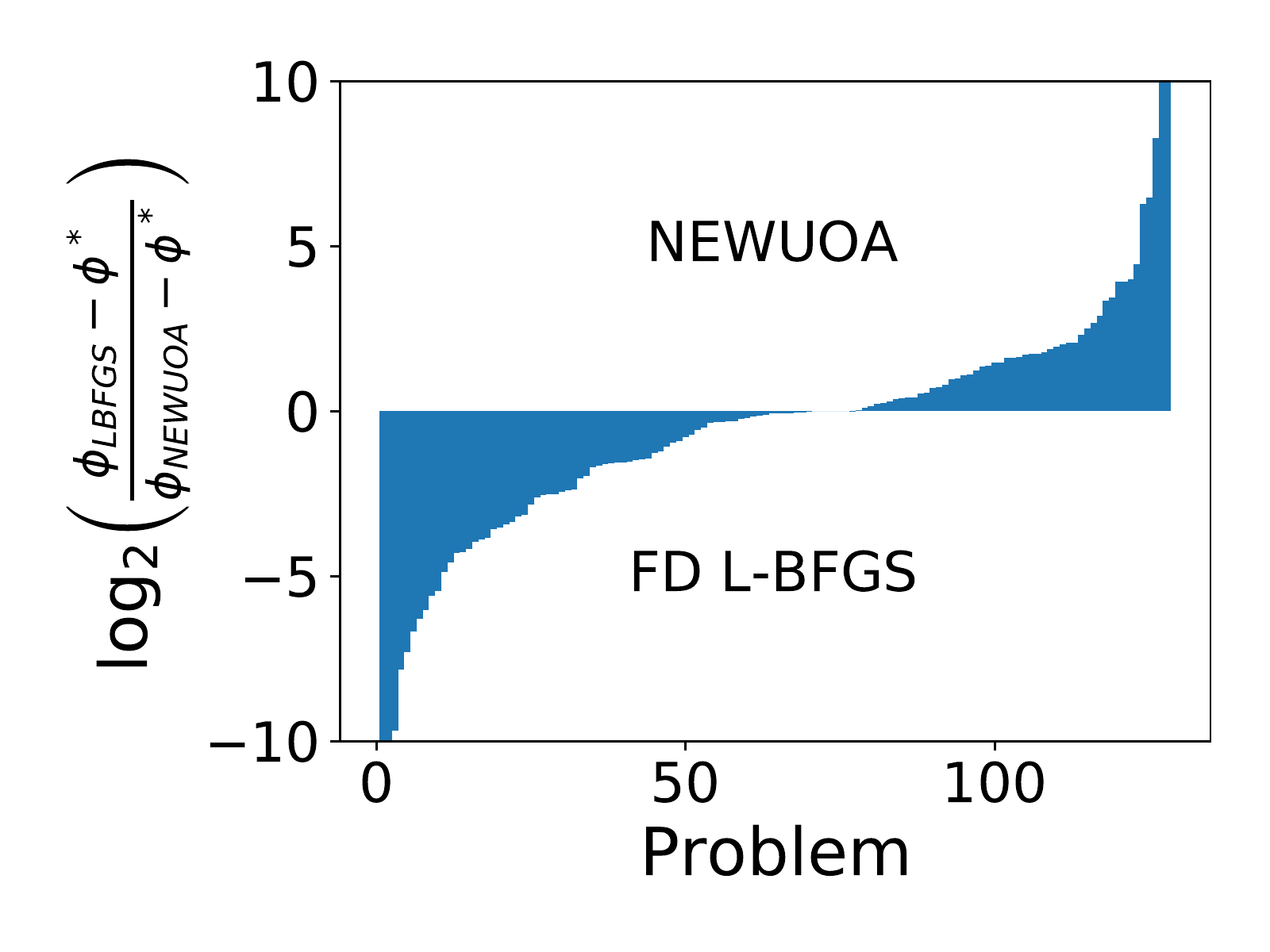}
\includegraphics[width=0.32\textwidth]{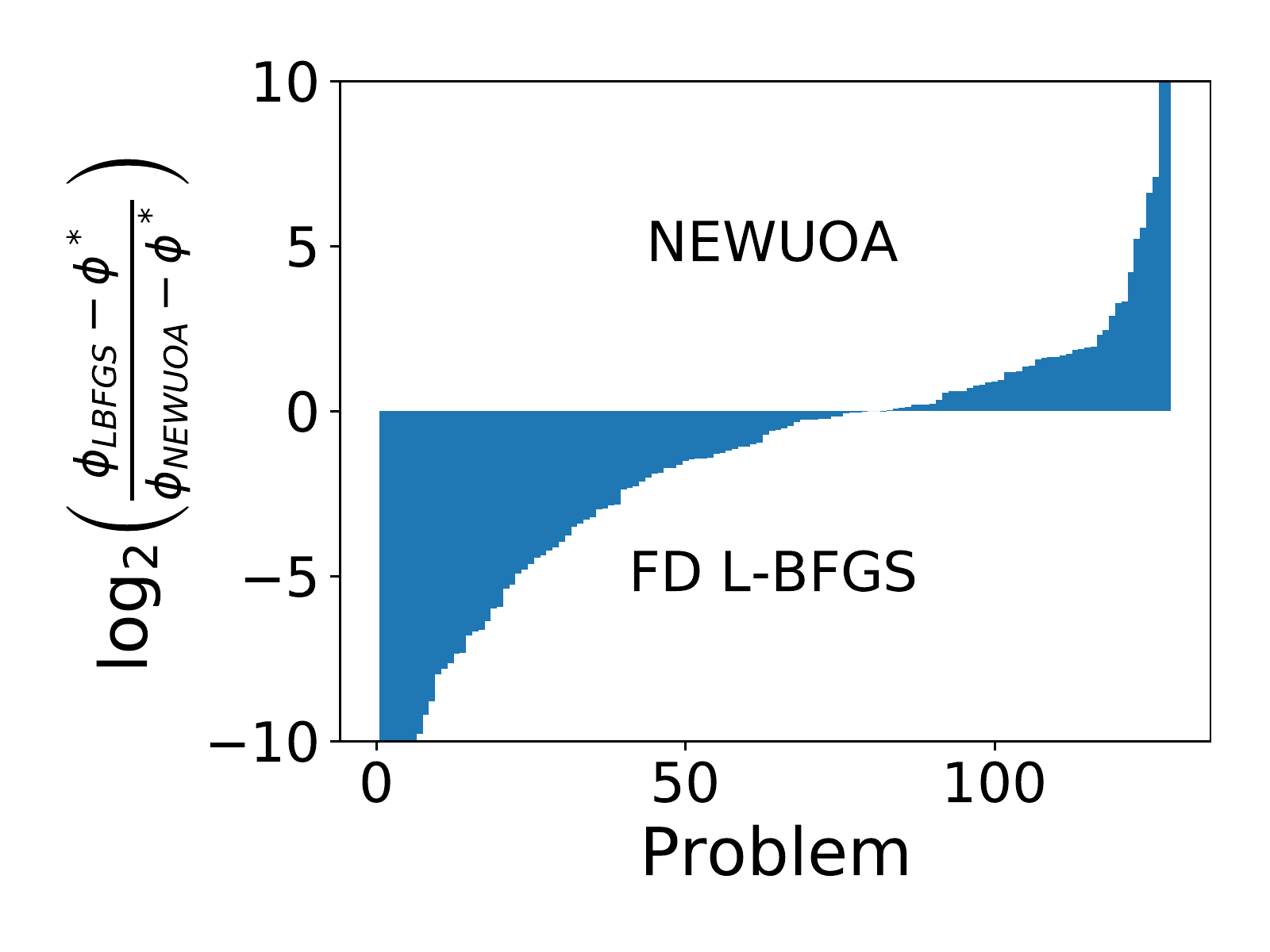} \\
\includegraphics[width=0.32\textwidth]{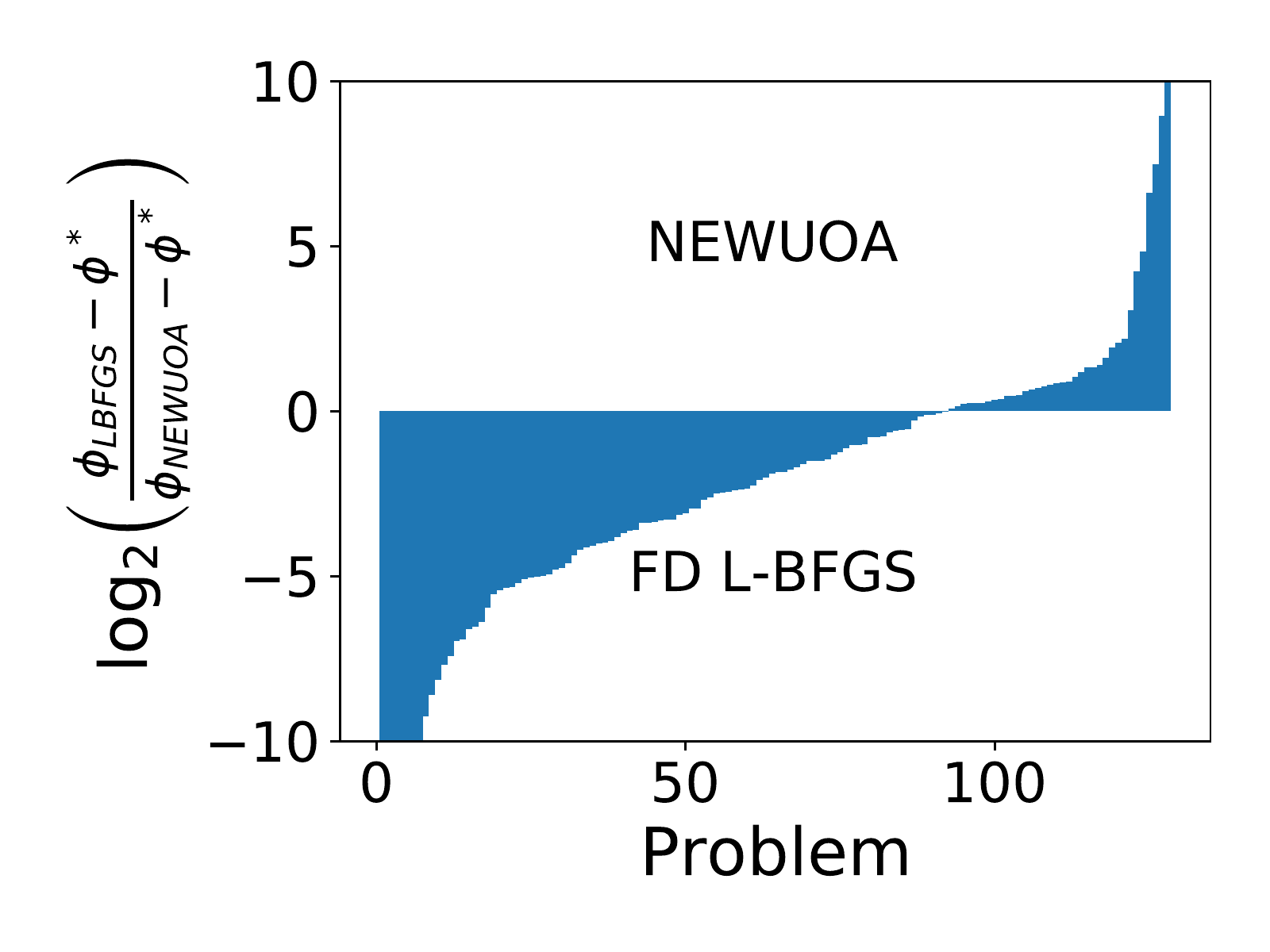}
\includegraphics[width=0.32\textwidth]{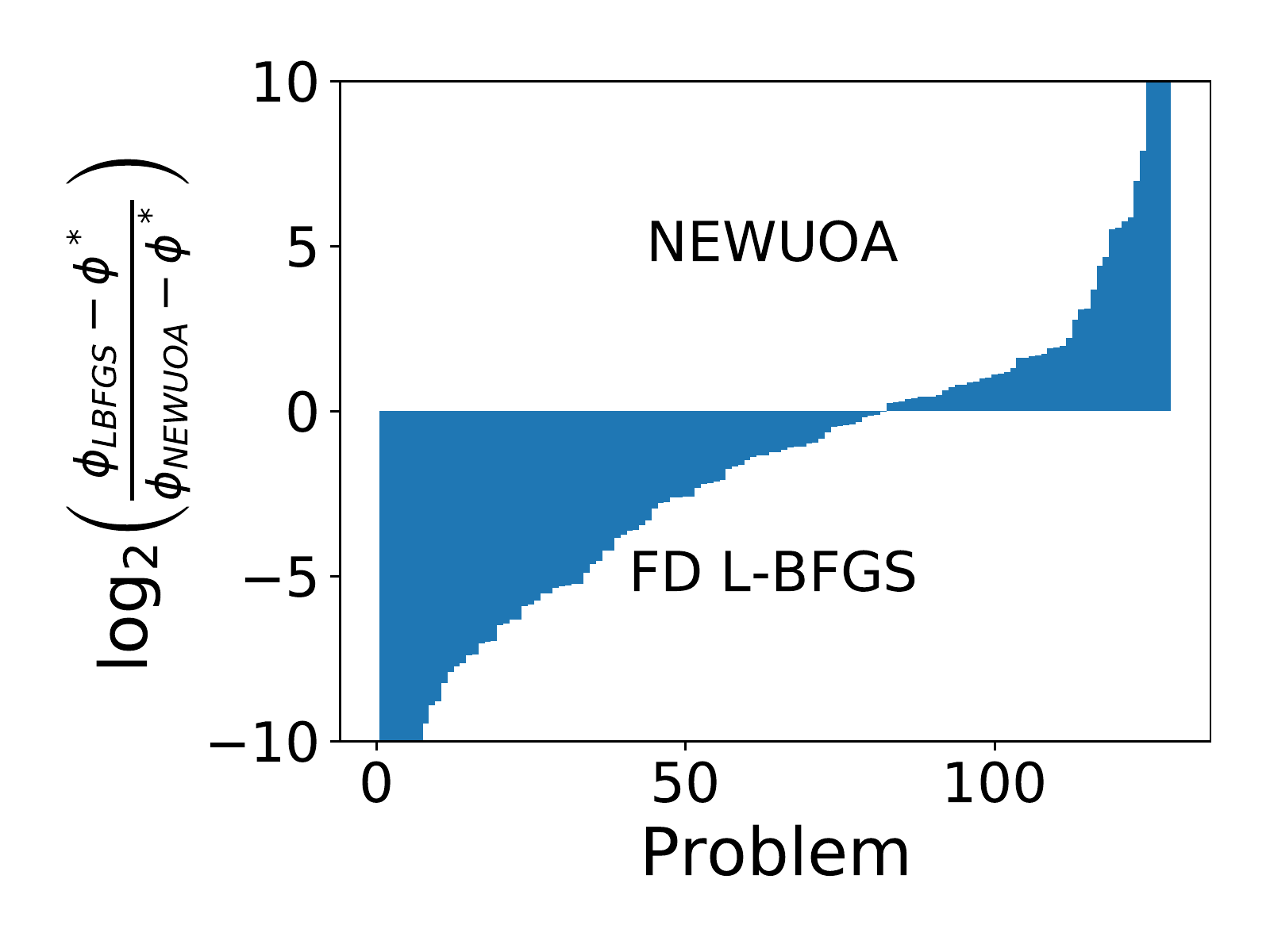}
\caption{Noisy Case. Log-ratio optimality gap profiles comparing forward difference {\sc l-bfgs} against {\sc newuoa} with $p = n + 2$ points. The noise levels are $\sigma_f = 10^{-1}$ (top left), $\sigma_f = 10^{-3}$ (top right), $10^{-5}$ (bottom left), and $10^{-7}$ (bottom right).}
\label{fig:new n+2 vs lbfgs}
\end{figure}

\begin{figure}[ht]
\centering
\includegraphics[width=0.32\textwidth]{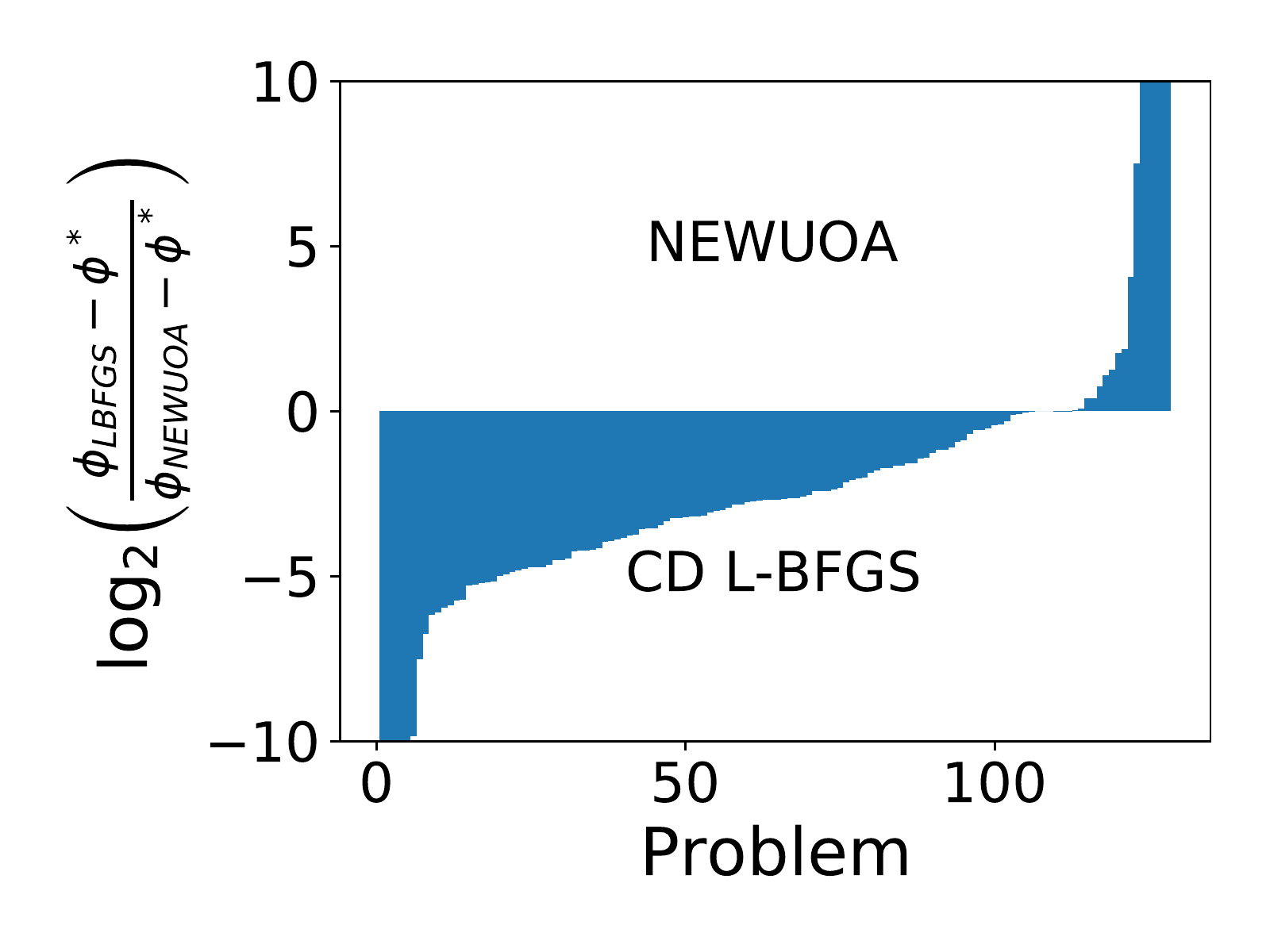}
\includegraphics[width=0.32\textwidth]{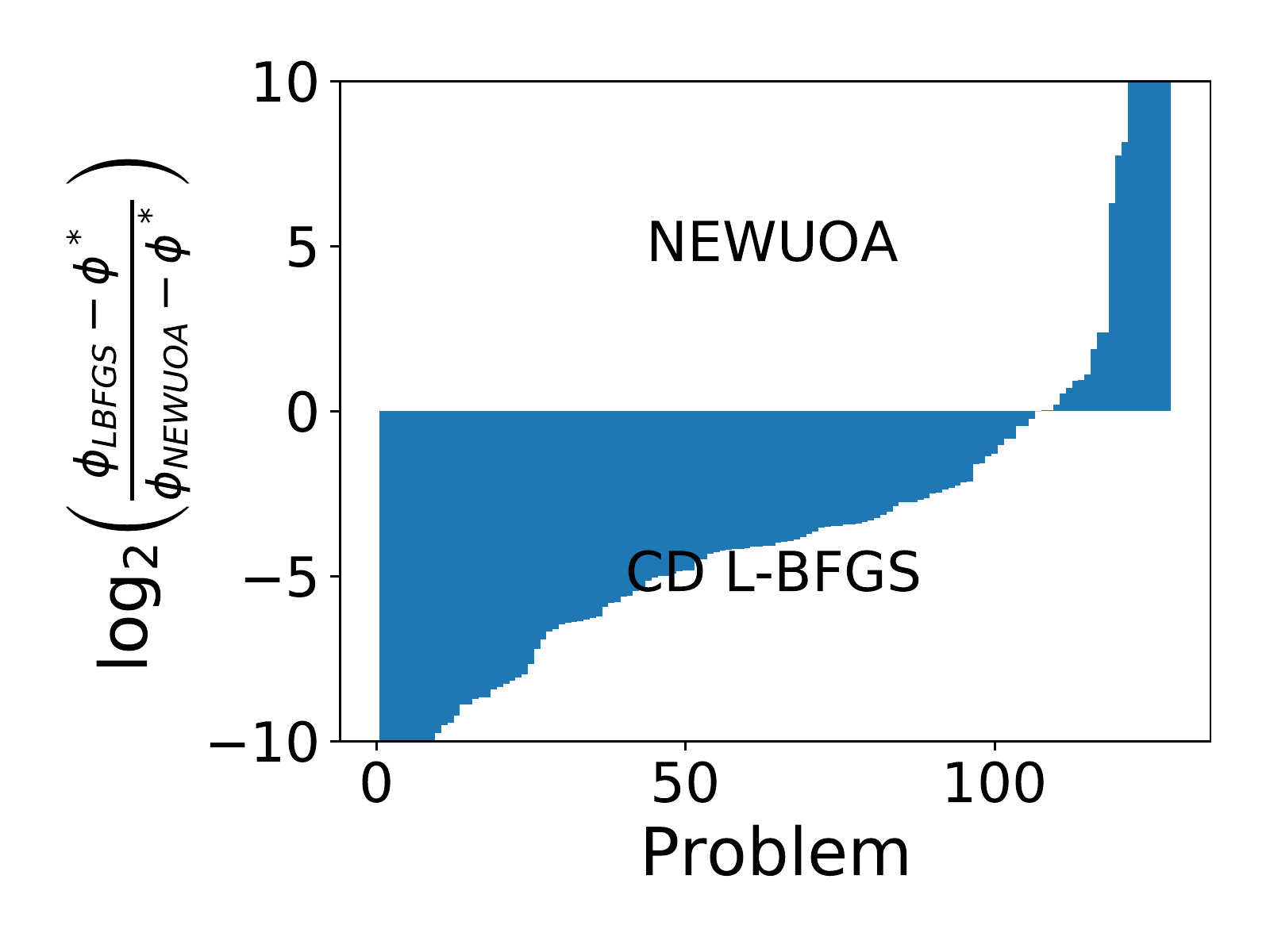} \\
\includegraphics[width=0.32\textwidth]{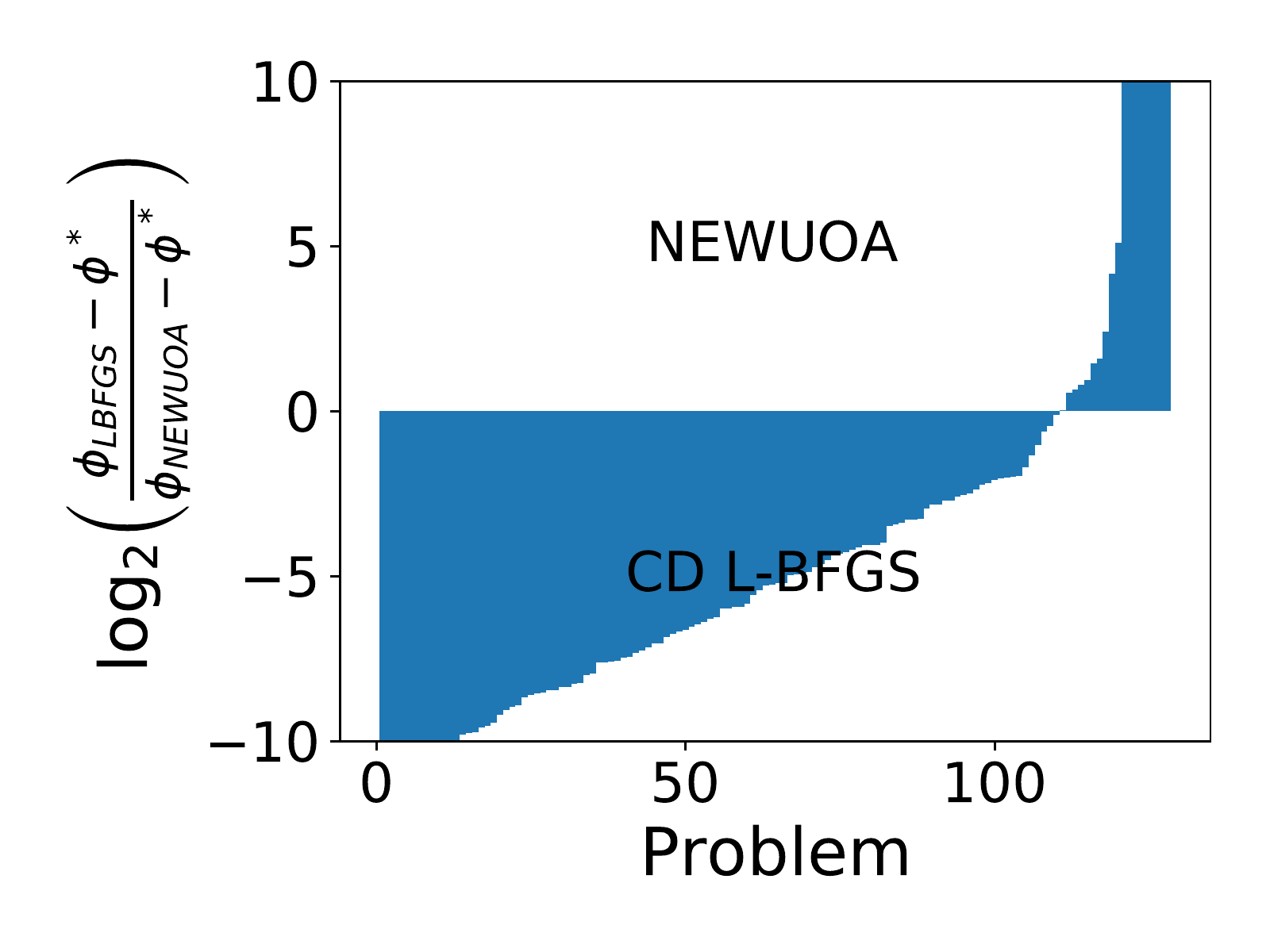}
\includegraphics[width=0.32\textwidth]{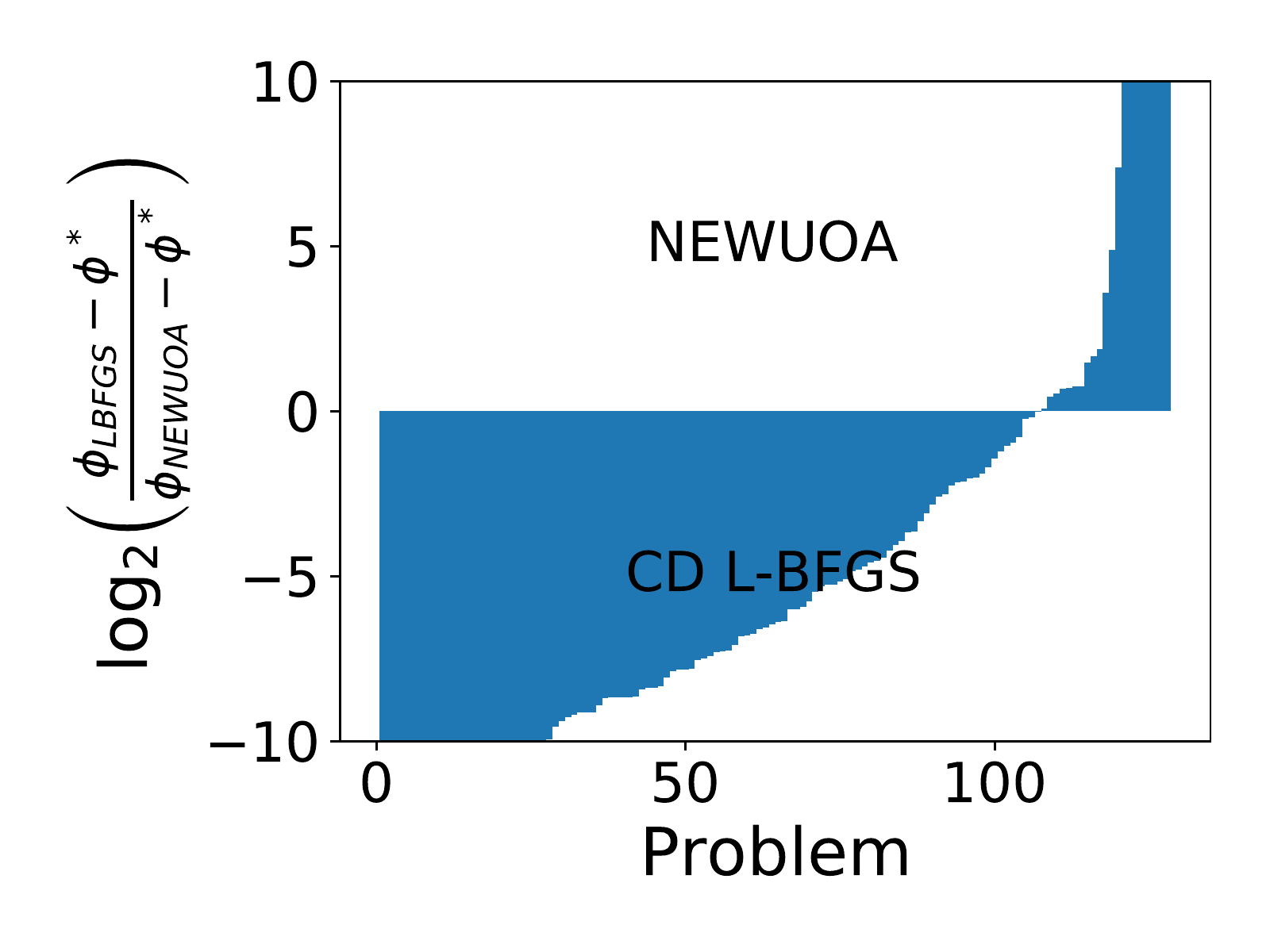}
\caption{Noisy Case. Log-ratio optimality gap profiles comparing central difference {\sc l-bfgs} against {\sc newuoa} with $p = 3n + 1$ points. The noise levels are $\sigma_f = 10^{-1}$ (top left), $\sigma_f = 10^{-3}$ (top right), $10^{-5}$ (bottom left), and $10^{-7}$ (bottom right).}
\label{fig:new 3n+1 vs lbfgs}
\end{figure}

When we compare {\sc newuoa} with $p = n + 2$ points against forward difference {\sc l-bfgs}, we see that {\sc l-bfgs} is now competitive against {\sc newuoa}, unlike when {\sc newuoa} employed $p = 2n + 1$ points; see Figure \ref{fig:new n+2 vs lbfgs}. However, when we increase the number of points to $p = 3n + 1$, {\sc newuoa} is still not competitive against central differencing, as seen in Figure \ref{fig:new 3n+1 vs lbfgs}.

\subsection{Trust Region Radius}

\begin{figure}[ht]
\centering
\includegraphics[width=0.32\textwidth]{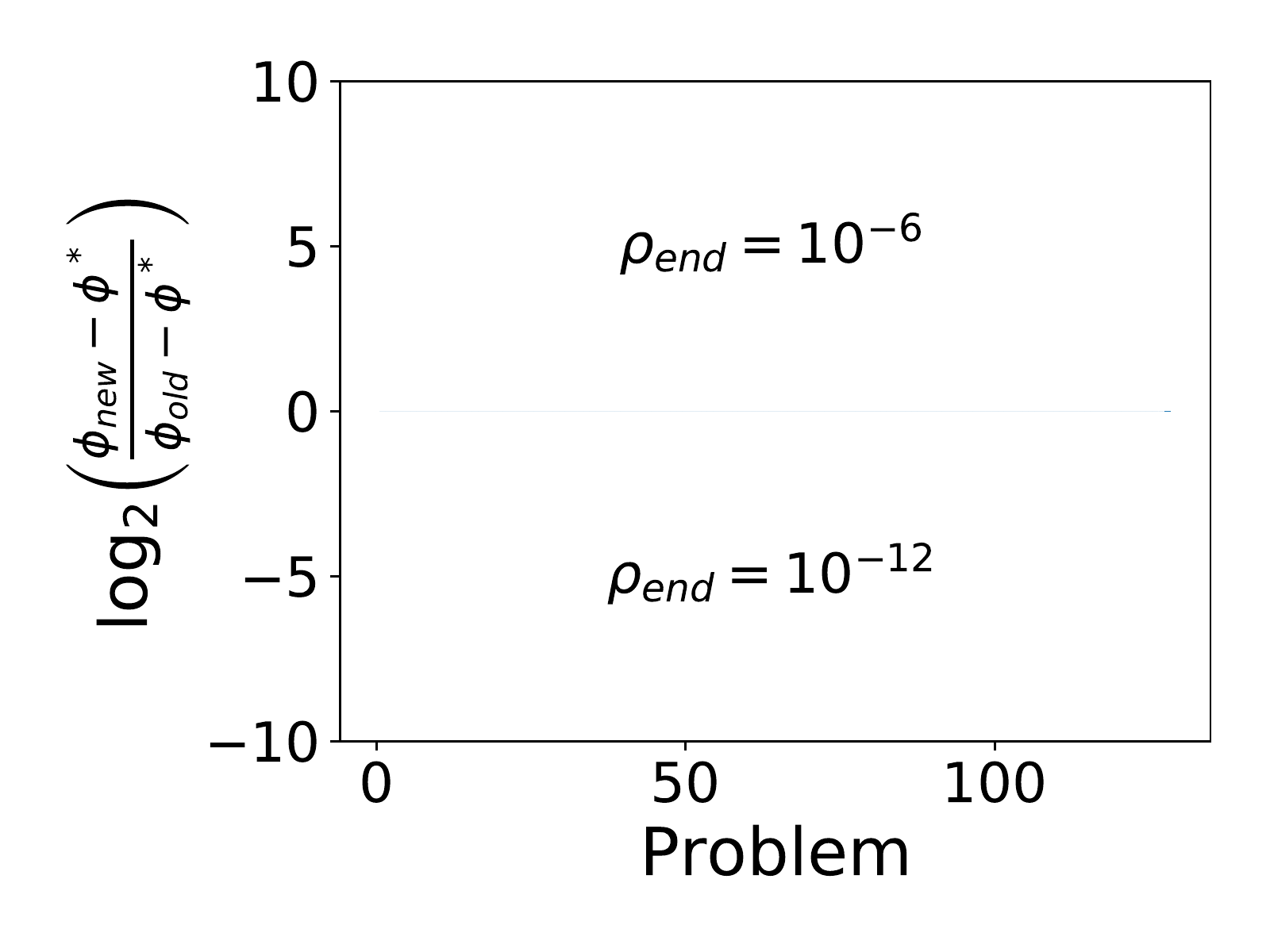}
\includegraphics[width=0.32\textwidth]{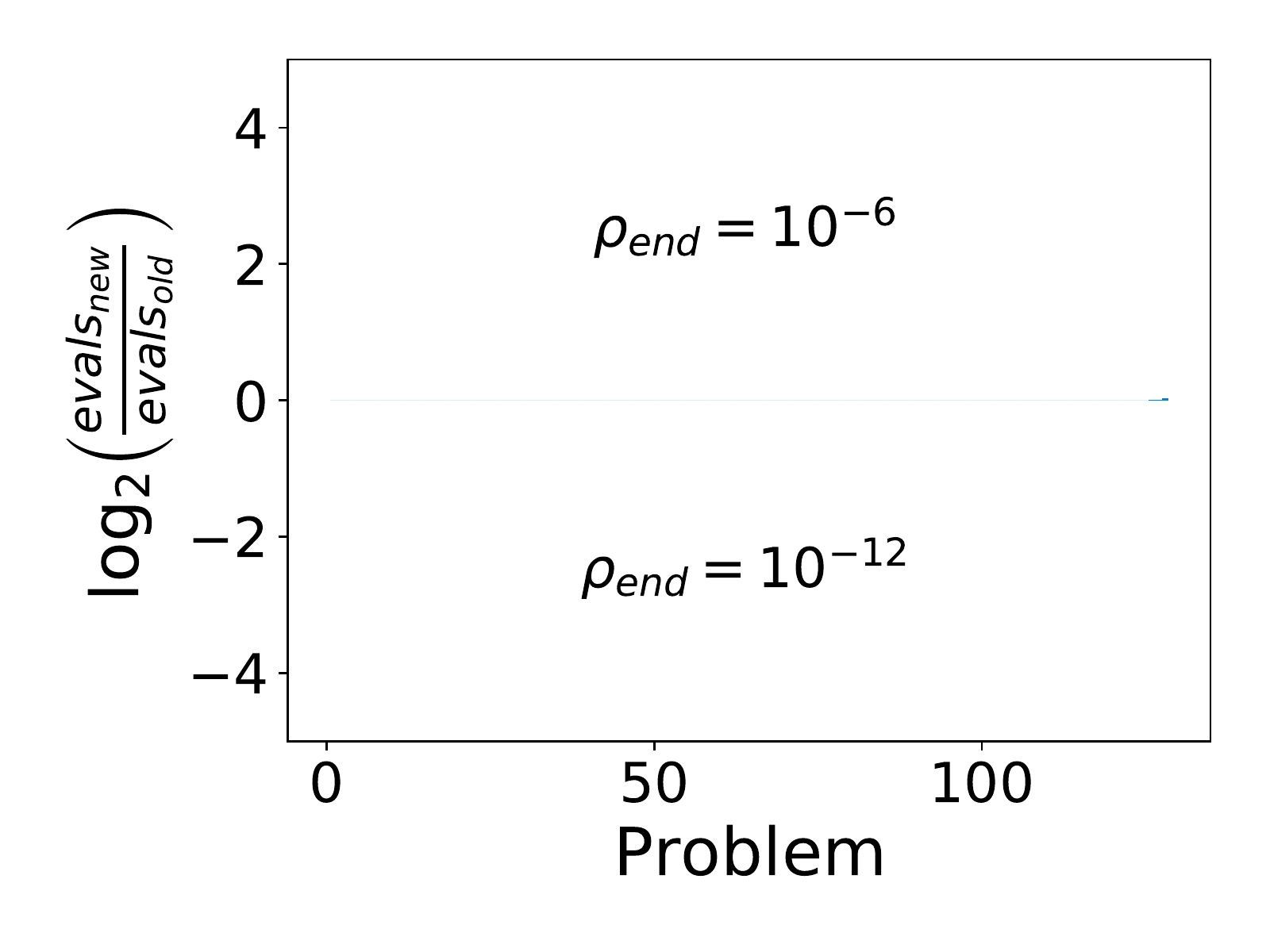}
\caption{Noisy Case with $\sigma_f = 10^{-5}$. Log-ratio optimality gap profiles comparing {\sc newuoa} with $\rho_{\text{end}} = 10^{-6}$ and $10^{-12}$.}
\label{fig:new rho}
\end{figure}

One can also ask if decreasing the final trust region radius could allow {\sc newuoa} to converge to a better solution. To test this, we change $\rho_{\text{end}} = 10^{-12}$ and compare against the default $\rho_{\text{end}} = 10^{-6}$ in Figure \ref{fig:new rho}. From our experiments, changing the final trust region radius makes almost no difference on the solution quality and efficiency.

\section{Complete Numerical Results}
\label{sec:ext_num}
In this appendix, we present our complete numerical results.

\subsection{Unconstrained Optimization}

\subsubsection{Noiseless Functions}
\label{app:smooth, det}

In this section, we list the number of function evaluations (\#feval), ratio between number of function evaluations to the number of function evaluations taken by {\sc newuoa} (RATIO), CPU time (CPU), and optimality gap ($\phi(x) - \phi^*$) for each problem instance in Tables \ref{tab:no-noise1}-\ref{tab:no-noise5}. Function evaluations marked with a $^*$ denote cases where the algorithm reached the maximum number of function evaluations. Instances where {\sc newuoa} samples a point that satisfies \eqref{eq:term crit} within the initial $2 n + 1$ evaluations are denoted by $^{**}$. Problems marked with a $^\ddagger$ denote cases where {\sc newuoa} and {\sc l-bfgs} converge to different minimizers. Optimality gaps marked with a $\dagger$ denote failures of the algorithm to converge to a valid solution satisfying \eqref{eq:term crit}.

\begin{landscape}

\begin{table}[htp]
    \centering
    \footnotesize
    \begin{tabular}{ l  r | r r r r | r r r r | r r r r}
        \toprule
        \multicolumn{2}{c}{ } &  \multicolumn{4}{|c|}{NEWUOA} & \multicolumn{4}{c|}{FD L-BFGS} & \multicolumn{4}{c}{CD L-BFGS} \\
        \midrule
        Problem & $n$ & \#feval & RATIO & CPU & $\phi(x) - \phi^*$ & \#feval & RATIO & CPU & $\phi(x) - \phi^*$ & \#feval & RATIO & CPU & $\phi(x) - \phi^*$ \\
        \midrule
        AIRCRFTB & $5$ & $661$ & \num{ 1.000} & $0.1$ & \num{ 6.963e-07} & $296$ & \num{ 0.448} & $0.02$ & \num{ 6.518e-07} & $536$ & \num{ 0.811} & $0.03$ & \num{ 8.288e-07} \\
        ALLINITU & $4$ & $41$ & \num{ 1.000} & $0.09$ & \num{ 3.073e-06} & $58$ & \num{ 1.415} & $0.01$ & \num{ 4.240e-06} & $103$ & \num{ 2.512} & $0.01$ & \num{ 4.240e-06} \\
        ARWHEAD & $100$ & $201^{**}$ & \num{ 1.000} & $0.09$ & \num{ 0.000e+00} & $1130$ & \num{ 5.622} & $0.06$ & \num{ 2.474e-07} & $2240$ & \num{ 11.144} & $0.07$ & \num{ 2.476e-07} \\
        BARD & $3$ & $69$ & \num{ 1.000} & $0.09$ & \num{ 7.032e-08} & $102$ & \num{ 1.478} & $0.01$ & \num{ 5.244e-07} & $176$ & \num{ 2.551} & $0.01$ & \num{ 5.246e-07} \\
        BDQRTIC & $100$ & $3682$ & \num{ 1.000} & $1.85$ & \num{ 3.731e-04} & $3178$ & \num{ 0.863} & $0.14$ & \num{ 2.923e-04} & $6308$ & \num{ 1.713} & $0.18$ & \num{ 2.914e-04} \\
        BIGGS3 & $3$ & $73$ & \num{ 1.000} & $0.09$ & \num{ 4.599e-07} & $85$ & \num{ 1.164} & $0.01$ & \num{ 5.422e-07} & $152$ & \num{ 2.082} & $0.01$ & \num{ 5.423e-07} \\
        BIGGS5$^\ddagger$ & $5$ & $103$ & \num{ 1.000} & $0.08$ & \num{-3.722e-05} & $416$ & \num{ 4.039} & $0.04$ & \num{ 6.651e-07} & $767$ & \num{ 7.447} & $0.05$ & \num{ 4.321e-07} \\
        BIGGS6$^\ddagger$ & $6$ & $338$ & \num{ 1.000} & $0.08$ & \num{-1.444e-04} & $323$ & \num{ 0.956} & $0.03$ & \num{ 5.357e-07} & $593$ & \num{ 1.754} & $0.03$ & \num{ 5.654e-07} \\
        BOX2 & $2$ & $2^{**}$ & \num{ 1.000} & $0.08$ & \num{ 5.847e-19} & $57$ & \num{ 28.500} & $0.01$ & \num{ 9.947e-07} & $95$ & \num{ 47.500} & $0.01$ & \num{ 9.947e-07} \\
        BOX3 & $3$ & $2^{**}$ & \num{ 1.000} & $0.08$ & \num{ 5.847e-19} & $52$ & \num{ 26.000} & $0.01$ & \num{ 7.371e-07} & $91$ & \num{ 45.500} & $0.01$ & \num{ 7.371e-07} \\
        BRKMCC & $2$ & $10$ & \num{ 1.000} & $0.08$ & \num{ 1.440e-07} & $32$ & \num{ 3.200} & $0.0$ & \num{ 4.066e-09} & $51$ & \num{ 5.100} & $0.0$ & \num{ 4.088e-09} \\
        BROWNAL & $10$ & $917$ & \num{ 1.000} & $0.11$ & \num{ 9.998e-07} & $177$ & \num{ 0.193} & $0.01$ & \num{ 1.707e-07} & $331$ & \num{ 0.361} & $0.01$ & \num{ 1.707e-07} \\
        BROWNAL & $100$ & $18069$ & \num{ 1.000} & $9.53$ & \num{ 9.993e-07} & $1039$ & \num{ 0.058} & $0.07$ & \num{ 1.040e-07} & $2051$ & \num{ 0.114} & $0.1$ & \num{ 1.047e-07} \\
        BROWNAL & $200$ & $100000^*$ & \num{ 1.000} & $207.31$ & \num{ 3.184e-06} & $823$ & \num{ 0.008} & $0.13$ & \num{ 8.254e-07} & $1626$ & \num{ 0.016} & $0.16$ & \num{ 8.243e-07} \\
        BROWNDEN & $4$ & $140$ & \num{ 1.000} & $0.09$ & \num{ 2.900e-02} & $131$ & \num{ 0.936} & $0.01$ & \num{ 1.448e-02} & $225$ & \num{ 1.607} & $0.01$ & \num{ 1.449e-02} \\
        CLIFF & $2$ & $84$ & \num{ 1.000} & $0.09$ & \num{ 5.220e-07} & $240$ & \num{ 2.857} & $0.02$ & \num{ 3.101e-07} & $365$ & \num{ 4.345} & $0.02$ & \num{ 2.076e-07} \\
        CRAGGLVY & $4$ & $182$ & \num{ 1.000} & $0.09$ & \num{ 8.761e-07} & $135$ & \num{ 0.742} & $0.01$ & \num{ 5.854e-07} & $244$ & \num{ 1.341} & $0.01$ & \num{ 5.854e-07} \\
        CRAGGLVY & $10$ & $511$ & \num{ 1.000} & $0.1$ & \num{ 1.650e-06} & $482$ & \num{ 0.943} & $0.03$ & \num{ 1.176e-06} & $910$ & \num{ 1.781} & $0.04$ & \num{ 1.168e-06} \\
        CRAGGLVY & $50$ & $2670$ & \num{ 1.000} & $0.45$ & \num{ 1.523e-05} & $2301$ & \num{ 0.862} & $0.1$ & \num{ 1.289e-05} & $4544$ & \num{ 1.702} & $0.14$ & \num{ 1.292e-05} \\
        CRAGGLVY & $100$ & $5679$ & \num{ 1.000} & $2.92$ & \num{ 3.219e-05} & $4503$ & \num{ 0.793} & $0.21$ & \num{ 1.932e-05} & $8946$ & \num{ 1.575} & $0.29$ & \num{ 1.934e-05} \\
        CUBE & $2$ & $173$ & \num{ 1.000} & $0.08$ & \num{ 8.157e-07} & $114$ & \num{ 0.659} & $0.01$ & \num{ 3.107e-07} & $190$ & \num{ 1.098} & $0.01$ & \num{ 3.036e-07} \\
        DENSCHNA & $2$ & $22$ & \num{ 1.000} & $0.09$ & \num{ 2.893e-10} & $38$ & \num{ 1.727} & $0.01$ & \num{ 9.517e-10} & $64$ & \num{ 2.909} & $0.01$ & \num{ 9.521e-10} \\
        DENSCHNB & $2$ & $21$ & \num{ 1.000} & $0.09$ & \num{ 9.135e-08} & $25$ & \num{ 1.190} & $0.0$ & \num{ 3.867e-11} & $42$ & \num{ 2.000} & $0.0$ & \num{ 3.877e-11} \\
        DENSCHNC & $2$ & $66$ & \num{ 1.000} & $0.09$ & \num{ 7.333e-08} & $68$ & \num{ 1.030} & $0.01$ & \num{ 8.949e-08} & $111$ & \num{ 1.682} & $0.01$ & \num{ 8.951e-08} \\
        DENSCHND & $3$ & $312$ & \num{ 1.000} & $0.07$ & \num{ 9.794e-07} & $208$ & \num{ 0.667} & $0.02$ & \num{ 9.534e-07} & $356$ & \num{ 1.141} & $0.02$ & \num{ 9.537e-07} \\
        DENSCHNE & $3$ & $127$ & \num{ 1.000} & $0.09$ & \num{ 4.487e-07} & $141$ & \num{ 1.110} & $0.01$ & \num{ 1.599e-07} & $237$ & \num{ 1.866} & $0.01$ & \num{ 1.210e-07} \\
        DENSCHNF & $2$ & $28$ & \num{ 1.000} & $0.09$ & \num{ 5.079e-09} & $48$ & \num{ 1.714} & $0.01$ & \num{ 5.450e-09} & $77$ & \num{ 2.750} & $0.01$ & \num{ 5.468e-09} \\
        \bottomrule
    \end{tabular}
    \caption{Noiseless Unconstrained \texttt{CUTEst} Problems Tested. $n$ is the number of variables.}
    \label{tab:no-noise1}
\end{table}

\begin{table}[htp]
    \centering
    \footnotesize
    \begin{tabular}{ l  r | r r r r | r r r r | r r r r}
        \toprule
        \multicolumn{2}{c}{ } &  \multicolumn{4}{|c|}{NEWUOA} & \multicolumn{4}{c|}{FD L-BFGS} & \multicolumn{4}{c}{CD L-BFGS} \\
        \midrule
        Problem & $n$ & \#feval & RATIO & CPU & $\phi(x) - \phi^*$ & \#feval & RATIO & CPU & $\phi(x) - \phi^*$ & \#feval & RATIO & CPU & $\phi(x) - \phi^*$ \\
        \midrule
        DIXMAANA & $15$ & $596$ & \num{ 1.000} & $0.09$ & \num{ 9.909e-07} & $138$ & \num{ 0.232} & $0.01$ & \num{ 1.096e-07} & $265$ & \num{ 0.445} & $0.01$ & \num{ 1.096e-07} \\
        DIXMAANA & $90$ & $1594$ & \num{ 1.000} & $0.61$ & \num{ 9.626e-07} & $738$ & \num{ 0.463} & $0.03$ & \num{ 6.577e-07} & $1465$ & \num{ 0.919} & $0.04$ & \num{ 6.578e-07} \\
        DIXMAANA & $300$ & $8117$ & \num{ 1.000} & $37.53$ & \num{ 9.936e-07} & $2720$ & \num{ 0.335} & $0.15$ & \num{ 2.979e-09} & $5428$ & \num{ 0.669} & $0.21$ & \num{ 2.978e-09} \\
        DIXMAANB & $15$ & $395$ & \num{ 1.000} & $0.1$ & \num{ 8.586e-07} & $122$ & \num{ 0.309} & $0.01$ & \num{ 3.845e-09} & $233$ & \num{ 0.590} & $0.01$ & \num{ 3.845e-09} \\
        DIXMAANB & $90$ & $1491$ & \num{ 1.000} & $0.56$ & \num{ 9.724e-07} & $647$ & \num{ 0.434} & $0.03$ & \num{ 6.536e-10} & $1283$ & \num{ 0.860} & $0.04$ & \num{ 6.537e-10} \\
        DIXMAANB & $300$ & $6446$ & \num{ 1.000} & $27.48$ & \num{ 9.980e-07} & $2117$ & \num{ 0.328} & $0.15$ & \num{ 7.825e-10} & $4223$ & \num{ 0.655} & $0.15$ & \num{ 7.845e-10} \\
        DIXMAANC & $15$ & $197$ & \num{ 1.000} & $0.08$ & \num{ 9.747e-07} & $158$ & \num{ 0.802} & $0.01$ & \num{ 8.216e-07} & $302$ & \num{ 1.533} & $0.01$ & \num{ 8.216e-07} \\
        DIXMAANC & $90$ & $1611$ & \num{ 1.000} & $0.63$ & \num{ 9.914e-07} & $558$ & \num{ 0.346} & $0.03$ & \num{ 1.352e-07} & $1104$ & \num{ 0.685} & $0.03$ & \num{ 1.353e-07} \\
        DIXMAANC & $300$ & $4643$ & \num{ 1.000} & $18.87$ & \num{ 9.964e-07} & $2120$ & \num{ 0.457} & $0.15$ & \num{ 5.674e-07} & $4227$ & \num{ 0.910} & $0.15$ & \num{ 5.672e-07} \\
        DIXMAAND & $15$ & $251$ & \num{ 1.000} & $0.1$ & \num{ 9.036e-07} & $179$ & \num{ 0.713} & $0.01$ & \num{ 4.536e-07} & $340$ & \num{ 1.355} & $0.01$ & \num{ 4.536e-07} \\
        DIXMAAND & $90$ & $2090$ & \num{ 1.000} & $0.81$ & \num{ 9.552e-07} & $745$ & \num{ 0.356} & $0.04$ & \num{ 6.242e-07} & $1474$ & \num{ 0.705} & $0.04$ & \num{ 6.242e-07} \\
        DIXMAAND & $300$ & $7155$ & \num{ 1.000} & $33.18$ & \num{ 9.791e-07} & $2425$ & \num{ 0.339} & $0.18$ & \num{ 3.136e-07} & $4834$ & \num{ 0.676} & $0.18$ & \num{ 3.136e-07} \\
        DIXMAANE & $15$ & $377$ & \num{ 1.000} & $0.1$ & \num{ 7.225e-07} & $293$ & \num{ 0.777} & $0.02$ & \num{ 4.958e-07} & $564$ & \num{ 1.496} & $0.02$ & \num{ 4.959e-07} \\
        DIXMAANE & $90$ & $3081$ & \num{ 1.000} & $1.28$ & \num{ 9.840e-07} & $3499$ & \num{ 1.136} & $0.15$ & \num{ 8.197e-07} & $6956$ & \num{ 2.258} & $0.2$ & \num{ 8.202e-07} \\
        DIXMAANE & $300$ & $18666$ & \num{ 1.000} & $107.25$ & \num{ 9.993e-07} & $19633$ & \num{ 1.052} & $0.96$ & \num{ 5.643e-07} & $39197$ & \num{ 2.100} & $1.42$ & \num{ 6.031e-07} \\
        DIXMAANF & $15$ & $274$ & \num{ 1.000} & $0.09$ & \num{ 9.301e-07} & $226$ & \num{ 0.825} & $0.01$ & \num{ 7.285e-07} & $434$ & \num{ 1.584} & $0.01$ & \num{ 7.286e-07} \\
        DIXMAANF & $90$ & $3104$ & \num{ 1.000} & $1.33$ & \num{ 9.964e-07} & $2674$ & \num{ 0.861} & $0.11$ & \num{ 3.730e-07} & $5313$ & \num{ 1.712} & $0.15$ & \num{ 3.731e-07} \\
        DIXMAANF & $300$ & $16182$ & \num{ 1.000} & $92.5$ & \num{ 9.988e-07} & $14504$ & \num{ 0.896} & $0.69$ & \num{ 9.006e-07} & $28952$ & \num{ 1.789} & $1.06$ & \num{ 9.011e-07} \\
        DIXMAANG & $15$ & $286$ & \num{ 1.000} & $0.08$ & \num{ 8.435e-07} & $277$ & \num{ 0.969} & $0.01$ & \num{ 2.101e-07} & $533$ & \num{ 1.864} & $0.02$ & \num{ 2.102e-07} \\
        DIXMAANG & $90$ & $3429$ & \num{ 1.000} & $1.47$ & \num{ 9.922e-07} & $2857$ & \num{ 0.833} & $0.12$ & \num{ 6.520e-07} & $5678$ & \num{ 1.656} & $0.16$ & \num{ 6.522e-07} \\
        DIXMAANG & $300$ & $21081$ & \num{ 1.000} & $126.91$ & \num{ 9.999e-07} & $14806$ & \num{ 0.702} & $0.78$ & \num{ 9.894e-07} & $30158$ & \num{ 1.431} & $1.1$ & \num{ 8.871e-07} \\
        DIXMAANH & $15$ & $375$ & \num{ 1.000} & $0.1$ & \num{ 8.025e-07} & $280$ & \num{ 0.747} & $0.02$ & \num{ 1.341e-07} & $537$ & \num{ 1.432} & $0.02$ & \num{ 1.342e-07} \\
        DIXMAANH & $90$ & $3106$ & \num{ 1.000} & $1.27$ & \num{ 9.943e-07} & $2588$ & \num{ 0.833} & $0.11$ & \num{ 6.407e-07} & $5138$ & \num{ 1.654} & $0.15$ & \num{ 6.405e-07} \\
        DIXMAANH & $300$ & $35784$ & \num{ 1.000} & $212.15$ & \num{ 9.999e-07} & $14207$ & \num{ 0.397} & $0.69$ & \num{ 7.042e-07} & $28356$ & \num{ 0.792} & $1.05$ & \num{ 6.752e-07} \\
        DIXMAANI & $15$ & $586$ & \num{ 1.000} & $0.1$ & \num{ 9.553e-07} & $549$ & \num{ 0.937} & $0.03$ & \num{ 9.880e-07} & $1060$ & \num{ 1.809} & $0.04$ & \num{ 9.880e-07} \\
        DIXMAANI & $90$ & $13092$ & \num{ 1.000} & $6.05$ & \num{ 9.966e-07} & $15925$ & \num{ 1.216} & $0.62$ & \num{ 9.862e-07} & $33132$ & \num{ 2.531} & $0.91$ & \num{ 9.270e-07} \\
        DIXMAANI & $300$ & $121834$ & \num{ 1.000} & $774.33$ & \num{ 9.999e-07} & $150107^*$ & \num{ 1.232} & $6.82$ & \num{ 2.119e-06} & $150155^*$ & \num{ 1.232} & $5.3$ & \num{ 1.720e-05} \\
        DIXMAANJ & $15$ & $563$ & \num{ 1.000} & $0.1$ & \num{ 9.736e-07} & $479$ & \num{ 0.851} & $0.03$ & \num{ 1.809e-07} & $926$ & \num{ 1.645} & $0.03$ & \num{ 1.797e-07} \\
        DIXMAANJ & $90$ & $12417$ & \num{ 1.000} & $6.11$ & \num{ 9.923e-07} & $11508$ & \num{ 0.927} & $0.53$ & \num{ 9.628e-07} & $22883$ & \num{ 1.843} & $0.74$ & \num{ 9.706e-07} \\
        DIXMAANJ & $300$ & $120362$ & \num{ 1.000} & $718.6$ & \num{ 9.987e-07} & $77321$ & \num{ 0.642} & $3.61$ & \num{ 9.974e-07} & $150156^*$ & \num{ 1.248} & $5.37$ & \num{ 1.028e-06} \\
        DIXMAANK & $15$ & $559$ & \num{ 1.000} & $0.1$ & \num{ 9.181e-07} & $533$ & \num{ 0.953} & $0.03$ & \num{ 3.242e-07} & $1029$ & \num{ 1.841} & $0.04$ & \num{ 3.230e-07} \\
        DIXMAANK & $90$ & $12601$ & \num{ 1.000} & $5.72$ & \num{ 9.969e-07} & $7463$ & \num{ 0.592} & $0.3$ & \num{ 9.902e-07} & $15018$ & \num{ 1.192} & $0.41$ & \num{ 9.893e-07} \\
        DIXMAANK & $300$ & $80808$ & \num{ 1.000} & $455.59$ & \num{ 9.997e-07} & $23567$ & \num{ 0.292} & $1.12$ & \num{ 9.623e-07} & $47649$ & \num{ 0.590} & $1.74$ & \num{ 9.108e-07} \\
        DIXMAANL & $15$ & $752$ & \num{ 1.000} & $0.1$ & \num{ 9.740e-07} & $502$ & \num{ 0.668} & $0.03$ & \num{ 4.740e-07} & $967$ & \num{ 1.286} & $0.03$ & \num{ 4.714e-07} \\
        DIXMAANL & $90$ & $10608$ & \num{ 1.000} & $4.87$ & \num{ 9.970e-07} & $10323$ & \num{ 0.973} & $0.41$ & \num{ 9.519e-07} & $20517$ & \num{ 1.934} & $0.57$ & \num{ 8.891e-07} \\
        DIXMAANL & $300$ & $150000^*$ & \num{ 1.000} & $875.24$ & \num{ 3.192e-06} & $34747$ & \num{ 0.232} & $1.67$ & \num{ 9.946e-07} & $69364$ & \num{ 0.462} & $2.51$ & \num{ 9.926e-07} \\
        \bottomrule
    \end{tabular}
    \caption{Noiseless Unconstrained \texttt{CUTEst} Problems Tested. $n$ is the number of variables.}
    \label{tab:no-noise2}
\end{table}

\begin{table}[htp]
    \centering
    \footnotesize
    \begin{tabular}{ l  r | r r r r | r r r r | r r r r}
        \toprule
        \multicolumn{2}{c}{ } &  \multicolumn{4}{|c|}{NEWUOA} & \multicolumn{4}{c|}{FD L-BFGS} & \multicolumn{4}{c}{CD L-BFGS} \\
        \midrule
        Problem & $n$ & \#feval & RATIO & CPU & $\phi(x) - \phi^*$ & \#feval & RATIO & CPU & $\phi(x) - \phi^*$ & \#feval & RATIO & CPU & $\phi(x) - \phi^*$ \\
        \midrule
        DQRTIC & $10$ & $502$ & \num{ 1.000} & $0.1$ & \num{ 9.719e-07} & $258$ & \num{ 0.514} & $0.02$ & \num{ 4.439e-07} & $488$ & \num{ 0.972} & $0.02$ & \num{ 4.439e-07} \\
        DQRTIC & $50$ & $2847$ & \num{ 1.000} & $0.43$ & \num{ 8.852e-07} & $1467$ & \num{ 0.515} & $0.06$ & \num{ 7.890e-07} & $2894$ & \num{ 1.017} & $0.08$ & \num{ 7.888e-07} \\
        DQRTIC & $100$ & $6061$ & \num{ 1.000} & $2.81$ & \num{ 9.769e-07} & $3277$ & \num{ 0.541} & $0.13$ & \num{ 5.471e-07} & $6508$ & \num{ 1.074} & $0.17$ & \num{ 5.462e-07} \\
        EDENSCH & $36$ & $928$ & \num{ 1.000} & $0.14$ & \num{ 2.175e-04} & $577$ & \num{ 0.622} & $0.03$ & \num{ 7.341e-05} & $1131$ & \num{ 1.219} & $0.03$ & \num{ 7.341e-05} \\
        EIGENALS & $6$ & $5^{**}$ & \num{ 1.000} & $0.08$ & \num{ 0.000e+00} & $83$ & \num{ 16.600} & $0.01$ & \num{ 7.344e-09} & $152$ & \num{ 30.400} & $0.01$ & \num{ 7.099e-09} \\
        EIGENALS & $110$ & $55000^*$ & \num{ 1.000} & $40.44$ & \num{ 1.716e-03} & $42138$ & \num{ 0.766} & $2.12$ & \num{ 9.758e-07} & $55108^*$ & \num{ 1.002} & $2.13$ & \num{ 2.142e-05} \\
        EIGENBLS$^\ddagger$ & $6$ & $113$ & \num{ 1.000} & $0.09$ & \num{-3.033e-04} & $91$ & \num{ 0.805} & $0.01$ & \num{ 9.727e-09} & $167$ & \num{ 1.478} & $0.01$ & \num{ 9.145e-09} \\
        EIGENBLS & $110$ & $37505$ & \num{ 1.000} & $27.4$ & \num{ 9.993e-07} & $55010^*$ & \num{ 1.467} & $2.68$ & \num{ 9.280e-04} & $55088^*$ & \num{ 1.469} & $2.1$ & \num{ 4.002e-02} \\
        EIGENCLS & $30$ & $1630$ & \num{ 1.000} & $0.19$ & \num{ 9.383e-07} & $2634$ & \num{ 1.616} & $0.12$ & \num{ 9.353e-07} & $5174$ & \num{ 3.174} & $0.16$ & \num{ 8.934e-07} \\
        ENGVAL1 & $2$ & $32$ & \num{ 1.000} & $0.08$ & \num{ 6.651e-09} & $38$ & \num{ 1.188} & $0.0$ & \num{ 2.191e-07} & $60$ & \num{ 1.875} & $0.0$ & \num{ 2.191e-07} \\
        ENGVAL1 & $50$ & $2234$ & \num{ 1.000} & $0.38$ & \num{ 5.211e-05} & $585$ & \num{ 0.262} & $0.02$ & \num{ 3.074e-05} & $1148$ & \num{ 0.514} & $0.03$ & \num{ 3.074e-05} \\
        ENGVAL1 & $100$ & $1941$ & \num{ 1.000} & $0.92$ & \num{ 1.076e-04} & $1337$ & \num{ 0.689} & $0.07$ & \num{ 4.197e-06} & $2651$ & \num{ 1.366} & $0.07$ & \num{ 4.198e-06} \\
        EXPFIT & $2$ & $42$ & \num{ 1.000} & $0.09$ & \num{ 9.588e-07} & $46$ & \num{ 1.095} & $0.01$ & \num{ 1.427e-09} & $75$ & \num{ 1.786} & $0.0$ & \num{ 1.420e-09} \\
        FLETCBV3$^\ddagger$ & $10$ & $1106$ & \num{ 1.000} & $0.11$ & \num{ 1.063e-03} & $261$ & \num{ 0.236} & $0.01$ & \num{-1.664e-05} & $521$ & \num{ 0.471} & $0.02$ & \num{-4.710e-05} \\
        FLETCBV3$^\ddagger$ & $100$ & $50000^*$ & \num{ 1.000} & $28.66$ & \num{ 1.333e+05} & $19708$ & \num{ 0.394} & $0.87$ & \num{-3.486e-01} & $16125$ & \num{ 0.323} & $0.54$ & \num{-1.107e+01} \\
        FLETCHBV$^\ddagger$ & $10$ & $945$ & \num{ 1.000} & $0.09$ & \num{ 1.490e+05} & $263$ & \num{ 0.278} & $0.02$ & \num{-1.278e+01} & $278$ & \num{ 0.294} & $0.01$ & \num{-5.322e+03} \\
        FLETCHBV$^\ddagger$ & $100$ & $50000^*$ & \num{ 1.000} & $27.77$ & \num{ 1.210e+13} & $18068$ & \num{ 0.361} & $0.8$ & \num{-6.240e+06} & $20933$ & \num{ 0.419} & $0.67$ & \num{-7.016e+07} \\
        FREUROTH & $2$ & $55$ & \num{ 1.000} & $0.09$ & \num{ 2.214e-05} & $81$ & \num{ 1.473} & $0.01$ & \num{ 3.932e-05} & $134$ & \num{ 2.436} & $0.01$ & \num{ 3.932e-05} \\
        FREUROTH & $10$ & $367$ & \num{ 1.000} & $0.09$ & \num{ 6.524e-04} & $286$ & \num{ 0.779} & $0.02$ & \num{ 1.194e-06} & $538$ & \num{ 1.466} & $0.02$ & \num{ 1.194e-06} \\
        FREUROTH & $50$ & $5840$ & \num{ 1.000} & $0.86$ & \num{ 5.844e-03} & $1053$ & \num{ 0.180} & $0.04$ & \num{ 4.888e-04} & $2073$ & \num{ 0.355} & $0.06$ & \num{ 4.886e-04} \\
        FREUROTH & $100$ & $2881$ & \num{ 1.000} & $1.44$ & \num{ 1.188e-02} & $1950$ & \num{ 0.677} & $0.09$ & \num{ 2.601e-03} & $3868$ & \num{ 1.343} & $0.12$ & \num{ 2.589e-03} \\
        GENROSE & $5$ & $199$ & \num{ 1.000} & $0.08$ & \num{ 9.086e-07} & $226$ & \num{ 1.136} & $0.02$ & \num{ 9.527e-09} & $411$ & \num{ 2.065} & $0.02$ & \num{ 9.367e-09} \\
        GENROSE & $10$ & $645$ & \num{ 1.000} & $0.1$ & \num{ 8.369e-07} & $896$ & \num{ 1.389} & $0.06$ & \num{ 8.050e-07} & $1614$ & \num{ 2.502} & $0.06$ & \num{ 1.032e-07} \\
        GENROSE & $100$ & $17249$ & \num{ 1.000} & $9.45$ & \num{ 9.954e-07} & $25453$ & \num{ 1.476} & $1.01$ & \num{ 8.013e-07} & $50198^*$ & \num{ 2.910} & $1.4$ & \num{ 1.939e-05} \\
        GULF & $3$ & $876$ & \num{ 1.000} & $0.12$ & \num{ 6.954e-07} & $225$ & \num{ 0.257} & $0.03$ & \num{ 7.913e-08} & $392$ & \num{ 0.447} & $0.04$ & \num{ 6.048e-07} \\
        HAIRY & $2$ & $261$ & \num{ 1.000} & $0.09$ & \num{ 1.214e-06} & $86$ & \num{ 0.330} & $0.01$ & \num{ 2.209e-07} & $149$ & \num{ 0.571} & $0.01$ & \num{ 2.996e-07} \\
        HELIX & $3$ & $65$ & \num{ 1.000} & $0.09$ & \num{ 3.569e-08} & $149$ & \num{ 2.292} & $0.01$ & \num{ 6.812e-07} & $86$ & \num{ 1.323} & $0.0$ & \num{ 8.521e+02}$^\dagger$ \\
        JENSMP$^\ddagger$ & $2$ & $4^{**}$ & \num{ 1.000} & $0.09$ & \num{-5.469e+02} & $16$ & \num{ 4.000} & $0.0$ & \num{ 0.000e+00} & $21$ & \num{ 5.250} & $0.0$ & \num{ 0.000e+00} \\
        KOWOSB & $4$ & $150$ & \num{ 1.000} & $0.09$ & \num{ 7.089e-07} & $183$ & \num{ 1.220} & $0.02$ & \num{ 2.824e-07} & $332$ & \num{ 2.213} & $0.02$ & \num{ 2.825e-07} \\
        MEXHAT & $2$ & $46$ & \num{ 1.000} & $0.09$ & \num{ 1.036e-02}$^\dagger$ & $203$ & \num{ 4.413} & $0.02$ & \num{ 7.395e-07} & $339$ & \num{ 7.370} & $0.02$ & \num{ 4.707e-07} \\
        \bottomrule
    \end{tabular}
    \caption{Noiseless Unconstrained \texttt{CUTEst} Problems Tested. $n$ is the number of variables.}
    \label{tab:no-noise3}
\end{table}

\begin{table}[htp]
    \centering
    \footnotesize
    \begin{tabular}{ l  r | r r r r | r r r r | r r r r}
        \toprule
        \multicolumn{2}{c}{ } &  \multicolumn{4}{|c|}{NEWUOA} & \multicolumn{4}{c|}{FD L-BFGS} & \multicolumn{4}{c}{CD L-BFGS} \\
        \midrule
        Problem & $n$ & \#feval & RATIO & CPU & $\phi(x) - \phi^*$ & \#feval & RATIO & CPU & $\phi(x) - \phi^*$ & \#feval & RATIO & CPU & $\phi(x) - \phi^*$ \\
        \midrule
        MOREBV & $10$ & $379$ & \num{ 1.000} & $0.09$ & \num{ 9.458e-07} & $366$ & \num{ 0.966} & $0.02$ & \num{ 6.782e-07} & $695$ & \num{ 1.834} & $0.03$ & \num{ 6.843e-07} \\
        MOREBV & $50$ & $21488$ & \num{ 1.000} & $3.71$ & \num{ 9.988e-07} & $25029^*$ & \num{ 1.165} & $0.98$ & \num{ 1.368e-06} & $25035^*$ & \num{ 1.165} & $0.68$ & \num{ 6.088e-06} \\
        MOREBV & $100$ & $11143$ & \num{ 1.000} & $6.29$ & \num{ 9.999e-07} & $15512$ & \num{ 1.392} & $0.59$ & \num{ 9.985e-07} & $31067$ & \num{ 2.788} & $0.82$ & \num{ 9.992e-07} \\
        NCB20B & $21$ & $478$ & \num{ 1.000} & $0.11$ & \num{ 4.093e-05} & $169$ & \num{ 0.354} & $0.01$ & \num{-4.368e-11} & $322$ & \num{ 0.674} & $0.01$ & \num{-4.334e-11} \\
        NCB20B & $22$ & $627$ & \num{ 1.000} & $0.11$ & \num{ 4.361e-05} & $154$ & \num{ 0.246} & $0.01$ & \num{ 2.493e-05} & $292$ & \num{ 0.466} & $0.01$ & \num{ 2.493e-05} \\
        NCB20B & $50$ & $3345$ & \num{ 1.000} & $0.59$ & \num{ 9.936e-05} & $1465$ & \num{ 0.438} & $0.07$ & \num{ 9.906e-05} & $2892$ & \num{ 0.865} & $0.1$ & \num{ 9.912e-05} \\
        NCB20B & $100$ & $7309$ & \num{ 1.000} & $4.02$ & \num{ 1.963e-04} & $3376$ & \num{ 0.462} & $0.21$ & \num{ 1.794e-04} & $6708$ & \num{ 0.918} & $0.3$ & \num{ 1.797e-04} \\
        NCB20B & $180$ & $11313$ & \num{ 1.000} & $22.07$ & \num{ 3.470e-04} & $4375$ & \num{ 0.387} & $0.35$ & \num{ 3.384e-04} & $8718$ & \num{ 0.771} & $0.56$ & \num{ 3.379e-04} \\
        NONDIA & $10$ & $198$ & \num{ 1.000} & $0.1$ & \num{ 4.267e-07} & $167$ & \num{ 0.843} & $0.01$ & \num{ 1.459e-11} & $298$ & \num{ 1.505} & $0.01$ & \num{ 1.053e-11} \\
        NONDIA & $20$ & $426$ & \num{ 1.000} & $0.1$ & \num{ 8.903e-07} & $341$ & \num{ 0.800} & $0.02$ & \num{ 2.571e-10} & $655$ & \num{ 1.538} & $0.02$ & \num{ 2.929e-10} \\
        NONDIA & $30$ & $569$ & \num{ 1.000} & $0.12$ & \num{ 9.936e-07} & $492$ & \num{ 0.865} & $0.02$ & \num{ 3.712e-07} & $956$ & \num{ 1.680} & $0.03$ & \num{ 3.738e-07} \\
        NONDIA & $50$ & $799$ & \num{ 1.000} & $0.2$ & \num{ 9.829e-07} & $792$ & \num{ 0.991} & $0.03$ & \num{ 3.825e-08} & $1556$ & \num{ 1.947} & $0.04$ & \num{ 3.691e-08} \\
        NONDIA & $90$ & $1302$ & \num{ 1.000} & $0.58$ & \num{ 9.192e-07} & $1301$ & \num{ 0.999} & $0.06$ & \num{ 1.233e-07} & $2574$ & \num{ 1.977} & $0.07$ & \num{ 1.188e-07} \\
        NONDIA & $100$ & $1683$ & \num{ 1.000} & $0.89$ & \num{ 9.873e-07} & $1543$ & \num{ 0.917} & $0.07$ & \num{ 1.028e-10} & $3057$ & \num{ 1.816} & $0.09$ & \num{ 1.125e-10} \\
        NONDQUAR & $100$ & $50000^*$ & \num{ 1.000} & $29.28$ & \num{ 8.675e-06} & $50056^*$ & \num{ 1.001} & $1.88$ & \num{ 1.137e-05} & $50189^*$ & \num{ 1.004} & $1.32$ & \num{ 3.611e-05} \\
        OSBORNEA$^\ddagger$ & $5$ & $1094$ & \num{ 1.000} & $0.09$ & \num{ 1.204e-05} & $477$ & \num{ 0.436} & $0.05$ & \num{ 8.442e-07} & $888$ & \num{ 0.812} & $0.06$ & \num{ 8.390e-07} \\
        OSBORNEB & $11$ & $1529$ & \num{ 1.000} & $0.14$ & \num{ 9.678e-07} & $1044$ & \num{ 0.683} & $0.08$ & \num{ 9.551e-07} & $2216$ & \num{ 1.449} & $0.11$ & \num{ 6.312e-07} \\
        PENALTY1 & $4$ & $577$ & \num{ 1.000} & $0.09$ & \num{ 9.994e-07} & $868$ & \num{ 1.504} & $0.09$ & \num{ 9.134e-07} & $1425$ & \num{ 2.470} & $0.09$ & \num{ 9.856e-07} \\
        PENALTY1 & $10$ & $1332$ & \num{ 1.000} & $0.11$ & \num{ 9.944e-07} & $1523$ & \num{ 1.143} & $0.09$ & \num{ 9.903e-07} & $2942$ & \num{ 2.209} & $0.11$ & \num{ 8.873e-07} \\
        PENALTY1 & $50$ & $6082$ & \num{ 1.000} & $0.89$ & \num{ 9.985e-07} & $5773$ & \num{ 0.949} & $0.22$ & \num{ 9.994e-07} & $11797$ & \num{ 1.940} & $0.32$ & \num{ 9.652e-07} \\
        PENALTY1 & $100$ & $13215$ & \num{ 1.000} & $6.44$ & \num{ 9.969e-07} & $11175$ & \num{ 0.846} & $0.43$ & \num{ 7.672e-07} & $21381$ & \num{ 1.618} & $0.57$ & \num{ 9.919e-07} \\
        PFIT1LS$^\ddagger$ & $3$ & $417$ & \num{ 1.000} & $0.11$ & \num{ 2.892e-04} & $351$ & \num{ 0.842} & $0.03$ & \num{ 9.871e-07} & $403$ & \num{ 0.966} & $0.03$ & \num{ 6.795e-07} \\
        PFIT2LS$^\ddagger$ & $3$ & $767$ & \num{ 1.000} & $0.1$ & \num{ 1.243e-02} & $1502^*$ & \num{ 1.958} & $0.17$ & \num{ 2.818e-04} & $1509^*$ & \num{ 1.967} & $0.11$ & \num{ 1.528e-03} \\
        PFIT3LS$^\ddagger$ & $3$ & $907$ & \num{ 1.000} & $0.1$ & \num{ 8.228e-02} & $1502^*$ & \num{ 1.656} & $0.17$ & \num{ 1.504e-02} & $1501^*$ & \num{ 1.655} & $0.11$ & \num{ 3.622e-02} \\
        PFIT4LS$^\ddagger$ & $3$ & $1088$ & \num{ 1.000} & $0.1$ & \num{ 2.673e-01} & $1505^*$ & \num{ 1.383} & $0.15$ & \num{ 5.424e-02} & $1508^*$ & \num{ 1.386} & $0.1$ & \num{ 9.180e-02} \\
        QUARTC & $25$ & $1004$ & \num{ 1.000} & $0.13$ & \num{ 8.886e-07} & $765$ & \num{ 0.762} & $0.04$ & \num{ 4.947e-07} & $1492$ & \num{ 1.486} & $0.05$ & \num{ 4.946e-07} \\
        QUARTC & $100$ & $6061$ & \num{ 1.000} & $2.82$ & \num{ 9.769e-07} & $3277$ & \num{ 0.541} & $0.13$ & \num{ 5.471e-07} & $6508$ & \num{ 1.074} & $0.16$ & \num{ 5.462e-07} \\
        SINEVAL & $2$ & $254$ & \num{ 1.000} & $0.08$ & \num{ 6.706e-09} & $309$ & \num{ 1.217} & $0.04$ & \num{ 3.585e-09} & $522$ & \num{ 2.055} & $0.04$ & \num{ 1.305e-07} \\
        \bottomrule
    \end{tabular}
    \caption{Noiseless Unconstrained \texttt{CUTEst} Problems Tested. $n$ is the number of variables.}
    \label{tab:no-noise4}
\end{table}

\begin{table}[htp]
    \centering
    \footnotesize
    \begin{tabular}{ l  r | r r r r | r r r r | r r r r}
        \toprule
        \multicolumn{2}{c}{ } &  \multicolumn{4}{|c|}{NEWUOA} & \multicolumn{4}{c|}{FD L-BFGS} & \multicolumn{4}{c}{CD L-BFGS} \\
        \midrule
        Problem & $n$ & \#feval & RATIO & CPU & $\phi(x) - \phi^*$ & \#feval & RATIO & CPU & $\phi(x) - \phi^*$ & \#feval & RATIO & CPU & $\phi(x) - \phi^*$ \\
        \midrule
        SINQUAD & $5$ & $125$ & \num{ 1.000} & $0.09$ & \num{ 6.682e-06} & $61$ & \num{ 0.488} & $0.01$ & \num{ 3.061e-06} & $108$ & \num{ 0.864} & $0.01$ & \num{ 3.061e-06} \\
        SINQUAD & $50$ & $3325$ & \num{ 1.000} & $0.55$ & \num{ 1.132e-03} & $795$ & \num{ 0.239} & $0.04$ & \num{ 4.489e-06} & $1562$ & \num{ 0.470} & $0.05$ & \num{ 4.553e-06} \\
        SINQUAD & $100$ & $8318$ & \num{ 1.000} & $4.3$ & \num{ 3.985e-03} & $1545$ & \num{ 0.186} & $0.07$ & \num{ 2.755e-03} & $3062$ & \num{ 0.368} & $0.09$ & \num{ 2.758e-03} \\
        SISSER & $2$ & $24$ & \num{ 1.000} & $0.09$ & \num{ 8.489e-07} & $54$ & \num{ 2.250} & $0.01$ & \num{ 5.876e-07} & $92$ & \num{ 3.833} & $0.01$ & \num{ 5.876e-07} \\
        SPARSQUR & $10$ & $184$ & \num{ 1.000} & $0.09$ & \num{ 8.876e-07} & $184$ & \num{ 1.000} & $0.01$ & \num{ 6.336e-07} & $348$ & \num{ 1.891} & $0.01$ & \num{ 6.336e-07} \\
        SPARSQUR & $50$ & $1645$ & \num{ 1.000} & $0.29$ & \num{ 6.775e-07} & $1046$ & \num{ 0.636} & $0.04$ & \num{ 3.744e-07} & $2065$ & \num{ 1.255} & $0.06$ & \num{ 3.744e-07} \\
        SPARSQUR & $100$ & $2932$ & \num{ 1.000} & $1.37$ & \num{ 8.999e-07} & $2150$ & \num{ 0.733} & $0.09$ & \num{ 3.880e-07} & $4270$ & \num{ 1.456} & $0.12$ & \num{ 3.880e-07} \\
        TOINTGSS & $10$ & $234$ & \num{ 1.000} & $0.09$ & \num{ 1.038e-05} & $24$ & \num{ 0.103} & $0.0$ & \num{ 0.000e+00} & $45$ & \num{ 0.192} & $0.0$ & \num{ 0.000e+00} \\
        TOINTGSS & $50$ & $775$ & \num{ 1.000} & $0.18$ & \num{ 9.723e-06} & $104$ & \num{ 0.134} & $0.0$ & \num{-3.000e-09} & $205$ & \num{ 0.265} & $0.01$ & \num{-3.000e-09} \\
        TOINTGSS & $100$ & $1199$ & \num{ 1.000} & $0.55$ & \num{ 9.778e-06} & $204$ & \num{ 0.170} & $0.02$ & \num{ 4.000e-09} & $405$ & \num{ 0.338} & $0.01$ & \num{ 4.000e-09} \\
        TQUARTIC & $5$ & $92$ & \num{ 1.000} & $0.09$ & \num{ 5.444e-07} & $85$ & \num{ 0.924} & $0.01$ & \num{ 1.078e-08} & $150$ & \num{ 1.630} & $0.01$ & \num{ 1.074e-08} \\
        TQUARTIC & $10$ & $277$ & \num{ 1.000} & $0.09$ & \num{ 9.847e-07} & $172$ & \num{ 0.621} & $0.01$ & \num{ 5.228e-09} & $325$ & \num{ 1.173} & $0.01$ & \num{ 5.313e-09} \\
        TQUARTIC & $50$ & $5022$ & \num{ 1.000} & $0.75$ & \num{ 9.949e-07} & $646$ & \num{ 0.129} & $0.03$ & \num{ 1.681e-08} & $1265$ & \num{ 0.252} & $0.03$ & \num{ 1.725e-08} \\
        TQUARTIC & $100$ & $20515$ & \num{ 1.000} & $10.51$ & \num{ 9.991e-07} & $1538$ & \num{ 0.075} & $0.07$ & \num{ 6.818e-09} & $3053$ & \num{ 0.149} & $0.08$ & \num{ 6.223e-09} \\
        TRIDIA & $10$ & $186$ & \num{ 1.000} & $0.09$ & \num{ 9.019e-07} & $210$ & \num{ 1.129} & $0.01$ & \num{ 5.451e-08} & $396$ & \num{ 2.129} & $0.01$ & \num{ 5.446e-08} \\
        TRIDIA & $20$ & $446$ & \num{ 1.000} & $0.1$ & \num{ 8.062e-07} & $690$ & \num{ 1.547} & $0.04$ & \num{ 8.068e-07} & $1340$ & \num{ 3.004} & $0.04$ & \num{ 8.064e-07} \\
        TRIDIA & $30$ & $768$ & \num{ 1.000} & $0.13$ & \num{ 9.962e-07} & $1481$ & \num{ 1.928} & $0.07$ & \num{ 6.309e-07} & $2906$ & \num{ 3.784} & $0.08$ & \num{ 6.315e-07} \\
        TRIDIA & $50$ & $1446$ & \num{ 1.000} & $0.3$ & \num{ 8.806e-07} & $3128$ & \num{ 2.163} & $0.15$ & \num{ 9.300e-07} & $6187$ & \num{ 4.279} & $0.19$ & \num{ 9.292e-07} \\
        TRIDIA & $100$ & $3450$ & \num{ 1.000} & $2.11$ & \num{ 9.968e-07} & $8783$ & \num{ 2.546} & $0.32$ & \num{ 9.677e-07} & $17468$ & \num{ 5.063} & $0.43$ & \num{ 9.208e-07} \\
        WATSON & $12$ & $4966$ & \num{ 1.000} & $0.21$ & \num{ 9.982e-07} & $1324$ & \num{ 0.267} & $0.08$ & \num{ 9.773e-07} & $2708$ & \num{ 0.545} & $0.11$ & \num{ 9.950e-07} \\
        WATSON & $31$ & $15500^*$ & \num{ 1.000} & $1.11$ & \num{ 6.438e-05} & $13875$ & \num{ 0.895} & $0.65$ & \num{ 9.957e-07} & $15530^*$ & \num{ 1.002} & $0.51$ & \num{ 5.664e-06} \\
        WOODS & $4$ & $497$ & \num{ 1.000} & $0.08$ & \num{ 1.434e-07} & $178$ & \num{ 0.358} & $0.02$ & \num{ 4.447e-09} & $312$ & \num{ 0.628} & $0.02$ & \num{ 4.353e-09} \\
        WOODS & $100$ & $50000^*$ & \num{ 1.000} & $28.69$ & \num{ 1.371e-04} & $2972$ & \num{ 0.059} & $0.13$ & \num{ 7.208e-07} & $5902$ & \num{ 0.118} & $0.16$ & \num{ 6.252e-07} \\
        ZANGWIL2 & $2$ & $12$ & \num{ 1.000} & $0.09$ & \num{ 5.085e-07} & $11$ & \num{ 0.917} & $0.0$ & \num{-1.000e-10} & $19$ & \num{ 1.583} & $0.0$ & \num{-9.999e-11} \\
        \bottomrule
    \end{tabular}
    \caption{Noiseless Unconstrained \texttt{CUTEst} Problems Tested. $n$ is the number of variables.}
    \label{tab:no-noise5}
\end{table}

\end{landscape}

\subsubsection{Noisy Functions}
\label{app:smooth, noisy}

In this section, we list the best optimality gap ($\phi(x) - \phi^*$) and the number of function evaluations (\#feval) needed to achieve this for each problem instance, varying the noise level $\sigma_f \in \{10^{-1}, 10^{-3}, 10^{-5}, 10^{-7}\}$, in Tables \ref{tab:noise1}-\ref{tab:noise13}. Function evaluations marked with a $^*$ denote cases where the algorithm reached the maximum number of function evaluations. Instances where {\sc newuoa} samples a point that satisfies \eqref{eq:term crit} within the initial $2 d + 1$ evaluations are denoted by $^{**}$. Problems marked with a $^\ddagger$ denote cases where {\sc newuoa} and {\sc l-bfgs} converge to different minimizers.

\begin{table}[!htp]
    \centering
    \footnotesize
\begin{tabular}{ l r r | r r | r r | r r}
        \toprule
        \multicolumn{3}{c|}{ } &  \multicolumn{2}{c|}{NEWUOA} & \multicolumn{2}{c|}{FD L-BFGS} & \multicolumn{2}{c}{CD L-BFGS}\\
        \midrule
        Problem & $n$ & $\sigma_f$ & \#feval & $\phi(x) - \phi^*$ & \#feval & $\phi(x) - \phi^*$ & \#feval & $\phi(x) - \phi^*$ \\
        \midrule
        AIRCRFTB & $5$ & \num{1e-01} & $43$ & \num{ 7.257e-01} & $53$ & \num{ 7.277e-01} & $231$ & \num{ 2.623e-01} \\
        AIRCRFTB & $5$ & \num{1e-03} & $70$ & \num{ 1.772e-01} & $715$ & \num{ 9.938e-02} & $904$ & \num{ 6.648e-05} \\
        AIRCRFTB & $5$ & \num{1e-05} & $159$ & \num{ 1.006e-02} & $1221$ & \num{ 3.355e-04} & $927$ & \num{ 1.042e-05} \\
        AIRCRFTB & $5$ & \num{1e-07} & $317$ & \num{ 2.893e-04} & $1176$ & \num{ 4.119e-06} & $1106$ & \num{ 7.139e-12} \\
        ALLINITU & $4$ & \num{1e-01} & $26$ & \num{ 1.020e-01} & $24$ & \num{ 2.881e+00} & $108$ & \num{ 1.980e-03} \\
        ALLINITU & $4$ & \num{1e-03} & $38$ & \num{ 9.901e-04} & $97$ & \num{ 2.614e-04} & $189$ & \num{ 7.467e-07} \\
\bottomrule
    \end{tabular}
    \caption{Noisy Unconstrained \texttt{CUTEst} Problems Tested. $n$ is the number of variables. $\sigma_f$ is the standard deviation of the noise.}
    \label{tab:noise1}
\end{table}

\begin{table}[!htp]
    \centering
    \footnotesize
    \begin{tabular}{ l r r | r r | r r | r r}
        \toprule
        \multicolumn{3}{c|}{ } &  \multicolumn{2}{c|}{NEWUOA} & \multicolumn{2}{c|}{FD L-BFGS} & \multicolumn{2}{c}{CD L-BFGS}\\
        \midrule
        Problem & $n$ & $\sigma_f$ & \#feval & $\phi(x) - \phi^*$ & \#feval & $\phi(x) - \phi^*$ & \#feval & $\phi(x) - \phi^*$ \\
        \midrule
        BROWNAL & $10$ & \num{1e-01} & $46$ & \num{ 9.212e-02} & $53$ & \num{ 2.719e-02} & $384$ & \num{ 4.202e-05} \\
        BROWNAL & $10$ & \num{1e-03} & $123$ & \num{ 1.232e-02} & $65$ & \num{ 5.759e-05} & $223$ & \num{ 2.969e-05} \\
        BROWNAL & $10$ & \num{1e-05} & $248$ & \num{ 3.218e-04} & $77$ & \num{ 3.062e-05} & $740$ & \num{ 1.092e-05} \\
        BROWNAL & $10$ & \num{1e-07} & $458$ & \num{ 6.915e-06} & $77$ & \num{ 2.974e-05} & $821$ & \num{ 2.971e-05} \\
        BROWNAL & $100$ & \num{1e-01} & $610$ & \num{ 1.164e+00} & $725$ & \num{ 5.055e-02} & $2664$ & \num{ 1.297e-05} \\
        BROWNAL & $100$ & \num{1e-03} & $4084$ & \num{ 2.877e-02} & $725$ & \num{ 5.541e-04} & $3677$ & \num{ 1.287e-05} \\
\bottomrule
    \end{tabular}
    \caption{Noisy Unconstrained \texttt{CUTEst} Problems Tested. $n$ is the number of variables. $\sigma_f$ is the standard deviation of the noise.}
    \label{tab:noise2}
\end{table}

\begin{table}[!htp]
    \centering
    \footnotesize
    \begin{tabular}{ l r r | r r | r r | r r}
        \toprule
        \multicolumn{3}{c|}{ } &  \multicolumn{2}{c|}{NEWUOA} & \multicolumn{2}{c|}{FD L-BFGS} & \multicolumn{2}{c}{CD L-BFGS}\\
        \midrule
        Problem & $n$ & $\sigma_f$ & \#feval & $\phi(x) - \phi^*$ & \#feval & $\phi(x) - \phi^*$ & \#feval & $\phi(x) - \phi^*$ \\
        \midrule
        DENSCHNB & $2$ & \num{1e-01} & $14$ & \num{ 2.282e-02} & $23$ & \num{ 2.369e-01} & $293$ & \num{ 4.982e-04} \\
        DENSCHNB & $2$ & \num{1e-03} & $17$ & \num{ 1.807e-05} & $40$ & \num{ 8.803e-04} & $387$ & \num{ 1.432e-06} \\
        DENSCHNB & $2$ & \num{1e-05} & $27$ & \num{ 2.693e-07} & $61$ & \num{ 1.017e-06} & $42$ & \num{ 5.983e-11} \\
        DENSCHNB & $2$ & \num{1e-07} & $26$ & \num{ 7.582e-09} & $32$ & \num{ 2.194e-08} & $136$ & \num{ 9.498e-12} \\
        DENSCHNC & $2$ & \num{1e-01} & $30$ & \num{-7.755e-02} & $18$ & \num{ 5.033e-02} & $281$ & \num{ 3.323e-02} \\
        DENSCHNC & $2$ & \num{1e-03} & $62$ & \num{ 2.780e-05} & $80$ & \num{ 1.775e-02} & $219$ & \num{ 2.208e-05} \\
\bottomrule
    \end{tabular}
    \caption{Noisy Unconstrained \texttt{CUTEst} Problems Tested. $n$ is the number of variables. $\sigma_f$ is the standard deviation of the noise.}
    \label{tab:noise3}
\end{table}

\begin{table}[!htp]
    \centering
    \footnotesize
    \begin{tabular}{ l r r | r r | r r | r r}
        \toprule
        \multicolumn{3}{c|}{ } &  \multicolumn{2}{c|}{NEWUOA} & \multicolumn{2}{c|}{FD L-BFGS} & \multicolumn{2}{c}{CD L-BFGS}\\
        \midrule
        Problem & $n$ & $\sigma_f$ & \#feval & $\phi(x) - \phi^*$ & \#feval & $\phi(x) - \phi^*$ & \#feval & $\phi(x) - \phi^*$ \\
        \midrule
        DIXMAANC & $15$ & \num{1e-01} & $123$ & \num{ 3.678e-01} & $696$ & \num{ 9.417e-01} & $413$ & \num{ 4.277e-03} \\
        DIXMAANC & $15$ & \num{1e-03} & $214$ & \num{ 1.694e-03} & $784$ & \num{ 9.657e-03} & $627$ & \num{ 7.834e-05} \\
        DIXMAANC & $15$ & \num{1e-05} & $252$ & \num{ 1.523e-04} & $621$ & \num{ 3.658e-05} & $413$ & \num{ 2.695e-08} \\
        DIXMAANC & $15$ & \num{1e-07} & $249$ & \num{ 2.681e-07} & $528$ & \num{ 4.150e-08} & $627$ & \num{ 1.283e-10} \\
        DIXMAANC & $90$ & \num{1e-01} & $853$ & \num{ 1.678e+00} & $2077$ & \num{ 3.626e+00} & $2649$ & \num{ 1.112e-02} \\
        DIXMAANC & $90$ & \num{1e-03} & $2546$ & \num{ 1.886e-02} & $1016$ & \num{ 6.657e-02} & $1880$ & \num{ 4.973e-05} \\
\bottomrule
    \end{tabular}
    \caption{Noisy Unconstrained \texttt{CUTEst} Problems Tested. $n$ is the number of variables. $\sigma_f$ is the standard deviation of the noise.}
    \label{tab:noise4}
\end{table}

\begin{table}[!htp]
    \centering
    \footnotesize
    \begin{tabular}{ l r r | r r | r r | r r}
        \toprule
        \multicolumn{3}{c|}{ } &  \multicolumn{2}{c|}{NEWUOA} & \multicolumn{2}{c|}{FD L-BFGS} & \multicolumn{2}{c}{CD L-BFGS}\\
        \midrule
        Problem & $n$ & $\sigma_f$ & \#feval & $\phi(x) - \phi^*$ & \#feval & $\phi(x) - \phi^*$ & \#feval & $\phi(x) - \phi^*$ \\
        \midrule
        DIXMAANG & $15$ & \num{1e-01} & $118$ & \num{ 1.301e+00} & $53$ & \num{ 8.760e-01} & $1730$ & \num{ 2.446e-02} \\
        DIXMAANG & $15$ & \num{1e-03} & $219$ & \num{ 3.590e-03} & $532$ & \num{ 8.830e-03} & $1104$ & \num{ 4.679e-05} \\
        DIXMAANG & $15$ & \num{1e-05} & $259$ & \num{ 1.628e-05} & $310$ & \num{ 9.965e-04} & $1730$ & \num{ 9.830e-08} \\
        DIXMAANG & $15$ & \num{1e-07} & $354$ & \num{ 1.419e-07} & $1235$ & \num{ 2.106e-06} & $1730$ & \num{ 3.410e-10} \\
        DIXMAANG & $90$ & \num{1e-01} & $1253$ & \num{ 1.630e+00} & $2284$ & \num{ 1.752e+00} & $4496$ & \num{ 9.061e-02} \\
        DIXMAANG & $90$ & \num{1e-03} & $3690$ & \num{ 9.785e-02} & $2011$ & \num{ 3.834e-02} & $6194$ & \num{ 1.981e-03} \\
\bottomrule
    \end{tabular}
    \caption{Noisy Unconstrained \texttt{CUTEst} Problems Tested. $n$ is the number of variables. $\sigma_f$ is the standard deviation of the noise.}
    \label{tab:noise5}
\end{table}

\begin{table}[!htp]
    \centering
    \footnotesize
    \begin{tabular}{ l r r | r r | r r | r r}
        \toprule
        \multicolumn{3}{c|}{ } &  \multicolumn{2}{c|}{NEWUOA} & \multicolumn{2}{c|}{FD L-BFGS} & \multicolumn{2}{c}{CD L-BFGS}\\
        \midrule
        Problem & $n$ & $\sigma_f$ & \#feval & $\phi(x) - \phi^*$ & \#feval & $\phi(x) - \phi^*$ & \#feval & $\phi(x) - \phi^*$ \\
        \midrule
        DIXMAANK & $15$ & \num{1e-01} & $123$ & \num{ 2.036e-01} & $70$ & \num{ 3.310e-01} & $270$ & \num{ 5.345e-02} \\
        DIXMAANK & $15$ & \num{1e-03} & $226$ & \num{ 4.148e-02} & $1229$ & \num{ 1.660e-02} & $1024$ & \num{ 8.734e-04} \\
        DIXMAANK & $15$ & \num{1e-05} & $533$ & \num{ 4.054e-04} & $759$ & \num{ 7.894e-05} & $1730$ & \num{ 3.452e-07} \\
        DIXMAANK & $15$ & \num{1e-07} & $710$ & \num{ 1.111e-06} & $724$ & \num{ 3.219e-07} & $1595$ & \num{ 4.578e-09} \\
        DIXMAANK & $90$ & \num{1e-01} & $962$ & \num{ 2.115e+00} & $961$ & \num{ 9.536e-01} & $1880$ & \num{ 3.049e-01} \\
        DIXMAANK & $90$ & \num{1e-03} & $1921$ & \num{ 4.371e-02} & $2718$ & \num{ 9.449e-02} & $10211$ & \num{ 3.178e-03} \\
\bottomrule
    \end{tabular}
    \caption{Noisy Unconstrained \texttt{CUTEst} Problems Tested. $n$ is the number of variables. $\sigma_f$ is the standard deviation of the noise.}
    \label{tab:noise6}
\end{table}

\begin{table}[!htp]
    \centering
    \footnotesize
    \begin{tabular}{ l r r | r r | r r | r r}
        \toprule
        \multicolumn{3}{c|}{ } &  \multicolumn{2}{c|}{NEWUOA} & \multicolumn{2}{c|}{FD L-BFGS} & \multicolumn{2}{c}{CD L-BFGS}\\
        \midrule
        Problem & $n$ & $\sigma_f$ & \#feval & $\phi(x) - \phi^*$ & \#feval & $\phi(x) - \phi^*$ & \#feval & $\phi(x) - \phi^*$ \\
        \midrule
        EIGENBLS$^\ddagger$ & $6$ & \num{1e-01} & $51$ & \num{ 1.421e-01} & $32$ & \num{ 1.096e+00} & $252$ & \num{ 3.986e-03} \\
        EIGENBLS$^\ddagger$ & $6$ & \num{1e-03} & $109$ & \num{ 3.843e-02} & $240$ & \num{ 1.223e-03} & $376$ & \num{ 4.336e-05} \\
        EIGENBLS$^\ddagger$ & $6$ & \num{1e-05} & $217$ & \num{-1.840e-01} & $247$ & \num{ 2.792e-05} & $299$ & \num{ 1.340e-07} \\
        EIGENBLS$^\ddagger$ & $6$ & \num{1e-07} & $307$ & \num{-1.849e-01} & $190$ & \num{ 9.576e-08} & $256$ & \num{ 5.404e-10} \\
        EIGENBLS & $110$ & \num{1e-01} & $705$ & \num{ 3.389e+00} & $1761$ & \num{ 1.122e+01} & $4023$ & \num{ 1.707e+00} \\
        EIGENBLS & $110$ & \num{1e-03} & $5295$ & \num{ 1.452e+00} & $2754$ & \num{ 1.596e+00} & $7292$ & \num{ 1.552e+00} \\
\bottomrule
    \end{tabular}
    \caption{Noisy Unconstrained \texttt{CUTEst} Problems Tested. $n$ is the number of variables. $\sigma_f$ is the standard deviation of the noise.}
    \label{tab:noise7}
\end{table}

\begin{table}[!htp]
    \centering
    \footnotesize
    \begin{tabular}{ l r r | r r | r r | r r}
        \toprule
        \multicolumn{3}{c|}{ } &  \multicolumn{2}{c|}{NEWUOA} & \multicolumn{2}{c|}{FD L-BFGS} & \multicolumn{2}{c}{CD L-BFGS}\\
        \midrule
        Problem & $n$ & $\sigma_f$ & \#feval & $\phi(x) - \phi^*$ & \#feval & $\phi(x) - \phi^*$ & \#feval & $\phi(x) - \phi^*$ \\
        \midrule
        FREUROTH & $2$ & \num{1e-01} & $38$ & \num{ 6.304e-02} & $35$ & \num{ 5.015e+01} & $368$ & \num{ 3.877e-04} \\
        FREUROTH & $2$ & \num{1e-03} & $35$ & \num{ 5.507e-01} & $347$ & \num{ 1.660e-03} & $141$ & \num{ 1.496e-05} \\
        FREUROTH & $2$ & \num{1e-05} & $70$ & \num{ 8.659e-07} & $83$ & \num{ 9.381e-05} & $387$ & \num{ 3.302e-08} \\
        FREUROTH & $2$ & \num{1e-07} & $73$ & \num{ 4.594e-10} & $138$ & \num{ 1.112e-06} & $261$ & \num{ 4.996e-10} \\
        FREUROTH & $10$ & \num{1e-01} & $125$ & \num{ 4.555e+01} & $137$ & \num{ 5.548e+01} & $1466$ & \num{ 8.708e-03} \\
        FREUROTH & $10$ & \num{1e-03} & $371$ & \num{ 3.964e-03} & $467$ & \num{ 1.341e-01} & $562$ & \num{ 4.686e-05} \\
\bottomrule
    \end{tabular}
    \caption{Noisy Unconstrained \texttt{CUTEst} Problems Tested. $n$ is the number of variables. $\sigma_f$ is the standard deviation of the noise.}
    \label{tab:noise8}
\end{table}

\begin{table}[!htp]
    \centering
    \footnotesize
    \begin{tabular}{ l r r | r r | r r | r r}
        \toprule
        \multicolumn{3}{c|}{ } &  \multicolumn{2}{c|}{NEWUOA} & \multicolumn{2}{c|}{FD L-BFGS} & \multicolumn{2}{c}{CD L-BFGS}\\
        \midrule
        Problem & $n$ & $\sigma_f$ & \#feval & $\phi(x) - \phi^*$ & \#feval & $\phi(x) - \phi^*$ & \#feval & $\phi(x) - \phi^*$ \\
        \midrule
        MEXHAT & $2$ & \num{1e-01} & $30$ & \num{ 1.166e-02} & $98$ & \num{ 4.477e-01} & $146$ & \num{ 4.231e-01} \\
        MEXHAT & $2$ & \num{1e-03} & $43$ & \num{ 1.102e-02} & $80$ & \num{ 8.450e-01} & $709$ & \num{ 4.227e-01} \\
        MEXHAT & $2$ & \num{1e-05} & $45$ & \num{ 1.110e-02} & $203$ & \num{ 4.173e-01} & $232$ & \num{ 1.568e-03} \\
        MEXHAT & $2$ & \num{1e-07} & $44$ & \num{ 1.036e-02} & $210$ & \num{ 4.226e-01} & $660$ & \num{ 1.383e-04} \\
        MOREBV & $10$ & \num{1e-01} & $1^{**}$ & \num{ 1.599e-02} & $19$ & \num{ 1.599e-02} & $21$ & \num{ 1.599e-02} \\
        MOREBV & $10$ & \num{1e-03} & $81$ & \num{ 1.262e-02} & $131$ & \num{ 1.451e-02} & $937$ & \num{ 1.506e-03} \\
\bottomrule
    \end{tabular}
    \caption{Noisy Unconstrained \texttt{CUTEst} Problems Tested. $n$ is the number of variables. $\sigma_f$ is the standard deviation of the noise.}
    \label{tab:noise9}
\end{table}

\begin{table}[!htp]
    \centering
    \footnotesize
    \begin{tabular}{ l r r | r r | r r | r r}
        \toprule
        \multicolumn{3}{c|}{ } &  \multicolumn{2}{c|}{NEWUOA} & \multicolumn{2}{c|}{FD L-BFGS} & \multicolumn{2}{c}{CD L-BFGS}\\
        \midrule
        Problem & $n$ & $\sigma_f$ & \#feval & $\phi(x) - \phi^*$ & \#feval & $\phi(x) - \phi^*$ & \#feval & $\phi(x) - \phi^*$ \\
        \midrule
        NONDIA & $50$ & \num{1e-01} & $317$ & \num{ 8.665e-01} & $269$ & \num{ 1.505e+00} & $423$ & \num{ 1.520e+00} \\
        NONDIA & $50$ & \num{1e-03} & $635$ & \num{ 8.814e-03} & $375$ & \num{ 1.498e+00} & $1659$ & \num{ 2.467e-03} \\
        NONDIA & $50$ & \num{1e-05} & $851$ & \num{ 7.108e-05} & $3599$ & \num{ 4.640e-04} & $1659$ & \num{ 6.908e-06} \\
        NONDIA & $50$ & \num{1e-07} & $1132$ & \num{ 1.517e-06} & $1731$ & \num{ 1.971e-05} & $1556$ & \num{ 1.127e-08} \\
        NONDIA & $90$ & \num{1e-01} & $573$ & \num{ 3.893e-01} & $1087$ & \num{ 6.154e-01} & $744$ & \num{ 5.843e-01} \\
        NONDIA & $90$ & \num{1e-03} & $1374$ & \num{ 3.301e-01} & $563$ & \num{ 5.766e-01} & $2052$ & \num{ 9.232e-03} \\
\bottomrule
    \end{tabular}
    \caption{Noisy Unconstrained \texttt{CUTEst} Problems Tested. $n$ is the number of variables. $\sigma_f$ is the standard deviation of the noise.}
    \label{tab:noise10}
\end{table}

\begin{table}[!htp]
    \centering
    \footnotesize
    \begin{tabular}{ l r r | r r | r r | r r}
        \toprule
        \multicolumn{3}{c|}{ } &  \multicolumn{2}{c|}{NEWUOA} & \multicolumn{2}{c|}{FD L-BFGS} & \multicolumn{2}{c}{CD L-BFGS}\\
        \midrule
        Problem & $n$ & $\sigma_f$ & \#feval & $\phi(x) - \phi^*$ & \#feval & $\phi(x) - \phi^*$ & \#feval & $\phi(x) - \phi^*$ \\
        \midrule
        PFIT3LS$^\ddagger$ & $3$ & \num{1e-01} & $162$ & \num{ 1.868e+02} & $390$ & \num{ 7.944e+01} & $817$ & \num{ 3.682e-01} \\
        PFIT3LS$^\ddagger$ & $3$ & \num{1e-03} & $546$ & \num{ 7.561e+00} & $640$ & \num{ 9.159e-01} & $612$ & \num{ 3.276e-02} \\
        PFIT3LS$^\ddagger$ & $3$ & \num{1e-05} & $1012$ & \num{ 9.384e-02} & $615$ & \num{ 4.587e-02} & $702$ & \num{ 3.154e-02} \\
        PFIT3LS$^\ddagger$ & $3$ & \num{1e-07} & $1012$ & \num{ 8.229e-02} & $444$ & \num{ 4.374e-02} & $1504^*$ & \num{ 2.782e-02} \\
        PFIT4LS$^\ddagger$ & $3$ & \num{1e-01} & $234$ & \num{ 2.464e+02} & $78$ & \num{ 2.144e+03} & $723$ & \num{ 1.842e-01} \\
        PFIT4LS$^\ddagger$ & $3$ & \num{1e-03} & $612$ & \num{ 2.736e+01} & $758$ & \num{ 6.570e+00} & $715$ & \num{ 1.234e-01} \\
\bottomrule
    \end{tabular}
    \caption{Noisy Unconstrained \texttt{CUTEst} Problems Tested. $n$ is the number of variables. $\sigma_f$ is the standard deviation of the noise.}
    \label{tab:noise11}
\end{table}

\begin{table}[!htp]
    \centering
    \footnotesize
    \begin{tabular}{ l r r | r r | r r | r r}
        \toprule
        \multicolumn{3}{c|}{ } &  \multicolumn{2}{c|}{NEWUOA} & \multicolumn{2}{c|}{FD L-BFGS} & \multicolumn{2}{c}{CD L-BFGS}\\
        \midrule
        Problem & $n$ & $\sigma_f$ & \#feval & $\phi(x) - \phi^*$ & \#feval & $\phi(x) - \phi^*$ & \#feval & $\phi(x) - \phi^*$ \\
        \midrule
        TOINTGSS & $10$ & \num{1e-01} & $65$ & \num{ 4.850e-01} & $759$ & \num{ 2.072e+00} & $101$ & \num{ 1.831e-02} \\
        TOINTGSS & $10$ & \num{1e-03} & $164$ & \num{ 1.448e-02} & $301$ & \num{ 3.304e-02} & $165$ & \num{ 8.167e-06} \\
        TOINTGSS & $10$ & \num{1e-05} & $211$ & \num{ 2.845e-05} & $124$ & \num{ 3.624e-04} & $200$ & \num{ 2.744e-08} \\
        TOINTGSS & $10$ & \num{1e-07} & $277$ & \num{ 2.113e-07} & $124$ & \num{ 3.539e-06} & $200$ & \num{ 5.802e-11} \\
        TOINTGSS & $50$ & \num{1e-01} & $388$ & \num{ 1.033e+00} & $563$ & \num{ 1.183e+01} & $724$ & \num{ 6.155e-02} \\
        TOINTGSS & $50$ & \num{1e-03} & $878$ & \num{ 1.276e-02} & $415$ & \num{ 4.651e-02} & $2132$ & \num{ 1.554e-05} \\
\bottomrule
    \end{tabular}
    \caption{Noisy Unconstrained \texttt{CUTEst} Problems Tested. $n$ is the number of variables. $\sigma_f$ is the standard deviation of the noise.}
    \label{tab:noise12}
\end{table}

% \newpage

\begin{table}[!htp]
    \centering
    \footnotesize
    \begin{tabular}{ l r r | r r | r r | r r}
        \toprule
        \multicolumn{3}{c|}{ } &  \multicolumn{2}{c|}{NEWUOA} & \multicolumn{2}{c|}{FD L-BFGS} & \multicolumn{2}{c}{CD L-BFGS}\\
        \midrule
        Problem & $n$ & $\sigma_f$ & \#feval & $\phi(x) - \phi^*$ & \#feval & $\phi(x) - \phi^*$ & \#feval & $\phi(x) - \phi^*$ \\
        \midrule
        WATSON & $12$ & \num{1e-01} & $91$ & \num{ 1.602e-01} & $260$ & \num{ 2.747e-01} & $407$ & \num{ 9.362e-02} \\
        WATSON & $12$ & \num{1e-03} & $170$ & \num{ 7.112e-02} & $752$ & \num{ 7.265e-03} & $737$ & \num{ 8.831e-03} \\
        WATSON & $12$ & \num{1e-05} & $708$ & \num{ 6.531e-04} & $1487$ & \num{ 2.403e-03} & $1470$ & \num{ 1.808e-04} \\
        WATSON & $12$ & \num{1e-07} & $1664$ & \num{ 2.355e-04} & $1077$ & \num{ 1.868e-04} & $1559$ & \num{ 1.367e-05} \\
        WATSON & $31$ & \num{1e-01} & $283$ & \num{ 1.444e+00} & $567$ & \num{ 2.596e+00} & $1379$ & \num{ 2.727e-01} \\
        WATSON & $31$ & \num{1e-03} & $982$ & \num{ 1.556e-01} & $1465$ & \num{ 1.648e-01} & $4316$ & \num{ 6.367e-03} \\
\bottomrule
    \end{tabular}
    \caption{Noisy Unconstrained \texttt{CUTEst} Problems Tested. $n$ is the number of variables. $\sigma_f$ is the standard deviation of the noise.}
    \label{tab:noise13}
\end{table}

% \newpage

\subsection{Nonlinear Least Squares Problems}
\label{app: ls}

\subsubsection{Noiseless Functions}
\label{app:ls det}

In this section, we list the number of function evaluations (\#feval), CPU time (CPU), and optimality gap ($\phi(x) - \phi^*$) for each problem instance in Tables \ref{tab:ls_no-noise}. Function evaluations marked with a $^*$ denote cases where the algorithm reached the maximum number of function evaluations. 

\begin{table}[!htp] %DETERMINISTIC
	\centering
    \footnotesize
    \begin{tabular}{ l  r | r r   | r r}
        \toprule
        \multicolumn{2}{c}{ } &  \multicolumn{2}{|c|}{LMDER} & \multicolumn{2}{c}{DFOLS} \\
        \midrule
        Problem & $(n,m)$ & \#feval  & $\phi(x) - \phi^*$ &  \#feval   & $\phi(x) - \phi^*$ \\
        \midrule
LINEAR(FR) & $ (9,45)$ & $21$    & \num{ 1.421e-14} & $45$    & \num{-1.421e-14} \\
LINEAR(FR) & $ (9,45)$ & $21$    & \num{ 1.279e-13} & $49$    & \num{ 0.000e+00} \\
LINEAR(R1) & $ (7,35)$ & $16$    & \num{-3.099e-07} & $61$    & \num{-3.099e-07} \\
LINEAR(R1) & $ (7,35)$ & $16$    & \num{-3.099e-07} & $56$    & \num{-3.099e-07} \\
LINEAR(R10RC) & $ (7,35)$ & $16$    & \num{ 1.493e-08} & $56$    & \num{ 1.493e-08} \\
LINEAR(R10RC) & $ (7,35)$ & $16$    & \num{ 1.493e-08} & $62$    & \num{ 1.493e-08} \\
ROSENBR & $ (2,2)$ & $39$    & \num{ 0.000e+00} & $50$    & \num{ 3.788e-24} \\
ROSENBR & $ (2,2)$ & $17$    & \num{ 0.000e+00} & $14$    & \num{ 1.498e-20} \\
HELIX & $ (3,3)$ & $49$    & \num{ 3.941e-59} & $44$    & \num{ 2.532e-28} \\
HELIX & $ (3,3)$ & $57$    & \num{ 5.921e-47} & $80$    & \num{ 5.476e-20} \\
POWELLSG & $ (4,4)$ & $331$    & \num{ 6.610e-35} & $54$    & \num{ 1.074e-18} \\
POWELLSG & $ (4,4)$ & $346$    & \num{ 7.106e-35} & $47$    & \num{ 1.467e-14} \\
FREUROTH & $ (2,2)$ & $60$    & \num{ 3.688e-06} & $122$    & \num{ 3.679e-06} \\
FREUROTH & $ (2,2)$ & $87$    & \num{ 3.793e-06} & $107$    & \num{ 3.679e-06} \\
BARD & $ (3,15)$ & $21$    & \num{ 3.066e-10} & $32$    & \num{ 3.066e-10} \\
BARD & $ (3,15)$ & $169$    & \num{ 1.742e+01} & $107$    & \num{ 1.066e-01} \\
KOWOSB & $ (4,11)$ & $99$    & \num{ 4.429e-12} & $77$    & \num{ 3.849e-12} \\
MEYER & $ (3,16)$ & $456$    & \num{-4.829e-06} & $1500^*$    & \num{ 3.715e+04} \\
WATSON & $ (6,31)$ & $57$    & \num{ 5.359e-11} & $81$    & \num{ 5.355e-11} \\
WATSON & $ (6,31)$ & $85$    & \num{ 5.357e-11} & $105$    & \num{ 5.355e-11} \\
WATSON & $ (9,31)$ & $144$    & \num{ 1.393e-13} & $409$    & \num{ 1.468e-13} \\
WATSON & $ (9,31)$ & $134$    & \num{ 1.382e-13} & $499$    & \num{ 1.454e-13} \\
WATSON & $ (12,31)$ & $276$    & \num{ 1.198e-15} & $340$    & \num{ 2.214e-09} \\
WATSON & $ (12,31)$ & $181$    & \num{ 3.636e-16} & $1013$    & \num{ 2.246e-09} \\
BOX3D & $ (3,10)$ & $25$    & \num{ 3.390e-32} & $52$    & \num{ 1.921e-27} \\
JENSMP & $ (2,10)$ & $46$    & \num{-1.764e-05} & $70$    & \num{-1.764e-05} \\
BROWNDEN & $ (4,20)$ & $94$    & \num{ 1.684e-03} & $167$    & \num{ 1.626e-03} \\
BROWNDEN & $ (4,20)$ & $96$    & \num{ 1.733e-03} & $189$    & \num{ 1.626e-03} \\
CHEBYQUAD & $ (6,6)$ & $59$    & \num{ 4.434e-32} & $54$    & \num{ 1.239e-23} \\
CHEBYQUAD & $ (7,7)$ & $58$    & \num{ 6.915e-32} & $49$    & \num{ 5.120e-29} \\
CHEBYQUAD & $ (8,8)$ & $168$    & \num{-2.725e-10} & $258$    & \num{-2.743e-10} \\
CHEBYQUAD & $ (9,9)$ & $94$    & \num{ 4.143e-32} & $87$    & \num{ 4.535e-28} \\
CHEBYQUAD & $ (10,10)$ & $108$    & \num{ 1.731e-03} & $191$    & \num{-3.036e-10} \\
CHEBYQUAD & $ (11,11)$ & $320$    & \num{-4.428e-10} & $352$    & \num{-4.481e-10} \\
BROWNAL & $ (10,10)$ & $79$    & \num{ 2.840e-29} & $52$    & \num{ 7.676e-19} \\
OSBORNE1 & $ (5,33)$ & $157$    & \num{-3.025e-12} & $783$    & \num{ 3.737e-12} \\
OSBORNE2 & $ (11,65)$ & $101$    & \num{-3.705e-09} & $150$    & \num{-3.706e-09} \\
OSBORNE2 & $ (11,65)$ & $148$    & \num{ 1.750e+00} & $138$    & \num{ 1.750e+00} \\
BDQRTIC & $ (8,8)$ & $508$    & \num{ 3.759e-06} & $336$    & \num{ 8.111e-06} \\
BDQRTIC & $ (10,12)$ & $662$    & \num{ 3.152e-06} & $457$    & \num{ 3.874e-06} \\
BDQRTIC & $ (11,14)$ & $1046$    & \num{ 3.296e-06} & $554$    & \num{ 2.046e-05} \\
BDQRTIC & $ (12,16)$ & $1119$    & \num{-1.268e-06} & $561$    & \num{ 5.775e-05} \\
    \bottomrule
    \end{tabular}
    \caption{Benchmarking Problems Tested. $n$ is the number of variables and $m$ is the dimension of the residual vector.}
    \label{tab:ls_no-noise}
\end{table}

\begin{table}[!htp] %DETERMINISTIC
    \centering
    \footnotesize
    \begin{tabular}{ l  r | r r   | r r}
        \toprule
        \multicolumn{2}{c}{ } &  \multicolumn{2}{|c|}{LMDER} & \multicolumn{2}{c}{DFOLS} \\
        \midrule
        Problem & $(n,m)$ & \#feval  & $\phi(x) - \phi^*$ &  \#feval   & $\phi(x) - \phi^*$ \\
        \midrule
CUBE & $ (5,5)$ & $388$    & \num{ 0.000e+00} & $295$    & \num{ 3.119e-18} \\
CUBE & $ (6,6)$ & $1024$    & \num{ 0.000e+00} & $714$    & \num{ 4.127e-24} \\
CUBE & $ (8,8)$ & $4008^*$    & \num{ 2.090e-08} & $4000^*$    & \num{ 4.753e-11} \\
MANCINO & $ (5,5)$ & $19$    & \num{ 4.224e-22} & $25$    & \num{ 6.176e-20} \\
MANCINO & $ (5,5)$ & $25$    & \num{ 4.686e-22} & $27$    & \num{ 8.536e-18} \\
MANCINO & $ (8,8)$ & $28$    & \num{ 5.795e-22} & $23$    & \num{ 3.152e-12} \\
MANCINO & $ (10,10)$ & $34$    & \num{ 8.218e-22} & $25$    & \num{ 3.069e-11} \\
MANCINO & $ (12,12)$ & $53$    & \num{ 1.322e-22} & $32$    & \num{ 2.454e-17} \\
MANCINO & $ (12,12)$ & $53$    & \num{ 5.201e-22} & $35$    & \num{ 2.798e-09} \\
HEART8LS & $ (8,8)$ & $37$    & \num{ 3.402e-30} & $46$    & \num{ 1.082e-23} \\
HEART8LS & $ (8,8)$ & $109$    & \num{ 3.402e-30} & $2293$    & \num{ 1.001e-10} \\
    \bottomrule
    \end{tabular}
    \caption{Benchmarking Problems Tested. $n$ is the number of variables and $m$ is the dimension of the residual vector.}
    \label{tab:ls_no-noise2}
\end{table}

% \newpage
\subsubsection{Noisy Functions}
\label{app: ls noisy}

In this section, we list the best optimality gap ($\phi(x) - \phi^*$) and the number of function evaluations (\#feval) needed to achieve this for each problem instance, varying the noise level $\sigma_f \in \{10^{-1}, 10^{-3}, 10^{-5}, 10^{-7}\}$, in Tables \ref{tab:ls_noise1}-\ref{tab:ls_noise4}. Function evaluations marked with a $^*$ denote cases where the algorithm reached the maximum number of function evaluations. 

\begin{table}[!htp] %NOISY CASE 1
	\centering
    \footnotesize
\begin{tabular}{ l r  r | r r   | r r}
        \toprule
        \multicolumn{3}{c}{ } &  \multicolumn{2}{|c|}{LMDER} & \multicolumn{2}{c}{DFOLS} \\
        \midrule
        Problem & $(n,m)$ &$\sigma_f$ & \#feval  & $\phi(x) - \phi^*$ &  \#feval   & $\phi(x) - \phi^*$ \\
        \midrule
LINEAR(FR) & $ (9,45)$ & \num{1e-01} & $65$    & \num{ 3.565e+01} & $70$    & \num{ 3.641e+01} \\
LINEAR(FR) & $ (9,45)$ & \num{1e-03} & $126$    & \num{ 5.008e-02} & $97$    & \num{ 6.606e-02} \\
LINEAR(FR) & $ (9,45)$ & \num{1e-05} & $115$    & \num{ 5.170e-04} & $81$    & \num{ 2.889e-05} \\
LINEAR(FR) & $ (9,45)$ & \num{1e-07} & $114$    & \num{ 5.523e-06} & $76$    & \num{ 2.296e-06} \\
LINEAR(FR) & $ (9,45)$ & \num{1e-01} & $229$    & \num{ 2.822e+02} & $77$    & \num{ 1.112e+03} \\
LINEAR(FR) & $ (9,45)$ & \num{1e-03} & $115$    & \num{ 6.086e-02} & $89$    & \num{ 2.179e-02} \\
\bottomrule
    \end{tabular}
    \caption{Benchmarking Problems Tested. $n$ is the number of variables and $m$ is the dimension of the residual vector.}
    \label{tab:ls_noise1}
\end{table}

\begin{table}[!htp] %NOISY CASE2
	\centering
    \footnotesize
    \begin{tabular}{ l r  r | r r   | r r}
        \toprule
        \multicolumn{3}{c}{ } &  \multicolumn{2}{|c|}{LMDER} & \multicolumn{2}{c}{DFOLS} \\
        \midrule
        Problem & $(n,m)$ &$\sigma_f$ & \#feval  & $\phi(x) - \phi^*$ &  \#feval   & $\phi(x) - \phi^*$ \\
        \midrule
FREUROTH & $ (2,2)$ & \num{1e-01} & $42$    & \num{ 1.724e+00} & $68$    & \num{ 4.634e-02} \\
FREUROTH & $ (2,2)$ & \num{1e-03} & $48$    & \num{ 7.582e-03} & $75$    & \num{ 4.555e-06} \\
FREUROTH & $ (2,2)$ & \num{1e-05} & $59$    & \num{ 7.704e-05} & $113$    & \num{ 3.717e-06} \\
FREUROTH & $ (2,2)$ & \num{1e-07} & $64$    & \num{ 4.523e-06} & $91$    & \num{ 3.679e-06} \\
FREUROTH & $ (2,2)$ & \num{1e-01} & $55$    & \num{ 2.546e+00} & $83$    & \num{ 1.310e-02} \\
FREUROTH & $ (2,2)$ & \num{1e-03} & $66$    & \num{ 1.265e-02} & $78$    & \num{ 1.331e-05} \\
\bottomrule
    \end{tabular}
    \caption{Benchmarking Problems Tested. $n$ is the number of variables and $m$ is the dimension of the residual vector.}
    \label{tab:ls_noise2}
\end{table}

\begin{table}[!htp] %NOISY CASE3
	\centering
    \footnotesize
    \begin{tabular}{ l r  r | r r   | r r}
        \toprule
        \multicolumn{3}{c}{ } &  \multicolumn{2}{|c|}{LMDER} & \multicolumn{2}{c}{DFOLS} \\
        \midrule
        Problem & $(n,m)$ &$\sigma_f$ & \#feval  & $\phi(x) - \phi^*$ &  \#feval   & $\phi(x) - \phi^*$ \\
        \midrule
BOX3D & $ (3,10)$ & \num{1e-01} & $44$    & \num{ 1.674e-01} & $53$    & \num{ 1.181e-01} \\
BOX3D & $ (3,10)$ & \num{1e-03} & $63$    & \num{ 1.532e-03} & $43$    & \num{ 5.394e-06} \\
BOX3D & $ (3,10)$ & \num{1e-05} & $56$    & \num{ 1.738e-09} & $46$    & \num{ 5.818e-12} \\
BOX3D & $ (3,10)$ & \num{1e-07} & $55$    & \num{ 1.758e-13} & $57$    & \num{ 2.248e-13} \\
JENSMP & $ (2,10)$ & \num{1e-01} & $46$    & \num{ 3.728e+03} & $47$    & \num{ 2.684e-01} \\
JENSMP & $ (2,10)$ & \num{1e-03} & $44$    & \num{ 4.102e-02} & $51$    & \num{ 1.660e-02} \\
\bottomrule
    \end{tabular}
    \caption{Benchmarking Problems Tested. $n$ is the number of variables and $m$ is the dimension of the residual vector.}
    \label{tab:ls_noise3}
\end{table}

\begin{table}[!htp] %NOISY CASE4
	\centering
    \footnotesize
    \begin{tabular}{ l r  r | r r   | r r}
        \toprule
        \multicolumn{3}{c}{ } &  \multicolumn{2}{|c|}{LMDER} & \multicolumn{2}{c}{DFOLS} \\
        \midrule
        Problem & $(n,m)$ &$\sigma_f$ & \#feval  & $\phi(x) - \phi^*$ &  \#feval   & $\phi(x) - \phi^*$ \\
        \midrule
BROWNAL & $ (10,10)$ & \num{1e-01} & $127$    & \num{ 1.419e+00} & $89$    & \num{ 1.655e+02} \\
BROWNAL & $ (10,10)$ & \num{1e-03} & $147$    & \num{ 5.220e-07} & $99$    & \num{ 3.284e-05} \\
BROWNAL & $ (10,10)$ & \num{1e-05} & $155$    & \num{ 4.102e-10} & $89$    & \num{ 9.707e-10} \\
BROWNAL & $ (10,10)$ & \num{1e-07} & $131$    & \num{ 1.084e-13} & $75$    & \num{ 1.308e-14} \\
OSBORNE1 & $ (5,33)$ & \num{1e-01} & $53$    & \num{ 1.279e+00} & $53$    & \num{ 5.737e+00} \\
OSBORNE1 & $ (5,33)$ & \num{1e-03} & $75$    & \num{ 3.424e-03} & $86$    & \num{ 4.512e-01} \\
\bottomrule
    \end{tabular}
    \caption{Benchmarking Problems Tested. $n$ is the number of variables and $m$ is the dimension of the residual vector.}
    \label{tab:ls_noise4}
\end{table}

\begin{table}[!htp] %NOISY CASE5
    \centering
    \footnotesize
    \begin{tabular}{ l r  r | r r   | r r}
        \toprule
        \multicolumn{3}{c}{ } &  \multicolumn{2}{|c|}{LMDER} & \multicolumn{2}{c}{DFOLS} \\
        \midrule
        Problem & $(n,m)$ &$\sigma_f$ & \#feval  & $\phi(x) - \phi^*$ &  \#feval   & $\phi(x) - \phi^*$ \\
        \midrule
MANCINO & $ (5,5)$ & \num{1e-01} & $47$    & \num{ 7.329e-02} & $48$    & \num{ 1.334e-02} \\
MANCINO & $ (5,5)$ & \num{1e-03} & $45$    & \num{ 7.341e-06} & $44$    & \num{ 5.405e-08} \\
MANCINO & $ (5,5)$ & \num{1e-05} & $39$    & \num{ 7.341e-10} & $33$    & \num{ 1.657e-12} \\
MANCINO & $ (5,5)$ & \num{1e-07} & $39$    & \num{ 7.342e-14} & $25$    & \num{ 3.807e-15} \\
MANCINO & $ (5,5)$ & \num{1e-01} & $47$    & \num{ 7.435e-02} & $58$    & \num{ 2.132e-02} \\
MANCINO & $ (5,5)$ & \num{1e-03} & $45$    & \num{ 7.511e-06} & $44$    & \num{ 1.457e-05} \\
\bottomrule
    \end{tabular}
    \caption{Benchmarking Problems Tested. $n$ is the number of variables and $m$ is the dimension of the residual vector.}
    \label{tab:ls_noise5}
\end{table}

% \newpage

\subsection{Constrained Optimization}
\label{app: cons}

\subsubsection{Test Problem Summary} 
\label{app: cons property}

%In this section, we list the properties of the functions tested in the constrained section.
\begin{table}[!htp]
\centering
\footnotesize
\begin{tabular}{l l | r r r | r r r}
\toprule
	&	&	Number of	&		&		&	Number of	&		&		\\
	&	Objective	&	Variables 	&	One-sided	&	Two-sided	&	Constraints 	&	Linear 	&  Nonlinear 	\\
Problem & Type & ($n$) & Bounds & Bounds & ($m$) & Ineq. & Ineq.\\
\midrule
BURKEHAN	&	quadratic	&	1	&	1	&	0	&	1	&	0	&	1	\\
CANTILVR	&	linear	&	5	&	5	&	0	&	1	&	0	&	1	\\
CB2	&	linear	&	3	&	0	&	0	&	3	&	0	&	3	\\
CB3	&	linear	&	3	&	0	&	0	&	3	&	0	&	3	\\
CHACONN1	&	linear	&	3	&	0	&	0	&	3	&	0	&	3	\\
\bottomrule
\end{tabular}
\caption{Properties of small-scale \texttt{CUTEst} problems.}
\label{tab:smallscaleprobs}
\end{table}
\normalsize

\begin{table}[!htp]
\centering
\footnotesize
\begin{tabular}{l l | l r r}
\toprule
Problem & SIF Parameter & Testing Name & $d$ & $m$ \\
\midrule
SVANBERG & $N=10$ & SVANBERGN10 & 10 & 10\\
SVANBERG & $N=50$ & SVANBERGN50 & 50 & 50\\
SVANBERG & $N=100$ & SVANBERGN100 & 100 & 100\\
SVANBERG & $N=500$ & SVANBERGN500 & 500 & 500\\
\midrule
READING4 & $N=2$ & READING4N2 & 2 & 2\\
READING4 & $N=50$ & READING4N50 & 50 & 50\\
READING4 & $N=100$ & READING4N100 & 100 & 100\\
READING4 & $N=500$ & READING4N500& 500 & 500\\
\midrule
COSHFUN & $M=3$ & COSHFUNM3 & 10 & 3\\
COSHFUN & $M=8$ & COSHFUNM8 & 25 & 8\\
COSHFUN & $M=14$ & COSHFUNM14 & 43 & 14\\
COSHFUN & $M=20$ & COSHFUNM20 & 61 & 20\\
\bottomrule
\end{tabular}
\caption{Instances of variable-size \texttt{CUTEst} problems.}
\label{tab:varsizeprobs}
\end{table}
\normalsize

%\begin{center}
%\begin{figure}[h]
%\label{fig:progress}
%\caption{Progress of algorithms on problem HS92}
%\centering
%%\includegraphics[scale=0.65]{figures/constrained/plot_HS104}
%\includegraphics[scale=0.65]{figures/constrained/plot_HS92}
%%\includegraphics[scale=0.65]{figures/constrained/plot_HALDMADS}
%\end{figure}
%\end{center}

% \newpage

\subsubsection{Noiseless Functions}
\label{app: cons det}

In this section, we list the final objective value($\phi(x)$), number of function evaluations (\#feval), CPU time (CPU), and feasibility error(feaserr) for each problem instance in Tables \ref{table:detA}-\ref{table:detAtol}.

\begin{table}[!htp]
\centering
\footnotesize
\begin{tabular}{r | l l l l | l l l l}
\toprule
\multicolumn{1}{c|}{ } &  \multicolumn{4}{c|}{KNITRO} & \multicolumn{4}{c}{COBYLA}\\
\midrule
Problem	&	$\phi(x)$	&	\#feval	&	CPU	&	feaserr	&	$\phi(x)$	&	\#fevals	&	CPU	&	feaserr	\\
\midrule
CB3	&	2.0000	&	40	&	0.023	&	2.22e-16	&	2.0000	&	44	&	0.157	&	0.00e+00	\\
CHACONN2	&	2.0000	&	40	&	0.021	&	0.00e+00	&	2.0000	&	48	&	0.156	&	0.00e+00	\\
GIGOMEZ3	&	2.0000	&	40	&	0.022	&	0.00e+00	&	2.0000	&	45	&	0.160	&	0.00e+00	\\
HS100	&	680.6301	&	360	&	0.160	&	0.00e+00	&	680.6300	&	454	&	0.272	&	0.00e+00	\\
DIPIGRI	&	680.6301	&	389	&	0.152	&	2.84e-14	&	680.6300	&	458	&	0.288	&	0.00e+00	\\
HS93	&	135.0760	&	2247	&	1.040	&	4.44e-16	&	135.0760	&	3000*	&	0.873	&	0.00e+00	\\
HS64	&	6299.8424	&	144	&	0.093	&	0.00e+00	&	6299.8400	&	409	&	0.244	&	0.00e+00	\\
POLAK3	&	5.9330	&	902	&	0.414	&	3.89e-16	&	5.9330	&	3139	&	0.071	&	0.00e+00	\\
POLAK1	&	2.7183	&	60	&	0.034	&	0.00e+00	&	2.7183	&	326	&	0.232	&	0.00e+00	\\
HS104	&	3.9512	&	718	&	0.282	&	0.00e+00	&	3.9512	&	3977	&	0.131	&	0.00e+00	\\
HS100MOD	&	678.6796	&	477	&	0.207	&	0.00e+00	&	678.6790	&	385	&	0.240	&	0.00e+00	\\
TWOBARS	&	1.5087	&	63	&	0.037	&	0.00e+00	&	1.5087	&	83	&	0.172	&	0.00e+00	\\
HS85	&	-2.2156	&	117	&	0.078	&	7.11e-15	&	-2.2156	&	299	&	0.242	&	0.00e+00	\\
CB2	&	1.9522	&	72	&	0.041	&	1.11e-16	&	1.9522	&	111	&	0.173	&	0.00e+00	\\
CHACONN1	&	1.9522	&	55	&	0.030	&	0.00e+00	&	1.9522	&	120	&	0.175	&	0.00e+00	\\
MADSEN	&	0.6164	&	73	&	0.042	&	0.00e+00	&	0.6164	&	112	&	0.175	&	0.00e+00	\\
HS66	&	0.5182	&	43	&	0.026	&	4.44e-16	&	0.5182	&	101	&	0.170	&	0.00e+00	\\
MINMAXBD	&	115.7064	&	420	&	0.335	&	6.17e-14	&	115.7060	&	924	&	0.395	&	0.00e+00	\\
CANTILVR	&	1.3400	&	170	&	0.078	&	3.61e-16	&	1.3400	&	281	&	0.227	&	0.00e+00	\\
HS92	&	1.3627	&	976	&	0.988	&	4.60e-17	&	1.3627	&	3000*	&	0.993	&	0.00e+00	\\
HALDMADS	&	0.0001	&	113	&	0.073	&	4.44e-16	&	0.0001	&	669	&	0.341	&	0.00e+00	\\
HS34	&	-0.8340	&	59	&	0.033	&	4.44e-16	&	-0.8340	&	47	&	0.149	&	0.00e+00	\\
HS90	&	1.3627	&	570	&	0.475	&	0.00e+00	&	1.3629	&	2000*	&	0.066	&	0.00e+00	\\
SNAKE	&	0.0000	&	25	&	0.029	&	1.02e-16	&	0.0000	&	72	&	0.165	&	0.00e+00	\\
HS72	&	727.6794	&	108	&	0.103	&	1.73e-18	&	727.6794	&	1002	&	0.424	&	4.77e-18	\\
GIGOMEZ2	&	1.9522	&	79	&	0.048	&	2.22e-16	&	1.9522	&	111	&	0.169	&	1.11e-16	\\
SYNTHES1	&	0.7593	&	96	&	0.041	&	0.00e+00	&	0.7593	&	174	&	0.186	&	1.54e-23	\\
HS88	&	1.3627	&	209	&	0.156	&	0.00e+00	&	1.3627	&	167	&	0.203	&	2.09e-17	\\
HS13	&	0.9999	&	140	&	0.155	&	2.61e-13	&	1.0000	&	86	&	0.168	&	1.07e-27	\\
MATRIX2	&	0.0000	&	275	&	0.131	&	3.60e-17	&	0.0000	&	216	&	0.199	&	5.20e-18	\\
WOMFLET	&	0.0000	&	103	&	0.070	&	5.89e-10	&	0.0000	&	117	&	0.181	&	8.99e-18	\\
S365	&	0.0000	&	63	&	0.032	&	1.72e-07	&	0.0000	&	210	&	0.215	&	4.35e-17	\\
\midrule
%MGH17SLS	&	1.0224	&	1053	&	0.388	&	0.00E+00	&	1.1037	&	2500*	&	0.237	&	0.00E+00	\\
%DMN15102LS	&	11456.7766	&	33000*	&	112.470	&	0.00E+00	&	10303.7000	&	33000*	&	29.662	&	0.00E+00	\\
HS67	&	-1162.1187	&	285	&	0.198	&	0.00e+00	&	-980.1600	&	7000*	&	0.885	&	0.00e+00	\\
CRESC50	&	0.5952	&	16704	&	17.912	&	0.00e+00	&	7.9917	&	50000*	&	4.863	&	0.00e+00	\\
CRESC4	&	0.8719	&	688	&	0.363	&	0.00e+00	&	44.5901	&	4000*	&	0.318	&	0.00e+00	\\
\midrule
HS91	&	1.3627	&	1229	&	1.057	&	9.96e-18	&	1.1199	&	2500*	&	0.583	&	6.50e-04	\\
DEMBO7	&	174.7870	&	3225	&	1.594	&	3.61e-16	&	250.0720	&	10000*	&	0.238	&	7.29e-02	\\
POLAK5	&	50.0000	&	65	&	0.043	&	0.00e+00	&	34.4929	&	1500*	&	0.571	&	1.55e+01	\\
HS101	&	1809.7691	&	3500*	&	1.908	&	1.38e-13	&	3000.1900	&	3500*	&	0.077	&	1.99e-01	\\
HS89	&	1.3627	&	315	&	0.257	&	2.17e-17	&	0.6520	&	1500*	&	0.742	&	5.12e-02	\\
HS102	&	911.8816	&	3500*	&	1.965	&	4.67e-12	&	3000.3500	&	3500*	&	0.096	&	3.60e-01	\\
HS103	&	543.6680	&	3500*	&	1.858	&	2.57e-16	&	3000.1700	&	3500*	&	0.091	&	1.76e-01	\\
\midrule
SPIRAL	&	0.0195	&	1500*	&	0.796	&	1.34e-05	&	0.0958	&	1500*	&	0.508	&	0.00e+00	\\
POLAK6	&	-78.6745	&	474	&	1.311	&	8.55e+01	&	0.0000	&	2500*	&	0.765	&	0.00e+00	\\
\midrule
S365MOD	&	0.2500	&	162	&	0.204	&	1.25e+00	&	0.0301	&	247	&	0.211	&	3.26e-01	\\
BURKEHAN	&	10.0000	&	13	&	0.007	&	1.01e+02	&	0.0000	&	54	&	0.175	&	1.00e+00	\\
POLAK2	&	29.9570	&	2641	&	5.119	&	5.39e+01	&	-259.2500	&	5500*	&	0.794	&	3.14e+02	\\
\bottomrule
\end{tabular}
\caption{Summary of the results for small-scale noiseless constrained \texttt{CUTEst} problems. The horizontal bars divide cases (i), (ii), (iii), and (iv).}
\label{table:detA}
\normalsize
\end{table}

We now comment on outcomes (ii)-(iv) mentioned in Section~\ref{cons-noiseless}.

\textit{(ii) The solvers converged to feasible points with different objective function values.}  There are three such problems, for all of which COBYLA  terminated due to the limit on the number of function evaluations.  We removed the limit on the number of evaluations for COBYLA to see if it converges to a different local solution.  For HS67, COBYLA terminates with a final objective function value of $-1116.415$ after $32606$ evaluations, and for CRESC4 it converges to a feasible solution with $f=1.03523$ after $4706634$ evaluations. 

\smallskip
\textit{(iii) One of the solvers terminated at an infeasible point.}  There are seven problems for which {\sc knitro} terminated at a feasible point but {\sc cobyla} at an infeasible one.  For all of these problems, {\sc cobyla} hits the limit on the number of function evaluations; it can make further progress in feasibility if the budget of evaluations is increased.  For problems HS101, HS102, and HS103, {\sc knitro} also hits the evaluation limit; however, it gets a better solution --a feasible one with a lower objective value-- for all three problems.

For two problems, {\sc cobyla} terminated at a feasible point but {\sc knitro} at an infeasible one.  For problem POLAK6, {\sc knitro} stalls at an infeasible point.  For SPIRAL, it gets close to feasibility but runs out its evaluation budget before it can reduce the feasibility error further.

\smallskip
\textit{(iv) Both solvers terminated at infeasible points.}   For all three of those problems, {\sc knitro} declares convergence to infeasible stationary points.  For problems S365 and BURKEHAN, {\sc cobyla} terminates due to the condition on the final trust region radius whereas for POLAK2 it stops due to the function evaluation limit.

% \newpage

\begin{landscape}
\begin{table}[!htp]
\tiny
\begin{tabular}{r| r| r|  r l | r l| r| r l | r l| r| r l | r l}
\toprule
\multicolumn{2}{c|}{ } & \multicolumn{5}{c|}{$TOL=10^{-6}$} & \multicolumn{5}{c|}{$TOL=10^{-3}$} & \multicolumn{5}{c}{$TOL=10^{-1}$}\\
\midrule
\multicolumn{2}{c|}{ } & & \multicolumn{2}{|c|}{KNITRO} &  \multicolumn{2}{|c|}{COBYLA} &  & \multicolumn{2}{|c|}{KNITRO} &  \multicolumn{2}{|c|}{COBYLA} & & \multicolumn{2}{|c|}{KNITRO} & \multicolumn{2}{c}{COBYLA}\\
\midrule
problem & $\phi^\ast$ & target $\phi$ & $\phi(x)$ & \#fevals & $\phi(x)$ & \#fevals & target $\phi$ & $\phi(x)$ & \#fevals & $\phi(x)$ & \#fevals & target $\phi$ & $\phi(x)$ & \#fevals & $\phi(x)$ & \#fevals \\
\midrule
CB3	&	2.000000	&	2.00000	&	2.00000	&	30	&	2.00000	&	25	&	2.00200	&	2.00000	&	30	&	2.00000	&	25	&	2.20000	&	2.00000	&	30	&	2.00351	&	18	\\
CHACONN2	&	2.000000	&	2.00000	&	2.00000	&	30	&	2.00000	&	32	&	2.00200	&	2.00000	&	30	&	2.00021	&	24	&	2.20000	&	2.00000	&	30	&	2.06187	&	19	\\
GIGOMEZ3	&	2.000000	&	2.00000	&	2.00000	&	30	&	2.00000	&	39	&	2.00200	&	2.00000	&	30	&	2.00000	&	39	&	2.20000	&	2.00000	&	30	&	2.00000	&	39	\\
HS100	&	680.630057	&	680.63074	&	680.63013	&	124	&	680.63000	&	247	&	681.31069	&	680.75511	&	75	&	680.63000	&	247	&	748.69306	&	714.00000	&	1	&	714.00000	&	1	\\
DIPIGRI	&	680.630057	&	680.63074	&	680.63013	&	124	&	680.63000	&	410	&	681.31069	&	680.75511	&	75	&	680.63100	&	141	&	748.69306	&	714.00000	&	1	&	714.00000	&	1	\\
HS93	&	135.075964	&	135.07610	&	135.07610	&	871	&	135.07600	&	2416	&	135.21104	&	135.14058	&	229	&	135.21000	&	270	&	148.58356	&	137.06640	&	1	&	137.06600	&	1	\\
HS64	&	6299.842428	&	6299.84873	&	6299.73981	&	64	&	6299.84000	&	335	&	6306.14227	&	6302.48265	&	54	&	6305.80000	&	316	&	6929.82667	&	6687.88225	&	30	&	6684.54000	&	239	\\
POLAK3	&	5.933003	&	5.93301	&	5.93301	&	469	&	5.93300	&	2637	&	5.93894	&	5.93353	&	349	&	5.93644	&	2350	&	6.52630	&	5.93353	&	349	&	5.94016	&	2328	\\
POLAK1	&	2.718282	&	2.71828	&	2.71828	&	50	&	2.71828	&	265	&	2.72100	&	2.71828	&	50	&	2.71834	&	253	&	2.99011	&	2.71828	&	50	&	2.83811	&	250	\\
HS104	&	3.951163	&	3.95117	&	3.95117	&	290	&	3.95116	&	1764	&	3.95511	&	3.95148	&	231	&	3.95381	&	730	&	4.34628	&	3.95946	&	193	&	4.15766	&	173	\\
HS100MOD	&	678.679638	&	678.68032	&	678.67969	&	173	&	678.68000	&	169	&	679.35832	&	679.08174	&	63	&	679.18500	&	50	&	746.54760	&	714.00000	&	1	&	714.00000	&	1	\\
TWOBARS	&	1.508652	&	1.50865	&	1.50865	&	50	&	1.50865	&	33	&	1.51016	&	1.50865	&	50	&	1.50892	&	21	&	1.65952	&	1.50865	&	50	&	1.65447	&	10	\\
HS85	&	-2.215605	&	-2.21560	&	-2.21560	&	110	&	-2.21560	&	299	&	-2.21339	&	-2.21560	&	110	&	-2.21500	&	247	&	-1.99404	&	-2.21560	&	110	&	-2.02180	&	158	\\
CB2	&	1.952224	&	1.95223	&	1.95222	&	52	&	1.95222	&	57	&	1.95418	&	1.95222	&	52	&	1.95230	&	41	&	2.14745	&	1.95222	&	52	&	1.95230	&	41	\\
CHACONN1	&	1.952224	&	1.95223	&	1.95222	&	40	&	1.95222	&	55	&	1.95418	&	1.95222	&	40	&	1.95222	&	55	&	2.14745	&	1.95222	&	40	&	2.00282	&	13	\\
MADSEN	&	0.616432	&	0.61643	&	0.61643	&	58	&	0.61643	&	65	&	0.61743	&	0.61643	&	58	&	0.61643	&	65	&	0.71643	&	0.61643	&	58	&	0.65433	&	22	\\
HS66	&	0.518163	&	0.51816	&	0.51816	&	33	&	0.51816	&	37	&	0.51916	&	0.51866	&	21	&	0.51817	&	28	&	0.61816	&	0.58000	&	1	&	0.58000	&	1	\\
MINMAXBD	&	115.706440	&	115.70656	&	115.70644	&	343	&	115.70600	&	800	&	115.82215	&	115.70644	&	343	&	115.72100	&	766	&	127.27708	&	118.74188	&	246	&	115.72100	&	766	\\
CANTILVR	&	1.339956	&	1.33996	&	1.33994	&	107	&	1.33995	&	189	&	1.34130	&	1.33994	&	107	&	1.33997	&	130	&	1.47395	&	1.33994	&	107	&	1.34758	&	107	\\
HS92	&	1.362657	&	1.36266	&	1.36266	&	844	&	1.36265	&	2121	&	1.36402	&	1.36327	&	405	&	1.36310	&	203	&	1.49892	&	1.36327	&	405	&	1.36310	&	203	\\
HALDMADS	&	0.000122	&	0.00012	&	0.00012	&	97	&	0.00013	&	669	&	0.00112	&	0.00012	&	97	&	0.00057	&	51	&	0.10012	&	0.00012	&	97	&	0.00057	&	51	\\
HS34	&	-0.834032	&	-0.83403	&	-0.83403	&	54	&	-0.83400	&	40	&	-0.83303	&	-0.83403	&	54	&	-0.83400	&	25	&	-0.73403	&	-0.83186	&	44	&	-0.83400	&	25	\\
HS90	&	1.362657	&	1.36266	&	1.36185	&	259	&	1.36287	&	2000	&	1.36402	&	1.36185	&	259	&	1.36311	&	159	&	1.49892	&	1.36185	&	259	&	1.36311	&	159	\\
SNAKE	&	0.000000	&	0.00000	&	0.00000	&	17	&	0.00000	&	67	&	0.00100	&	0.00000	&	17	&	0.00009	&	61	&	0.10000	&	0.00000	&	17	&	0.00009	&	61	\\
HS72	&	727.679358	&	727.68009	&	727.66249	&	72	&	727.67973	&	877	&	728.40704	&	727.66249	&	72	&	728.36304	&	856	&	800.44729	&	727.66249	&	72	&	738.77714	&	832	\\
GIGOMEZ2	&	1.952224	&	1.95223	&	1.95222	&	64	&	1.95223	&	47	&	1.95418	&	1.95222	&	64	&	1.95266	&	32	&	2.14745	&	1.95222	&	64	&	2.08642	&	13	\\
SYNTHES1	&	0.759284	&	0.75929	&	0.75928	&	72	&	0.75928	&	79	&	0.76028	&	0.75928	&	72	&	0.75928	&	79	&	0.85928	&	0.75928	&	72	&	0.79837	&	17	\\
HS88	&	1.362657	&	1.36266	&	1.36249	&	189	&	1.36266	&	138	&	1.36402	&	1.36249	&	189	&	1.36387	&	123	&	1.49892	&	1.36249	&	189	&	1.36387	&	123	\\
HS13	&	0.999872	&	0.99987	&	0.99987	&	116	&	1.00000	&	86	&	1.00087	&	1.00060	&	88	&	1.00085	&	41	&	1.09987	&	1.07956	&	40	&	1.07745	&	22	\\
MATRIX2	&	0.000000	&	0.00000	&	0.00000	&	131	&	0.00000	&	187	&	0.00100	&	0.00001	&	107	&	0.00000	&	187	&	0.10000	&	0.00001	&	107	&	0.09622	&	46	\\
WOMFLET	&	0.000000	&	0.00000	&	0.00000	&	78	&	0.00000	&	104	&	0.00100	&	0.00000	&	78	&	0.00000	&	104	&	0.10000	&	0.00000	&	78	&	0.00000	&	104	\\
S365	&	0.000000	&	0.00000	&	0.00000	&	63	&	0.00000	&	179	&	0.00100	&	0.00000	&	63	&	0.00000	&	179	&	0.10000	&	0.00000	&	63	&	0.00000	&	179	\\
\bottomrule
\end{tabular}
\caption{Function evaluations to obtain $\phi_k \leq \phi^\ast + TOL \cdot \max\{1.0,|\phi^\ast|\}$ for small scale noiseless constrained problems.}
\label{table:detAtol}
\end{table}
\end{landscape}
%We next question whether the larger performance gap for increased TOL values is related to the use of a quadratic term in the {\sc knitro} algorithm (although this is a very weak quadratic term since it is constructed by the LBFGS formulae with memory size 1).  To assess the potential help of the quadratic term, we increase the memory size for the LBFGS matrix to 10 and observe the resulting improvement in the performance of {\sc knitro} (again, in terms of function evaluations) as reflected by Figure~\ref{fig:efficiency_mem}.  Interestingly, it looks like the additional memory helps at early iterations of the algorithm more than it does in improving the solution precision. 

%\includegraphics[scale=0.31]{figures/constrained/profile3_TOL=1e-6}
%\includegraphics[scale=0.31]{figures/constrained/profile3_TOL=1e-3}
%\includegraphics[scale=0.31]{figures/constrained/profile3_TOL=1e-1}

\begin{figure}
	\includegraphics[width=0.32\textwidth]{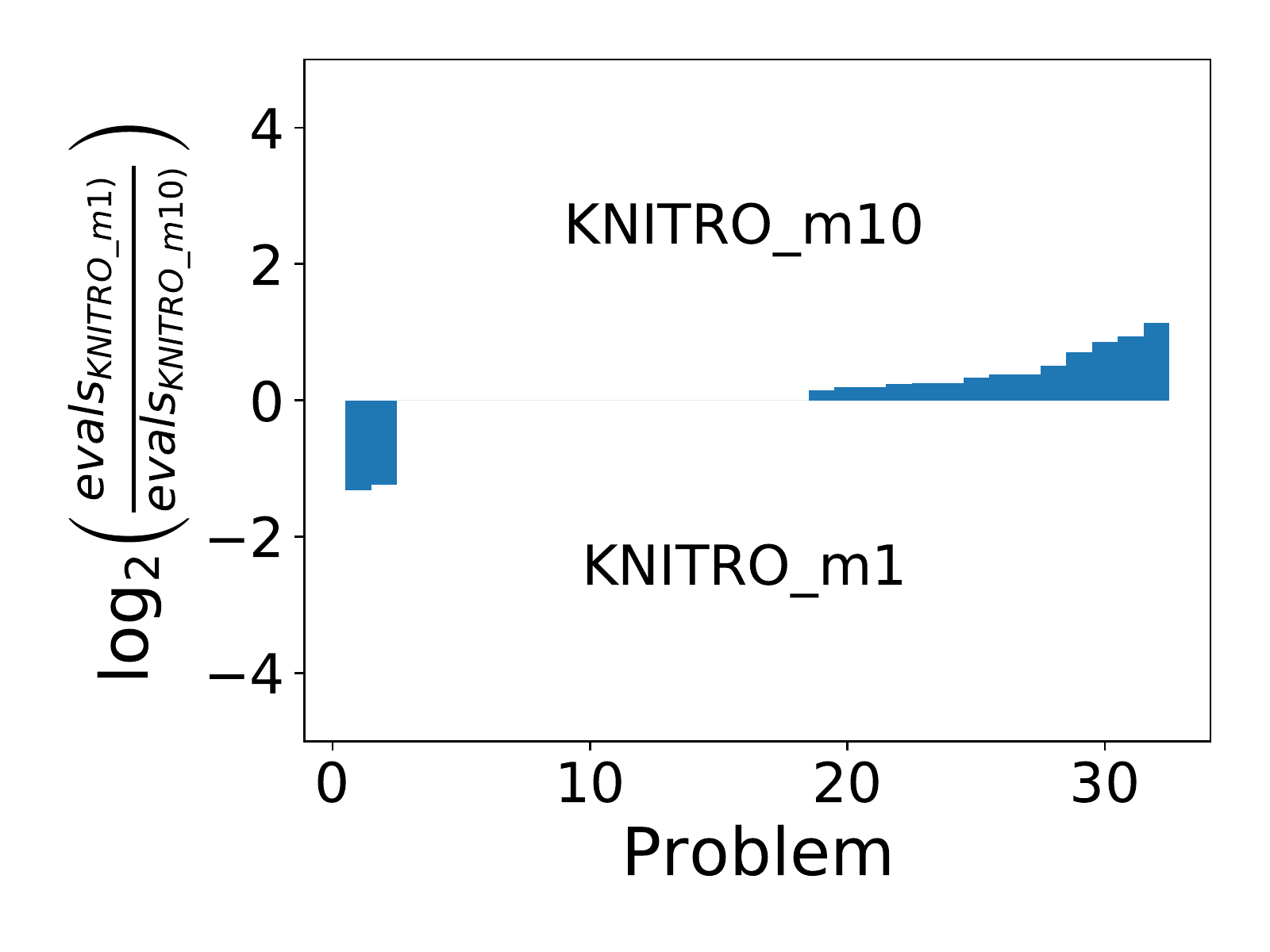}
	\includegraphics[width=0.32\textwidth]{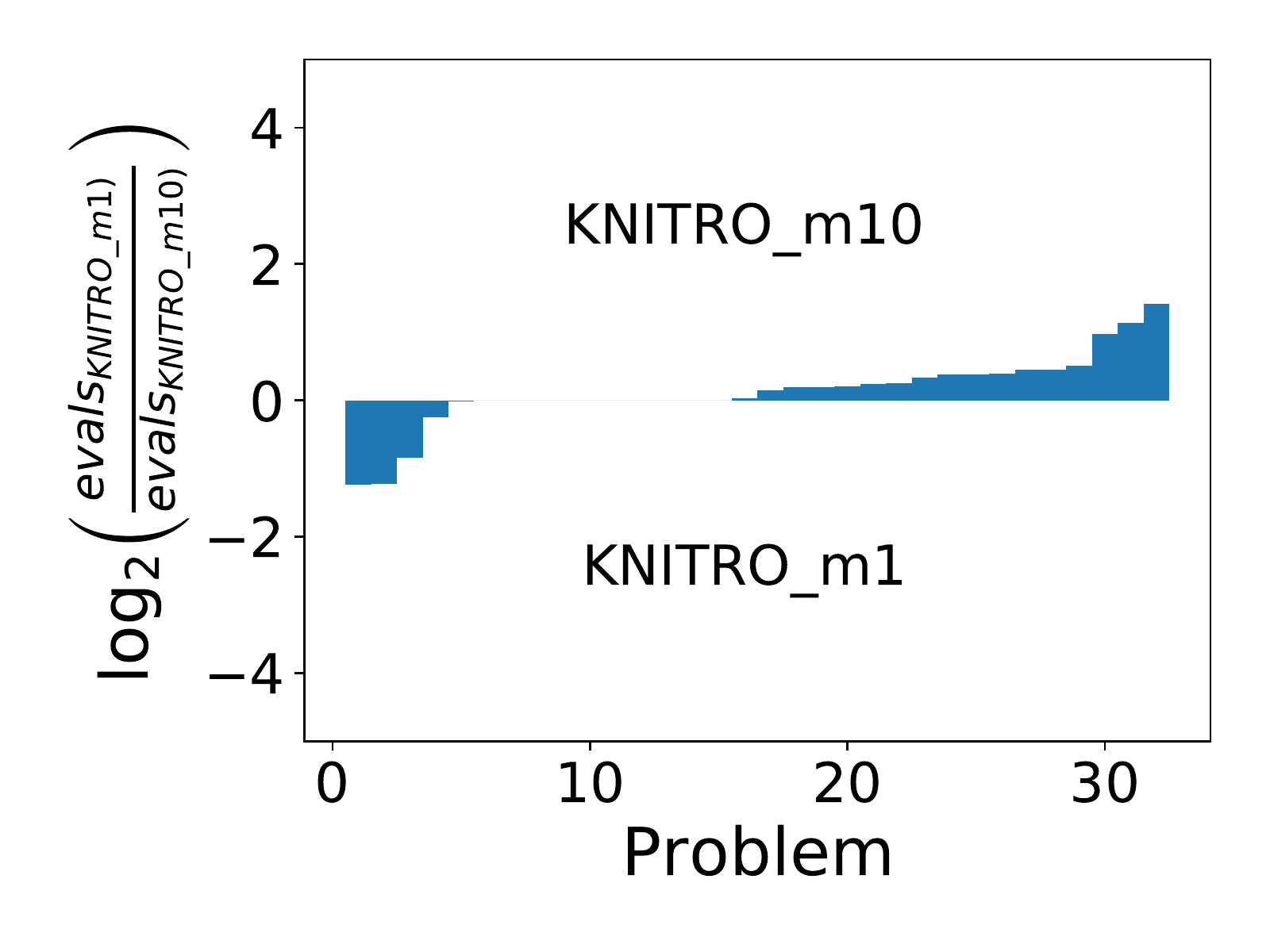}	
	\includegraphics[width=0.32\textwidth]{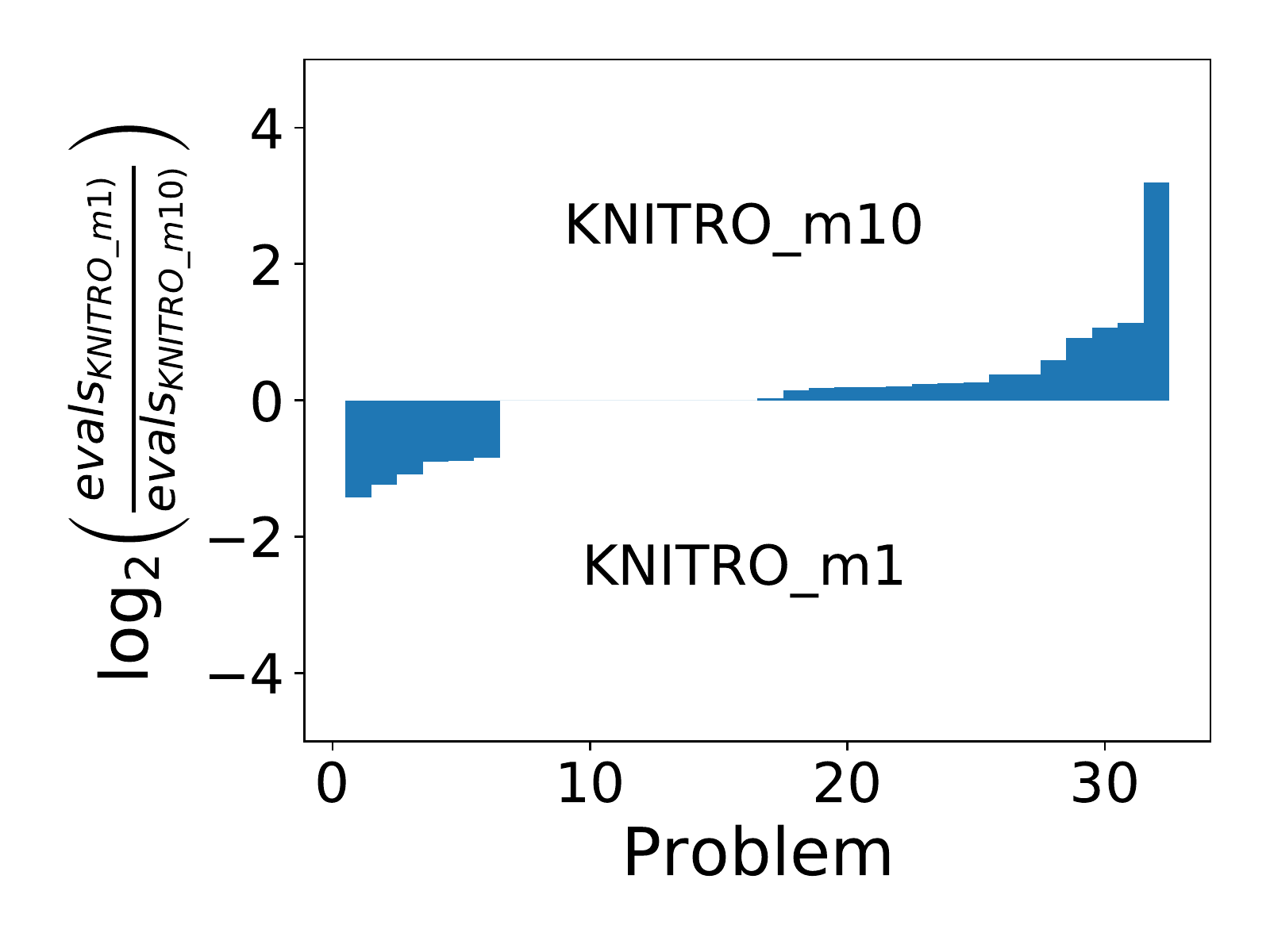}
    \caption{{\em Efficiency, Noiseless Case.} Log-ratio profiles comparing {\sc knitro} with an {\sc l-bfgs} Hessian approximation of memory one and memory 10 when $\epsilon(x) = 0$. The figures measure number of function evaluations to satisfy \eqref{eq:term} for $\tau = 10^{-1}$ (left), $10^{-3}$ (middle), $10^{-6}$ (right). } 
    \label{fig:efficiency_mem}
\end{figure}

% \newpage
\subsubsection{Noisy Functions}
\label{app: cons noisy}

In this section, we list the final objective value($\phi(x)$), feasibility error(feaserr) and the number of function evaluations (\#feval) needed to achieve this for each problem instance, varying the noise level $\sigma_f \in \{10^{-1}, 10^{-3}, 10^{-5}, 10^{-7}\}$, in Tables \ref{table:noisyobj}-\ref{table:noisyfeas5}. Function evaluations marked with a $^*$ denote cases where the algorithm reached the maximum number of function evaluations. An `f' marks the cases where the feasibility error is large compared to the noise level, i.e., $\Vert \max \{ \psi(x),0  \} \Vert_\infty \geq \sqrt{3 \sigma_f}$.

%\begin{landscape}
\begin{table}[!htp]
\tiny
\begin{tabular}{r |l l l l l | l l l l l}
\toprule
 & \multicolumn{5}{c}{KNITRO} & \multicolumn{5}{c}{COBYLA} \\
\midrule 
$\sigma_f$ &	$0.0$ &	1e-07	&	1e-05	&	1e-03	&	1e-01	&	$0.0$ &	1.00e-07	&	1.00e-05	&	1.00e-03	&	1.00e-01	\\
\midrule
CB3	&	4.92e-12	&	1.00e-08	&	2.00e-06	&	1.35e-04	&	1.51e-01	&	8.36e-12	&	1.00e-08	&	6.00e-06	&	1.10e-03	&	2.97e-01	\\
CHACONN2	&	1.72e-12	&	1.00e-08	&	1.60e-05	&	3.80e-04	&	8.27e-02	&	1.72e-12	&	1.00e-08	&	3.00e-06	&	1.01e-03	&	f	\\
GIGOMEZ3	&	9.80e-13	&	1.00e-08	&	2.00e-06	&	1.35e-04	&	1.62e-01	&	9.80e-13	&	1.00e-08	&	1.00e-05	&	9.30e-05	&	f	\\
HS100	&	2.05e-12	&	2.37e-06	&	2.14e-04	&	1.34e-02	&	2.16e+00	&	5.87e-07	&	1.63e-06	&	9.62e-04	&	3.13e-02	&	3.65e-02	\\
DIPIGRI	&	4.09e-12	&	3.24e-05	&	1.87e-04	&	1.25e-02	&	1.98e+00	&	2.05e-12	&	1.63e-06	&	4.23e-04	&	3.21e-02	&	9.14e-01	\\
HS93	&	3.90e-09	&	4.34e-03	&	6.21e-02	&	2.49e-01	&	2.69e+00	&	3.96e-05	&	1.52e-01	&	4.88e-01	&	1.79e+00	&	1.99e+00	\\
\bottomrule
\end{tabular}
\caption{Optimality gap $|\phi(x_k)-\phi^\ast|$ for small scale noisy constrained \texttt{CUTEst} problems.}
\label{table:noisyobj}
\end{table}
%\end{landscape}

% \newpage

\begin{table}[!htp]
\tiny
\begin{tabular}{r |l l l l l | l l l l l}
\toprule
 & \multicolumn{5}{c}{KNITRO} & \multicolumn{5}{c}{COBYLA} \\
\midrule 
	$\sigma_f$ & $0.0$ &	1e-07	&	1e-05	&	1e-03	&	1e-01	&	$0.0$ &	1.00e-07	&	1.00e-05	&	1.00e-03	&	1.00e-01	\\
\midrule
CB3	&	2.22e-16	&	0.00e+00	&	0.00e+00	&	7.50e-05	&	0.00e+00	&	0.00e+00	&	0.00e+00	&	1.20e-05	&	0.00e+00	&	3.71e-02	\\
CHACONN2	&	0.00e+00	&	0.00e+00	&	0.00e+00	&	1.22e-04	&	2.54e-01	&	0.00e+00	&	0.00e+00	&	5.00e-06	&	0.00e+00	&	1.53e+00	\\
GIGOMEZ3	&	0.00e+00	&	0.00e+00	&	0.00e+00	&	7.50e-05	&	0.00e+00	&	0.00e+00	&	0.00e+00	&	0.00e+00	&	3.32e-04	&	5.89e-01	\\
HS100	&	0.00e+00	&	5.00e-06	&	2.10e-05	&	6.18e-04	&	1.01e-01	&	0.00e+00	&	0.00e+00	&	0.00e+00	&	0.00e+00	&	1.67e-01	\\
DIPIGRI	&	2.84e-14	&	3.10e-05	&	3.00e-06	&	5.53e-04	&	1.38e-01	&	0.00e+00	&	0.00e+00	&	0.00e+00	&	0.00e+00	&	0.00e+00	\\
HS93	&	4.44e-16	&	0.00e+00	&	1.00e-05	&	1.03e-03	&	1.13e-01	&	0.00e+00	&	0.00e+00	&	0.00e+00	&	0.00e+00	&	0.00e+00	\\
\bottomrule
\end{tabular}
\caption{Feasibility error $\|\max\{\psi(x_k),0\}\|_\infty$ for small scale noisy constrained \texttt{CUTEst} problems.}
\label{table:noisyfeas1}
\end{table}

% \newpage

%\begin{landscape}
\begin{table}[!htp]
\tiny
\begin{tabular}{r |l l l| l l l | l l l | l l l}
\toprule
& \multicolumn{6}{|c|}{$\tau=10^{-6}$} & \multicolumn{6}{|c}{$\tau=10^{-2}$}\\
\midrule
& \multicolumn{3}{|c|}{KNITRO} & \multicolumn{3}{c|}{COBYLA} & \multicolumn{3}{|c|}{KNITRO} & \multicolumn{3}{c}{COBYLA} \\
problem	&	$\phi(x)$	&	feaserr	&	\#fevals	&	$\phi(x)$	&	feaserr	&	\#fevals	&		$\phi(x)$	&	feaserr	&	\#fevals	&	$\phi(x)$	&	feaserr	&	\#fevals	\\
\midrule
CB3	&	2.000000	&	0.00e+00	&	35	&	2.000000	&	0.00e+00	&	30	&	1.9999980	&	6.00e-06	&	30	&	2.000030	&	1.20e-05	&	20	\\
CHACONN2	&	1.999998	&	6.00e-06	&	30	&	2.000000	&	0.00e+00	&	29	&	1.9999980	&	6.00e-06	&	30	&	1.999631	&	5.01e-04	&	22	\\
GIGOMEZ3	&	2.000000	&	0.00e+00	&	35	&	2.000000	&	0.00e+00	&	24	&	2.0000000	&	0.00e+00	&	35	&	2.000000	&	0.00e+00	&	24	\\
HS100	&	680.630087	&	8.00e-06	&	173	&	680.630055	&	1.15e-04	&	125	&	680.7482460	&	0.00e+00	&	75	&	680.707367	&	0.00e+00	&	87	\\
DIPIGRI	&	680.630025	&	3.10e-05	&	189	&	680.630055	&	2.80e-05	&	205	&	680.7513980	&	0.00e+00	&	75	&	680.638226	&	2.42e-04	&	114	\\
\bottomrule
\end{tabular}
\caption{Number of function evaluations to obtain \eqref{eq:term_eps} for noise level $\sigma_f = 10^{-7}$.}
\label{table:noisyfeas2}
\end{table}
%\end{landscape}

% \newpage

%\begin{landscape}
\begin{table}[!htp]
\tiny
\begin{tabular}{r |l l l| l l l | l l l | l l l}
\toprule
& \multicolumn{6}{|c|}{$\tau=10^{-6}$} & \multicolumn{6}{|c}{$\tau=10^{-2}$}\\
\midrule
& \multicolumn{3}{|c|}{KNITRO} & \multicolumn{3}{c|}{COBYLA} & \multicolumn{3}{|c|}{KNITRO} & \multicolumn{3}{c}{COBYLA} \\
problem	&	$\phi(x)$	&	feaserr	&	\#fevals	&	$\phi(x)$	&	feaserr	&	\#fevals	&		$\phi(x)$	&	feaserr	&	\#fevals	&	$\phi(x)$	&	feaserr	&	\#fevals	\\
\midrule
CB3	&	2.000003	&	0.00e+00	&	47	&	2.000000	&	8.00e-06	&	32	&	1.999981	&	2.50e-05	&	30	&	2.000031	&	1.10e-05	&	20	\\
CHACONN2	&	2.000006	&	4.00e-06	&	45	&	2.000001	&	2.79e-04	&	31	&	1.999999	&	1.90e-05	&	30	&	1.999624	&	5.08e-04	&	22	\\
GIGOMEZ3	&	2.000003	&	0.00e+00	&	47	&	2.000003	&	1.26e-03	&	25	&	2.000003	&	0.00e+00	&	47	&	2.000003	&	1.26e-03	&	25	\\
HS100	&	680.630284	&	9.00e-06	&	160	&	680.630275	&	1.47e-04	&	152	&	680.849165	&	0.00e+00	&	46	&	680.757842	&	1.82e-03	&	56	\\
DIPIGRI	&	680.630221	&	7.00e-06	&	149	&	680.630243	&	5.00e-05	&	164	&	680.817892	&	1.54e-03	&	56	&	680.718994	&	3.36e-04	&	58	\\
\bottomrule
\end{tabular}
\caption{Number of function evaluations to obtain \eqref{eq:term_eps} for noise level $\sigma_f = 10^{-5}$}
\label{table:noisyfeas3}
\end{table}
%\end{landscape}

% \newpage

%\begin{landscape}
\begin{table}[!htp]
\tiny
\begin{tabular}{r |l l l| l l l | l l l | l l l}
\toprule
& \multicolumn{6}{|c|}{$\tau=10^{-6}$} & \multicolumn{6}{|c}{$\tau=10^{-2}$}\\
\midrule
& \multicolumn{3}{|c|}{KNITRO} & \multicolumn{3}{c|}{COBYLA} & \multicolumn{3}{|c|}{KNITRO} & \multicolumn{3}{c}{COBYLA} \\
problem	&	$\phi(x)$	&	feaserr	&	\#fevals	&	$\phi(x)$	&	feaserr	&	\#fevals	&		$\phi(x)$	&	feaserr	&	\#fevals	&	$\phi(x)$	&	feaserr	&	\#fevals	\\
\midrule
CB3	&	2.000135	&	7.50e-05	&	57	&	2.001105	&	0.00e+00	&	f	&	1.991529	&	1.86e-02	&	25	&	1.999922	&	9.37e-04	&	17	\\
CHACONN2	&	2.000382	&	1.20e-04	&	30	&	2.000360	&	1.84e-03	&	f	&	1.989816	&	2.59e-02	&	25	&	2.007665	&	3.71e-02	&	19	\\
GIGOMEZ3	&	2.000135	&	7.50e-05	&	f	&	1.999908	&	3.32e-04	&	35	&	2.000135	&	7.50e-05	&	f	&	1.999908	&	3.32e-04	&	35	\\
HS100	&	680.643454	&	6.18e-04	&	118	&	680.653293	&	2.80e-03	&	f	&	680.862278	&	2.56e-03	&	47	&	680.751163	&	1.19e-02	&	46	\\
DIPIGRI	&	680.642602	&	5.53e-04	&	146	&	680.653238	&	2.48e-03	&	f	&	680.858601	&	0.00e+00	&	47	&	680.938844	&	1.37e-02	&	46	\\
\bottomrule
\end{tabular}
\caption{Number of function evaluations to obtain \eqref{eq:term_eps} for noise level $\sigma_f = 10^{-3}$.}
\label{table:noisyfeas4}
\end{table}
%\end{landscape}

% \newpage

%\begin{landscape}
\begin{table}[!htp]
\tiny
\begin{tabular}{r |l l l| l l l | l l l | l l l}
\toprule
& \multicolumn{6}{|c|}{$\tau=10^{-6}$} & \multicolumn{6}{|c}{$\tau=10^{-2}$}\\
\midrule
& \multicolumn{3}{|c|}{KNITRO} & \multicolumn{3}{c|}{COBYLA} & \multicolumn{3}{|c|}{KNITRO} & \multicolumn{3}{c}{COBYLA} \\
problem	&	$\phi(x)$	&	feaserr	&	\#fevals	&	$\phi(x)$	&	feaserr	&	\#fevals	&		$\phi(x)$	&	feaserr	&	\#fevals	&	$\phi(x)$	&	feaserr	&	\#fevals	\\
\midrule
CB3	&	2.151366	&	0.00e+00	&	200	&	2.297095	&	3.71e-02	&	f	&	2.151799	&	0.00e+00	&	183	&	2.297095	&	3.71e-02	&	f	\\
CHACONN2	&	1.917274	&	2.54e-01	&	f	&	0.544480	&	1.53e+00	&	f	&	1.917274	&	2.54e-01	&	f	&	0.544480	&	1.53e+00	&	f	\\
GIGOMEZ3	&	2.162155	&	0.00e+00	&	f	&	1.558092	&	5.89e-01	&	f	&	2.162155	&	0.00e+00	&	f	&	1.554036	&	5.32e-01	&	20	\\
HS100	&	682.785696	&	1.01e-01	&	f	&	680.545842	&	2.09e-01	&	f	&	682.785696	&	1.01e-01	&	f	&	680.679716	&	3.84e-01	&	25	\\
DIPIGRI	&	682.612679	&	1.38e-01	&	f	&	681.544064	&	0.00e+00	&	49	&	682.612679	&	1.38e-01	&	f	&	681.282921	&	7.16e-02	&	39	\\
\bottomrule
\end{tabular}
\caption{Number of function evaluations to obtain \eqref{eq:term_eps} for noise level $\sigma_f = 10^{-1}$.}
\label{table:noisyfeas5}
\end{table}

\end{document}